\documentclass[12pt,reqno]{amsart}
\usepackage{amsmath}
\usepackage{amssymb}
\usepackage{amstext}
\usepackage{mathrsfs}
\usepackage{a4wide}
\usepackage{graphicx}
\usepackage{bbm}
\usepackage{verbatim}
\usepackage{bm}
\usepackage{graphicx}
\usepackage{tikz}
\usepackage{pgfplots}

\allowdisplaybreaks \numberwithin{equation}{section}
\usepackage{color}
\usepackage{cases}
\usepackage{hyperref}
\hypersetup{hypertex=true,
	colorlinks=true,
	linkcolor=blue,
	anchorcolor=blue,
	citecolor=blue }

\numberwithin{equation}{section}

\newtheorem{theorem}{Theorem}[section]
\newtheorem{proposition}[theorem]{Proposition}

\newtheorem{lemma}[theorem]{Lemma}

\theoremstyle{definition}

\theoremstyle{remark}
\newtheorem{remark}[theorem]{Remark}

\newcommand{\ep}{\varepsilon}

\newcommand{\R}{\mathbb{R}}

\begin{document}

	\title
	[Steady vortex rings of small cross-section]{Existence, uniqueness and stability of steady vortex rings of small cross-section}
	
	\author{Daomin Cao, Guolin Qin, Weilin Yu, Weicheng Zhan, Changjun Zou}
	
	\address{Institute of Applied Mathematics, Chinese Academy of Sciences, Beijing 100190, and University of Chinese Academy of Sciences, Beijing 100049,  P.R. China}
	\email{dmcao@amt.ac.cn}
	
	\address{School of Mathematical Sciences, Peking University, Beijing 100871, P.R. China}
	\email{qinguolin18@mails.ucas.ac.cn}
	       
	\address{Yanqi Lake Beijing Institute of Mathematical Sciences and Applications, Beijing 101408, P.R. China}
	\email{weilinyu@amss.ac.cn}
	
	\address{School of Mathematical Sciences, Xiamen University, Xiamen, Fujian, 361005, P.R. China}
	\email{zhanweicheng@amss.ac.cn}	
	
	\address{Department of Mathematics, Sichuan University, Chengdu, Sichuan, 610064, P.R. China}
	\email{zouchangjun@amss.ac.cn}
	

	\begin{abstract}           
    This paper is concerned with steady vortex rings in an ideal fluid of uniform density, which are special global axi-symmetric solutions of the three-dimensional incompressible Euler equation. We systematically establish the existence, uniqueness and nonlinear orbital stability of steady vortex rings of small cross-section for which the potential vorticity is constant throughout the core. The latter two answer a long-standing question since the pioneering work of Fraenkel and Berger \cite{BF1} (Acta Math., 1974). To achieve our goal, we rescale the Stokes stream function of vortex ring by its cross-section radius, and expand it at the well-known Rankine vortex using Taylor's formula, where the estimates for coefficients are obtained by a decomposition according to Green's function and local Pohozaev identities. The main observations are: The stream function is even and has a translational invariance in $z$-direction; the zero point for the $r$-coefficient of linear term in its expansion determines the asymptotic location of vortex ring, which appears as the condition to eliminate the degenerate direction in Lyapunov-Schmidt reduction argument for existence; the non-vanishing condition of the second order $r$-coefficient at foresaid zero point is verified as one of the essential factors for uniqueness, while the negativity means that these vortex rings maximizes the functional composed of kinetic energy and impulse. By applying the Arnol'd's dual variational principle together with the uniqueness result, we are finally able to prove the nonlinear orbital stability of thin vortex rings. This result gives a large class of stable vortex rings supported on topological tori, which is different from Hill's spherical vortex discussed by Choi \cite{Choi20} (Comm. Pure Appl. Math., 2023).  Our approach also provides a general theory for well-posedness of localized steady vortex model, including vortex rings of small cross-section and regularization for planar point vortices. 
	\end{abstract}
	
	\maketitle{\small{\bf Keywords:} The 3D Euler equation, Steady vortex rings, Existence, Uniqueness, Nonlinear stability   \\
		
		}
	
	\tableofcontents
	
	\section{Introduction and main results}
	The motion of particles in an ideal fluid in $\mathbb R^3$ is described by the following Euler equation
	\begin{align}\label{1-1}
		\begin{cases}
			\partial_t\mathbf v+(\mathbf v\cdot \nabla)\mathbf v=-\nabla P,
			\\
			\nabla\cdot\mathbf v=0,
		\end{cases}
	\end{align}
 where  $\mathbf{v}(\boldsymbol x,t)$ is the velocity field	and $P(\boldsymbol x,t)$ is the pressure. Corresponding to $\mathbf{v}$ is its vorticity vector defined by $\pmb{\omega}:=\nabla\times\mathbf{v}$. Taking curl of the first equation in Euler equation \eqref{1-1}, H. Helmholtz obtained the equation for vorticity
	\begin{align}\label{1-2}
		\begin{cases}
			\partial_t \pmb{\omega}+(\mathbf{v}\cdot\nabla)\pmb{\omega}=(\pmb{\omega}\cdot\nabla)\mathbf{v},
			\\
			\mathbf{v}=\nabla\times(-\Delta)^{-1}\pmb{\omega}.
		\end{cases}
	\end{align}
	We refer to \cite{Che, MB} for more details about this system.

	We are interested in solutions of the Euler equation whose vorticities are large and uniformly concentrated near an evolving smooth curve embedded in entire $\mathbb R^3$. This type of solutions, \emph{vortex filaments}, have been a subject of active studies for a long time. By the first Helmholtz theorem, in $\mathbb R^3$ a vortex must form a loop with compact support. The simplest vortex loop is a circular \emph{vortex ring}, whose analysis traces back to the works of Helmholtz \cite{Hel} in 1858 and Lord Kelvin \cite{Tho} in 1867. Vortex rings are an intriguing marvel of fluid dynamics that can be easily observed experimentally, e.g. when smoke is ejected from a tube, a bubble rises in a liquid, or an ink is dropped in another fluid, and so on. We refer the reader to \cite{Akh, MGT, Sa92} for some good historical reviews of the achievements in experimental, analytical, and numerical studies of vortex rings.
	
	Helmholtz detected that vortex rings have an approximately steady form and travel with a large constant velocity along the axis of the ring. In 1970, Fraenkel \cite{Fra1} (see also \cite{Fra2}) provided a first constructive proof for the existence of a vortex ring concentrated around a torus with fixed radius $r^*$ with a small, nearly singular cross-section $\epsilon>0$, traveling with constant speed $\sim|\ln\ep|$, rigorously establishing the behavior predicted by Helmholtz. Indeed, Lord Kelvin and Hicks showed that such a vortex ring would approximately move at the velocity (see \cite{Lam,Tho})
	\begin{equation}\label{KH}
		\frac{\kappa}{4\pi r^*}\left(\ln\frac{8r^*}{\epsilon}-\frac{1}{4}\right),
	\end{equation}
	where $\kappa$ denotes its circulation. Fraenkel's result is consistent with the Kelvin--Hicks formula \eqref{KH}.

	Roughly speaking, vortex rings can be characterized simply as an axi-symmetric flow with a (`thin' or `fat') toroidal vortex tube. Here the word `toroidal' means topologically equivalent to a torus. In the usual cylindrical coordinate frame $\{\mathbf{e}_r, \mathbf{e}_\theta, \mathbf{e}_z\}$, the velocity field $\mathbf{v}$ of an axi-symmetric flow can be expressed in the following way
	\begin{equation*}
		\mathbf{v}=v^r(r,z)\mathbf{e}_r+v^\theta(r,z)\mathbf{e}_\theta+v^z(r,z)\mathbf{e}_z.
	\end{equation*}
	The component $v^\theta$ in the $\mathbf{e}_\theta$ direction is usually called the swirl velocity. If an axi-symmetric flow is non-swirling (i.e., $v^\theta \equiv 0$), then the vorticity admits its angular component $\omega^\theta$ only, namely, $\pmb{\omega}=\omega^\theta \mathbf{e}_\theta$. Let $\zeta=\omega^\theta/r$ be the potential vorticity. Then the vorticity equation \eqref{1-2} is reduced to an active scalar equation for $\zeta$
	\begin{equation}\label{1-3}
		\partial_t \zeta+\mathbf{v}\cdot\nabla\zeta=0,\ \ \ \mathbf{v}=\nabla\times(-\Delta)^{-1}\left(r\zeta \right).
	\end{equation}
	We shall refer to an axi-symmetric non-swirling flow as `\emph{vortex ring}' if there is a toroidal region inside of which $\boldsymbol{\omega}\not = 0$ (the core), and outside of which $\boldsymbol{\omega}= 0$.
	By a \emph{steady vortex ring} we mean a vortex ring that moves vertically at a constant speed forever without changing its shape or size. In other words, a steady vortex ring is of the form
	\begin{equation}\label{1-4}
		\zeta(\boldsymbol x,t)=\zeta(\boldsymbol x+t\mathbf{v}_\infty),
	\end{equation}
	where $\mathbf{v}_\infty=-W\mathbf{e}_z$ is a constant propagation speed. Substituting \eqref{1-4} into \eqref{1-3}, we arrive at a stationary equation
	\begin{equation}\label{1-5}
		\left(\mathbf{v}_\infty+\mathbf{v}\right)\cdot\nabla\zeta=0,\ \ \ \mathbf{v}=\nabla\times(-\Delta)^{-1}\left(r\zeta \right).
	\end{equation}
	In 1894, Hill \cite{Hil} found an explicit solution of \eqref{1-5} supported in a sphere. In 1972, Norbury \cite{Nor72} provided a constructive proof for the existence of steady vortex rings with constant $\zeta$ that are close to Hill's vortex but are homeomorphic to a solid torus; and he also presented some numerical results for the existence of a family of steady vortex rings of small cross-section \cite{Nor73}. General existence results of steady vortex rings with a given vorticity function was first established by Fraenkel--Berger \cite{BF1} in 1974. Following these pioneering works, the existence and abundance of steady vortex rings has been rigorously established; see \cite{AS, AT88, Bad, CWZ, DW2,FT,MGT, Ni, VS, YJ} and the references therein.

	Compared with the results on the existence, rather limited work has been done on the uniqueness of steady vortex rings. In 1986, Amick--Fraenkel \cite{AF} proved that Hill's vortex is the unique solution when viewed in a natural weak formulation by the method of moving planes; and they (1988) \cite{AF88} also established local uniqueness for Norbury's nearly spherical vortex. However, to the best of our current knowledge, the uniqueness of steady vortex rings of small cross-section is still open. The first goal of this paper is to give an answer to this question.

	The stability for steady flows is a classical subject of study in fluid dynamics. In \cite{GS1,GS2,GS3}, the authors used the Arnol'd's variational principle to study the stability viscous or inviscid vorticies. Very recently, Choi \cite{Choi20} established the orbital stability of Hill's spherical vortex. We would like to mention that Hill's vortex and vortex ring of small cross-section are two extremal cases of a one-parameter family (see \cite{Nor73}), and it is still not clear whether there exists a class of stable steady vortex ring whose support is on a topological torus (a stable vortex ring in Norbury's class is given in \cite{CQZZ3}). Recent numerical computations in \cite{Pro} revealed that while `thin' vortex rings remain neutrally stable to axi-symmetric perturbations, they become linearly unstable to such perturbations when they are sufficiently `fat'. By virtue of our local uniqueness result, we will establish orbital stability of a family of steady vortex rings of small cross-section, which is also the second main goal of this paper.
	
	To explain the main idea of our approach, let us first consider the regularization for steady planer $k$-point vortices (the construction of a series of $k$ compact-supported vortices approximating the $k$-point vortices) discussed in \cite{CPY,SV,T}. In \cite{Cao4,CQZZ2,CYZ}, the authors proved that these $k$ vortices must locate near the critical point of the Kirchhoff-Routh function $\mathcal K_k$  defined as
	\begin{equation*}
		\mathcal{K}_k(\boldsymbol x_1,\boldsymbol x_2,..,\boldsymbol x_k):=\frac{1}{2}\sum\limits_{i\neq j}^k \kappa_i\kappa_jG(\boldsymbol x_i,\boldsymbol x_j)+\frac{1}{2}\sum\limits_{i=1}^k\kappa_i^2h(\boldsymbol x_i,\boldsymbol x_i)
	\end{equation*}
	with $\boldsymbol x_i\in\mathbb R^2$ the locations, $\kappa_i\in\mathbb R$ the circulation of each vortex, $G$ the Green's function for $-\Delta$ and $h$ its regular part, see \cite{Lin}. This phenomenon is not a miracle since the $k$-point vortex system is dominated by 
	$$\kappa_i \dot {\boldsymbol x}_i=\nabla^{\perp}_{\boldsymbol x_i}\mathcal K_k, \quad 1\le i\le k,$$
	and in particular, $(\boldsymbol x_1,\boldsymbol x_2,..,\boldsymbol x_k)$ must be a critical point of $\mathcal{K}_k$ when the system constitutes an equilibrium. Actually, if one rescales the stream function $\Psi_\varepsilon$ for steady patch type regularization series (The vortices are characteristic functions of $k$ small domains) near the $i^{\mathrm{th}}$ vortex by its radius parameter $\varepsilon$ and expand it by Taylor's formula at the limiting function, then for the point $\varepsilon\boldsymbol y+\boldsymbol x_i$ near the $i^{\mathrm{th}}$ vortex boundary we have
	\begin{equation*}
		\begin{split}
	          \tilde\Psi_{\varepsilon,i}(\boldsymbol y ):&=\Psi_\varepsilon(\varepsilon\boldsymbol y+\boldsymbol x_i)\\
	          &=\frac{\kappa_i}{\pi}\cdot w(\boldsymbol y)+\tilde\phi_\varepsilon+\frac{\kappa_i}{2\pi}\ln\frac{1}{\varepsilon}+\frac{\nabla_{\boldsymbol{x}_i}\mathcal K_k}{\kappa_i}\cdot\varepsilon \boldsymbol y+\varepsilon\boldsymbol y\cdot\frac{\nabla^2_{\boldsymbol {x}_i}\mathcal K_k}{2\kappa_i}\cdot\varepsilon \boldsymbol y^\perp+O(\varepsilon^3),
	    \end{split}
	\end{equation*}
	with 
		\begin{equation*}
		w(\boldsymbol y)=\left\{
		\begin{array}{lll}
			\frac{1}{4}(1-|\boldsymbol y|^2), \ \ \ &\mathrm{if} \ |\boldsymbol y|\le 1,\\
			\\
			\frac{1}{2}\ln\frac{1}{|\boldsymbol y|}, &\mathrm{if} \ |\boldsymbol y|\ge 1
		\end{array}
		\right.
	\end{equation*}
	the stream function of well-known Rankine vortex, and $\tilde\phi_\varepsilon$ an $O(\varepsilon^2)$-perturbation term induced by the second order term. Now the condition $(\boldsymbol x_1,\boldsymbol x_2,..,\boldsymbol x_k)$ being a critical point of $\mathcal{K}_k$ is equal to the linear term of above expansion being nearly zero for $i=1,\cdots,k$. Moreover, the authors in \cite{Cao4} applied the assumption $\det(\nabla^2_{\boldsymbol x_i}\mathcal K_k)\neq 0$ for $i=1,\cdots,k$ together with a non-degenerate condition of limiting function (The kernel of the linearized operator $\mathbb Lv=-\Delta v+v\cdot\boldsymbol \delta_{|\boldsymbol y=1|}$ is a two-dimensional space spanned by $\partial_1 w$ and $\partial_2 w$) to prove the local uniqueness of regularization, which means there is only one regularization series for prescribed location, circulation and vorticity distribution of each vortex. On the other hand, taking advantage of integration by part, the kinetic energy $\int |\nabla \Psi_\varepsilon|^2$ of flow is equal to $\varepsilon^{-2}\int_{D_\varepsilon}\Psi_\varepsilon$ with $D_\varepsilon$ the union of $k$ vortex domain. Thus for $i=1,\cdots,k$ all eigenvalues of matrix $\nabla^2_{\boldsymbol x_i}\mathcal K_k$ being negative means that the system attains its strict local maximum, which is the main observation used in \cite{CYZ} to derive the nonlinear stability for the regularization series. 
	
	In this paper, we will also expand the Stokes stream function for vortex rings of small cross-section at the limiting function, and show that these phenomena of concentrated vorticity follow the same pattern, which involves a lot of hard work on estimate for coefficients. However, since the problem has a $z$-directional translation invariance in view of the symmetry property obtained in Appendix \ref{appA}, we can assume the stream function is even with respect to $r$-axis, and hence simply the discussion for coefficients in $z$-direction. To be more precise, let $\psi_0(r,z)=Wr^2|\ln \varepsilon|/2$ be the Stokes stream function of translational velocity $W|\ln\varepsilon|\mathbf e_z$, $\boldsymbol q_\varepsilon=(q_\varepsilon,0)$ be the asymptotic location of cross-section center, and $s_\varepsilon$ of order $O(\varepsilon)$ be the asymptotic radius for vortex cross-section. Then for a point $(r,z)=s_\varepsilon\boldsymbol y+\boldsymbol q_\varepsilon=(s_\varepsilon y_1+q_\varepsilon,s_\varepsilon y_2)$ near the cross-section boundary, we have the expansion for rescaled stream function.
	\begin{equation}\label{exp}
		\begin{split}
			\tilde\psi_\varepsilon(\boldsymbol y)-\tilde\psi_0(\boldsymbol y):&=\psi_\varepsilon(s_\varepsilon\boldsymbol y+\boldsymbol q_\varepsilon)-\frac{W}{2}\cdot(s_\varepsilon y_1+q_\varepsilon)^2\ln\frac{1}{\varepsilon}\\
			&=V_{\boldsymbol q_\varepsilon,\varepsilon}(\boldsymbol y)+\phi_\varepsilon(s_\varepsilon\boldsymbol y+\boldsymbol q_\varepsilon)+\frac{q_\varepsilon\kappa}{2\pi}\ln\frac{1}{\varepsilon}-\left(Wq_\varepsilon\ln\frac{1}{\varepsilon}-\frac{\kappa}{4\pi}\ln\frac{8q_\varepsilon}{s_\varepsilon}+\frac{\kappa}{16\pi}\right)\cdot s_\varepsilon y_1\\
			&\quad-\frac{3\kappa}{16\pi q_\varepsilon}\ln\frac{1}{\varepsilon}\cdot (s_\varepsilon |\boldsymbol y|)^2-\frac{W}{2}\ln\frac{1}{\varepsilon}\cdot(s_\varepsilon y_1)^2+O(\varepsilon^2)\\
			&=V_{\boldsymbol q_\varepsilon,\varepsilon}(\boldsymbol y)+\phi_\varepsilon(s_\varepsilon\boldsymbol y+\boldsymbol q_\varepsilon)+\frac{q_\varepsilon\kappa}{2\pi}\ln\frac{1}{\varepsilon}-\mathbf c_1\cdot s_\varepsilon \tilde r\\
			&\quad-\mathbf c_{2,\varepsilon}\cdot\ln\frac{1}{\varepsilon}\cdot (s_\varepsilon |\boldsymbol y|)^2-\mathbf c_{2,W}\cdot\ln\frac{1}{\varepsilon}\cdot(s_\varepsilon y_1)^2+O(\varepsilon^2)
		\end{split}
	\end{equation}
    with $V_{\boldsymbol q_\varepsilon,\varepsilon}(\boldsymbol y)$ the limit function defined in \eqref{approxsolu}, $\phi_\varepsilon$ an $O(\varepsilon^2|\ln\varepsilon|)$-perturbation term whose estimate is given in \eqref{phiestimate}, and $\mathbf c_0,\mathbf c_1,\mathbf c_{2,\varepsilon},\mathbf c_{2,W}$ four positive coefficients, see also Lemma \ref{B1} and Remark \ref{remarkB} in Appendix \ref{appB}. We conclude the main observations as three criteria in the follows:
	
	\begin{itemize}
		\item [(a).]  \emph{The $r$-coordinate $q_\varepsilon$ of vortex ring location is the unique zero point of coefficient $\mathbf c_1$ in \eqref{exp}, which corresponds to the Kelvin-Hicks Formula \eqref{KH} exactly.}
		\item [(b).]  \emph{The coefficient $\mathbf c_{2,\varepsilon}\neq 0$ in \eqref{exp} is one of the necessary condition for uniqueness of vortex rings. (The other is the non-degeneracy condition on the limiting function $w$.)}
		\item [(c).]  \emph{The condition $\mathbf c_{2,\varepsilon}>0$ and $\mathbf c_{2,W}>0$ in \eqref{exp} means the vortex maximizing the functional composed of kinetic energy and impulse, which leads to the nonlinear orbital stability of vortex rings.}
	\end{itemize}
	
\bigskip
	
	We shall focus on steady vortex rings for which $\zeta$ is a constant throughout the core. As remarked by Fraenkel \cite{Fra2}, this simplest of all admissible vorticity distributions has been a favorite for over a century, including the well-known Hill's spherical vortex and Norbury's nearly spherical vortex. Now, we turn to state our main results. To this end, we need to introduce some notations. We shall say that a scalar function $\vartheta:\mathbb R^3\to \mathbb R$ is axi-symmetric if it has the form of $\vartheta(\boldsymbol x)=\vartheta(r,z)$, and a subset $\Omega\subset \mathbb R^3$ is axi-symmetric if its characteristic function $\boldsymbol 1_\Omega$ is axi-symmetric. The cross-section parameter $\sigma$ of an axi-symmetric set $\Omega\subset \mathbb R^3$ is defined by
	\begin{equation*}
		\sigma(\Omega):=\frac{1}{2}\cdot\sup\left\{\boldsymbol \delta_{z}(\boldsymbol x,\boldsymbol y)\,\,|\,\,\boldsymbol x,\boldsymbol y\in \Omega \right\},
	\end{equation*}
	where the axi-symmetric distance $\boldsymbol \delta_z$ is given by
	\begin{equation*}
		\boldsymbol \delta_z(\boldsymbol x,\boldsymbol y):=\inf\left\{|\boldsymbol x-Q(\boldsymbol y)|\,\,\,\,| \,\, \ Q \ \text{is a rotation around}\ \mathbf{e}_z\right\}.
	\end{equation*}
	Let $\mathcal{C}_r=\{\boldsymbol x\in\mathbb{R}^3~|x_1^2+x_2^2=r^2,x_3=0\}$ be a circle of radius $r$ on the plane perpendicular to $\mathbf{e}_z$. For an axi-symmetric set $\Omega\subset \mathbb{R}^3$, we define the axi-symmetric distance between $\Omega$ and $\mathcal C_r$ as follows
	\begin{equation*}
		\text{dist}_{\mathcal{C}_r}(\Omega)=\sup_{\boldsymbol x\in \Omega}\inf_{\boldsymbol x'\in{\mathcal{C}_r}}|\boldsymbol x-\boldsymbol x'|.
	\end{equation*}
	The circulation of a steady vortex ring $\zeta$ is given by
	\begin{equation*}
		\frac{1}{2\pi}\int_{\mathbb{R}^3}\zeta(\boldsymbol x)d\boldsymbol x.
	\end{equation*}
	A steady vortex ring $\zeta$ is said to be \emph{centralized} if $\zeta$ is symmetric non-increasing in $z$, namely,
	\begin{equation*}
		\begin{split}
			& \zeta(r,z)=\zeta(r,-z),\ \ \text{and} \\
			& \zeta(r,z)\ \text{is a non-increasing function of}\ z\ \text{for}\ z>0,\,\,\,\text{for each fixed}\ r>0.
		\end{split}
	\end{equation*}
	
	Our first main result is on the existence of steady vortex rings of small cross-section for which $\zeta$ is constant throughout the core. The existence for such kind of solutions was proved in \cite{CWZ,Fra2,FT} by different methods. However, we will construct steady vortex rings from a new perspective of Stokes stream function, where we use the scaled Rankine vortex as the back ground solution and transform the construction into a fix point problem of perturbation term. This process also casts a profound light on our decomposition of stream function in Section \ref{sec3} for uniqueness.
	
	\begin{theorem}[Existence]\label{thm1}
		Let $\kappa$ and $W$ be two positive numbers. Then there exists a small number $\varepsilon_0>0$ such that, for every $\varepsilon\in (0,\varepsilon_0]$ there is a centralized steady vortex ring $\zeta_\varepsilon$ with fixed circulation $\kappa$ and translational velocity $W\ln \varepsilon\,\mathbf{e}_z$. Moreover,
		\begin{itemize}
			\item [(i)]$\zeta_\varepsilon=\varepsilon^{-2}\boldsymbol 1_{\Omega_\varepsilon}$ for some axi-symmetric topological torus $\Omega_\varepsilon\subset \mathbb R^3$.
			\item [(ii)]It holds $C_1\varepsilon \le \sigma\left(\Omega_\varepsilon\right)<C_2\varepsilon$ for some constants $0<C_1<C_2$.
			\item [(iii)]As $\varepsilon\to 0$, $\mathrm{dist}_{\mathcal C_{r^*}}( \Omega_\varepsilon)\to0$, where $r^*={\kappa}/{4\pi W}$.
		\end{itemize}
	\end{theorem}
	
	Our result on the existence is established by an improved Lyapunov-Schmidt reduction argument on planar vortex patch problem in \cite{CPY}, where the degenerate direction is eliminated by vortex ring location consistent with foresaid condition (a). Compared with the method used in \cite{CPY}, our approach in the present paper is the modulation argument being used to deal with a non-uniform elliptic operator and Heaviside nonlinearity for the first time. To obtain desired estimates, we use an equivalent integral formulation of the problem, and introduce a weighted $L^\infty$ norm to handle the degeneracy at infinity and singularity near $z$-axis. Another difficulty in our construction is the lack of compactness, which arises from whole-space $\mathbb R^3$. To overcome it, we will use a few techniques, so that versions of Ascoli--Arzel\`a theorem can be applied to recover the compactness.
	
	There are similar results on the existence for different types of steady vortex rings in the works \cite{AS, Bad,CWZ,DV, Fra1, Fra2, FT}. For instance, de Valeriola et al. \cite{DV} constructed vortex rings with $C^{1,\alpha}$ regularity by mountain pass theorem, and recently Cao et al. \cite{CWZ} studied desingularization of vortex rings by solving variational problems for the potential vorticity $\zeta$. However, in the absence of a comprehensive uniqueness theory, the problem that whether the solutions with fixed vorticity distributions obtained by various methods coincide remains unsolved. Our second main result is to address this question.
	
	\begin{theorem}[Uniqueness]\label{thm2}
		Let $\kappa$ and $W$ be two positive numbers. Let $\{\zeta^{(1)}_\varepsilon\}_{\varepsilon>0}$ and $\{\zeta^{(2)}_\varepsilon\}_{\varepsilon>0}$ be two families of centralized steady vortex rings with the same circulation $\kappa$ and translational velocity $W\ln \varepsilon\,\mathbf{e}_z$. If, in addition,
		\begin{itemize}
			\item [(i)]$\zeta^{(1)}_\varepsilon=\varepsilon^{-2}\boldsymbol 1_{\Omega^{(1)}_\varepsilon}$ and $\zeta^{(2)}_\varepsilon=\varepsilon^{-2}\boldsymbol 1_{\Omega^{(2)}_\varepsilon}$ for certain axi-symmetric topological tori $\Omega^{(1)}_\varepsilon$, $\Omega^{(2)}_\varepsilon\subset \mathbb R^3$.
			\item [(ii)]As $\varepsilon\to 0$, $\sigma\left(\Omega^{(i)}_\varepsilon\right)<L\varepsilon/2$ for $i=1,2$ with a constant $L>0$.
			\item [(iii)]There exists a $\delta_0>0$ such that $\Omega^{(1)}_\varepsilon \cup \Omega^{(2)}_\varepsilon\subset \left\{\boldsymbol x\in \mathbb{R}^3\mid \sqrt{x_1^2+x_2^2}\ge \delta_0 \right\}$  for all $\varepsilon>0$.
		\end{itemize}
		Then there exists a small $\varepsilon_0>0$ such that $\zeta^{(1)}_\varepsilon\equiv\zeta^{(2)}_\varepsilon$ for all $\varepsilon \in (0,\varepsilon_0]$.
	\end{theorem}
     
     \begin{remark}
     	Although our uniqueness theorem is concerned with the case where the potential vorticity is constant throughout the core, the approach in this paper can be easily applied to a more general situation. As we will discuss in the following sections, if we let the Stokes stream function $\psi$ satisfy a semilinear equation 
     	$$\mathcal L \psi=\lambda (\psi-\frac{1}{2}\mathscr{W}r^2-\mu)_+^\gamma$$ 
     	with $\gamma\in[0,\infty)$, then the potential vorticity $\zeta=\mathcal L \psi$ will form a thin vortex ring automatically. Theorem \ref{thm2} deals with $\gamma=0$ in particular, which is the most difficult case due to the discontinuity of potential vorticity, and the 2-dimensional kernel of linearized operator plays an important role in the proof of uniqueness. In Section 3, the stream function of Rankine vortex
     	is proved to be the limiting function of $\psi$ after scaling by $\sigma_\varepsilon$, and a series of measure type estimates are carried out to obtain uniqueness.  While for $\gamma\in (0,1)\cup(1,\infty)$, the uniqueness can be proved similarly: the limiting function $w_\gamma$ turns to satisfy
     	\begin{equation*}
     		w_\gamma(\boldsymbol y):=\left\{
     		\begin{array}{lll}
     			v_\gamma(|\boldsymbol y|), \ \ \ & \mathrm{if} \ |\boldsymbol y|\le 1,\\
     			|v'_\gamma(1)|\ln\frac{1}{|\boldsymbol y|}, & \mathrm{if} \ |\boldsymbol y|>1,
     		\end{array}
     		\right.
     	\end{equation*}
     	with $v_\gamma(\boldsymbol y)=v_\gamma(|\boldsymbol y|)$ the unique radial solution of 
     	\begin{equation*}
     		-\Delta v=v^\gamma, \ \ v\in H^1_0(B_{1}(\boldsymbol 0)), \ \ v>0 \ \ \text{in} \ B_{1}(\boldsymbol 0),
     	\end{equation*}
        which is easier to handle than $\gamma=0$ (the former Rankine vortex case) since the potential vorticity is continuous. (The limiting function for case $\gamma=0$ in $B_{1}(\boldsymbol 0)$ is the first eigenfunction of eigenvalue problem $-\Delta v=\lambda v^*$, which will be the same as in \cite{CYZ}.) On the other hand, for $\gamma\in(0,\infty)$ the estimate in Lemma \ref{C4} will be different according to the value of $\gamma$ (this does not lead to essential difficulty), which yields minor difference in the second main term of Kelvin-Hicks formula \eqref{KH}, see Theorem 3.1 in \cite{Fra1}. 
     	
     	With this general uniqueness result in hand, a dual variational principle can be used to derive the nonlinear stability for a large class of thin vortex rings following the foresaid condition (c) for \eqref{exp} and strategy in Section \ref{sec4}.  
     \end{remark}
	
	To obtain the uniqueness, we first give a rough estimate for the Stokes stream function of vortex rings by blow up analysis and Taylor's expansion near Rankine vortex. Then we improve the estimate step by step, and obtain a more accurate version of Kelvin--Hicks formula \eqref{KH} in Proposition \ref{prop3-2}. Actually, our result is slightly stronger than Fraenkel's in \cite{Fra1} by a careful study of vortex boundary and a bootstrap procedure. With a delicate enough estimate in hand, a local Pohozaev identity can be used to obtain the information of the second order term $\mathbf c_{2,\varepsilon}$ in the Taylor's expansion \eqref{exp} at limiting function, and derive a contradiction if there are two different vortex rings satisfying assumptions in Theorem \ref{thm2}. It is notable that the methods in \cite{AF,AF88} depend strongly on specific distribution of vorticity in cross-section.  While our method has much broader applicability. We also conjecture that the uniqueness stated in theorem \ref{thm2} will hold for vortex rings of small cross-section with general vorticity distribution by condition (b) for \eqref{exp}, provided the corresponding limiting function satisfies a non-degeneracy condition.
	
	Using the uniqueness result in Theorem \ref{thm2}, we can further show that the solutions constructed in Theorem \ref{thm1} are orbitally stable in the Lyapunov sense. Recalling \eqref{1-3}, for an axi-symmetric flow without swirl, the vorticity equation \eqref{1-2} can be reduced to the active scalar equation for the potential vorticity $\zeta=\omega^\theta/r$:
	\begin{equation}\label{1-7}
		\begin{cases}
			\partial_t \zeta+\mathbf{v}\cdot\nabla\zeta=0,\,\qquad\,\,\,\,\, \, \,  \,  \boldsymbol x\in \mathbb R^3,\ \ t>0, \\
			\mathbf{v}=\nabla\times(-\Delta)^{-1}\left(r\zeta \right),\,\,\,\,\,\,\,\boldsymbol x\in \mathbb R^3,\ \ t>0,\\
			\zeta|_{t=0} =\zeta_0,\,\qquad\qquad\qquad\,\,  \boldsymbol x\in \mathbb R^3.
		\end{cases}
	\end{equation}
	The existence and uniqueness of solutions $\zeta(x,t)$ can be studied analogously as the two-dimensional case. We refer to \cite{B96,Choi20, MB, Nobi, Sain, Ukh} for some discussion in this direction. Let $BC([0,\infty);X)$ denote the space of all bounded continuous functions from $[0,\infty)$ into a Banach space $X$. Define the weighted space $L^1_\text{w}(\mathbb R^3)$ by 
	$$L^1_\text{w}(\mathbb R^3)=\{\vartheta: \mathbb R^3 \to \mathbb R\ \text{measurable} \mid r^2\vartheta\in L^1(\mathbb R^3)\}.$$ 
	We introduce the kinetic energy of the fluid
	\begin{equation*}
		E[\zeta]:=\frac{1}{2}\int_{\mathbb{R}^3}|\mathbf{v}(\boldsymbol x)|^2d\boldsymbol x,\ \ \ \mathbf{v}=\nabla\times(-\Delta)^{-1}\left(r\zeta \right),
	\end{equation*}
	and its impulse
	\begin{equation*}
		\mathcal{P}[\zeta]=\frac{1}{2}\int_{\mathbb{R}^3}r^2\zeta(\boldsymbol x)d \boldsymbol x=\pi\int_\Pi r^3 \zeta drdz.
	\end{equation*}
	The following result has been established, see e.g. Lemma 3.4 in \cite{Choi20}.
	\begin{proposition}\label{Pro1}
		For any non-negative axi-symmetric function $\zeta_0\in L^1\cap L^\infty\cap L^1_\mathrm{w}(\mathbb R^3)$ satisfying $r\zeta_0\in L^\infty(\mathbb R^3)$, there exists a unique weak solution $\zeta\in BC([0,\infty);L^1\cap L^\infty\cap L^1_\mathrm{w}(\mathbb R^3))$ of \eqref{1-7} for the initial data $\zeta_0$ such that $\zeta(\cdot,t)\ge 0$ and is axi-symmetric,
		\begin{equation*}
			\begin{array}{ll}
				\|\zeta(\cdot,t)\|_{L^p(\mathbb{R}^3)} =\|\zeta_0\|_{L^p(\mathbb{R}^3)},\ \ 1\le p\le \infty, &\\
				\mathcal{P}[\zeta(\cdot,t)] = \mathcal{P}[\zeta_0], \,\,\,E[\zeta(\cdot,t)] =E[\zeta_0],\ \ \ \text{for all}\ t>0,&
			\end{array}
		\end{equation*}
and, for any $0<\upsilon_1<\upsilon_2<\infty$ and for all $t>0$,
		\begin{equation*}
			\int_{\{\boldsymbol x\in\mathbb{R}^3\,\,\mid\,\, \upsilon_1<\zeta(\boldsymbol x,t)<\upsilon_2\}}\zeta(\boldsymbol x,t)d\boldsymbol x=\int_{\{\boldsymbol x\in\mathbb{R}^3\,\,\mid\,\, \upsilon_1<\zeta_0(\boldsymbol x)<\upsilon_2\}}\zeta_0(\boldsymbol x)d\boldsymbol x.
		\end{equation*}
	\end{proposition}
	
	\bigskip
	
	By denoting the space of admissible functions as
	\begin{equation*}
		\mathcal{F}_\varepsilon:=\left\{\zeta\in L^\infty(\mathbb R^3)\mid \zeta: \text{axi-symmetric},\ 0\le  \zeta\le 1/\varepsilon^2, \ \|\zeta\|_{L^1(\mathbb R^3)}\le2\pi\kappa \right\},
	\end{equation*}
	we will see that the vortex ring $\zeta_\varepsilon$ is a maximizer of the maximization problem
	\begin{equation*}
		\mathcal{E}_{\varepsilon}=\sup_{\zeta\in \mathcal{F}_\varepsilon}\left(E[\zeta]-W\ln\frac{1}{\varepsilon}\mathcal{P}[\zeta]   \right),
	\end{equation*}
	which is the consequence of foresaid condition (c) and Riesz rearrangement inequality, see \cite{Bad,CWZ}. Then we can apply Theorem \ref{thm2} to claim that the maximizer $\zeta_\varepsilon$ is unique modulo a $z$-directional translation, and use the concentration-compactness principle to verify the nonlinear orbital stability as follows.
	\begin{theorem}[Stability]\label{thm4}
		The steady vortex ring $\zeta_\varepsilon$ in Theorem \ref{thm1} is stable up to  translations in the following sense:
		
		For any $\eta>0$, there exists $\delta_1>0$ such that for any non-negative axi-symmetric function $\zeta_0$ satisfying $\zeta_0, r\zeta_0\in L^\infty(\mathbb R^3)$ and
		\begin{equation*}
			\|\zeta_0-\zeta_\varepsilon\|_{L^1\cap L^2(\mathbb R^3)}+\|r^2(\zeta_0-\zeta_\varepsilon)\|_{L^1(\mathbb R^3)}\le \delta_1,
		\end{equation*}
		the corresponding solution $\zeta(\boldsymbol x,t)$ of \eqref{1-7} for the initial data $\zeta_0$ satisfies
		\begin{equation*}
			\inf_{\tau\in \mathbb R}\left\{\|\zeta(\cdot-\tau \mathbf{e}_z,t)-\zeta_\varepsilon\|_{L^1\cap L^2(\mathbb R^3)}+\|r^2(\zeta(\cdot-\tau \mathbf{e}_z,t)-\zeta_\varepsilon\|_{L^1(\mathbb R^3)} \right\}\le \eta,\,\,\,\,\text{for all}\,\,t>0.
		\end{equation*}
Here, $\|\cdot\|_{L^1\cap L^2(\mathbb R^3)}$ means $\|\cdot\|_{L^1(\mathbb R^3)}+\|\cdot\|_{L^2(\mathbb R^3)}$.
	\end{theorem}
	
	Except for bringing a new perspective of expansion for Stokes stream function near the scaled Rankine vortex as in Appendix \ref{appB}, our another innovation is that a modified Pohozaev identity is used to obtain the delicate estimate for the coefficients in this expansion, see Appendix \ref{appC}. There are essential difficulties in this paper compared with \cite{Cao4,CYZ} since the elliptic operator related to axi-symmetric flow is no longer Laplacian. Fortunately, inspired by the construction in Section \ref{sec2}, we will decompose the stream function of vortex ring into two parts, where one is corresponding to a scaled Rankine vortex with its mirror vortex (to make the boundary condition consistent), and the other is the remaining regular part. Then we can use a modified Pohozaev identity \eqref{C-1} to derive our desired estimate. In fact, the concrete form of expansion for stream function is much more easier to obtain compared with estimates for their coefficients. With suitable tools in hand to handle these estimates, we believe that a similar approach can be applied to verify the uniqueness and nonlinear stability of localized vorticity discussed in \cite{EWZ} for water wave equation, and in \cite{Ao} for gSQG (generalized surface quasi-geostrophic) equation. 
	
	The paper is organized as follows. In Section \ref{sec2}, we construct vortex rings of small cross-section by a Lyapunov--Schmidt reduction argument. In Section \ref{sec3}, we expand the Stokes stream function of vortex ring near the scaled Rankine vortex, and prove the uniqueness result stated in Theorem \ref{thm2}. A revised Kelvin-Hicks formula is obtained as a byproduct in Proposition \ref{prop3-2}. Then the nonlinear orbital stability is proved in Section \ref{sec4} based on the uniqueness and variational method. In Appendix \ref{appA}, we use the method of moving planes to obtain the $z$-directional symmetry of the problem. The expansion for stream function is calculated in Appendix \ref{appB}, where we also give a precise description for the cross-section boundary. In Appendix \ref{appC}, a modified Pohozaev identity is established to derive delicate estimates for coefficients in the expansion, which is the key element in the proof of uniqueness in Section \ref{sec3}.
	
	\bigskip
	\bigskip
	
	\section{Existence}\label{sec2}
	
	The main objective of this paper is to study steady vortex rings, which are actually traveling-wave solutions for \eqref{1-7}.
	Thanks to the continuity equation in \eqref{1-1}, we can find a Stokes stream function $\Psi$ such that
	\begin{equation*}
		\mathbf{v}=\frac{1}{r}\left(-\frac{\partial\Psi}{\partial z}\mathbf{e}_r+\frac{\partial\Psi}{\partial r}\mathbf{e}_z\right).
	\end{equation*}
	In terms of the Stokes stream function $\Psi$, the problem of steady vortex rings can be reduced to a steady problem on the meridional half plane $\Pi=\{(r,z)\mid r>0\}$ of the form:
	\begin{numcases}
		{ }
		\label{2-1} \mathcal{L}\Psi =0, \,\ \, \ \ \ \ \ \ \  \ \, &\text{in}\  $\Pi \setminus A$,  \\
		\label{2-2} \mathcal{L}\Psi=\lambda f_0(\Psi), \ \ \ \  &\text{in}\  $A$,\\
		\label{2-3}  \Psi(0,z)=-\mu \le 0, \ \ z\in\mathbb R&\\
		\label{2-4} \Psi=0, \ \ \ \  &\text{on}\  $\partial A$,\\
		\label{2-5} \frac{1}{r}\frac{\partial \Psi}{\partial r} \to -\mathscr{W}\ \text{and}\ \frac{1}{r}\frac{\partial \Psi}{\partial z} \to 0,\ \ \text{as}\ \ r^2+z^2\to \infty,&
	\end{numcases}
	where
	\begin{equation*}
		\mathcal{L}:=-\frac{1}{r}\frac{\partial}{\partial r}\Big(\frac{1}{r}\frac{\partial}{\partial r}\Big)-\frac{1}{r^2}\frac{\partial^2}{\partial z^2}.
	\end{equation*}
	Here the positive vorticity function $f_0$ and the vortex-strength parameter $\lambda>0$ are prescribed; $A$ is the (a priori unknown) cross-section of the vortex ring; $\mu$ is called the flux constant measuring the flow rate between the $z$-axis and $\partial A$; The constant $\mathscr{W}>0$ is the translational speed, and the condition \eqref{2-5} means that the limit of the velocity field $\mathbf{v}$ at infinity is $-\mathscr{W} \mathbf{e}_z$. For a detailed derivation of this system, we refer to \cite{AF, Choi20, BF1} and the references therein.
	
	By the maximum principle, we see that $\Psi>0$ in $A$ and $\Psi<0$ in $\Pi\backslash \bar{A}$. Therefore the cross-section $A$ is given by
	\begin{equation*}
		A=\left\{(r,z)\in \Pi \mid \Psi(r,z)>0 \right\}.
	\end{equation*}
	It is convenient to write
	\begin{equation*}
		\Psi(r,z)=\psi(r,z)-\frac{1}{2}\mathscr{W}r^2-\mu,
	\end{equation*}
	where $\psi$ is the stream function due to vorticity. In addition, it is also convenient to define
	\[
	f(\tau)=\left\{
	\begin{array}{ll}
		0, &    \tau\le 0, \\
		f_0(\tau),                  & \tau>0,
	\end{array}
	\right.
	\]
	so that $\lambda f(\Psi)$ is exactly the potential vorticity $\zeta$. We now can rewrite \eqref{2-1}-\eqref{2-5} as
	\begin{numcases}
		{(\mathscr{P})\ \ \  }  
		\label{2-6} \mathcal{L}\psi =\lambda f(\psi-\frac{1}{2}\mathscr{W}r^2-\mu),\ \ \ \text{in}\ \Pi,  &  \\
		\label{2-7}  \psi(0,z)=0, &\\
		\label{2-8} \psi,\ \ {|\nabla \psi|}/{r} \to 0 \ \ \text{as}\ \ r^2+z^2\to \infty.&
	\end{numcases}
	In the following, we will focus on the construction of $\psi$ satisfying $(\mathscr{P})$.
	
	\bigskip
	
	\subsection{Formulation of the problem}
	
	In order to simplify notations, in Section \ref{sec2} and \ref{sec3} we will use  $\mathbb R^2_+=\{\boldsymbol x=(x_1,x_2) \ | \ x_1>0\}$
	to substitute the meridional half plane $\Pi$, and abbreviate the elliptic operator $\mathcal L$ as
	\begin{equation}\label{delta*}
		\Delta^*:=\frac{1}{x_1}\text{div}\left(\frac{1}{x_1}\nabla\right).
	\end{equation}
	We will use $\varepsilon:=\lambda^{-1/2}$ as the  parameter instead of $\lambda$ in the rest of this paper. Since we are concerned with steady vortex rings for which $\zeta$ is a constant throughout the core, we will choose the vorticity function $f$ in \eqref{2-6} having the following form
	\[
	f(\tau)=\left\{
	\begin{array}{ll}
		0, &    \tau\le0, \\
		1,                  & \tau>0,
	\end{array}
	\right.
	\]
	and the cross-section of the vortex ring is
	$$A_\varepsilon=\left\{\boldsymbol x\in \mathbb R^2_+ \ \big| \ \psi_\varepsilon-\frac{W}{2}x_1^2\ln\frac{1}{\varepsilon}>\mu_\varepsilon\right\}$$ for some flux constant $\mu_\varepsilon>0$.
	Here we let $\mathscr{W}$ equal $W|\ln \varepsilon|$ according to Kelvin--Hicks formula \eqref{1-3}. The fact that $\mu_\varepsilon>0$ means $A_\varepsilon$ will not touch the $x_2$-axis. Thus we can rewrite $(\mathscr{P})$ as
	\begin{equation}\label{2-9}
		\begin{cases}
			-\varepsilon^2{\Delta^*}\psi=
			\boldsymbol1_{\left\{\psi-\frac{W}{2}x_1^2\ln\frac{1}{\varepsilon}>\mu_\varepsilon\right\}}, & \text{in} \ \mathbb R^2_+,
			\\
			\psi=0, & \text{when} \ x_1=0,
			\\
			\psi, \ |\nabla\psi|/x_1\to0, &\text{as} \ |\boldsymbol x |\to \infty.
		\end{cases}
	\end{equation}
	Since \eqref{2-9} is invariant in $x_2$-direction, we may assume
	\begin{equation}\label{2-10}
		\psi(x_1,x_2)=\psi(x_1,-x_2)
	\end{equation}
	due to the method of moving planes in Appendix \ref{appA} (see also Lemma 2.1 in \cite{AF88}), which means the steady vortex ring $\zeta_\varepsilon$ corresponding to $\psi_\varepsilon$ is centralized.
	
	The existence result in Theorem \ref{thm1} can be deduced from the following proposition.
	\begin{proposition}\label{prop2-1}
		For every $\kappa>0$ and $W>0$, there exists an $\varepsilon_0>0$ such that for each $\varepsilon\in (0,\varepsilon_0]$, problem \eqref{2-9} has a solution $\psi_\varepsilon$ satisfying \eqref{2-10}. Moreover,
		\begin{itemize}
			\item[(i)] The cross-section $A_\varepsilon$ is a convex domain, and satisfies
			\begin{equation*}
				B_{\sqrt{\frac{\kappa}{q_\varepsilon\pi}}\varepsilon(1-L_1\varepsilon|\ln\varepsilon|)}(\boldsymbol q_\varepsilon)\subset A_\varepsilon\subset B_{\sqrt{\frac{\kappa}{q_\varepsilon\pi}}\varepsilon(1+L_2\varepsilon|\ln\varepsilon|)}(\boldsymbol q_\varepsilon),
			\end{equation*}
			where $L_1$, $L_2$ are two positive constants independent of $\varepsilon$, and $\boldsymbol q=(q_\varepsilon,0)$ is on $x_1$-axis with the estimate
			$$q_{\varepsilon}=\frac{\kappa}{4\pi W}+O\left(\frac{1}{|\ln\varepsilon|}\right).$$
			\item[(ii)]
			As $\varepsilon\to 0$, it holds
			\begin{equation*}
				\kappa_\varepsilon:=\varepsilon^{-2}\int_{A_\varepsilon}x_1d\boldsymbol x\to \kappa.
			\end{equation*}
		\end{itemize}
	\end{proposition}
	\begin{remark}
		Notice that in Proposition \ref{prop2-1}, the circulation parameter $\kappa_\varepsilon$ is not fixed, which only has the limiting behavior described in property (ii). To obtain a family of vortex rings with fixed circulation $\kappa$ as in Theorem \ref{thm1}, we can rescale $\psi_\varepsilon$ as follows
		$$\bar \psi_\varepsilon(\boldsymbol x):=\frac{\kappa_\varepsilon^2}{\kappa^2}\cdot \psi_\varepsilon\left(\frac{\kappa}{\kappa_\varepsilon}\cdot \boldsymbol x\right).$$
		Then $\bar \psi_\varepsilon(\boldsymbol x)$ satisfies
		\begin{equation*}
			-\bar\varepsilon^2\Delta^*\bar\psi_\varepsilon=\boldsymbol 1_{\{\bar\psi_\varepsilon-\frac{W}{2}x_1^2\ln\frac{1}{\varepsilon}>\bar\mu_\varepsilon\}},
		\end{equation*}
		where
		$$\bar\varepsilon=\frac{\kappa_\varepsilon}{\kappa}\cdot \varepsilon, \quad \text{and} \quad \bar\mu_\varepsilon=\frac{\kappa^2}{\kappa_\varepsilon^2}\cdot \mu_\varepsilon.$$
		It is easy to verify that
		$$\frac{1}{\bar \varepsilon^2}\int_{\mathbb R^2_+} x_1\boldsymbol 1_{\{\bar\psi_\varepsilon-\frac{W}{2}x_1^2\ln\frac{1}{\varepsilon}>\bar\mu_\varepsilon\}}d\boldsymbol x=\kappa,$$
		and the vortex ring $\bar\zeta_\varepsilon$ corresponding to $\bar\psi_\varepsilon$ satisfies all assumptions in Theorem \ref{thm1}.
	\end{remark}

	For the study of steady vortex rings of small cross-section, our main tool is the Green's representation of Stokes stream function $\psi_\varepsilon$. To be more rigorous, $\psi_\varepsilon$ satisfies the integral equation
	\begin{equation}\label{2-11}
		\psi_\varepsilon(\boldsymbol x)=\frac{1}{\varepsilon^2}\int_{\mathbb R^2_+}G_*(\boldsymbol x,\boldsymbol x') \boldsymbol 1_{\{\boldsymbol x'=(x'_1,x'_2)\in \mathbb R^2_+ \,\, |\, \, \psi_\varepsilon(\boldsymbol x')-\frac{W}{2}(x'_1)^2\ln\frac{1}{\varepsilon}>\mu_\varepsilon\}}(\boldsymbol x')d\boldsymbol x',
	\end{equation}
	where ${G_*}(\boldsymbol x,\boldsymbol x')$ is the Green's function for $-{\Delta^*}$ with boundary condition in \eqref{2-9}. Using Biot--Savart law in $\mathbb R^3$ and a coordinate transformation, we can derive an explicit formula of ${G_*}(\boldsymbol x,\boldsymbol x')$ as
	\begin{equation*}
		{G_*}(\boldsymbol x,\boldsymbol x')=\frac{x_1x_1'^2}{4\pi}\int_{-\pi}^\pi\frac{\cos\theta d\theta}{\left[(x_2-x_2')^2+x_1^2+x_1'^2-2x_1 x_1'\cos\theta\right]^{\frac{1}{2}}}.
	\end{equation*}
	Then, by denoting the distance
	\begin{equation}\label{rho}
		\rho(\boldsymbol x, \boldsymbol x')=\frac{(x_1-x_1')^2+(x_2-x_2')^2}{x_1x_1'},
	\end{equation}
	we have the following asymptotic estimates
	\begin{equation}\label{2-12}
		\begin{split}
		{G_*}(\boldsymbol x,\boldsymbol x^\prime)&=
		\frac{x_1^{1/2}x_1^{\prime 3/2}}{4\pi}\left(\ln\left(\frac{1}{\rho}\right)
		+2\ln 8-4\right)\\
		&\quad+\frac{x_1^{1/2}x_1^{\prime 3/2}}{4\pi}\left(\frac{3}{16}\cdot\rho\ln\frac{1}{\rho}+\left(\frac{3\ln8}{8}-\frac{1}{8}\right)\rho+o(\rho)\right),\quad \text{as} \ \rho\to 0,
		\end{split}
	\end{equation}
	and
	\begin{equation}\label{2-13}
		{G_*}(\boldsymbol x,\boldsymbol x^\prime) =\frac{x_1^{1/2}x_1^{\prime 3/2}}{4}\left(\frac{1}{\rho^{3/2}}+O(\rho^{-5/2})\right), \quad \text{as} \  \rho\to \infty,
	\end{equation}
	which can be found in \cite{Fen,Fra2,GS3,Lam,Sve}. Actually, the theory of elliptic integrals can be used to obtain a more precise expansion of ${G_*}$ in $\rho$.
	
	To simplify integral equation \eqref{2-11}, we let $\boldsymbol q_\varepsilon=(q_\varepsilon,0)$ with $q_\varepsilon>0$ bounded and determined later near the limiting position of the cross-section, and split ${G_*}$ as
	\begin{equation*}
		{G_*}(\boldsymbol x,\boldsymbol x')=q_\varepsilon^2G(\boldsymbol x,\boldsymbol x')+H(\boldsymbol x,\boldsymbol x'),
	\end{equation*}
	where
	\begin{equation*}
		G(\boldsymbol x,\boldsymbol x')=\frac{1}{4\pi}\ln\frac{(x_1+x_1')^2+(x_2-x_2')^2}{(x_1-x_1')^2+(x_2-x_2')^2},
	\end{equation*}
	is the Green's function for $-\Delta$ in right half plane, and $H(\boldsymbol x,\boldsymbol x')$ is a relatively regular function. Note that the regular part is also dependent on our choice of $\boldsymbol q_\varepsilon$ in this decomposition. However, we leave out the parameter $q_\varepsilon$ in $H(\boldsymbol x,\boldsymbol x')$ for simplicity in this paper, which is because the estimates of $H(\boldsymbol x,\boldsymbol x')$ used in this section and next (for the seek of existence and uniqueness) are uniform for $\boldsymbol q_\varepsilon$ in the cross-section $A_\varepsilon$. For instance, by using Lemma 2.1 in \cite{CQZ} on a more general lake equation and let the depth function $b(\boldsymbol x)=1/x_1$, we can conclude that $H(\boldsymbol x,\boldsymbol x')\in C_{\text{loc}}^\alpha(\mathbb R^2_+\times \mathbb R^2_+)$ for every $\alpha\in (0,1)$, and $H(\boldsymbol x,\boldsymbol x')$ is quasi-Lipschitz near a fixed $\boldsymbol x'\in \mathbb R^2_+$, namely, for any $\boldsymbol x^{(1)},\boldsymbol x^{(2)}$ in $B_\delta(\boldsymbol x')\subset \mathbb R^2_+$ with $\delta$ a small constant, there exists a $C_\delta$ such that
	$$|H(\boldsymbol x^{(1)},\boldsymbol q_\varepsilon)-H(\boldsymbol x^{(2)},\boldsymbol q_\varepsilon)|\le C_\delta\cdot|\boldsymbol x^{(1)}-\boldsymbol x^{(2)}|(1+\ln|\boldsymbol x^{(1)}-\boldsymbol x^{(2)}|).$$
	
	Our argument of constructing  solutions to \eqref{2-9} is divided into several steps, which is known as the Lyapunov--Schmidt reduction or modulation. We will first give a series of approximate solutions of $\psi_\varepsilon$, so that \eqref{2-9} is transformed into a semilinear equation of the error term $\phi_\varepsilon$. Then, we establish the linear theory of the corresponding projected problem. The existence and limiting behavior of $\psi_\varepsilon$ will be obtained by contraction mapping theorem and one-dimensional reduction in the last part of our proof.

	\subsection{Approximate solutions}
	
	To give suitable approximate solutions to \eqref{2-9} and \eqref{2-10}, let us consider the following problem
	\begin{equation*}
		\begin{cases}
			-\varepsilon^2\Delta V_{\boldsymbol q_\varepsilon,\varepsilon}(\boldsymbol x)=q_\varepsilon^2\boldsymbol{1}_{B_{s_\varepsilon}(\boldsymbol q_\varepsilon)},  \ \ \ & \text{in} \ \mathbb R^2,\\
			V_{\boldsymbol q_\varepsilon,\varepsilon}(\boldsymbol x)=\frac{a_\varepsilon}{2\pi}\ln\frac{1}{\varepsilon}, &\text{on} \ \partial B_{s_\varepsilon}(\boldsymbol q_\varepsilon),
		\end{cases}
	\end{equation*}
	with $\boldsymbol q_\varepsilon=(q_\varepsilon,0)\in \mathbb R^2$, $q_\varepsilon\neq 0$ and $a_\varepsilon,s_\varepsilon$ are two parameters to be determined later. Actually, $V_{\boldsymbol q_\varepsilon,\varepsilon}(\boldsymbol x)$ corresponds to the stream function of a planar Rankine vortex after scaling, and we can write $V_{\boldsymbol q_\varepsilon,\varepsilon}(\boldsymbol x)$ explicitly as
	\begin{equation}\label{approxsolu}
		V_{\boldsymbol q_\varepsilon,\varepsilon}(\boldsymbol x)=\left\{
		\begin{array}{lll}
			\frac{a_\varepsilon}{2\pi}\ln\frac{1}{\varepsilon}+\frac{q_\varepsilon^2}{4\varepsilon^2}(s_\varepsilon^2-|\boldsymbol x-\boldsymbol q_\varepsilon|^2), \ \ \ &\mathrm{if} \ |\boldsymbol x-\boldsymbol q_\varepsilon|\le s_\varepsilon,\\
			\frac{a_\varepsilon}{2\pi}\ln\frac{1}{\varepsilon}\cdot\frac{\ln|\boldsymbol x-\boldsymbol q_\varepsilon|}{\ln s_\varepsilon},&\mathrm{if} \ |\boldsymbol x-\boldsymbol q_\varepsilon|\ge s_\varepsilon.
		\end{array}
		\right.
	\end{equation}
    To make $V_{\boldsymbol q_\varepsilon,\varepsilon}(\boldsymbol x)\in C^1(\mathbb R^2)$, it will hold the following relationship
	\begin{equation}\label{2-14}
		\mathcal N_\varepsilon:=\frac{a_\varepsilon}{2\pi}\ln\frac{1}{\varepsilon}\cdot\frac{1}{s_\varepsilon|\ln s_\varepsilon|}=\frac{s_\varepsilon}{2\varepsilon^2}\cdot q_\varepsilon^2,
	\end{equation}
	where $\mathcal N_\varepsilon$ is the value of gradient $|\nabla V_{\boldsymbol q_\varepsilon,\varepsilon}|$ on $\partial B_{s_\varepsilon}(\boldsymbol q_\varepsilon)$. From \eqref{2-14}, we see that $s_\varepsilon$ is asymptotically linearly dependent on $\varepsilon$ as
	$$s_\varepsilon=\left(\sqrt{\frac{a_\varepsilon}{\pi q_\varepsilon^2}}+o_\varepsilon(1)\right)\varepsilon.$$

	In our construction, $V_{\boldsymbol q_\varepsilon,\varepsilon}(\boldsymbol x)$ will be used as the building block of approximate solutions. To make it agree with the zero boundary value on $\partial \mathbb R^2_+$, for general $\boldsymbol x=(x_1,x_2)\in \mathbb R^2_+$ we denote $\boldsymbol {\bar x}=(-x_1, x_2)$ as the reflection of $\boldsymbol x$ with respect to $x_2$-axis, and let
	\begin{equation*}
		\begin{split}
			&\quad\mathcal V_{\boldsymbol q_\varepsilon,\varepsilon}(\boldsymbol x):=V_{\boldsymbol q_\varepsilon,\varepsilon}(\boldsymbol x)-V_{\boldsymbol{\bar q_\varepsilon},\varepsilon}(\boldsymbol x)\\
			&=\frac{1}{2\pi\varepsilon^2}\int_{\mathbb R^2_+}q_\varepsilon^2\ln\left(\frac{1}{|\boldsymbol x-\boldsymbol x'|}\right)\boldsymbol{1}_{B_{s_\varepsilon}(\boldsymbol q_\varepsilon)}(\boldsymbol x')d\boldsymbol x'- \frac{1}{2\pi\varepsilon^2}\int_{\mathbb R^2_+}q_\varepsilon^2\ln\left(\frac{1}{|\boldsymbol x-\boldsymbol{\bar x}'|}\right)\boldsymbol{1}_{B_{s_\varepsilon}(\boldsymbol q_\varepsilon)}(\boldsymbol x')d\boldsymbol x'\\
			&=\frac{q_\varepsilon^2}{\varepsilon^2}\int_{\mathbb R^2_+}G(\boldsymbol x,\boldsymbol x') \boldsymbol 1_{B_{s_\varepsilon}(\boldsymbol q_\varepsilon)}(\boldsymbol x')d\boldsymbol x',
		\end{split}
	\end{equation*}
	where $V_{\boldsymbol{\bar q_\varepsilon},\varepsilon}(\boldsymbol x)$ is the stream function of a mirror vortex on the left half plane, and $G(\boldsymbol x,\boldsymbol x')$ is the Green's function for $-\Delta$ on $\mathbb R^2_+$. Then $\mathcal V_{\boldsymbol q_\varepsilon,\varepsilon}(\boldsymbol x)$ is the unique solution to the following problem
	\begin{equation*}
		\begin{cases}
			-\varepsilon^2\Delta \mathcal{V}_{\boldsymbol q_\varepsilon,\varepsilon}(\boldsymbol x)=q_\varepsilon^2\boldsymbol{1}_{B_{s_\varepsilon}(\boldsymbol q_\varepsilon)},  \ \ \ & \text{on} \ \mathbb R^2_+,\\
			\mathcal{V}_{\boldsymbol q_\varepsilon,\varepsilon}=0, &\text{if} \ x_1=0,\\
			\mathcal{V}_{\boldsymbol q_\varepsilon,\varepsilon}, \ |\nabla\mathcal{V}_{\boldsymbol q_\varepsilon,\varepsilon}|/x_1\to0, &\text{as} \ |\boldsymbol x|\to \infty.
		\end{cases}
	\end{equation*}
	To approximate the remaining part of $\psi_\varepsilon$, let
	\begin{equation*}
		\mathcal H_{\boldsymbol q_\varepsilon,\varepsilon}(\boldsymbol x)=\frac{1}{\varepsilon^2}\int_{\mathbb R^2_+}H(\boldsymbol x,\boldsymbol x')\boldsymbol{1}_{B_{s_\varepsilon}(\boldsymbol q_\varepsilon)}(\boldsymbol x')d\boldsymbol x'.
	\end{equation*}
	According to the definition of $H(\boldsymbol x,\boldsymbol x')$, it is obvious that $\mathcal H_{\boldsymbol q_\varepsilon,\varepsilon}(\boldsymbol x)$ solves
	\begin{equation*}
		\begin{cases}
			-\varepsilon^2{\Delta^*}\left(\mathcal V_{\boldsymbol q_\varepsilon,\varepsilon}+\mathcal H_{\boldsymbol q_\varepsilon,\varepsilon}\right)=\boldsymbol{1}_{B_{s_\varepsilon}(\boldsymbol q_\varepsilon)}, \ \ \ &\text{on} \ \mathbb R^2_+,\\
			\mathcal{H}_{\boldsymbol q_\varepsilon,\varepsilon}=0, &\text{if} \ x_1=0,\\
			\mathcal{H}_{\boldsymbol q_\varepsilon,\varepsilon}, \ |\nabla\mathcal{H}_{\boldsymbol q_\varepsilon,\varepsilon}|/x_1\to0, &\text{as} \ |\boldsymbol x|\to \infty.
		\end{cases}
	\end{equation*}
	Moreover, using the definition of $H(\boldsymbol x,\boldsymbol x')$ and computing directly (as the calculation for regular part in Lemma \ref{B1}), we have
	\begin{equation*}
		\mathcal H_{\boldsymbol q_\varepsilon,\varepsilon}(\boldsymbol x)-\frac{s_\varepsilon^2\pi}{\varepsilon^2}H(\boldsymbol x,\boldsymbol q_\varepsilon)=\frac{1}{\varepsilon^2}\int_{\mathbb R^2_+}\left(H(\boldsymbol x,\boldsymbol x')-H(\boldsymbol x,\boldsymbol q_\varepsilon)\right)\boldsymbol{1}_{B_{s_\varepsilon}(\boldsymbol q_\varepsilon)}(\boldsymbol x')d\boldsymbol x'=O(\varepsilon|\ln\varepsilon|),
	\end{equation*}
	and
	\begin{equation*}
		\partial_1\mathcal H_{\boldsymbol q_\varepsilon,\varepsilon}(\boldsymbol x)=\frac{1}{\varepsilon^2}\int_{\mathbb R^2_+}\partial_{x_1}H(\boldsymbol x,\boldsymbol x')\boldsymbol{1}_{B_{s_\varepsilon}(\boldsymbol q_\varepsilon)}(\boldsymbol x')d\boldsymbol x'=O(|\ln\varepsilon|).
	\end{equation*}
	
	After all these preparations, we write a solution $\psi_\varepsilon$ to \eqref{2-9} as
	\begin{equation*}
		\psi_\varepsilon(\boldsymbol x)=\mathcal V_{\boldsymbol q_\varepsilon,\varepsilon}+\mathcal H_{\boldsymbol q_\varepsilon,\varepsilon}+\phi_\varepsilon,
	\end{equation*}
	where $\phi_\varepsilon(\boldsymbol x)$ is a small error term with the boundary condition
	\begin{equation*}
		\begin{cases}
			\phi_\varepsilon=0, & \text{when}\ x_1=0,
			\\
			\phi_\varepsilon, \ |\nabla\phi_\varepsilon|/x_1\to0, &\text{as} \ |\boldsymbol x |\to \infty,
		\end{cases}
	\end{equation*}
	and symmetry condition (by Appendix \ref{appA})
	\begin{equation*}
		\phi_\varepsilon(x_1,x_2)=\phi_\varepsilon(x_1,-x_2),\,\,\,\, \text{for all}\,\,x_1,x_2>0.
	\end{equation*}
	Then we can derive the semi-linear equation for $\phi_\varepsilon$ by direct computations
	\begin{equation*}
		\begin{split}
			0&=-x_1\varepsilon^2{\Delta^*}\left(\mathcal V_{\boldsymbol q_\varepsilon,\varepsilon}+\mathcal H_{\boldsymbol q_\varepsilon,\varepsilon}+\phi_\varepsilon\right)-x_1\boldsymbol1_{\{\psi_\varepsilon-\frac{W}{2}x_1^2\ln\frac{1}{\varepsilon}>\mu_\varepsilon\}}\\
			&=x_1\left(-\varepsilon^2{\Delta^*}(\mathcal V_{\boldsymbol q_\varepsilon,\varepsilon}+\mathcal H_{\boldsymbol q_\varepsilon,\varepsilon})-\boldsymbol1_{\{V_{\boldsymbol q_\varepsilon,\varepsilon}>\frac{a_\varepsilon}{2\pi}\ln\frac{1}{\varepsilon}\}}\right)\\
			& \ \ \ +\varepsilon^2\left(-x_1{\Delta^*}\phi_\varepsilon-\frac{2}{s_\varepsilon q_\varepsilon}\phi_\varepsilon(r,\theta)\boldsymbol\delta_{|\boldsymbol x-\boldsymbol q_\varepsilon|=s_\varepsilon}\right)\\
			& \ \ \ -\bigg(x_1\boldsymbol1_{\{\psi_\varepsilon-\frac{W}{2}x_1^2\ln\frac{1}{\varepsilon}>\mu_\varepsilon\}}-x_1\boldsymbol1_{\{V_{\boldsymbol q_\varepsilon,\varepsilon}>\frac{a_\varepsilon}{2\pi}\ln\frac{1}{\varepsilon}\}}-\frac{2}{s_\varepsilon q_\varepsilon}\phi_\varepsilon(s_\varepsilon,\theta)\boldsymbol\delta_{|\boldsymbol x-\boldsymbol q_\varepsilon|=s_\varepsilon}\bigg)\\
			&=\varepsilon^2\mathbb L_\varepsilon\phi_\varepsilon-\varepsilon^2R_\varepsilon(\phi_\varepsilon),
		\end{split}
	\end{equation*}
	where $\mathbb L_\varepsilon$ is a linear operator defined by
	\begin{equation}\label{2-15}
		\mathbb L_\varepsilon\phi =-x_1{\Delta^*}\phi -\frac{2}{s_\varepsilon q_\varepsilon}\phi (r,\theta)\boldsymbol\delta_{|\boldsymbol x-\boldsymbol q_\varepsilon|=s_\varepsilon}
	\end{equation}
	with the notation $\phi(r,\theta)=\phi(q_\varepsilon+r\cos\theta,r\sin\theta)$, and
	\begin{equation*}
		 R_\varepsilon(\phi)=\frac{1}{\varepsilon^2}\bigg(x_1\boldsymbol1_{\{\psi_\varepsilon-\frac{W}{2}x_1^2\ln\frac{1}{\varepsilon}>\mu_\varepsilon\}}-x_1\boldsymbol1_{\{V_{\boldsymbol q_\varepsilon,\varepsilon}>\frac{a_\varepsilon}{2\pi}\ln\frac{1}{\varepsilon}\}}-\frac{2}{s_\varepsilon q_\varepsilon}\phi(r,\theta)\boldsymbol\delta_{|\boldsymbol x-\boldsymbol q_\varepsilon|=s_\varepsilon}\bigg)
	\end{equation*}
	is the nonlinear perturbation.
	
	To make $R_\varepsilon(\phi_\varepsilon)$ as small as possible so as to apply the contraction mapping theorem to find solutions of $\mathbb L_\varepsilon\phi_\varepsilon=R_\varepsilon(\phi_\varepsilon)$, we choose the parameters $a_\varepsilon$ such that 
	\begin{equation}\label{2-16}
		\frac{a_\varepsilon}{2\pi}\ln\frac{1}{\varepsilon}=\mu_\varepsilon+\frac{W}{2}q_\varepsilon^2\ln\frac{1}{\varepsilon}-\mathcal H_{\boldsymbol q_\varepsilon,\varepsilon}(\boldsymbol q_\varepsilon)+V_{\boldsymbol{\bar q_\varepsilon},\varepsilon}(\boldsymbol q_\varepsilon),
	\end{equation}
	which is always possible since the last two terms in \eqref{2-16} is uniformly bounded. Moreover, we can obtain following asymptotic estimate
	$$a_\varepsilon=q_\varepsilon\kappa+O(1/|\ln\varepsilon|).$$
	For simplicity of further discussion later and in Appendix \ref{appB}, we will denote
	\begin{equation}\label{Udef}
		\begin{split}
		\mathbf U_{\boldsymbol q_\varepsilon,\varepsilon}(\boldsymbol x)&=\psi_\varepsilon(\boldsymbol x)-\phi_\varepsilon(\boldsymbol x)-\frac{W}{2}x_1^2\ln\frac{1}{\varepsilon}-\mu_\varepsilon\\
		&=\mathcal V_{\boldsymbol q_\varepsilon,\varepsilon}(\boldsymbol x)+\mathcal H_{\boldsymbol q_\varepsilon,\varepsilon}(\boldsymbol x)-\frac{W}{2}x_1^2\ln\frac{1}{\varepsilon}-\mu_\varepsilon.
		\end{split}
	\end{equation}
    as a first approximation of $V_{\boldsymbol q_\varepsilon,\varepsilon}(\boldsymbol x)-\frac{a_\varepsilon}{2\pi}\ln\frac{1}{\varepsilon}$. Then we see that $R_\varepsilon(\phi_\varepsilon)$ contains two main parts: one is nonlinear term on $\phi_\varepsilon$, and the other is induced by the difference of $\boldsymbol1_{\{U_{\boldsymbol q_\varepsilon,\varepsilon}>0\}}$ and $\boldsymbol1_{\{V_{\boldsymbol q_\varepsilon,\varepsilon}>\frac{a_\varepsilon}{2\pi}\ln\frac{1}{\varepsilon}\}}$.

	Since the parameters $s_\varepsilon,\mu_\varepsilon, a_\varepsilon$ all depend on $\boldsymbol q_\varepsilon$, problem \eqref{2-9} and \eqref{2-10} is transformed into finding the pairs $(q_\varepsilon,\phi_\varepsilon)$ for each $\varepsilon>0$ small, such that
	\begin{equation}\label{Eqforperturbation}
		\begin{cases}
			\mathbb L_\varepsilon\phi_\varepsilon=R_\varepsilon(\phi_\varepsilon), & \text{in} \ \mathbb R^2_+,
			\\
			\phi_\varepsilon=0, & \text{if} \ x_1=0,
			\\
			\phi_\varepsilon, \ |\nabla\phi_\varepsilon|/x_1\to0, &\text{as} \ |\boldsymbol x |\to \infty.
		\end{cases}
	\end{equation}
    If $\mathbb L_\varepsilon$ is invertible, then one can apply the contraction mapping theorem to solve \eqref{Eqforperturbation} by giving $R_\varepsilon(\phi_\varepsilon)$ a delicate estimate. However, this is not the case, and we will deal with the problem in the next part of construction.

    \bigskip
	
	\subsection{The linear theory}
	
	The purpose of this part is to discuss the properties of linearized operator $\mathbb L_\varepsilon$, and solve \eqref{Eqforperturbation} in a space modulo the degenerate direction. Although we know little about $\mathbb L_\varepsilon$, we see it will tend to the linear operator $\mathbb L^*_\varepsilon$ as $\varepsilon\to 0$, which is defined by
	\begin{equation}\label{linear}
		\mathbb L^*_\varepsilon\phi:=-\frac{1}{q_\varepsilon}\Delta\phi-\frac{2}{s_\varepsilon q_\varepsilon}\phi(r,\theta)\boldsymbol\delta_{|\boldsymbol x-\boldsymbol q_\varepsilon|=s_\varepsilon}.
	\end{equation}
	In \cite{CPY}, the authors dealt with a eigenvalue problem for Laplace–Beltrami operator on $\mathbb S^1$ to show that 
	\begin{equation*}
		\text{ker}(\mathbb L^*_\varepsilon)=\text{span} \left\{\frac{\partial V_{\boldsymbol q_\varepsilon,\varepsilon}}{\partial x_1},\frac{\partial V_{\boldsymbol q_\varepsilon,\varepsilon}}{\partial x_2}\right\},
	\end{equation*}
	where for $m=1,2$,
	\begin{equation*}
		\frac{\partial V_{\boldsymbol q_\varepsilon,\varepsilon}}{\partial x_m}=\left\{
		\begin{array}{lll}
			-\frac{q_\varepsilon^2}{2\varepsilon^2}(x_m-q_m), \ \ \ & \mathrm{if} \ |\boldsymbol x-\boldsymbol q_\varepsilon|\le s_\varepsilon,\\
			-\frac{a_\varepsilon|\ln \varepsilon|}{2\pi|\ln s_\varepsilon|}\frac{x_m-q_m}{|\boldsymbol x-\boldsymbol q_\varepsilon|^2}, & \mathrm{if} \ |\boldsymbol x-\boldsymbol q_\varepsilon|\ge s_\varepsilon
		\end{array}
		\right.
	\end{equation*}
	with $q_1=q_\varepsilon$ and $q_2=0$. 
	
	Recall that $\mathbb L_\varepsilon$ is defined on $\mathbb R^2_+$ and $\phi_\varepsilon$ is even symmetric with respect to $x_1$-axis. When $\varepsilon$ is chosen sufficiently small, $\mathbb L_\varepsilon$ will be very close to $\mathbb L_\varepsilon^*$, and the kernel of $\mathbb L_\varepsilon$ restricted on spaces composed of functions even in $x_2$ can be approximated by the space spanned by
	\begin{equation*}
		Z_{\boldsymbol q_\varepsilon, \varepsilon}=\chi_\varepsilon\cdot\frac{\partial V_{\boldsymbol q_\varepsilon,\varepsilon}}{\partial x_1}
	\end{equation*}
	(the part related to $\partial V_{\boldsymbol q_\varepsilon,\varepsilon}/\partial x_2$ is eliminated by the symmetry assumption), where $\chi_\varepsilon$ is a smooth truncation function satisfying
	\begin{equation}\label{truncation}
		\chi_\varepsilon(\boldsymbol x)=\left\{
		\begin{array}{lll}
			1, \ \ \  & \mathrm{if} \ |\boldsymbol x-\boldsymbol q_\varepsilon|< \delta_\varepsilon,\\
			0, & \mathrm{if} \ |\boldsymbol x-\boldsymbol q_\varepsilon|\ge 2\delta_\varepsilon
		\end{array}
		\right.
	\end{equation}
	for $\delta_\varepsilon=\varepsilon|\ln\varepsilon|$. Moreover, we assume that $\chi_\varepsilon$ are radially symmetric with respect to $\boldsymbol q_\varepsilon$ and
	\begin{equation*}
		|\nabla \chi_\varepsilon|\le \frac{2}{\delta_\varepsilon}, \quad\quad	|\nabla^2 \chi_\varepsilon|\le \frac{2}{\delta_\varepsilon^2}.
	\end{equation*}
    The purpose of adding a truncation $\chi_\varepsilon$ here is to make $Z_{\boldsymbol q_\varepsilon, \varepsilon}$ satisfy a suitable boundary condition so that we can integrate by part in following argument. To solve \eqref{Eqforperturbation}, we will first consider the following projected problem
	\begin{equation}\label{2-17}
		\begin{cases}
			\mathbb L_\varepsilon\phi=\mathbf h(\boldsymbol x)-\Lambda x_1 {\Delta^*}Z_{\boldsymbol q_\varepsilon, \varepsilon}, \ \ &\text{in} \ \mathbb R^2_+,\\
			\int_{\mathbb R^2_+} \frac{1}{x_1} \nabla\phi\cdot\nabla Z_{\boldsymbol q_\varepsilon, \varepsilon}d\boldsymbol x=0,\\
			\phi= 0, \ \ &\text{on} \ x_1=0,\\
			\phi,\ |\nabla\phi|/x_1\to0, \ \ &\text{as} \ |\boldsymbol x|\to \infty,
		\end{cases}
	\end{equation}
	where $\phi$ is even with respect to $x_1$-axis, $\text{supp} \, \mathbf h\subset B_{2s_\varepsilon}(\boldsymbol q_\varepsilon)$, and $\Lambda$ is the projection coefficient such that
	$$\int_{\mathbb R^2_+}Z_{\boldsymbol q_\varepsilon, \varepsilon}(\mathbb L_\varepsilon\phi-\mathbf h+\Lambda x_1 {\Delta^*}Z_{\boldsymbol q_\varepsilon, \varepsilon})d\boldsymbol x=0.$$
	Since $\mathbf h$ has a compact support, the behavior of $\phi$ away from $B_{2s_\varepsilon}(\boldsymbol q_\varepsilon)$ will be the same as the Green's function $G_*(\boldsymbol x,\boldsymbol q_\varepsilon)$. To match the asymptotic behavior of $\phi$ near $\partial \mathbb R^2_+$, we let
	\begin{equation}\label{rho2}
		\rho_1(\boldsymbol x):=\frac{(1+|\boldsymbol x-\boldsymbol q_\varepsilon|^2)^{\frac{3}{2}}}{1+x_1^2} \ \ \ \text{and} \ \ \	\rho_2(\boldsymbol x):=\left(\frac{1}{x_1}+1\right),
	\end{equation}
	and define the weighted $L^\infty$ norm of $\phi$ by
	\begin{equation}\label{2-18}
		||\phi||_*:=\sup_{\boldsymbol x\in\mathbb R^2_+}\rho_1(\boldsymbol x)\rho_2(\boldsymbol x)|\phi(\boldsymbol x)|.
	\end{equation}
	
    \bigskip
	
	To show the projection we make in \eqref{2-17} is reasonable, we will give a coercive estimate for the problem in the following lemma.
	\begin{lemma}\label{lem2-2}
		Assume that $\mathbf h$ satisfies $\mathrm{supp}\, \mathbf h\subset B_{2s_\varepsilon}(\boldsymbol q_\varepsilon)$ and $$\varepsilon^{1-\frac{2}{p}}\|\mathbf h\|_{W^{-1,p}(B_{Ls_\varepsilon}(\boldsymbol q_\varepsilon))}<\infty$$
		with $p\in (2,+\infty]$, then there exists a small $\varepsilon_0>0$ and a positive constant $c_0$ such that for any $\varepsilon\in(0,\varepsilon_0]$ and solution pair $(\phi,\Lambda)$ to \eqref{2-17}, one has
		\begin{equation}\label{2-19}
			\|\phi\|_*+\varepsilon^{1-\frac{2}{p}}\|\nabla\phi\|_{L^p(B_{Ls_\varepsilon}(\boldsymbol q_\varepsilon))}\le c_0\varepsilon^{1-\frac{2}{p}}\|\mathbf h\|_{W^{-1,p}(B_{Ls_\varepsilon}(\boldsymbol q_\varepsilon))},
		\end{equation}
		and
		\begin{equation}\label{2-20}
		   |\Lambda|\le c_0\varepsilon^{2-\frac{2}{p}}\|\mathbf h\|_{W^{-1,p}(B_{Ls_\varepsilon}(\boldsymbol q_\varepsilon))}.
		\end{equation}
	\end{lemma}
	\begin{proof}
		First we are to obtain an estimate for coefficient $\Lambda$. To proceed an energy method, we multiply the first equation in \eqref{2-17} by $Z_{\boldsymbol q_\varepsilon, \varepsilon}$. By integrations by parts we obtain
		\begin{equation}\label{2-21}
			\Lambda\int_{\mathbb R^2_+} \frac{1}{x_1}\nabla Z_{\boldsymbol q_\varepsilon, \varepsilon}\cdot\nabla Z_{\boldsymbol q_\varepsilon, \varepsilon}d\boldsymbol x=\int_{\mathbb R^2_+}Z_{\boldsymbol q_\varepsilon, \varepsilon}\mathbb L_\varepsilon\phi d\boldsymbol x-\int_{\mathbb R^2_+} Z_{\boldsymbol q_\varepsilon, \varepsilon}\mathbf h d\boldsymbol x.
		\end{equation}
		Recall the definition of $Z_{\boldsymbol q_\varepsilon, \varepsilon}$. For the integral in the left hand side of \eqref{2-21}, we have
		\begin{equation*}
			\begin{split}
				&\quad\int_{\mathbb R^2_+} \frac{1}{x_1}\nabla\left(\chi_\varepsilon\cdot\frac{\partial V_{\boldsymbol q_\varepsilon,\varepsilon}}{\partial x_1}\right)\cdot\nabla \left(\chi_\varepsilon\cdot\frac{\partial V_{\boldsymbol q_\varepsilon,\varepsilon}}{\partial x_1}\right)d\boldsymbol x\\
				&=\int_{\mathbb R^2_+} \frac{\chi_\varepsilon^2}{q_\varepsilon}\cdot\left(\nabla\frac{\partial V_{\boldsymbol q_\varepsilon,\varepsilon}}{\partial x_1}\right)^2d\boldsymbol x+\int_{\mathbb R^2_+} \frac{2\chi_\varepsilon\nabla\chi_\varepsilon}{q_\varepsilon}\cdot\left(\nabla\frac{\partial V_{\boldsymbol q_\varepsilon,\varepsilon}}{\partial x_1}\right)\cdot\frac{\partial V_{\boldsymbol q_\varepsilon,\varepsilon}}{\partial x_1}d\boldsymbol x\\
				&\quad+\int_{\mathbb R^2_+} \frac{(\nabla\chi_\varepsilon)^2}{q_\varepsilon}\cdot\left(\frac{\partial V_{\boldsymbol q_\varepsilon,\varepsilon}}{\partial x_1}\right)^2d\boldsymbol x+\frac{C}{\varepsilon^2}\cdot \delta_\varepsilon\\
				&=\frac{C_Z}{\varepsilon^2}\cdot (1+o_\varepsilon(1)),
			\end{split}
		\end{equation*}
		where $C_Z>0$ is some constant independent of $\varepsilon$. Let $\chi^*(\boldsymbol x)$ be a smooth truncation function taking the value $1$ in $B_{2s_\varepsilon}(\boldsymbol q_\varepsilon)$, and $0$ in $\mathbb R^2_+\setminus B_{Ls_\varepsilon}(\boldsymbol q_\varepsilon)$, which means $\chi^*\chi_\varepsilon=\chi_*$ by $\delta_\varepsilon=\varepsilon|\ln\varepsilon|$ and $s_\varepsilon=O(\varepsilon)$. Then it holds 
		\begin{equation*}
			\begin{array}{ll}
				\left\|\nabla\left(\chi^*\chi_\varepsilon\cdot\frac{\partial V_{\boldsymbol q_\varepsilon,\varepsilon}}{\partial x_1}\right)\right\|_{L^{p'}(B_{Ls_\varepsilon}(\boldsymbol q_\varepsilon))}
				&\le \left\|\left(\nabla \chi^*\right)\cdot\frac{\partial V_{\boldsymbol q_\varepsilon,\varepsilon}}{\partial x_1}\right\|_{L^{p'}(B_{Ls_\varepsilon}(\boldsymbol q_\varepsilon))}+\left\|\chi^*\cdot\left(\nabla\frac{\partial V_{\boldsymbol q_\varepsilon,\varepsilon}}{\partial x_1}\right)\right\|_{L^{p'}(B_{Ls_\varepsilon}(\boldsymbol q_\varepsilon))} \\
				 &\le\frac{C}{\varepsilon}\left(\int_{2s_\varepsilon}^{Ls_\varepsilon}\frac{\tau}{\tau^{p'}}d\tau\right)^{\frac{1}{p'}}+\frac{C}{\varepsilon^2}\left(\int_0^{s_\varepsilon}\tau d\tau\right)^{\frac{1}{p'}}+\left(\int_{s_\varepsilon}^{Ls_\varepsilon}\frac{\tau}{\tau^{2p'}}d\tau\right)^{\frac{1}{p'}}\\
				 &=C\varepsilon^{\frac{2}{p'}-2}.
			\end{array}
		\end{equation*}
	     Since $\text{supp}\, \mathbf h \subset B_{2s_\varepsilon}(\boldsymbol q_\varepsilon)$, for the second term in the right hand side of \eqref{2-21}, we have
	    \begin{equation*}
	    	\begin{split}
	    		\left|\int_{\mathbb R^2_+}\chi_\varepsilon\cdot\frac{\partial V_{\boldsymbol q_\varepsilon,\varepsilon}}{\partial x_1}\cdot\mathbf h d\boldsymbol x\right|&=\left|\int_{\mathbb R^2_+}\chi^*\chi_\varepsilon\cdot\frac{\partial V_{\boldsymbol q_\varepsilon,\varepsilon}}{\partial x_1} \cdot\mathbf h d\boldsymbol x \right|\\
	    		&\le C\|\mathbf h\|_{W^{-1,p}(B_{Ls_\varepsilon}(\boldsymbol q_\varepsilon))}\left\|\nabla\left(\chi^*\chi_\varepsilon\cdot\frac{\partial V_{\boldsymbol q_\varepsilon,\varepsilon}}{\partial x_1}\right)\right\|_{L^{p'}(B_{Ls_\varepsilon}(\boldsymbol q_\varepsilon))}\\
	    		&\le C\varepsilon^{\frac{2}{p'}-2} \|\mathbf h\|_{W^{-1,p}(B_{Ls_\varepsilon}(\boldsymbol q_\varepsilon))},
	    	\end{split}
	    \end{equation*}
	    where a Poincar\'{e} inequality
	    $$\left\|\chi^*\chi_\varepsilon\cdot\frac{\partial V_{\boldsymbol q_\varepsilon,\varepsilon}}{\partial x_1}\right\|_{L^{p'}(B_{Ls_\varepsilon}(\boldsymbol q_\varepsilon))}\le C\varepsilon\left\|\nabla\left(\chi^*\chi_\varepsilon\cdot\frac{\partial V_{\boldsymbol q_\varepsilon,\varepsilon}}{\partial x_1}\right)\right\|_{L^{p'}(B_{Ls_\varepsilon}(\boldsymbol q_\varepsilon))}$$
	    is used.
		For the first term in the right hand side of \eqref{2-21}, it holds
		\begin{equation*}
			\begin{split}
				&\quad\int_{\mathbb R^2_+}\chi_\varepsilon\cdot\frac{\partial V_{\boldsymbol q_\varepsilon,\varepsilon}}{\partial x_1}\cdot\mathbb L_\varepsilon\phi d\boldsymbol x=\int_{\mathbb R^2_+} \phi\cdot\mathbb L_\varepsilon\left(\chi_\varepsilon\cdot\frac{\partial V_{\boldsymbol q_\varepsilon,\varepsilon}}{\partial x_1}\right)d\boldsymbol x\\
				&=\int_{\mathbb R^2_+}\frac{1}{x_1}	\nabla\phi\cdot\nabla\left(\chi_\varepsilon\cdot\frac{\partial V_{\boldsymbol q_\varepsilon,\varepsilon}}{\partial x_1}\right)d\boldsymbol x-\frac{2}{s_\varepsilon q_\varepsilon}\int_{|\boldsymbol x-\boldsymbol q_\varepsilon|=s_\varepsilon}\phi\cdot\frac{\partial V_{\boldsymbol q_\varepsilon,\varepsilon}}{\partial x_1}\\
				&=-\int_{\mathbb R^2_+}\phi\cdot\nabla\left(\frac{1}{x_1}\right)\cdot\nabla\left(\chi_\varepsilon\cdot\frac{\partial V_{\boldsymbol q_\varepsilon,\varepsilon}}{\partial x_1}\right)d\boldsymbol x-\int_{\mathbb R^2_+}\phi\left(\frac{1}{x_1}-\frac{1}{q_\varepsilon}\right)\Delta\left(\frac{\partial V_{\boldsymbol q_\varepsilon,\varepsilon}}{\partial x_1}\right)d\boldsymbol x\\
				& \ \ \ -\int_{\mathbb R^2_+}\frac{\phi}{x_1}\cdot\left(2\nabla\chi_\varepsilon\cdot \nabla\left(\frac{\partial V_{\boldsymbol q_\varepsilon,\varepsilon}}{\partial x_1}\right)+(\Delta\chi_\varepsilon)\frac{\partial V_{\boldsymbol q_\varepsilon,\varepsilon}}{\partial x_1}\right)d\boldsymbol x,
			\end{split}
		\end{equation*}
		where we have used the fact that $\partial V_{\boldsymbol q_\varepsilon,\varepsilon}/\partial x_1$ is in the kernel of $\mathbb L^*_\varepsilon$ given in \eqref{linear}. Notice that for terms in above identity we have the following estimates
		\begin{equation*}
			\int_{\mathbb R^2_+}\left|\nabla\left(\chi_\varepsilon\cdot\frac{\partial V_{\boldsymbol q_\varepsilon,\varepsilon}}{\partial x_1}\right)\right|d\boldsymbol x\le C|\ln\varepsilon|,
		\end{equation*}
	    \begin{equation*}
	    	\int_{\mathbb R^2_+}\left|\left(\frac{1}{x_1}-\frac{1}{q_\varepsilon}\right)\Delta\left(\frac{\partial V_{\boldsymbol q_\varepsilon,\varepsilon}}{\partial x_1}\right)\right|d\boldsymbol x\le s_\varepsilon\cdot 2\pi s_\varepsilon\cdot \frac{1}{s_\varepsilon^2}= 2\pi,
	    \end{equation*}
		\begin{equation*}
			\int_{\mathbb R^2_+}\left|\nabla\chi_\varepsilon\cdot\left(\nabla\frac{\partial V_{\boldsymbol q_\varepsilon,\varepsilon}}{\partial x_1}\right)\right|d\boldsymbol x\le \frac{C}{\delta_\varepsilon}\cdot\int_{\delta_\varepsilon}^{2\delta_\varepsilon}\frac{1}{\tau}d\tau\le \frac{C}{\delta_\varepsilon},
		\end{equation*}
		\begin{equation*}
			\int_{\mathbb R^2_+}\left|(\Delta\chi_\varepsilon)\cdot\frac{\partial V_{\boldsymbol q_\varepsilon,\varepsilon}}{\partial x_1}\right|d\boldsymbol x\le \frac{C}{\delta_\varepsilon^2}\cdot\int_{\delta_\varepsilon}^{2\delta_\varepsilon}d\tau\le \frac{C}{\delta_\varepsilon}.
		\end{equation*}
		As a result, it holds
		\begin{equation*}
			\begin{split}
			\left|\int_{\mathbb R^2_+}\chi_\varepsilon\cdot\frac{\partial V_{\boldsymbol q_\varepsilon,\varepsilon}}{\partial x_1}\cdot \mathbb L_\varepsilon\phi d\boldsymbol x\right|&\le (|\ln\varepsilon|+\delta_\varepsilon^{-1})\|\phi\|_{L^\infty(B_{2\delta_\varepsilon}(\boldsymbol q_\varepsilon))}\\
			&\le (|\ln\varepsilon|+\delta_\varepsilon^{-1})\|\phi\|_*.
			\end{split}
		\end{equation*}
		Then combining all above estimates for \eqref{2-21}, we derive
		\begin{equation}\label{2-24}
			|\Lambda|\le C\varepsilon^2(|\ln\varepsilon|+\delta_\varepsilon^{-1})\cdot\|\phi\|_*+C\varepsilon^{\frac{2}{p'}}\|\mathbf h\|_{W^{-1,p}(B_{Ls_\varepsilon}(\boldsymbol q_\varepsilon))},
		\end{equation}
		By the explicit formulation of $Z_{\boldsymbol q_\varepsilon, \varepsilon}$ in $B_{Ls_\varepsilon}(\boldsymbol q_\varepsilon)$, it holds
		$$||x_1{\Delta^*} Z_{\boldsymbol q_\varepsilon, \varepsilon}||_{W^{-1,p}(B_{Ls_\varepsilon}(\boldsymbol q_\varepsilon))}\le C\|\nabla Z_{\boldsymbol q_\varepsilon, \varepsilon}\|_{L^{p}(B_{Ls_\varepsilon}(\boldsymbol q_\varepsilon))}=C\varepsilon^{\frac{2}{p}-2}.$$
		So we finally deduce from \eqref{2-24} that
		\begin{equation*}
			\begin{split}
				||\Lambda x_1{\Delta^*} Z_{\boldsymbol q_\varepsilon, \varepsilon}||_{W^{-1,p}(B_{Ls_\varepsilon}(\boldsymbol q_\varepsilon))}&\le C|\Lambda|\cdot \varepsilon^{\frac{2}{p}-2}\\
				&=C\varepsilon^{\frac{2}{p}}(|\ln\varepsilon|+\delta_\varepsilon^{-1})\cdot\|\phi\|_*+C\|\mathbf h\|_{W^{-1,p}(B_{Ls_\varepsilon}(\boldsymbol q_\varepsilon))}.
			\end{split}
		\end{equation*}
		
		Now we are to prove \eqref{2-19}. Suppose not, then there exist a sequence $\{\varepsilon_n\}$ tending to $0$ and $\{\phi_n\}$ such that
		\begin{equation}\label{2-25}
			\|\phi_n\|_*+\varepsilon_n^{1-\frac{2}{p}}\|\nabla\phi_n\|_{B_{Ls_{\varepsilon_n}}(\boldsymbol q_{\varepsilon_n}))}=1
		\end{equation}
		and
		\begin{equation*}
			\varepsilon_n^{1-\frac{2}{p}}\|\mathbf h\|_{W^{-1,p}(B_{Ls_{\varepsilon_n}}(\boldsymbol q_{\varepsilon_n}))}\le \frac{1}{n}.
		\end{equation*}
		Let
		\begin{equation*}
			\begin{split}
				-\text{div}\left(\frac{1}{x_1}\nabla\phi_n(\boldsymbol x)\right)&=\frac{2}{s_{\varepsilon_n} q_{\varepsilon_n}}\boldsymbol{\delta}_{|\boldsymbol x-\boldsymbol q_{\varepsilon_n}|=s_{\varepsilon_n}}\phi_n(r,\theta)+\mathbf h-\Lambda x_1{\Delta^*} Z_{\boldsymbol q_{\varepsilon_n}, \varepsilon_n}\\
				&=\frac{2}{s_{\varepsilon_n} q_{\varepsilon_n}}\boldsymbol{\delta}_{|\boldsymbol x-\boldsymbol q_{\varepsilon_n}|=s_{\varepsilon_n}}\phi_n(r,\theta)+f_n
			\end{split}
		\end{equation*}
		with $\text{supp}\, f_n\subset B_{2\delta_{\varepsilon_n}}(\boldsymbol q_{\varepsilon_n})$. For a general function $v$, we define its rescaled version centered at $\boldsymbol q_{\varepsilon_n}$ as:
		\begin{equation*}
			\tilde v(\boldsymbol y):=v(s_{\varepsilon_n}\boldsymbol y+\boldsymbol q_{\varepsilon_n}).
		\end{equation*}	
		Notice that parameter $s$ also depends on $\varepsilon_n$. Denoting $D_n=\{\boldsymbol y \ | \ s_{\varepsilon_n}\boldsymbol y+\boldsymbol q_{\varepsilon_n}\in\mathbb R^2_+\}$, then we obtain
		\begin{equation*}
			\int_{D_n}\frac{1}{s_{\varepsilon_n} y_1+q_{\varepsilon_n}}\cdot\nabla\tilde\phi_n\cdot\nabla\varphi d\boldsymbol y=2\int_{|\boldsymbol y|=1}\frac{1}{q_{\varepsilon_n}}\tilde\phi_n\varphi+\langle\tilde f_n,\varphi\rangle, \quad \forall \, \varphi\in C_0^\infty(D_n),
		\end{equation*}
		where for each $p\in(2,\infty]$, it holds
		\begin{equation*}
			\|\tilde f_n\|_{W^{-1,p}(B_{L}(\boldsymbol 0))}\le C\varepsilon_n^{1-\frac{2}{p}}\left(\varepsilon_n^{\frac{2}{p}}(|\ln\varepsilon_n|+\delta_{\varepsilon_n}^{-1})\cdot\|\phi_n\|_*+\|\mathbf h\|_{W^{-1,p}(B_{Ls_{\varepsilon_n}}(\boldsymbol q_{\varepsilon_n}))}\right)= o_n(1).
		\end{equation*}
		Hence $\tilde\phi_n$ is bounded in $C_{\text{loc}}^\alpha(\mathbb{R}^2)$ for some $\alpha>0$, and $\tilde\phi_n$ converges uniformly in any fixed compact set of $\mathbb{R}^2$ to $\phi^*\in L^\infty(\mathbb R^2)\cap C(\mathbb R^2)$, which satisfies
		\begin{equation*}
			-\Delta \phi^*=2\phi^*(r,\theta)\boldsymbol{\delta}_{|\boldsymbol y|=1}, \ \ \  \text{in} \ \mathbb R^2,
		\end{equation*}
		and $\phi^*$ can be written as
		$$\phi^*=C_1\frac{\partial w}{\partial y_1}+C_2\frac{\partial w}{\partial y_2}$$
		with
		\begin{equation*}
			w(\boldsymbol y)=\left\{
			\begin{array}{lll}
				\frac{1}{4}(1-|\boldsymbol y|^2), \ \ \ \ \ &\mathrm{if} \ |\boldsymbol y|\le 1,\\
				\\
				\frac{1}{2}\ln\frac{1}{|\boldsymbol y|}, &\mathrm{if} \ |\boldsymbol y|\ge 1.
			\end{array}
			\right.
		\end{equation*}
		Since $\phi^*$ is even with respect to $x_1$-axis, it holds $C_2=0$. Then, from the second equation in \eqref{2-17}, we have
		\begin{equation*}
			\int_{\mathbb R^2}\nabla\phi^*\nabla\frac{\partial w}{\partial x_1}=0.
		\end{equation*}
		Thus we get $C_1=0$, and $\phi_n\to 0$ in $B_{Ls_{\varepsilon_n}}(\boldsymbol q_{\varepsilon_n})$ as $n\to\infty$.
		
		To derive the estimate for $||\phi_n||_*$, we will use a comparison principle. We see that $\phi_n$ satisfy
		\begin{equation*}
			\begin{cases}
				\phi_n(\boldsymbol x)=0, & \text{on}\ x_1=0,
				\\
				\phi_n, \ |\nabla\phi_n|/x_1\to0, &\text{as} \ |\boldsymbol x |\to \infty.
			\end{cases}
		\end{equation*}
		Moreover, $\phi_n\to 0$ in $B_{Ls_{\varepsilon_n}}(\boldsymbol q_{\varepsilon_n})$ as $n\to\infty$, and $x_1{\Delta^*}\phi_n=0$ in $\mathbb R^2_+\setminus B_{Ls_{\varepsilon_n}}(\boldsymbol q_{\varepsilon_n})$. By letting
		$$\bar \phi_n(\boldsymbol x):= ||\phi_n||_{L^\infty(B_{Ls_{\varepsilon_n}}(\boldsymbol q_{\varepsilon_n}))}\cdot {G_*}(\boldsymbol x,\boldsymbol q_{\varepsilon_n})$$
		as a family of barrier functions, we have
		\begin{equation*}
			\begin{cases}
				\bar \phi_n-\phi_n\ge0, & \text{on}\ x_1=0,
				\\
				\bar \phi_n-\phi_n\ge0, &\text{as} \ |\boldsymbol x |\to \infty,
			\end{cases}
		\end{equation*}
		and
		\begin{equation*}
			x_1^2{\Delta^*}\bar \phi_n-x_1^2{\Delta^*}\phi_n=\Delta(\bar \phi_n-\phi_n)+x_1\nabla\left(\frac{1}{x_1}\right)\cdot\nabla(\bar \phi_n-\phi_n)=0, \ \text{in} \ \mathbb R^2_+\setminus B_{Ls_{\varepsilon_n}}(\boldsymbol q_{\varepsilon_n}).
		\end{equation*}
		Since $x_1\nabla(1/x_1)$ is locally bounded on $\mathbb R^2_+\setminus B_{Ls_{\varepsilon_n}}(\boldsymbol q_{\varepsilon_n})$, we can use the strong maximum principle to deduce $\phi_n\le\bar\phi_n$ on $\mathbb R^2_+\setminus B_{Ls_{\varepsilon_n}}(\boldsymbol q_{\varepsilon_n})$, and hence $|\phi_n|\le\bar\phi_n$ on $\mathbb R^2_+\setminus B_{Ls_{\varepsilon_n}}(\boldsymbol q_{\varepsilon_n})$. By the definition of $\bar\phi_n(\boldsymbol x)$, we have actually shown that
		\begin{equation}\label{2-26}
			||\phi_n||_*\le ||\phi_n||_{L^\infty(B_{Ls_{\varepsilon_n}}(\boldsymbol q_{\varepsilon_n}))}=o_n(1).
		\end{equation}
		
		On the other hand, for any $\tilde\varphi\in C_0^\infty(D_n)$ it holds
		\begin{equation*}
			\begin{split}
				&\left|\int_{D_n}\frac{1}{s_{\varepsilon_n}y_1+q_{\varepsilon_n}}\cdot\nabla\tilde\phi_n\cdot\nabla\tilde\varphi d\boldsymbol y\right|=\left|2\int_{|\boldsymbol y|=1}\frac{1}{q_{\varepsilon_n} }\tilde\phi_n\tilde\varphi+\langle\tilde f_n,\tilde\varphi\rangle \right|\\
				& \ \ \ \ \ \ \ \ \ =o_n(1)\cdot\|\tilde\varphi\|_{W^{1,1}(B_L(\boldsymbol 0))}+o_n(1)\cdot\|\tilde\varphi\|_{W^{1,p'}(B_L(\boldsymbol 0))}\\
				& \ \ \ \ \ \ \ \ \ =o_n(1)\cdot\left(\int_{B_L(\boldsymbol 0)}|\nabla\tilde\varphi|^{p'}\right)^{\frac{1}{p'}},
			\end{split}
		\end{equation*}
		which leads to
		\begin{equation}\label{2-27}
			\varepsilon_n^{1-\frac{2}{p}}\|\nabla\phi_n\|_{L^p(B_{Ls_{\varepsilon_n}(\boldsymbol q_{\varepsilon_n}}))}\le C||\nabla\tilde\phi_n||_{L^p(B_L(\boldsymbol 0))}=o_n(1).
		\end{equation}
		Combining \eqref{2-26} and \eqref{2-27}, we get a contradiction to \eqref{2-25}. Hence \eqref{2-19} holds, and \eqref{2-20} is a consequence of \eqref{2-19} and \eqref{2-24}.
	\end{proof}
	
	\bigskip
	
	Using this coercive estimate, we can obtain the following result on the solvability of projective problem \eqref{2-17}.
	\begin{lemma}\label{lem2-3}
		Suppose that $\mathrm{supp}\, \mathbf h\subset B_{2s_\varepsilon}(\boldsymbol q_\varepsilon)$ and $$\varepsilon^{1-\frac{2}{p}}\|\mathbf h\|_{W^{-1,p}(B_{Ls_\varepsilon}(\boldsymbol q_\varepsilon))}<\infty$$
		with $p\in(2,+\infty]$. Then there exists a small $\varepsilon_0>0$ such that for any $\varepsilon\in(0,\varepsilon_0]$, \eqref{2-17} has a unique solution $\phi_\varepsilon=\mathcal T_\varepsilon \, \mathbf h$, where $\mathcal T_\varepsilon$ is a linear operator. Moreover, there exists a constant $c_0>0$ independent of $\varepsilon$, such that
		\begin{equation}\label{2-28}
			\|\phi_\varepsilon\|_*+\varepsilon^{1-\frac{2}{p}}\|\nabla\phi_\varepsilon\|_{L^p(B_{Ls_\varepsilon}(\boldsymbol q_\varepsilon))}\le c_0\varepsilon^{1-\frac{2}{p}}\|\mathbf h\|_{W^{-1,p}(B_{Ls_\varepsilon}(\boldsymbol q_\varepsilon))},
		\end{equation}
	    where $L>0$ is a large constant.
	\end{lemma}
	\begin{proof}
		Let $H_a(\mathbb R^2_+)$ be the Hilbert space consists of functions $u$ satisfying the boundary condition
		\begin{equation*}
			\begin{cases}
				u=0, & \text{on}\  x_1=0,
				\\
				u, \ |\nabla u|/x_1\to0, &\text{as} \ |\boldsymbol x |\to \infty,
			\end{cases}
		\end{equation*}
		and endowed with the inner product
		\begin{equation*}
			[u,v]_{H_a(\mathbb R^2_+)}=\int_{\mathbb R^2_+} \frac{1}{x_1}\nabla u\cdot\nabla v d\boldsymbol x.
		\end{equation*}
		To yield the compactness of operator in $\mathbb R^2_+$, we need to introduce another weighted $L^\infty$ norm as
		\begin{equation*}
			||\phi||_{*,\nu}:=\sup_{\boldsymbol x\in\mathbb R^2}\rho_1(\boldsymbol x)^{1-\nu}\rho_2(\boldsymbol x)^{1-\nu}|\phi(\boldsymbol x)|,
		\end{equation*}
		where $0<\nu<1/4$ is a small number, and $\rho_1,\rho_2$ are defined in \eqref{rho2}. We introduce
		two spaces. The first one is
		\begin{equation*}
			E_\varepsilon:=\left\{u\in H_a(\mathbb R^2_+)\,\,\, \big|\, \,\, ||u||_{*,\nu}<\infty, \ u(x_1,x_2)=u(x_1,-x_2), \ \int_{\mathbb R^2_+}\frac{1}{x_1}\nabla u\cdot\nabla Z_{\boldsymbol q_\varepsilon,\varepsilon}d\boldsymbol x=0\right\}
		\end{equation*}
		with norm $||\cdot||_{*,\nu}$, and the second one is
		\begin{equation*}
			F_\varepsilon:=\left\{\mathbf h^* \in W^{-1,p}(B_{Ls_\varepsilon}(\boldsymbol q_\varepsilon)) \,\,\, \big| \,\,\, \ \mathbf h^*(x_1,x_2)=\mathbf h^*(x_1,-x_2), \ \int_{\mathbb R^2_+}Z_{\boldsymbol q_\varepsilon,\varepsilon}\mathbf h^*d\boldsymbol x=0\right\}
		\end{equation*}
		with $p\in (2,+\infty]$. Then for $\phi_\varepsilon\in E_\varepsilon$, problem \eqref{2-17} has an equivalent operation form
		\begin{equation*}
			\begin{split}
				 \phi_\varepsilon&=(-x_1{\Delta^*})^{-1}P_\varepsilon\left(\frac{1}{s_\varepsilon q_\varepsilon}\phi_\varepsilon(s_\varepsilon,\varepsilon)\boldsymbol{\delta}_{|\boldsymbol x-\boldsymbol q_\varepsilon|=s_\varepsilon}\right)+(-x_1{\Delta^*})^{-1} P_\varepsilon \mathbf h\\
				&=\mathscr K\phi_\varepsilon+(-x_1{\Delta^*})^{-1} P_\varepsilon \mathbf h,
			\end{split}
		\end{equation*}
		where
		\begin{equation*}
			(-x_1{\Delta^*})^{-1}u:=\int_{\mathbb R^2_+}{G_*}(\boldsymbol x,\boldsymbol x')x_1'^{-1}u(\boldsymbol x')d\boldsymbol x',
		\end{equation*}
		and $P_\varepsilon$ is the projection operator to $F_\varepsilon$. Since $Z_{\boldsymbol q_\varepsilon,\varepsilon}$ has a compact support due to the truncation \eqref{truncation}, by the definition of ${G_*}(\boldsymbol x,\boldsymbol x')$, we see that $\mathscr K$ maps $E_\varepsilon$ to $E_\varepsilon$.
		
		To show that $\mathscr K$ is a compact operator, we let $K_n:=\{\boldsymbol x\in\mathbb R^2\, | \,\, 1/n<x_1<n, \ |x_2|<n \}$ with $n\in N^*$. It is obvious that $K_n\to \mathbb R^2_+$ as $n\to +\infty$. Recall the asymptotic estimate for the Green's function ${G_*}$ given in \eqref{2-12} and \eqref{2-13}. For any small $\epsilon>0$, we can find an $N$ sufficiently large such that if $n>N$, then it holds
		$$\rho_1(\boldsymbol x)^{1-\nu}\rho_2(\boldsymbol x)^{1-\nu}|\mathscr Ku(\boldsymbol x)|<\epsilon, \ \ \ u\in E_\varepsilon, \ \ \ \boldsymbol x\in \mathbb R^2_+\setminus K_n.$$
		While for $\boldsymbol x\in K_n$, standard elliptic estimates shows that the $C^\alpha$ norm of $\mathscr Ku(\boldsymbol x)$ is bounded, and hence $\mathscr Ku(\boldsymbol x)$ is uniformly bounded and equi-continuous in $K_n$. By the Ascoli--Arzela theorem, we conclude that $\mathscr K$ is indeed a compact operator. (It is also noteworthy that this approach of recovering compactness is generally applicable in `gluing method', see the discussion in \cite{DW,DW2}.)
		
		Using the Fredholm alternative, \eqref{2-17} has a unique solution if the homogeneous equation
		\begin{equation*}
			\phi_\varepsilon=\mathscr K\phi_\varepsilon
		\end{equation*}
		has only trivial solution in $E_\varepsilon$, which can be obtained from Lemma \ref{lem2-2}. Now we let
		$$\mathcal T_\varepsilon:=(\text{Id}-\mathscr K)^{-1}(-x_1{\Delta^*})^{-1} P_\varepsilon,$$
		and the estimate \eqref{2-28} holds by Lemma \ref{lem2-2}. The proof is thus complete.
	\end{proof}

    \bigskip
	
	\subsection{The reduction and one-dimensional problem}
	
	Recall that our aim is to solve \eqref{Eqforperturbation}. However, since the linear operator $\mathbb L_\varepsilon$ has a nontrivial kernel, we have to settle for second best, and first deal with the projected problem in the space $E_\varepsilon$. Using the linear operator $\mathcal T_\varepsilon$ given in Lemma \ref{lem2-3}, we are to consider
	\begin{equation}\label{2-29}
		\phi_\varepsilon=\mathcal T_\varepsilon R_\varepsilon(\phi_\varepsilon)
	\end{equation}
	with
	\begin{equation*}
		 R_\varepsilon(\phi)=\frac{1}{\varepsilon^2}\bigg(x_1\boldsymbol1_{\{\psi_\varepsilon-\frac{W}{2}x_1^2\ln\frac{1}{\varepsilon}>\mu_\varepsilon\}}-x_1\boldsymbol1_{\{V_{\boldsymbol q_\varepsilon,\varepsilon}>\frac{a}{2\pi}\ln\frac{1}{\varepsilon}\}}-\frac{2}{s_\varepsilon q_\varepsilon}\phi(s_\varepsilon,\theta)\boldsymbol\delta_{|\boldsymbol x-\boldsymbol q_\varepsilon|=s_\varepsilon}\bigg)
	\end{equation*}
	for each small $\varepsilon\in (0,\varepsilon_0]$. In the following lemma, we will give a delicate estimate for the error term $R_\varepsilon(\phi_\varepsilon)$, so that a contraction mapping theorem can be applied to obtain the existence of $\phi_\varepsilon$ in $E_\varepsilon$.
	\begin{lemma}\label{lem2-4}
		There exists a small $\varepsilon_0>0$ such that for any $\varepsilon\in(0,\varepsilon_0]$, there is a unique solution $\phi_\varepsilon\in E_\varepsilon$ to \eqref{2-29}. Moreover  $\phi_\varepsilon$ satisfies
		\begin{equation}\label{2-30}
			\|\phi_\varepsilon\|_*+\varepsilon^{1-\frac{2}{p}}\|\nabla\phi_\varepsilon\|_{L^p(B_{Ls_\varepsilon}(\boldsymbol q_\varepsilon))}= O(\varepsilon|\ln\varepsilon|)
		\end{equation}
		with the norm $|| \cdot  ||_*$ defined in \eqref{2-11}, $p\in(2,+\infty]$.
	\end{lemma}
	\begin{proof}
		Denote $\mathcal G_\varepsilon:=\mathcal T_\varepsilon R_\varepsilon$, and a neighborhood of origin in $E_\varepsilon$ as
		\begin{equation*}
			\mathcal B_\varepsilon:=E_\varepsilon\cap \left\{\phi \ | \ \|\phi\|_*+\varepsilon^{1-\frac{2}{p}}\|\nabla\phi\|_{L^p(B_{Ls_\varepsilon}(\boldsymbol q_\varepsilon))}\le c^*\cdot\varepsilon|\ln\varepsilon|, \ p\in(2,\infty]\right\}
		\end{equation*}
		with $c^*$ a large positive constant. We will show that $\mathcal G_\varepsilon$ is a contraction map from $\mathcal B_\varepsilon$ to $\mathcal B_\varepsilon$, so that a unique fixed point $\phi_\varepsilon$ can be obtained by the contraction mapping theorem. Actually, letting $\mathbf h=R_\varepsilon(\phi)$ for $\phi\in\mathcal B_\varepsilon$, and noticing that $R_\varepsilon(\phi)$ satisfies assumptions for $\mathbf h$ in Lemma \ref{lem2-3} by Appendix \ref{appB}, we hence have
		\begin{equation*}
			\|\mathcal T_\varepsilon R_\varepsilon(\phi)\|_*+\varepsilon^{1-\frac{2}{p}}\|\nabla\mathcal T_\varepsilon R_\varepsilon(\phi)\|_{L^p(B_{Ls_\varepsilon}(\boldsymbol q_\varepsilon))}\le c_0\varepsilon^{1-\frac{2}{p}}\|R_\varepsilon(\phi) \|_{W^{-1,p}(B_{Ls_\varepsilon}(\boldsymbol q_\varepsilon))}.
		\end{equation*}
		
		To begin with, we are to show that $\mathcal G_\varepsilon$ maps $\mathcal B_\varepsilon$ continuously into itself. We use $\tilde v(\boldsymbol y)$ to denote $v(s_\varepsilon\boldsymbol y+\boldsymbol q_\varepsilon)$. For each $\varphi(\boldsymbol x)\in C_0^\infty(B_{Ls_\varepsilon}(\boldsymbol q_\varepsilon))$, in view of Lemma \ref{B2} and Lemma \ref{B3} in Appendix \ref{appB}, we have
		\begin{equation*}
			\begin{split}
				\langle R_\varepsilon(\phi),\varphi \rangle &=\frac{s_\varepsilon^2}{\varepsilon^2}\int_{B_L(\boldsymbol 0)}(s_\varepsilon y_1+q_\varepsilon)\left(\boldsymbol1_{\{\tilde\psi_\varepsilon-\frac{W}{2}(s_\varepsilon y_1+q_\varepsilon)^2\ln\frac{1}{\varepsilon}>\mu_\varepsilon\}}-\boldsymbol 1_{\{\tilde V_{\boldsymbol q_\varepsilon,\varepsilon}>\frac{a_\varepsilon}{2\pi}\ln\frac{1}{\varepsilon}\}}\right)\tilde\varphi d\boldsymbol y\\
				& \ \ \ -\frac{2}{q_\varepsilon}\int_0^{2\pi}\tilde\phi\tilde\varphi(1,\theta)d\theta\\
				&=(1+O(\varepsilon))\cdot q_\varepsilon \cdot\frac{s_\varepsilon^2}{\varepsilon^2}\int_0^{2\pi}\int_1^{1+t_\varepsilon(\theta)}t\tilde\varphi(t,\theta)dtd\theta-\frac{2}{q_\varepsilon}\int_0^{2\pi}\tilde\phi\tilde\varphi(1,\theta)d\theta\\
				&=\frac{s_\varepsilon^2}{\varepsilon^2}\cdot q_\varepsilon\int_0^{2\pi}\int_1^{1+t_\varepsilon(\theta)}t\tilde\varphi(1,\theta)dtd\theta-\frac{2}{q_\varepsilon}\int_0^{2\pi}\tilde\phi\tilde\varphi(1,\theta)d\theta\\
				& \ \ \ +\frac{s_\varepsilon^2}{\varepsilon^2}\cdot q_\varepsilon\int_0^{2\pi}\int_1^{1+t_\varepsilon(\theta)}t(\tilde\varphi(t,\theta)-\tilde\varphi(1,\theta))dtd\theta+O(\varepsilon)\cdot\int_0^{2\pi}|\tilde\varphi|d\theta\\
				&=\frac{s_\varepsilon^2}{\varepsilon^2}\cdot q_\varepsilon\int_0^{2\pi}\left(\frac{\tilde\phi(1,\theta)}{s_\varepsilon\mathcal N_\varepsilon}+O(\varepsilon|\ln\varepsilon|)\right)\tilde\varphi(1,\theta)d\theta+O(\varepsilon)\cdot\int_0^{2\pi}|\tilde\varphi|d\theta\\
				& \ \ \ +\frac{s_\varepsilon^2}{\varepsilon^2}\cdot q_\varepsilon\int_0^{2\pi}\int_1^{1+t_\varepsilon(\theta)}t\int_1^t\frac{\partial \tilde\varphi(r,\theta)}{\partial r}drdtd\theta-\frac{2}{q_\varepsilon}\int_0^{2\pi}\tilde\phi\tilde\varphi(1,\theta)d\theta\\
				&=\frac{s_\varepsilon^2}{\varepsilon^2}\cdot q_\varepsilon\int_0^{2\pi}|t_{\varepsilon,\mathcal W}+t_{\varepsilon,\tilde\phi}|\int_1^{1+t_\varepsilon(\theta)}\left|\frac{\partial \tilde\varphi(r,\theta)}{\partial r}\right|drd\theta+O(\varepsilon|\ln\varepsilon| )\cdot\|\tilde\varphi\|_{W^{1,p'}(B_L(\boldsymbol 0))}\\
				&=O(\varepsilon|\ln\varepsilon| )\cdot\|\tilde\varphi\|_{W^{1,p'}(B_L(\boldsymbol 0))},
			\end{split}
		\end{equation*}
		where we have used the definition of $\mathcal N_\varepsilon$ in \eqref{2-14}. Thus we have
		\begin{equation*}
			\varepsilon^{1-\frac{2}{p}}\|R_\varepsilon(\phi)\|_{W^{-1,p}(B_{Ls_\varepsilon}(\boldsymbol q_\varepsilon))}=O(\varepsilon|\ln\varepsilon|),
		\end{equation*}
		which yields
		\begin{equation*}
			\|\mathcal T_\varepsilon R_\varepsilon(\phi)\|_*+\varepsilon^{1-\frac{2}{p}}\|\nabla\mathcal T_\varepsilon R_\varepsilon(\phi)\|_{L^p(B_{Ls_\varepsilon}(\boldsymbol q_\varepsilon))}=O(\varepsilon|\ln\varepsilon|)
		\end{equation*}
		by Lemma \ref{lem2-3}. Arguing in a same way, we can deduce
		\begin{equation*}
			\varepsilon\|\nabla\phi\|_{L^\infty(B_{Ls_\varepsilon}(\boldsymbol q_\varepsilon))}=O(\varepsilon|\ln\varepsilon|)
		\end{equation*}
		from the estimate
		\begin{equation*}
			\varepsilon\|R_\varepsilon(\phi) \|_{W^{-1,\infty}(B_{Ls_\varepsilon}(\boldsymbol q_\varepsilon))}=O(\varepsilon|\ln\varepsilon|).
		\end{equation*}
		Thus operator $\mathcal G_\varepsilon$ indeed maps $\mathcal B_\varepsilon$ to $\mathcal B_\varepsilon$ continuously.
		
		In the next step, we are to verify that $\mathcal G_\varepsilon$ is a contraction mapping under the norm
		\begin{equation*}
			\|\cdot\|_{\mathcal G_\varepsilon}=\|\cdot\|_*+\varepsilon^{1-\frac{2}{p}}\|\cdot\|_{W^{1,p}(B_{Ls_\varepsilon}(\boldsymbol q_\varepsilon))}, \quad p\in(2,+\infty].
		\end{equation*}
		We already know that $\mathcal B_\varepsilon$ is close under this norm. Let $\phi_1$ and $\phi_2$ be two functions in $\mathcal B_\varepsilon$. From Lemma \ref{lem2-3}, it holds
		\begin{equation}\label{2-32}
			\|\mathcal G_\varepsilon\phi_1-\mathcal G_\varepsilon\phi_2\|_{\mathcal G_\varepsilon}\le C\varepsilon^{1-\frac{2}{p}}\|R_\varepsilon(\phi_1)-R_\varepsilon(\phi_2) \|_{W^{-1,p}(B_{Ls_\varepsilon}(\boldsymbol q_\varepsilon))},
		\end{equation}
		where
		\begin{equation*}
			\begin{split}
				R_\varepsilon(\phi_1)&-R_\varepsilon(\phi_2)\\
				&=\frac{1}{\varepsilon^2}\bigg(x_1\boldsymbol1_{\{\mathbf U_{\boldsymbol q_\varepsilon,\varepsilon}+\phi_1>0\}}-x_1\boldsymbol1_{\{\mathbf U_{\boldsymbol q_\varepsilon,\varepsilon}+\phi_2>0\}}-\frac{2}{s_\varepsilon q_\varepsilon}(\phi_1(r,\theta)-\phi_2(r,\theta))\boldsymbol\delta_{|\boldsymbol x-\boldsymbol q_\varepsilon|=s_\varepsilon}\bigg)
			\end{split}
		\end{equation*}
		with $\mathbf U_{\boldsymbol q_\varepsilon,\varepsilon}(\boldsymbol x)$ in \eqref{Udef} close to $V_{\boldsymbol q_\varepsilon,\varepsilon}(\boldsymbol x)-\frac{a_\varepsilon}{2\pi}\ln\frac{1}{\varepsilon}$. For $m=1,2$, let
		\begin{equation*}
			S_{m1}:=\{\boldsymbol y\,\,|\,\, \tilde{\mathbf U}_{\boldsymbol q_\varepsilon,\varepsilon}+\tilde\phi_m>0\}\cap B_{L}(\boldsymbol 0),
		\end{equation*}
		and
		\begin{equation*}
			S_{m2}:=\{\boldsymbol y\,\,|\,\, \tilde{\mathbf U}_{\boldsymbol q_\varepsilon,\varepsilon}+\tilde\phi_m<0\}\cap B_{L}(\boldsymbol 0).
		\end{equation*}
		Then it holds
		\begin{equation*}
			\boldsymbol1_{\{\tilde{\mathbf U}_{\boldsymbol q_\varepsilon,\varepsilon}+\tilde\phi_1>0\}}-\boldsymbol1_{\{\tilde{\mathbf U}_{\boldsymbol q_\varepsilon,\varepsilon}+\tilde\phi_2>0\}}=0, \ \ \  \text{in}\ \ (S_{11}\cap S_{21})\cup(S_{12}\cap S_{22}).
		\end{equation*}
	    Let
		\begin{equation*}
			\begin{split}
				\boldsymbol y_{\varepsilon,m}:&=(1+t_{\varepsilon,m}(\theta))(\cos\theta,\sin\theta)\\
				&=(1+t_{\varepsilon,\mathcal W}(\theta)+t_{\varepsilon,\tilde\phi_m}(\theta)+O(\varepsilon^2|\ln\varepsilon|))(\cos\theta,\sin\theta)\\
				&\in\{\boldsymbol y \,\,\,|\,\, \ \tilde{\mathbf U}_{\boldsymbol q_\varepsilon,\varepsilon}(\boldsymbol y_{\varepsilon,m})+\tilde\phi_m(\boldsymbol y_{\varepsilon,m})=\mu_\varepsilon\}\cap B_{2L}(\boldsymbol 0).
			\end{split}
		\end{equation*}
	    According to Lemma \ref{B3}, for each $\tilde\varphi(\boldsymbol y)\in C_0^\infty(B_L(\boldsymbol 0))$, we have
		\begin{equation*}
			\begin{split}
				&\quad\,\, \frac{s_\varepsilon^2}{\varepsilon^2}\int_{B_L(\boldsymbol 0)}(s_\varepsilon y_1+q_\varepsilon)\left(\boldsymbol1_{\{\tilde{\mathbf U}_{\boldsymbol q_\varepsilon,\varepsilon}+\tilde\phi_1>0\}\}}-\boldsymbol1_{\{\tilde{\mathbf U}_{\boldsymbol q_\varepsilon,\varepsilon}+\tilde\phi_2>0\}}\right)\tilde\varphi d\boldsymbol y\\
				&=\frac{s_\varepsilon^2}{\varepsilon^2}\left(\int_{S_{11}\cap S_{22}}(s_\varepsilon y_1+q_\varepsilon)\tilde\varphi d\boldsymbol y-\int_{S_{12}\cap S_{21}}(s_\varepsilon y_1+q_\varepsilon)\tilde\varphi d\boldsymbol y\right)\\
				 &=\frac{s_\varepsilon^2}{\varepsilon^2}\int_0^{2\pi}\int_{1+t_{\varepsilon,2}(\theta)}^{1+t_{\varepsilon,1}(\theta)}(s_\varepsilon y_1+q_\varepsilon)t\tilde\varphi dtd\theta\\
				 &=\frac{s_\varepsilon^2}{\varepsilon^2}\int_0^{2\pi}(t_{\varepsilon,1}-t_{\varepsilon,2})(s_\varepsilon y_1+q_\varepsilon)\tilde\varphi(1,\theta)d\theta\\
				 &\quad+\frac{s_\varepsilon^2}{\varepsilon^2}\int_0^{2\pi}\int_{1+t_{\varepsilon,2}(\theta)}^{1+t_{\varepsilon,1}(\theta)}(s_\varepsilon y_1+q_\varepsilon)r(\tilde\varphi(r,\theta)-\tilde\varphi(1,\theta))drd\theta\\
				 &=\frac{s_\varepsilon^2}{\varepsilon^2}\int_0^{2\pi}(t_{\varepsilon,\tilde\phi_1}-t_{\varepsilon,\tilde\phi_2})(s_\varepsilon y_1+q_\varepsilon)\tilde\varphi(1,\theta)d\theta+O\left((\varepsilon|\ln\varepsilon|^2) ^{\frac{1}{p}}\right)\cdot\sup_{\theta\in(0,2\pi]}|t_{\varepsilon,\tilde\phi_1}-t_{\varepsilon,\tilde\phi_2}|\cdot\|\tilde\varphi\|_{W^{1,p'}(B_L(\boldsymbol 0))}\\
				 &=\frac{s_\varepsilon^2}{\varepsilon^2}\int_0^{2\pi}(t_{\varepsilon,\tilde\phi_1}-t_{\varepsilon,\tilde\phi_2})(sy_1+q_\varepsilon)\tilde\varphi(1,\theta)d\theta+o_\varepsilon(1)\cdot\|\tilde\phi_1-\tilde\phi_2\|_{L^\infty(B_L(\boldsymbol 0))}\|\tilde\varphi\|_{W^{1,p'}(B_L(\boldsymbol 0))},
			\end{split}
		\end{equation*}
		where we have used the fact
		\begin{equation*}
			|t_{\varepsilon,\tilde\phi_1}-t_{\varepsilon,\tilde\phi_2}|\le C\|\tilde\phi_1-\tilde\phi_2\|_{L^\infty(B_L(\boldsymbol 0))}.
		\end{equation*}
		To handle the term
		$$\frac{s_\varepsilon^2}{\varepsilon^2}\int_0^{2\pi}(t_{\varepsilon,\tilde\phi_1}-t_{\varepsilon,\tilde\phi_2})(s_\varepsilon y_1+q_\varepsilon)\tilde\varphi(1,\theta)d\theta,$$
		we denote $\phi_*:=\tilde\phi_1-\tilde\phi_2$. Then it holds
		\begin{equation*}
			\begin{split}
				&\quad \tilde{\mathbf U}_{\boldsymbol q_\varepsilon,\varepsilon}(\boldsymbol y_{\varepsilon,1})-\tilde{\mathbf U}_{\boldsymbol q_\varepsilon,\varepsilon}(\boldsymbol y_{\varepsilon,2})=\tilde\phi_2(\boldsymbol y_{\varepsilon,2})-\tilde\phi_1(\boldsymbol y_{\varepsilon,1})\\
				&=\tilde\phi_2(\boldsymbol y_{\varepsilon,1})-\tilde\phi_1(\boldsymbol y_{\varepsilon,1})+\int_{1+t_{\varepsilon,1}(\theta)}^{1+t_{\varepsilon,2}(\theta)}\frac{\partial\tilde\phi_2(r,\theta)}{\partial r}dr\\
				&=\phi_*(1,\theta)+\int_1^{1+t_{\varepsilon,2}(\theta)}\frac{\partial\tilde\phi_*(r,\theta)}{\partial r}dr+\int_{1+t_{\varepsilon,1}(\theta)}^{1+t_{\varepsilon,2}(\theta)}\frac{\partial\tilde\phi_2(r,\theta)}{\partial r}dr.
			\end{split}
		\end{equation*}
		By the expansion
		\begin{equation*}
			\tilde{\mathbf U}_{\boldsymbol q_\varepsilon,\varepsilon}(\boldsymbol y_{\varepsilon,1})-\tilde{\mathbf U}_{\boldsymbol q_\varepsilon,\varepsilon}(\boldsymbol y_{\varepsilon,2})=-\frac{1}{s_\varepsilon\mathcal N_\varepsilon}(\boldsymbol y_{\varepsilon,1}-\boldsymbol y_{\varepsilon,2})+O(\varepsilon|\ln\varepsilon|^2),
		\end{equation*}
		we have
		\begin{equation*}
			\begin{split}
				&\quad t_{\varepsilon,\tilde\phi_1}-t_{\varepsilon,\tilde\phi_2}=|\boldsymbol y_{\varepsilon,1}-\boldsymbol y_{\varepsilon,2}|\\
				&=-s_\varepsilon\mathcal N_\varepsilon(1+o_\varepsilon(1))\cdot\left(\phi_*(1,\theta)+\int_1^{1+t_{\varepsilon,1}(\theta)}\frac{\partial\tilde\phi_*(r,\theta)}{\partial r}dr+\int_{1+t_{\varepsilon,1}(\theta)}^{1+t_{\varepsilon,2}(\theta)}\frac{\partial\tilde\phi_2(r,\theta)}{\partial r}dr\right).
			\end{split}
		\end{equation*}
		Then using the definition of $\mathcal N_\varepsilon$ in \eqref{2-14}, one can deduce
		\begin{equation*}
			\begin{split}
				 &\quad\frac{s_\varepsilon^2}{\varepsilon^2}\int_0^{2\pi}(t_{\varepsilon,\tilde\phi_1}-t_{\varepsilon,\tilde\phi_2})(s_\varepsilon y_1+q_\varepsilon)\tilde\varphi(1,\theta)d\theta=\frac{2}{q_\varepsilon}(1+o_\varepsilon(1))\cdot\int_0^{2\pi}(\tilde\phi_1-\tilde\phi_2)\tilde\varphi(1,\theta)d\theta\\
				& \ \ \ -\frac{2}{q_\varepsilon}(1+o_\varepsilon(1))\cdot\left(\int_1^{1+t_{\varepsilon,1}(\theta)}\frac{\partial\tilde\phi_*(r,\theta)}{\partial r}dr+\int_{1+t_{\varepsilon,2}(\theta)}^{1+t_{\varepsilon,2}(\theta)}\frac{\partial\tilde\phi_2(r,\theta)}{\partial r}dr\right)\\
				 &=\frac{2}{q_\varepsilon}\int_0^{2\pi}(\tilde\phi_1-\tilde\phi_2)\tilde\varphi(1,\theta)d\theta+o_\varepsilon(1)\cdot\|\tilde\phi_1-\tilde\phi_2\|_{L^\infty(B_L(\boldsymbol 0))}\\
				& \ \ \ +\left(O\left((\varepsilon|\ln\varepsilon|^2)^{\frac{1}{p}}\right)+
				\|\tilde\phi_2\|_{W^{1,p}(B_L(\boldsymbol 0))}\right)\cdot\|\tilde\phi_1
				-\tilde\phi_2\|_{L^\infty(B_L(\boldsymbol 0))}\cdot\|\tilde\varphi\|_{W^{1,p'}(B_L(\boldsymbol 0))}.
			\end{split}
		\end{equation*}
		Finally, we conclude that
		\begin{equation*}
			\varepsilon^{1-\frac{2}{p}}\|R_\varepsilon(\phi_1)-R_\varepsilon(\phi_2) \|_{W^{-1,p}(B_{L}(\boldsymbol z))}=o_\varepsilon(1)\cdot \|\phi_1-\phi_2\|_{\mathcal G_\varepsilon},
		\end{equation*}
		which yields
		\begin{equation*}
			\|\mathcal G_\varepsilon\phi_1-\mathcal G_\varepsilon\phi_2\|_{\mathcal G_\varepsilon}=o_\varepsilon(1)\cdot \|\phi_1-\phi_2\|_{\mathcal G_\varepsilon}
		\end{equation*}
		from \eqref{2-32}. Hence we have shown that $\mathcal G_\varepsilon$ is a contraction map from $\mathcal B_\varepsilon$ into itself.
		
		By applying the contraction mapping theorem, we now can claim that there is a unique $\phi_\varepsilon\in \mathcal B_\varepsilon$ such that $\phi_\varepsilon=\mathcal G_\varepsilon\phi_\varepsilon$, which satisfies \eqref{2-30}. Since $\|\phi_\varepsilon\|_{\mathcal G_\varepsilon}$ is bounded by a constant $C$ independent of $\boldsymbol q_\varepsilon$, we conclude that $\phi_\varepsilon$ is continuous with respect to $\boldsymbol q_\varepsilon$ in the norm $\|\cdot\|_{\mathcal G_\varepsilon}$.
	\end{proof}
	
	\bigskip
	
	From the above lemma, the problem of solving \eqref{Eqforperturbation} is now transformed into a one-dimensional problem: Finding the sufficient condition to ensure $\Lambda=0,$ which will also determine the location of $\boldsymbol q_\varepsilon=(q_\varepsilon,0)$ as a crucial parameter in approximate solutions. In the next lemma, we will derive a condition equivalent to $\Lambda=0$, which enables us to prove the existence of $\psi_\varepsilon$ as the solution of \eqref{2-9} since the components $\mathcal V_{\boldsymbol q_\varepsilon,\varepsilon}$, $\mathcal H_{\boldsymbol q_\varepsilon,\varepsilon}$, $\phi_\varepsilon$ are all determined.
	\begin{lemma}\label{lem2-5}
		If $\boldsymbol q_\varepsilon=(q_\varepsilon, 0)$ satisfies
		\begin{equation}\label{2-35}
			\varepsilon^2\int_{\mathbb R^2_+}\frac{1}{x_1}\nabla\psi_\varepsilon\cdot\nabla Z_{\boldsymbol q_\varepsilon,\varepsilon}d\boldsymbol x
			-\int_{A_\varepsilon}x_1\cdot Z_{\boldsymbol q_\varepsilon,\varepsilon}d\boldsymbol x=0,
		\end{equation}
		then $\psi_\varepsilon$ is a solution to \eqref{2-9} and \eqref{2-10}.
	\end{lemma}
	\begin{proof}
		If the assumption \eqref{2-35} holds true, from \eqref{2-29} we will have
		\begin{equation*}   	
			\varepsilon^2\Lambda\int_{\mathbb R^2_+}\frac{1}{x_1}\nabla Z_{\boldsymbol q_\varepsilon,\varepsilon}\cdot\nabla Z_{\boldsymbol q_\varepsilon,\varepsilon}d\boldsymbol x=0.
		\end{equation*}
		Proceeding as in the proof of Lemma \ref{lem2-2}, we deduce
		\begin{equation*}
			\varepsilon^2\int_{\mathbb R^2_+}\frac{1}{x_1}\nabla Z_{\boldsymbol q_\varepsilon,\varepsilon}\cdot\nabla Z_{\boldsymbol q_\varepsilon,\varepsilon}d\boldsymbol x=C_Z+o_\varepsilon(1)
		\end{equation*}
		with $C_Z$ the positive constant given by the left hand side of \eqref{2-21}. Hence it holds $\Lambda=0$ when $\varepsilon$ is sufficiently small. This fact implies that $\psi_\varepsilon$ is a solution to \eqref{2-9} and \eqref{2-10}.
	\end{proof}
	
	\bigskip
	
	Taking advantage of the above characterization, we are now in the position to prove Proposition \ref{prop2-1}.
	
	{\bf Proof of Proposition \ref{prop2-1}:}
	We will show that we can find $\boldsymbol q_\varepsilon=(q_\varepsilon, 0)$ such that \eqref{2-35} holds and $q_\varepsilon$ satisfies
	\begin{equation}\label{2-36}
		q_\varepsilon=\frac{\kappa}{4\pi W}+O\left(\frac{1}{|\ln\varepsilon|}\right).
	\end{equation}
	Since $\phi_\varepsilon\in E_\varepsilon$, we have
	\begin{equation*}
		\int_{\mathbb R^2_+}\frac{1}{x_1}\nabla\phi_\varepsilon\cdot\nabla Z_{\boldsymbol q_\varepsilon,\varepsilon}d\boldsymbol x=0.
	\end{equation*}
	Hence it holds
	\begin{equation*}
		\begin{split}
			&\quad\varepsilon^2\int_{\mathbb R^2_+}\frac{1}{x_1}\nabla\psi_\varepsilon\cdot\nabla Z_{\boldsymbol q_\varepsilon,\varepsilon}d\boldsymbol x
			-\int_{A_\varepsilon}x_1\cdot Z_{\boldsymbol q_\varepsilon,\varepsilon}d\boldsymbol x\\
			&=\varepsilon^2\int_{\mathbb R^2_+}\frac{1}{x_1}\nabla(\mathcal V_{\boldsymbol q_\varepsilon,\varepsilon}+\mathcal H_{\boldsymbol q_\varepsilon,\varepsilon})\cdot\nabla Z_{\boldsymbol q_\varepsilon,\varepsilon}d\boldsymbol x
			-\int_{A_\varepsilon}x_1\cdot Z_{\boldsymbol q_\varepsilon,\varepsilon}d\boldsymbol x\\
			&=\int_{B_{Ls_\varepsilon}(\boldsymbol q_\varepsilon)}x_1(\boldsymbol 1_{\{V_{\boldsymbol q_\varepsilon,\varepsilon}>\frac{a_\varepsilon}{2\pi}\ln\frac{1}{\varepsilon}\}}-\boldsymbol 1_{\{\psi_\varepsilon-\frac{W}{2}x_1^2\ln\frac{1}{\varepsilon}>\mu_\varepsilon\}})\cdot Z_{\boldsymbol q_\varepsilon,\varepsilon}d\boldsymbol x.
		\end{split}
	\end{equation*}
	By denoting
	$$\tilde Z_{\boldsymbol q_\varepsilon,\varepsilon}=Z_{\boldsymbol q_\varepsilon,\varepsilon}(s_\varepsilon\boldsymbol y+\boldsymbol q_\varepsilon),$$
	direct computation yields
	$$\|\tilde Z_{\boldsymbol q_\varepsilon,\varepsilon}\|_{W^{1,p'}(B_L(\boldsymbol 0))}=O(\varepsilon^{-1}).$$
	Note that
	\begin{equation*}
		\begin{split}
			&\quad\frac{2}{q_\varepsilon}\int_0^{2\pi}\tilde\phi_\varepsilon(1,\theta)\tilde Z_{\boldsymbol q_\varepsilon,\varepsilon}d\theta\\
			&=\int\frac{1}{s_\varepsilon y_1+q_\varepsilon}\cdot\nabla\tilde\phi_\varepsilon\cdot\nabla \tilde Z_{\boldsymbol q_\varepsilon,\varepsilon}d\boldsymbol y+O_\varepsilon(1)\cdot\left(\|\phi_\varepsilon\|_*+\varepsilon\|\nabla\phi_\varepsilon\|_{L^\infty(B_{Ls_\varepsilon}(\boldsymbol q_\varepsilon))}\right)\\
			&=O_\varepsilon(1)\cdot\left(\|\phi_\varepsilon\|_*+\varepsilon\|\nabla\phi_\varepsilon\|_{L^\infty(B_{Ls_\varepsilon}(\boldsymbol q_\varepsilon))}\right),
		\end{split}
	\end{equation*}
	due to the nondegeneracy property of operator $\mathbb L^*_\varepsilon$ defined in \eqref{linear}. Then, similar to the proof of Lemma \ref{lem2-5}, we can deduce
	\begin{equation*}
		\begin{split}
			&\quad\frac{1}{\varepsilon^2}\int_{B_{Ls_\varepsilon}(\boldsymbol q_\varepsilon)}x_1(\boldsymbol 1_{\{V_{\boldsymbol q_\varepsilon,\varepsilon}>\frac{a_\varepsilon}{2\pi}\ln\frac{1}{\varepsilon}\}}-\boldsymbol 1_{\{\psi_\varepsilon-\frac{W}{2}x_1^2\ln\frac{1}{\varepsilon}>\mu_\varepsilon\}})\cdot Z_{\boldsymbol q_\varepsilon,\varepsilon}d\boldsymbol x\\
			&=-\frac{s_\varepsilon^2}{\varepsilon^2}\cdot q_\varepsilon\int_0^{2\pi}\left(\frac{\tilde\phi(1,\theta)}{s_\varepsilon\mathcal N_\varepsilon}+s_\varepsilon\cos\theta\cdot \left(\frac{s_\varepsilon^2}{4\varepsilon^2}\cdot q_\varepsilon\ln\frac{1}{\varepsilon}-Wq_\varepsilon\ln\frac{1}{\varepsilon}\right)\right)\tilde Z_{\boldsymbol q_\varepsilon,\varepsilon}d\theta+O_\varepsilon(1)\\
			&=\frac{\pi}{2}\cdot\frac{s_\varepsilon^4}{\varepsilon^4}\cdot q_\varepsilon^3 \left(\frac{s_\varepsilon^2}{4\varepsilon^2}\cdot q_\varepsilon\ln\frac{1}{\varepsilon}-Wq_\varepsilon\ln\frac{1}{\varepsilon}\right)+O_\varepsilon(1).
		\end{split}
	\end{equation*}
	Since it holds $s_\varepsilon^2\pi q_\varepsilon/\varepsilon^2=\kappa+O(1/|\ln\varepsilon|)$ by our choice of $a_\varepsilon$ in \eqref{2-16}, condition \eqref{2-35} yields
	$$\frac{\kappa}{4\pi}-Wq_\varepsilon=O\left(\frac{1}{|\ln\varepsilon|}\right).$$
	Then we can solve above equation on $q_\varepsilon$ and obtain at least one $q_\varepsilon$ satisfying \eqref{2-36}. In view of Lemma \ref{lem2-5}, we obtain the existence of $\psi_\varepsilon$ for every $\varepsilon\in(0,\varepsilon_0]$. Moreover, the estimates for $A_\varepsilon$ can be deduced from Lemma \ref{lem2-4} and Appendix \ref{appB}. Thus we have completed the proof of Proposition \ref{prop2-1}.
	\qed
	
	\bigskip
	
	Notice that \eqref{2-36} is the necessary condition for the linear term in Lemma \ref{B1} being nearly vanishing, which is consistent with our observation (a) for \eqref{exp} in the introduction. A more precise estimate for the location $\boldsymbol q_\varepsilon=(q_\varepsilon,0)$ will be given in the next section by a prior estimate on coefficients of the expansion.

	\bigskip
	
	\bigskip
	
\section{Uniqueness}\label{sec3}
	
	In this section, we will prove the local uniqueness of a vortex ring of small cross-section for which $\zeta$ is constant throughout the core. It follows from \eqref{1-1} and the definition of Stokes stream function (see also \cite{DV}) that the quantity
	$$P+\frac{1}{2}|\mathbf v|^2-\frac{1}{\varepsilon^2}\left(\psi_\varepsilon-\frac{W}{2}x_1^2\ln\frac{1}{\varepsilon}\right)$$
	is a constant in the the cross-section $A_\varepsilon$, while the value of $P+\frac{1}{2}|\mathbf v|^2$ is invariant in $\mathbb R_+^2\setminus A_\varepsilon$. Moreover, we see that $A_\varepsilon$ has a positive distance from $x_2$-axis and simply-connected in view of assumptions of Theorem \ref{thm2} on $\Omega_\varepsilon$.  By the continuity of $P+\frac{1}{2}|\mathbf v|^2$, $A_\varepsilon$ is given by
	$$A_\varepsilon=\left\{\boldsymbol x\in \mathbb R^2_+ \,\,\, \big| \,\, \psi_\varepsilon-\frac{W}{2}x_1^2\ln\frac{1}{\varepsilon}>\mu_\varepsilon\right\},$$
	where $\mu_\varepsilon>0$ has a positive lower bound independent of $\varepsilon$. Using notations in Section 2, the Stokes stream function $\psi_\varepsilon$ should satisfy
	\begin{equation}\label{3-1}
		\begin{cases}
			-\varepsilon^2{\Delta^*}\psi_\varepsilon=\boldsymbol1_{A_\varepsilon}, & \text{in} \ \mathbb R^2_+,
			\\
			\psi_\varepsilon=0, & \text{on} \ x_1=0,
			\\
			\psi_\varepsilon, \ |\nabla\psi_\varepsilon|/x_1\to0, &\text{as} \ |\boldsymbol x |\to \infty.
		\end{cases}
	\end{equation}
	To discuss the uniqueness of vortex rings of small cross-section, we will fix the circulation
	\begin{equation}\label{3-2}
		\kappa=\frac{1}{\varepsilon^2}\int_{A_\varepsilon}x_1d\boldsymbol x,
	\end{equation}
	and the parameter $W$ in translational velocity $W\ln \varepsilon\, \mathbf e_z$. Since $\psi_\varepsilon$ determines the vortex ring $\zeta_\varepsilon$ absolutely, the uniqueness result in Theorem \ref{thm2} can be concluded from following proposition.
	\begin{proposition}\label{prop3-1}
		Let $\kappa$ and $W$ be two fixed positive constants. Suppose that the cross-section $A_\varepsilon$ is simply-connected with a positive distance from $x_2$-axis, and satisfies
		$$\mathrm{diam}\, A_\varepsilon< L\varepsilon, \quad \mathrm{as} \ \varepsilon\to 0$$
		for some $L>0$. Then for each $\varepsilon\in (0,\varepsilon_0]$ with $\varepsilon_0>0$ sufficiently small, equation \eqref{3-1} together with \eqref{3-2} has a unique solution $\psi_\varepsilon$ up to translations in the $x_2$-direction.
	\end{proposition}
	
	\bigskip
	
	To study the local behavior of $\psi_\varepsilon$ near $A_\varepsilon$, we denote
	$$\sigma_\varepsilon:=\frac{1}{2}\text{diam}\, A_\varepsilon$$
	as the cross-section parameter. By our assumptions, it will hold $\sigma_\varepsilon<L\varepsilon$ as $\varepsilon\to 0$. Intuitively, the maximum point of $\psi_\varepsilon$ in $A_\varepsilon$ gives the exact location of cross-section. So we can choose a point $\boldsymbol p_\varepsilon\in \bar A_\varepsilon$ satisfying
	$$\psi_\varepsilon(\boldsymbol p_\varepsilon)=\max_{\boldsymbol x\in \bar A_\varepsilon}\psi_\varepsilon(\boldsymbol x),$$
	which is achievable by maximum principle of $-{\Delta^*}$. In view of Lemma \ref{A1} in Appendix \ref{appA}, the set $A_\varepsilon$ must be symmetric with respect to some horizontal line $x_2=h$. Using the translation invariance of \eqref{3-1} in $x_2$-direction, we may always assume that $A_\varepsilon$ is even symmetric with respect to $x_1$-axis (i.e. $(x_1,x_2)\in A_\varepsilon$ if and only if
    $(x_1,-x_2)\in A_\varepsilon$). Then, by the integral equation
	\begin{equation*}
		\psi_\varepsilon(\boldsymbol x)=\frac{1}{\varepsilon^2}\int_{\mathbb R^2_+}{G_*}(\boldsymbol x,\boldsymbol x') \boldsymbol 1_{A_\varepsilon}(\boldsymbol x')d\boldsymbol x',
	\end{equation*}
	we see that $\psi_\varepsilon$ attains its maximum on $x_1$-axis, and
	$$\psi_\varepsilon(\boldsymbol x)-\frac{W}{2}\ln\frac{1}{\varepsilon}x_1^2<0, \quad \text{as}\ \ x_1\to +\infty. $$
	Thus we may assume that $\boldsymbol p_\varepsilon=(p_\varepsilon,0)$, where $p_\varepsilon\in (c_1,c_2)$ with $c_1<c_2$ two positive constants.
	
	Intuitively, as the cross-section $A_\varepsilon$ shrinks, $\psi_\varepsilon(\boldsymbol x)$ will tend to Green's function $G_*(\boldsymbol x,\boldsymbol q_0)$ with $\boldsymbol q_0$ the limit position of vortex ring. Recall that in the previous section of existence we let 
	$$\boldsymbol q_\varepsilon=(q_\varepsilon,0), \quad z_1>0$$
	be the center of $V_{\boldsymbol q_\varepsilon,\varepsilon}(\boldsymbol x)$ in the approximate stream function $\Psi_\varepsilon(\boldsymbol x)=V_{\boldsymbol q_\varepsilon,\varepsilon}-V_{\bar{\boldsymbol q_\varepsilon},\varepsilon}+\mathcal H_{\boldsymbol q_\varepsilon,\varepsilon}$, where $V_{\boldsymbol q_\varepsilon,\varepsilon}-V_{\bar{\boldsymbol q_\varepsilon},\varepsilon}$ corresponds to the singular $\log$-part and its mirror image (to make the boundary value zero) in the Green's function, and $\mathcal H_{\boldsymbol q_\varepsilon,\varepsilon}$ corresponds to the remaining regular part. To give a precise approximation for $\psi_\varepsilon$, we can use a same strategy and decompose the Green's function for $-{\Delta^*}$ in boundary condition of \eqref{3-1} as
	\begin{equation*}
		{G_*}(\boldsymbol x,\boldsymbol x')=q_\varepsilon^2G(\boldsymbol x,\boldsymbol x')+H(\boldsymbol x,\boldsymbol x'),
	\end{equation*}
	where $q_\varepsilon$ is the first coordinate of $\boldsymbol q_\varepsilon$ in $V_{\boldsymbol q_\varepsilon,\varepsilon}(\boldsymbol x)$ for the purpose of constructing approximation and estimating the error in Appendix \ref{appB}, $G(\boldsymbol x,\boldsymbol x')$ is the Green's function of $-\Delta$ on the half plane, and $H(\boldsymbol x,\boldsymbol x')$ is the rest regular part (we also omit the parameter $q_\varepsilon$ in $H(\boldsymbol x,\boldsymbol x')$). According to the choice for the center of $V_{\boldsymbol q_\varepsilon,\varepsilon}$, the location of $\boldsymbol q_\varepsilon$ is near $\boldsymbol p_\varepsilon$ so that $|p_\varepsilon-q_\varepsilon|=O(\varepsilon)$. An accurate limiting behavior (with respect to $\varepsilon$) for $q_\varepsilon$ will be given in the second part of this section by a more delicate approximation procedure.
	
	Applying this decomposition of ${G_*}(\boldsymbol x,\boldsymbol x')$, we can split the stream function $\psi_\varepsilon$ into two parts $\psi_{1,\varepsilon}+\psi_{2,\varepsilon}$, where
	\begin{equation*}
		\psi_{1,\varepsilon}(\boldsymbol x)=\frac{q_\varepsilon^2}{\varepsilon^2}\int_{\mathbb R^2_+}G(\boldsymbol x,\boldsymbol x')\boldsymbol{1}_{A_\varepsilon}(\boldsymbol x')d\boldsymbol x',
	\end{equation*}
	and
	\begin{equation*}
		\psi_{2,\varepsilon}(\boldsymbol x)=\frac{1}{\varepsilon^2}\int_{\mathbb R^2_+}H(\boldsymbol x,\boldsymbol x')\boldsymbol{1}_{A_\varepsilon}(\boldsymbol x')d\boldsymbol x'.
	\end{equation*}
	According to \eqref{3-1}, $\psi_{1,\varepsilon}(\boldsymbol x)$ solves the problem
	\begin{equation*}
		\begin{cases}
			-\varepsilon^2\Delta \psi_{1,\varepsilon}(\boldsymbol x)=q_\varepsilon^2\boldsymbol{1}_{A_\varepsilon},  \ \ \ &\text{in} \ \mathbb R^2_+,\\
			\psi_{1,\varepsilon}=0, &\text{on} \ x_1=0,\\
			\psi_{1,\varepsilon}, \ |\nabla\psi_{1,\varepsilon}|/x_1\to0, &\text{as} \ |\boldsymbol x|\to \infty,
		\end{cases}
	\end{equation*}
	and $\psi_{2,\varepsilon}(\boldsymbol x)$ satisfies
	\begin{equation*}
		\begin{cases}
			-\varepsilon^2{\Delta^*}(\psi_{1,\varepsilon}(\boldsymbol x)+\psi_{2,\varepsilon}(\boldsymbol x))=\boldsymbol{1}_{A_\varepsilon}, \ \ \ &\text{in} \ \mathbb R^2_+,\\
			\psi_{2,\varepsilon}=0, &\text{on} \ x_1=0,\\
			\psi_{2,\varepsilon}, \ |\nabla\psi_{2,\varepsilon}|/x_1\to0, &\text{as} \ |\boldsymbol x|\to \infty,
		\end{cases}
	\end{equation*}
	Since the cross-section 
	$$A_\varepsilon=\left\{\boldsymbol x\in \mathbb R^2_+ \,\,\, \big| \,\, \psi_{1,\varepsilon}(\boldsymbol x)+\psi_{2,\varepsilon}(\boldsymbol x)-\frac{W}{2}x_1^2\ln\frac{1}{\varepsilon}>\mu_\varepsilon\right\},$$
	is dependent both on $\psi_{1,\varepsilon}$ and $\psi_{2,\varepsilon}$, the above two equations constitute a coupled system of $\psi_{1,\varepsilon}$ and $\psi_{2,\varepsilon}$, which seems more complicated than \eqref{3-1}. However, it should be noted that $\psi_{1,\varepsilon}$ is a solution to a semi-linear Laplace equation. While $\psi_{2,\varepsilon}$ is a more regular function than $\psi_{1,\varepsilon}$ with the $L^\infty$ norm bounded independent of $\varepsilon$ due to the regularity of $H$. These fine properties enable us to decouple $\psi_{1,\varepsilon}$ and $\psi_{2,\varepsilon}$ in the above two equations, and use the local Pohozaev identity in Appendix \ref{appC} to establish a relationship between the cross-section $A_\varepsilon$ and the coefficient $s_\varepsilon\mathcal W_\varepsilon(s_\varepsilon)$ of linear term we calculated in Appendix \ref{appB}.
	
	To prove the uniqueness, the key idea is to obtain an estimate for the expansion of $\psi_\varepsilon$ as precise as possible, which are to be obtained by several steps of approximation and bootstrap procedure. In this process, we can also obtain a fine estimate for the circulation $\kappa$, traveling speed $W|\ln\varepsilon|$, cross-section parameter $\sigma_\varepsilon$ and the vortex ring center location $\boldsymbol q_\varepsilon=(q_\varepsilon,0)$, namely, an accurate version of Kelvin--Hicks formula \eqref{KH}.
	\begin{proposition}\label{prop3-2}
		For steady vortex rings of small cross-section depicted in Proposition \ref{prop3-1}, the parameters $\kappa$, $W$, $\sigma_\varepsilon$, and $q_\varepsilon$ satisfy
		\begin{equation*}
			 Wq_\varepsilon\ln\frac{1}{\varepsilon}=\frac{\kappa}{4\pi}\left(\ln\frac{8q_\varepsilon}{\sigma_\varepsilon}-\frac{1}{4}\right)+O(\varepsilon^2|\ln\varepsilon|), \quad  \mathrm{as} \ \varepsilon \to 0.
		\end{equation*}
	\end{proposition}
	
	In \cite{Fra1}, Fraenkel obtained a slightly weaker form of the above estimate with the error term $O(\varepsilon^2|\ln\varepsilon|^2)$.  We reach a level of $O(\varepsilon^2|\ln\varepsilon|)$ since the center $\boldsymbol q_\varepsilon$ of $V_{\boldsymbol q_\varepsilon,\varepsilon}$ in the approximate solution is chosen to be the location of vortex ring center. Actually, if we use maximum point $\boldsymbol p_\varepsilon$ of $\psi_\varepsilon$ rather than $\boldsymbol q_\varepsilon$ in above formula, then the error term will be $O(\varepsilon^2|\ln\varepsilon|^2)$  the same as in \cite{Fra1}.
	
	Our approach for showing the uniqueness is divided into several parts. In the first part of our proof, we give a coarse estimate for $\psi_\varepsilon$ and $A_\varepsilon$ by a blow up argument and deriving the limiting function. Then we will give a series of approximate solutions, and improve the estimates by carefully dealing with the error term compared with $\psi_\varepsilon$. These two steps can be regarded as an inverse of procedure we have done in Section 2: not to construct or joint together the building blocks this time, but to expand $\psi_\varepsilon$ and calculate the coefficients. In the last step, the uniqueness for $\psi_\varepsilon$ is obtained by the non-vanishing property of the second order coefficient $\mathbf c_{2,\varepsilon}$ in the expansion \eqref{exp} for $\psi_\varepsilon$ at vortex location $q_\varepsilon$, which is achieved by using a Pohozaev identity and contradiction.
	
	\bigskip
	
	\subsection{Asymptotic estimates for vortex ring}
	
	The purpose of this part is to obtain a first asymptotic estimate for $\psi_\varepsilon$ and cross section $A_\varepsilon$, and to obtain the following necessary condition on dynamics of the vortex ring, which is a coarse version of Kelvin--Hicks formula stated in Proposition \ref{prop3-2}.
	\begin{proposition}\label{prop3-3}
		As $\varepsilon\to 0$, it holds
		\begin{equation*}
			 Wq_\varepsilon\ln\frac{1}{\varepsilon}-\frac{\kappa}{4\pi}\ln\frac{8q_\varepsilon}{\sigma_\varepsilon}+\frac{\kappa}{16\pi}=o_\varepsilon(1),
		\end{equation*}
		where $q_\varepsilon$ is the first coordinate of $\boldsymbol q_\varepsilon=(q_\varepsilon,0)$ as the center for approximate solutions in Appendix \ref{appB}.
	\end{proposition}
	
	 Although the behavior of $\psi_{1,\varepsilon}$ near cross-section $A_\varepsilon$ is not well-understood at the first stage, it is not difficult to derive an estimate for $\psi_{1,\varepsilon}$ away from $A_\varepsilon$ by the integral equation and diameter restriction $\mathrm{diam}\, A_\varepsilon< L\varepsilon$ as $\varepsilon\to 0$. 
	\begin{lemma}\label{lem3-4}
		For every $\boldsymbol x\in \mathbb R^2_+\setminus\{\boldsymbol x\ | \ \mathrm{dist}(\boldsymbol x,A_\varepsilon)\le L\varepsilon\}$, we have
		\begin{equation*}
			\psi_{1,\varepsilon}(\boldsymbol x)=\frac{\kappa}{2\pi}\cdot p_\varepsilon\ln \frac{|\boldsymbol x-\boldsymbol {\bar p_\varepsilon}|}{|\boldsymbol x-\boldsymbol p_\varepsilon|}+O\left(\frac{ \varepsilon}{|\boldsymbol x-\boldsymbol p_\varepsilon|}\right),
		\end{equation*}
		and
		\begin{equation*}
			\nabla \psi_{1,\varepsilon}(\boldsymbol x)=-\frac{\kappa}{2\pi}\cdot p_\varepsilon\frac{\boldsymbol x-\boldsymbol p_\varepsilon}{|\boldsymbol x-\boldsymbol p_\varepsilon|^2}+\frac{\kappa}{2\pi}\cdot p_\varepsilon\frac{\boldsymbol x-\boldsymbol {\bar p_\varepsilon}}{|\boldsymbol x-\boldsymbol {\bar p_\varepsilon}|^2}+O\left(\frac{ \varepsilon}{|\boldsymbol x-\boldsymbol p_\varepsilon|^2}\right).
		\end{equation*}
	\end{lemma}
	\begin{proof}
		For every $\boldsymbol x\in \mathbb R^2_+\setminus\{\boldsymbol x \ | \ \mathrm{dist}(\boldsymbol x,A_\varepsilon)\le L\varepsilon\}$, it holds $\boldsymbol x \notin A_\varepsilon$. Recall the notation $\bar{\boldsymbol x}=(-x_1,x_2)$. For each $\boldsymbol x'\in A_\varepsilon$ we have
		\begin{equation*}
			|\boldsymbol x-\boldsymbol x'|=|\boldsymbol x-\boldsymbol p_\varepsilon|-\langle\frac{\boldsymbol x-\boldsymbol p_\varepsilon}{|\boldsymbol x-\boldsymbol p_\varepsilon|},\boldsymbol x'-\boldsymbol p_\varepsilon\rangle+O\left(\frac{|\boldsymbol x'-\boldsymbol p_\varepsilon|^2}{|\boldsymbol x-\boldsymbol p_\varepsilon|}\right),
		\end{equation*}
		and
		\begin{equation*}
			|\boldsymbol x-\boldsymbol {\bar x}'|=|\boldsymbol x-\boldsymbol {\bar p_\varepsilon}|-\langle\frac{\boldsymbol x-\boldsymbol {\bar p_\varepsilon}}{|\boldsymbol x-\boldsymbol {\bar p_\varepsilon}|},\boldsymbol {\bar x}'-\boldsymbol {\bar p_\varepsilon}\rangle+O\left(\frac{|\boldsymbol x'-\boldsymbol p_\varepsilon|^2}{|\boldsymbol x-\boldsymbol {\bar p_\varepsilon}|}\right).
		\end{equation*}
		Hence we deduce
		\begin{equation*}
			\begin{split}
				\psi_{1,\varepsilon}&(\boldsymbol x)=\frac{q_\varepsilon^2}{2\pi\varepsilon^2}\int_{A_\varepsilon}\ln \frac{|\boldsymbol x-\boldsymbol {\bar x}'|}{|\boldsymbol x-\boldsymbol x'|}d\boldsymbol x'\\
				&=\frac{\kappa}{2\pi}\cdot p_\varepsilon\ln \frac{|\boldsymbol x-\boldsymbol {\bar p_\varepsilon}|}{|\boldsymbol x-\boldsymbol p_\varepsilon|}+\frac{p_\varepsilon^2}{2\pi\varepsilon^2}\int_{A_\varepsilon}\ln \frac{|\boldsymbol x-\boldsymbol p_\varepsilon|}{|\boldsymbol x-\boldsymbol x'|}d\boldsymbol x'-\frac{p_\varepsilon^2}{2\pi\varepsilon^2}\int_{A_\varepsilon}\ln \frac{|\boldsymbol x-\boldsymbol {\bar p_\varepsilon}|}{|\boldsymbol x-\boldsymbol {\bar x}'|}d\boldsymbol x'+O\left(\frac{ \varepsilon}{|\boldsymbol x-\boldsymbol p_\varepsilon|}\right)\\
				&=\frac{\kappa}{2\pi}\cdot p_\varepsilon\ln \frac{|\boldsymbol x-\boldsymbol {\bar p_\varepsilon}|}{|\boldsymbol x-\boldsymbol p_\varepsilon|}+O\left(\frac{\varepsilon}{|\boldsymbol x-\boldsymbol p_\varepsilon|}\right),
			\end{split}
		\end{equation*}
		where we use the circulation constraint \eqref{3-2}, $|\boldsymbol x-\boldsymbol p_\varepsilon|<|\boldsymbol x-\boldsymbol {\bar p_\varepsilon}|$ and $|p_\varepsilon-q_\varepsilon|=O(\varepsilon)$. Similarly, from the relation
		\begin{equation*}
			\frac{\boldsymbol x-\boldsymbol p_\varepsilon}{|\boldsymbol x-\boldsymbol p_\varepsilon|^2}-\frac{\boldsymbol x-\boldsymbol x'}{|\boldsymbol x-\boldsymbol x'|^2}=O\left(\frac{\varepsilon}{|\boldsymbol x-\boldsymbol p_\varepsilon|^2}\right),
		\end{equation*}
		and
		\begin{equation*}
			\frac{\boldsymbol x-\boldsymbol {\bar p_\varepsilon}}{|\boldsymbol x-\boldsymbol {\bar p_\varepsilon}|^2}-\frac{\boldsymbol x-\boldsymbol {\bar x}'}{|\boldsymbol x-\boldsymbol {\bar x}'|^2}=O\left(\frac{\varepsilon}{|\boldsymbol x-\boldsymbol {\bar p_\varepsilon}|^2}\right),
		\end{equation*}
		we can obtain
		\begin{equation*}
			\nabla \psi_{1,\varepsilon}(\boldsymbol x)=-\frac{\kappa}{2\pi}\cdot p_\varepsilon\frac{\boldsymbol x-\boldsymbol p_\varepsilon}{|\boldsymbol x-\boldsymbol p_\varepsilon|^2}+\frac{\kappa}{2\pi}\cdot p_\varepsilon\frac{\boldsymbol x-\boldsymbol {\bar p_\varepsilon}}{|\boldsymbol x-\boldsymbol {\bar p_\varepsilon}|^2}+O\left(\frac{\varepsilon}{|\boldsymbol x-\boldsymbol p_\varepsilon|^2}\right).
		\end{equation*}
		Thus the proof is complete.
	\end{proof}
	
	Compared with the main term $\psi_{1,\varepsilon}(\boldsymbol x)$, the $L^\infty$ norm of secondary term $\psi_{2,\varepsilon}(\boldsymbol x)$ is bounded due to the regularity of $H(\boldsymbol x,\boldsymbol x')$, as can be seen from the following estimate, and we can therefore obtain its estimates in the whole right half-plane.
	\begin{lemma}\label{lem3-5}
		For $\boldsymbol x\in \mathbb R^2_+$, it holds
		\begin{equation*}
			\psi_{2,\varepsilon}(\boldsymbol x)=\frac{\kappa}{p_\varepsilon} H(\boldsymbol x,\boldsymbol p_\varepsilon)+O(\varepsilon|\ln \varepsilon|).
		\end{equation*}
	\end{lemma}
	\begin{proof}
		Using the definition of $H(\boldsymbol x,\boldsymbol x')$, the circulation constraint \eqref{3-2} and computing directly as the proof of Lemma \ref{B1}, it holds
		\begin{equation*}
			\begin{split}
				\psi_{2,\varepsilon}(\boldsymbol x)-\frac{\kappa}{p_\varepsilon} H(\boldsymbol x,\boldsymbol p_\varepsilon)&=\frac{1}{\varepsilon^2}\int_{\mathbb R^2_+}\left(H(\boldsymbol x,\boldsymbol x')-H(\boldsymbol x,\boldsymbol p_\varepsilon)\right)\boldsymbol{1}_{A_\varepsilon}d\boldsymbol x'+O(\varepsilon)\\
				&=O(\varepsilon|\ln \varepsilon|),
			\end{split}
		\end{equation*}
		which is the desired result.
	\end{proof}

\bigskip
	
	Having made above preparations, we now turn to study the local behavior of $\psi_{1,\varepsilon}$ near $\boldsymbol p_\varepsilon$, and the main estimates of this part can be concluded from the following proposition.
	\begin{proposition}\label{prop3-6}
		$\psi_{1,\varepsilon}$ has the following asymptotic behavior as $\varepsilon\to 0$,
		\begin{equation*}
			\psi_{1,\varepsilon}(\boldsymbol x)=\frac{\sigma_\varepsilon^2}{\varepsilon^2}\cdot p_\varepsilon^2\left(w\left(\frac{\boldsymbol x-\boldsymbol p_\varepsilon}{\sigma_\varepsilon}\right)+o_\varepsilon(1)\right)+\mu_\varepsilon+\frac{W}{2}p_\varepsilon^2\ln\frac{1}{\varepsilon}-\frac{\kappa}{p_\varepsilon} H(\boldsymbol p_\varepsilon,\boldsymbol p_\varepsilon),\ \ \  \boldsymbol x\in B_{L\varepsilon}(\boldsymbol p_\varepsilon)
		\end{equation*}
		with
		\begin{equation*}
			w(\boldsymbol y)=\begin{cases}
				\frac{1}{4}(1-|\boldsymbol y|^2),\quad & \mathrm{if} \ |\boldsymbol y|\le 1,\\
				\frac{1}{2}\ln \frac{1}{|\boldsymbol y|},\quad &\mathrm{if} \ |\boldsymbol y|\geq 1,
			\end{cases}
		\end{equation*}
		where the flux constant $\mu_\varepsilon$ satisfies
		\begin{equation*}
			\frac{\kappa}{2\pi}\cdot p_\varepsilon\ln \left(\frac{1}{\sigma_\varepsilon}\right)-\frac{\kappa}{2\pi}\cdot p_\varepsilon\ln\frac{1}{2p_\varepsilon}+\frac{\kappa}{p_\varepsilon} H(\boldsymbol p_\varepsilon,\boldsymbol p_\varepsilon)-\frac{W}{2}p_\varepsilon^2\ln\frac{1}{\varepsilon}-\mu_\varepsilon=o_\varepsilon(1),
		\end{equation*}
		and the cross-section parameter $\sigma_\varepsilon$ satisfies
		\begin{equation*}
			\frac{|A_\varepsilon|}{\sigma_\varepsilon^2}\rightarrow \pi.
		\end{equation*}
	\end{proposition}
	
   According to the asymptotic behavior of $\psi_{\varepsilon,1}$ and $\psi_{\varepsilon,2}$ in Lemma \ref{lem3-5} and Proposition \ref{prop3-2}, we see that $\boldsymbol p_\varepsilon\in A_\varepsilon$ and is not on $\partial A_\varepsilon$.
	
	\bigskip
	
	The proof for Proposition \ref{prop3-2} will be divided into several Lemmas. By integration by part, it holds
	$$ \frac{1}{\varepsilon^2}\int_{A_\varepsilon}x_1\left(\psi_\varepsilon-\frac{W}{2}\ln\frac{1}{\varepsilon}x_1^2-\mu_\varepsilon\right)_+d\boldsymbol x=\int_{\Omega_\varepsilon}|\mathbf v-W\ln\frac{1}{\varepsilon}\mathbf e_z|^2,$$
	where the right hand side is the kinetic energy of rotation in vortex core with $\Omega_\varepsilon$ a topological torus in $\mathbb R^3$, and $\mathbf v$ the velocity field of the fluid. To apply the regularity theory for elliptic equation and obtain the $L^\infty$ estimate for $\psi_{1,\varepsilon}$ after blow up, we will prove the following lemma, which physically means this part of energy is bounded.
	\begin{lemma}\label{lem3-7}
		As $\varepsilon\to 0$, it holds
		\begin{equation*}
			 \frac{1}{\varepsilon^2}\int_{A_\varepsilon}x_1\left(\psi_\varepsilon-\frac{W}{2}\ln\frac{1}{\varepsilon}x_1^2-\mu_\varepsilon\right)_+d\boldsymbol x=O_\varepsilon(1).
		\end{equation*}
	\end{lemma}
	\begin{proof}
		We take $\psi_+=\left(\psi_\varepsilon-\frac{W}{2}\ln\frac{1}{\varepsilon}x_1^2-\mu_\varepsilon\right)_+$ as the upper truncation of $\psi_\varepsilon$ in $A_\varepsilon$. From equation \eqref{3-1}, it holds
		\begin{equation*}
			\begin{cases}
				-\varepsilon^2{\Delta^*} \psi_+(\boldsymbol x)=\boldsymbol{1}_{A_\varepsilon},\\
				\psi_+(\boldsymbol x)=0, \ \ \ \text{on} \ \partial A_\varepsilon.
			\end{cases}
		\end{equation*}
		Thus we can integrate by part to obtain
		\begin{equation*}
			\int_{A_\varepsilon}\frac{1}{x_1}|\nabla\psi_+|^2d\boldsymbol x=\frac{1}{\varepsilon^2}\int_{A_\varepsilon}x_1\psi_+d\boldsymbol x\le \frac{C|A_\varepsilon|^{1/2}}{\varepsilon^2}\left( \int_{A_\varepsilon}|\psi_+|^2d\boldsymbol x\right)^{1/2},
		\end{equation*}
		where we use the restriction $p_\varepsilon\in(c_1,c_2)$ and $A_\varepsilon$ is away from $x_2$-axis. By the Sobolev embedding, it holds
		\begin{equation*}
			\left( \int_{A_\varepsilon}|\psi_+|^2d\boldsymbol x\right)^{1/2}\le C\int_{A_\varepsilon}|\nabla\psi_+|d\boldsymbol x.
		\end{equation*}
		Hence we deduce
		\begin{equation*}
			\int_{A_\varepsilon}\frac{1}{x_1}|\nabla\psi_+|^2d\boldsymbol x\le\frac{C|A_\varepsilon|^{1/2}}{\varepsilon^2}\int_{A_\varepsilon}|\nabla\psi_+|d\boldsymbol x\le\frac{C|A_\varepsilon|}{\varepsilon^2}\left( \int_{A_\varepsilon}|\nabla\psi_+|^2d\boldsymbol x\right)^{1/2}.
		\end{equation*}
		Using the circulation constraint \eqref{3-2}, we finally obtain
		\begin{equation*}
			\frac{1}{\varepsilon^2}\int_{A_\varepsilon}x_1\psi_+d\boldsymbol x=\int_{A_\varepsilon}\frac{1}{x_1}|\nabla\psi_+|^2d\boldsymbol x=O_\varepsilon(1),
		\end{equation*}
		which is the estimate we need by the definition of $\psi_+$.
	\end{proof}
	
	Now we introduce a scaling version of $\psi_{1,\varepsilon}$ by letting
	\begin{equation*}
		w_\varepsilon(\boldsymbol y)=\frac{1}{p_\varepsilon^2}\cdot\frac{\varepsilon^2}{\sigma_\varepsilon^2}\left(\psi_{1,\varepsilon}(\sigma_\varepsilon \boldsymbol y+\boldsymbol p_\varepsilon)+\frac{\kappa}{p_\varepsilon} H(\boldsymbol p_\varepsilon,\boldsymbol p_\varepsilon)-\frac{W}{2}p_\varepsilon^2\ln\frac{1}{\varepsilon}-\mu_\varepsilon\right),
	\end{equation*}
	so that $w_\varepsilon$ is supported near a unit disc, and satisfies a Possion equation
	\begin{equation}\label{3-3}
		-\Delta w_\varepsilon=\boldsymbol 1_{\{w_\varepsilon>0\}}+f(\sigma_\varepsilon \boldsymbol y+\boldsymbol p_\varepsilon, w_\varepsilon), \ \ \  \text{in}\ \mathbb{R}^2,
	\end{equation}
	with
	\begin{equation*}
		f(\boldsymbol x,w):=\frac{q_\varepsilon^2}{p_\varepsilon^2}\cdot \boldsymbol 1_{\{\psi_\varepsilon(\boldsymbol x)-Wx_1^2\ln\frac{1}{\varepsilon}-\mu_\varepsilon>0\}}-\boldsymbol 1_{\{w>0\}}.
	\end{equation*}
	By the definition of $A_\varepsilon$ and estimate for $\psi_{2,\varepsilon}$ in Lemma \ref{lem3-5}, we see that
    $$w_\varepsilon(\boldsymbol y)=O(\varepsilon|\ln\varepsilon|), \quad \text{for} \ \ \sigma_\varepsilon \boldsymbol y+\boldsymbol p_\varepsilon\in \partial A_\varepsilon.$$

    Intuitively,  as $\varepsilon\to 0$, the equation satisfied by the limit function $w$ of $\{w_\varepsilon\}$ is $-\Delta w=\boldsymbol 1_{\{w>0\}}$ since $f(\sigma_\varepsilon \boldsymbol y+\boldsymbol p_\varepsilon, w_\varepsilon)$ goes to zero as we will show in Lemma \ref{lem3-9}. Then by method of moving plane in \cite{Fb}, $w(\boldsymbol x)$ must be radially symmetric and hence be the stream function of well-known Rankine vortex. To show the convergence, we are to give a uniform bound for the $L^\infty$ norm of $w_\varepsilon$ by regularity theory for elliptic equations.
	\begin{lemma}\label{lem3-8}
		There exists a constant $C_L>0$  independent of $\varepsilon$ such that $$||w_\varepsilon||_{L^\infty(B_L(\boldsymbol 0))}\leq C_L.$$
	\end{lemma}
	\begin{proof}
		It follows from Lemma \ref{lem3-7} and the assumption on $p_\varepsilon$ that
		\begin{align*}
			 O_\varepsilon(1)&=\frac{1}{\varepsilon^2}\int_{A_\varepsilon}x_1\left(\psi_\varepsilon-\frac{W}{2}\ln\frac{1}{\varepsilon}x_1^2-\mu_\varepsilon\right)_+d\boldsymbol x\\
			&=\frac{\sigma_\varepsilon^4}{\varepsilon^4}\cdot\left(p_\varepsilon^3+O(\varepsilon)\right)\cdot\int_{\mathbb R^2} (w_\varepsilon)_+ d\boldsymbol y+O(\varepsilon|\ln\varepsilon|).
		\end{align*}
		Notice that $\kappa=\varepsilon^{-2}\cdot p_\varepsilon|A_\varepsilon|+o_\varepsilon(1) \leq C\sigma_\varepsilon^2/\varepsilon^2$. Then we deduce
		$$\int_{\mathbb R^2} (w_\varepsilon)_+d\boldsymbol y\le C.$$
		Applying Morse iteration on \eqref{3-3} (Theorem 4.1 in \cite{Han}), we obtain
		$$||(w_\varepsilon)_+||_{L^\infty(B_L(\boldsymbol 0))}\leq C.$$
		To prove that the $L^\infty$ norm of $w_\varepsilon$ is bounded, we consider the following problem.
		\begin{equation*}
			\begin{cases}-\Delta w_1=\boldsymbol 1_{\{w_\varepsilon>0\}}+f(\sigma_\varepsilon \boldsymbol y+p_\varepsilon, w_\varepsilon),\quad  &\text{in} \ B_L(\boldsymbol 0),\\
				w_1=0,&\text{on} \ \partial B_L(\boldsymbol 0).\end{cases}
		\end{equation*}
		It is obvious that $|w_1|\le C$, and the bounds here and after are always dependent on $L$. Let $w_2:=w_\varepsilon-w_1$. Since $\sup_{B_L(\boldsymbol 0)} w_\varepsilon\ge 0$, function $w_2$ is harmonic in $B_L(\boldsymbol 0)$ and satisfies
		$$\sup_{B_L(\boldsymbol 0)} w_2\ge \sup_{B_L(\boldsymbol 0)} w_\varepsilon-C\geq -C.$$
		On the other hand, we infer from $||(w_\varepsilon)_+||_{L^\infty(B_L(\boldsymbol 0))}\leq C$ that
		$$\sup_{B_L(\boldsymbol 0)} w_2\leq \sup_{B_L(\boldsymbol 0)} w_\varepsilon+C\leq M,$$
		for some constant $M$. Hence $M-w_2$ is a positive harmonic function. Using the Harnack inequality, we have
		$$\sup_{B_L(\boldsymbol 0)}( M- w_2)\leq \tilde C\inf_{B_L(\boldsymbol 0)}(M- w_2)\leq \tilde C(M+\sup_{B_L(\boldsymbol 0)} w_2)\leq C$$
		with $\tilde C$ some positive constant. Since $\sup_{B_R(\boldsymbol 0)}( M- w_2)=M-\inf_{B_R(\boldsymbol 0)} w_2$, we deduce
		$$\inf_{B_R(\boldsymbol 0)} w_2\geq C,$$
		which implies the boundedness of $L^\infty$ norm for $w_\varepsilon$.
	\end{proof}

\bigskip
	
	The limiting function for $w_\varepsilon$ as $\varepsilon\to 0$ is established in the following lemma, which give the local behavior we need in Proposition \ref{prop3-2}.
	\begin{lemma}\label{lem3-9}
		As $\varepsilon\to0$, it holds
		$$w_\varepsilon\to w$$
		in $C_{\mathrm{loc}}^1(\mathbb{R}^2)$ for some radial function $w$.
	\end{lemma}
	\begin{proof}
		For $\boldsymbol y\in \mathbb R^2\setminus B_L(\boldsymbol 0)$, we infer from Lemma \ref{lem3-4} and Lemma \ref{lem3-5} that
		\begin{align*}
			w_\varepsilon(\boldsymbol y)&=\frac{1}{p_\varepsilon^2}\cdot\left(\psi_{1,\varepsilon}(\sigma_\varepsilon \boldsymbol y+\boldsymbol p_\varepsilon)+\frac{\kappa}{p_\varepsilon} H(\boldsymbol p_\varepsilon,\boldsymbol p_\varepsilon)-\frac{W}{2}p_\varepsilon^2\ln\frac{1}{\varepsilon}-\mu_\varepsilon\right)\\
			&=\frac{|A_\varepsilon|\cdot(1+O(\varepsilon))}{\sigma_\varepsilon^2}\cdot\left(\frac{1}{2\pi}\ln \left(\frac{1}{|\sigma_\varepsilon\boldsymbol y|}\right)-\frac{1}{2\pi}\ln\frac{1}{|\sigma_\varepsilon\boldsymbol y+\boldsymbol {\bar p_\varepsilon}-\boldsymbol p_\varepsilon|}\right.\\
			&\quad\left.+\frac{1}{p_\varepsilon^2} H(\boldsymbol p_\varepsilon,\boldsymbol p_\varepsilon)
			-\frac{W}{2\kappa}\cdot p_\varepsilon\ln\frac{1}{\varepsilon}-\frac{\mu_\varepsilon}{p_\varepsilon\kappa}+O\left(\frac{1}{L}\right)\right)\\
			&=\frac{|A_\varepsilon|\cdot(1+O(\varepsilon))}{\sigma_\varepsilon^2}\cdot\frac{1}{2\pi}\ln \frac{1}{|\boldsymbol y|}\\
			&\quad+\frac{|A_\varepsilon|\cdot(1+O(\varepsilon))}{\sigma_\varepsilon^2}\cdot\left(\frac{1}{2\pi}\ln \left(\frac{1}{\sigma_\varepsilon}\right)-\frac{1}{2\pi}\ln\frac{1}{|\sigma_\varepsilon\boldsymbol y+\boldsymbol {\bar p_\varepsilon}-\boldsymbol p_\varepsilon|}\right.\\
			&\quad\left.+\frac{1}{p^2_\varepsilon} H(\boldsymbol p_\varepsilon,\boldsymbol p_\varepsilon)-\frac{W}{2\kappa}\cdot p_\varepsilon\ln\frac{1}{\varepsilon}-\frac{\mu_\varepsilon}{p_\varepsilon\kappa}+O\left(\frac{1}{L}\right)\right).
		\end{align*}
		Since $|A_\varepsilon|/\sigma^2_\varepsilon\le C$ and  $||w_\varepsilon||_{L^\infty(B_L(\boldsymbol 0))}\leq C_L$ by Lemma \ref{lem3-8}, we may assume
		$$|A_\varepsilon|/\sigma^2_\varepsilon\to\boldsymbol t,$$
		and
		\begin{equation*}
			\begin{split}
				\frac{|A_\varepsilon|\cdot(1+O(\varepsilon))}{\sigma_\varepsilon^2}\cdot\bigg(\frac{1}{2\pi}\ln \left(\frac{1}{\sigma_\varepsilon}\right)&-\frac{1}{2\pi}\ln\frac{1}{|\sigma_\varepsilon\boldsymbol y+\boldsymbol {\bar p_\varepsilon}-\boldsymbol p_\varepsilon|}\\
				&+\frac{1}{p^2_\varepsilon} H(\boldsymbol p_\varepsilon,\boldsymbol p_\varepsilon)-\frac{W}{2\kappa}\cdot p_\varepsilon\ln\frac{1}{\varepsilon}-\frac{\mu_\varepsilon}{p_\varepsilon\kappa}\bigg)\to \tau,
			\end{split}
		\end{equation*}
		for some $\boldsymbol t\in[0,+\infty)$ and $\tau\in(-\infty,+\infty)$.
		By \eqref{3-3}, we may further assume that $w_\varepsilon\to w$ in $C_{\text{loc}}^1(\mathbb{R}^2)$ and $w$ satisfies
		\begin{equation*}
			\begin{cases}
				-\Delta w=\boldsymbol1_{\{w>0\}},\quad &\text{in}\,B_R(\boldsymbol 0),\\
				w=\frac{\boldsymbol t}{2\pi}\ln \frac{1}{|\boldsymbol y|}+\tau+O\left(\frac{1}{L}\right), &\text{in} \, B_R(\boldsymbol 0)\setminus B_L(\boldsymbol 0).
			\end{cases}	
		\end{equation*}
		Moreover, $w$ will satisfy the integral equation
		$$w(\boldsymbol y)=\frac{1}{2\pi}\int_{\mathbb{R}^2} \ln\left(\frac{1}{|\boldsymbol y-\boldsymbol y'|}\right) \boldsymbol 1_{\{w>0\}}(\boldsymbol y')d\boldsymbol y'+\tau.$$
		Then the method of moving planes shows that $w$ is radial and decreasing (See e.g. \cite{T}), which completes the proof of this lemma.
	\end{proof}
	
	\bigskip
	
	\noindent{\bf Proof of Proposition \ref{prop3-6}:}
	By the definition of $\sigma_\varepsilon$, there exists a $\boldsymbol y_\varepsilon$ with $|\boldsymbol y_\varepsilon|=1$ and $\sigma_\varepsilon \boldsymbol y_\varepsilon+\boldsymbol p_\varepsilon\in \partial A_\varepsilon$. Thus it holds
	\begin{equation*}
		w(\boldsymbol y)=\begin{cases}
			\frac{1}{4}(1-|\boldsymbol y|^2),\quad &|\boldsymbol y|\le 1,\\
			\frac{1}{2}\ln \frac{1}{|\boldsymbol y|},\quad &|\boldsymbol y|\ge 1.
		\end{cases}
	\end{equation*}
	We further have that $\boldsymbol t=\pi$ and $\tau+O(1/L)=0$. Since $\tau$ is not depend on $L$, while $O(1/L)\to 0$ as $L\to+\infty$, one must have $\tau=0$ and $O(1/L)=0$. The proof of Proposition \ref{prop3-6} is hence complete. \qed
	
	\bigskip
	
	To prove Proposition \ref{prop3-3}, we will construct a family of approximate solutions $\Phi_{\boldsymbol q_\varepsilon,\varepsilon}(\boldsymbol x)$, and apply the Pohozaev identity in Appendix \ref{appC} to $\psi_{1,\varepsilon}$. Let us recall the definition of functions $\mathcal V_{\boldsymbol q_\varepsilon,\varepsilon}$, $\mathcal H_{\boldsymbol q_\varepsilon,\varepsilon}$, whose properties are discussed in the second part of Section \ref{sec2} and in Appendix \ref{appB}. We choose the approximate solutions to \eqref{3-1} and \eqref{3-2} of the form
	\begin{equation*}
		\begin{split}
		\Phi_{\boldsymbol q_\varepsilon,\varepsilon}(\boldsymbol x):&=V_{\boldsymbol q_\varepsilon,\varepsilon}(\boldsymbol x)-V_{\bar{\boldsymbol q}_\varepsilon,\varepsilon}(\boldsymbol x)+\mathcal H_{\boldsymbol q_\varepsilon,\varepsilon}(\boldsymbol x)\\
		&=\mathcal V_{\boldsymbol q_\varepsilon,\varepsilon}(\boldsymbol x)+\mathcal H_{\boldsymbol q_\varepsilon,\varepsilon}(\boldsymbol x),
		\end{split}
	\end{equation*}
	where the parameters $q_\varepsilon$, $a_\varepsilon$ and $s_\varepsilon$ in $\Phi_{\boldsymbol q_\varepsilon,\varepsilon}(\boldsymbol x)$ are chosen so that 
	\begin{equation*}
		\partial_1\Phi_{\boldsymbol q_\varepsilon,\varepsilon}(\boldsymbol p_\varepsilon)=0,
	\end{equation*}
	\begin{equation*}
		\frac{a_\varepsilon}{2\pi}\ln\frac{1}{\varepsilon}=\mu_\varepsilon+\frac{W}{2}q_\varepsilon^2\ln\frac{1}{\varepsilon}-\mathcal H_{\boldsymbol q_\varepsilon,\varepsilon}(\boldsymbol q_\varepsilon)+V_{\boldsymbol{\bar q_\varepsilon},\varepsilon}(\boldsymbol q_\varepsilon),
	\end{equation*}
	and
	\begin{equation*}
		\frac{a_\varepsilon}{2\pi}\ln\frac{1}{\varepsilon}\cdot\frac{1}{s_\varepsilon|\ln s_\varepsilon|}=\frac{s_\varepsilon}{2\varepsilon^2}\cdot q_\varepsilon^2.
	\end{equation*}
	The same as \eqref{2-14} in Section 2, here we also denote
	\begin{equation*}
		\mathcal N_\varepsilon:=\frac{a_\varepsilon}{2\pi}\ln\frac{1}{\varepsilon}\cdot\frac{1}{s_\varepsilon|\ln s_\varepsilon|}=\frac{s_\varepsilon}{2\varepsilon^2}\cdot q_\varepsilon^2
	\end{equation*}
	as the value of $|\nabla V_{\boldsymbol q_\varepsilon,\varepsilon}|$ at $|\boldsymbol x-\boldsymbol q_\varepsilon|= s_\varepsilon$, which will blow up at the rate $1/\varepsilon$. Notice the first condition \eqref{3-4} is equivalent to
	\begin{equation*}
		\frac{|p_\varepsilon-q_\varepsilon|}{2\varepsilon^2}+O(\varepsilon)=\partial_1 V_{\bar{\boldsymbol q_\varepsilon},\varepsilon}(\boldsymbol p_\varepsilon)-\partial_1\mathcal H_{\boldsymbol q_\varepsilon,\varepsilon}(\boldsymbol p_\varepsilon)+O(\varepsilon).
	\end{equation*}
	By the asymptotic estimates given in Proposition \ref{prop3-6} and implicit function theorem, the system \eqref{3-4}-\eqref{3-6} is solvable. Moreover, we have estimates
	\begin{equation}\label{3-4}
		|p_\varepsilon-q_\varepsilon|=O(\varepsilon^2|\ln\varepsilon|),
	\end{equation}
	\begin{equation}\label{3-5}
		\frac{a_\varepsilon}{2\pi}\ln\frac{1}{\varepsilon}=\mu_\varepsilon+\frac{W}{2}q_\varepsilon^2\ln\frac{1}{\varepsilon}+O_\varepsilon(1),
	\end{equation}
	and
	\begin{equation}\label{3-6}
		|\sigma_\varepsilon-s_\varepsilon|=o(\varepsilon).
	\end{equation}
    We are now in the position to prove Proposition \ref{prop3-3}.
	
	\noindent{\bf Proof of Proposition \ref{prop3-3}:}
    We can apply the local Pohozaev identity \eqref{C-1} in Appendix \ref{appC} to $\psi_{1,\varepsilon}$ and obtain
	\begin{equation*}
		\begin{split}
			&\quad-\int_{\partial B_\delta(\boldsymbol q_\varepsilon)} \frac{\partial\psi_{1,\varepsilon}}{\partial \nu}\frac{\partial\psi_{1,\varepsilon}}{\partial x_1}dS+ \frac{1}{2}\int_{\partial B_\delta(\boldsymbol q_\varepsilon)} |\nabla\psi_{1,\varepsilon}|^2 \nu_1dS\\
			&=-\frac{q_\varepsilon^2}{\varepsilon^2}\int_{B_\delta(\boldsymbol q_\varepsilon)} \partial_1\psi_{2,\varepsilon}(\boldsymbol x)\cdot\boldsymbol 1_{A_\varepsilon}(\boldsymbol x)d\boldsymbol x+\frac{q_\varepsilon^2}{\varepsilon^2}\int_{B_\delta(\boldsymbol q_\varepsilon)} Wx_1\ln\frac{1}{\varepsilon}\cdot\boldsymbol 1_{A_\varepsilon}(\boldsymbol x)d\boldsymbol x,
		\end{split}
	\end{equation*}
	where $\delta$ is a small positive number. Since $|A_\varepsilon|/\sigma^2_\varepsilon\to \pi$ as $\varepsilon\to 0$, from the isoperimetric inequality and \eqref{3-4}-\eqref{3-6}, we see that $A_\varepsilon$ tends to a disc with radius $\sigma_\varepsilon\to s_\varepsilon^*:=(\frac{\kappa}{q_\varepsilon\pi})^{1/2}\varepsilon$ centered at $\boldsymbol p_\varepsilon$. Thus the crucial symmetry difference for $A_\varepsilon$ and $B_{s_\varepsilon^*}(\boldsymbol q_\varepsilon)$ needed in Appendix \ref{appC} can be estimated as $|A_\varepsilon \bigtriangleup B_{s_\varepsilon^*}(\boldsymbol q_\varepsilon)|=o(\varepsilon^2)$.
	
	Then according to Lemma \ref{C5}, we have
	\begin{equation*}
		 Wq_\varepsilon\ln\frac{1}{\varepsilon}-\frac{\kappa}{4\pi}\ln\frac{8q_\varepsilon}{\sigma_\varepsilon}+\frac{\kappa}{16\pi}=o_\varepsilon(1).
	\end{equation*}
	So we have finished the proof of Proposition \ref{prop3-3}. \qed
	
	\bigskip
	
	\subsection{Refined estimates and revised Kelvin--Hicks formula }
	
	For the uniqueness of $\psi_\varepsilon$, we need to improve the results in Propositions \ref{prop3-3} and \ref{prop3-6}. So we reconsider the problem \eqref{3-1}
	\begin{equation*}
		\begin{cases}
			-\varepsilon^2{\Delta^*}\psi_\varepsilon=\boldsymbol1_{\{\psi_\varepsilon-\frac{W}{2}x_1^2\ln\frac{1}{\varepsilon}>\mu_\varepsilon\}}, & \text{in} \ \mathbb R^2_+,
			\\
			\psi_\varepsilon=0, & \text{on} \ x_1=0,
			\\
			\psi_\varepsilon,  \ |\nabla\psi_\varepsilon|/x_1\to0, & \text{as} \ |\boldsymbol x |\to \infty,
		\end{cases}
	\end{equation*}
	together with circulation constraint \eqref{3-2}
	\begin{equation*}
		\frac{1}{\varepsilon^2}\int_{A_\varepsilon}x_1d\boldsymbol x=\kappa.
	\end{equation*}
	
	By a same approach in Section 2, we also denote the difference of $\psi_\varepsilon$ and $\Phi_{\boldsymbol q_\varepsilon,\varepsilon}$ as the error term
	\begin{equation*}
		\phi_\varepsilon(\boldsymbol x):=\psi_\varepsilon(\boldsymbol x)-\Phi_{\boldsymbol q_\varepsilon,\varepsilon}(\boldsymbol x).
	\end{equation*}
	Hence our task in this part is to improve the estimate for $\phi_\varepsilon$. Recall the definition of the norm in \eqref{2-18} as
    \begin{equation*}
    ||\phi||_*:=\sup_{\boldsymbol x\in\mathbb R^2_+}\rho_1(\boldsymbol x)\rho_2(\boldsymbol x)|\phi(\boldsymbol x)|
    \end{equation*}
    with the weight 
    \begin{equation*}
    	\rho_1(\boldsymbol x):=\frac{(1+|\boldsymbol x-\boldsymbol q_\varepsilon|^2)^{\frac{3}{2}}}{1+x_1^2} \ \ \ \text{and} \ \ \	\rho_2(\boldsymbol x):=\left(\frac{1}{x_1}+1\right),
    \end{equation*}
    According to the asymptotic estimates in Proposition \ref{prop3-6}, we have the following lemma concerning the behavior of $\phi_\varepsilon(\boldsymbol x)$ in $\mathbb R^2_+$.
	\begin{lemma}\label{lem3-10}
		As $\varepsilon\to 0$, it holds
		$$\|\phi_\varepsilon\|_*=o_\varepsilon(1).$$
	\end{lemma}
	\begin{proof}
		In view of Proposition \ref{prop3-6} and estimates \eqref{3-4}--\eqref{3-6}, it is obvious that
		\begin{equation*}
			||\phi_\varepsilon||_{L^\infty(B_{Ls_\varepsilon}(\boldsymbol q_\varepsilon))}=o_\varepsilon(1).
		\end{equation*}
		
		While for those $\boldsymbol x$ far away from $B_{Ls_\varepsilon}(\boldsymbol q_\varepsilon)$, it holds
		\begin{equation*}
			\phi_\varepsilon(\boldsymbol x)=\frac{1}{\varepsilon^2}\int_{\mathbb R^2_+}{G_*}(\boldsymbol x,\boldsymbol x')(\boldsymbol 1_{A_\varepsilon}(\boldsymbol x')-\boldsymbol 1_{B_{s_\varepsilon}(\boldsymbol q_\varepsilon)}(\boldsymbol x'))d\boldsymbol x'.
		\end{equation*}
		Since
		\begin{equation*}
			\frac{1}{\varepsilon^2}||\boldsymbol 1_{A_\varepsilon}-\boldsymbol 1_{B_{s_\varepsilon}(\boldsymbol q_\varepsilon)}||_{L^1(B_{Ls_\varepsilon}(\boldsymbol q_\varepsilon))}=o_\varepsilon(1),
		\end{equation*}
		we can use the expansion
		\begin{equation*}
			\left(\frac{1}{x_1}+1\right){G_*}(\boldsymbol x,\boldsymbol x')\le C\cdot \frac{1+x_1^2}{(1+|\boldsymbol x-\boldsymbol q_\varepsilon|^2)^{\frac{3}{2}}},
		\end{equation*}
		and Young inequality to derive
		\begin{equation*}
			||\phi_\varepsilon||_*=o_\varepsilon(1),
		\end{equation*}
		which yields the conclusion.
	\end{proof}
	
	\bigskip
	
	In view of linearization for \eqref{3-1} near the approximate solution $\Phi_{\boldsymbol q_\varepsilon,\varepsilon}(\boldsymbol x)$, the perturbation $\phi_\varepsilon$ will satisfy the equation
	\begin{equation}\label{phieq}
		\begin{cases}
			\mathbb L_\varepsilon\phi_\varepsilon=R_\varepsilon(\phi_\varepsilon), & \text{in} \ \mathbb R^2_+,
			\\
			\partial_1\phi_\varepsilon(\boldsymbol p_\varepsilon)=0,
			\\
			\phi_\varepsilon=0, & \text{on} \ x_1=0,\\
			\phi_\varepsilon, \ |\nabla\phi_\varepsilon|/x_1\to0, &\text{as} \ |\boldsymbol x |\to \infty,
		\end{cases}
	\end{equation}
	where $\mathbb L_\varepsilon$ is the linear operator defined by
	\begin{equation*}
		\mathbb L_\varepsilon\phi_\varepsilon=-x_1{\Delta^*}\phi_\varepsilon-\frac{2}{s_\varepsilon q_\varepsilon}\phi_\varepsilon(r,\theta)\boldsymbol\delta_{|\boldsymbol x-\boldsymbol q_\varepsilon|=s_\varepsilon},
	\end{equation*}
	and
	\begin{equation*}
		 R_\varepsilon(\phi_\varepsilon)=\frac{1}{\varepsilon^2}\bigg(x_1\boldsymbol1_{\{\psi_\varepsilon-\frac{W}{2}x_1^2\ln\frac{1}{\varepsilon}>\mu_\varepsilon\}}-x_1\boldsymbol1_{\{V_{\boldsymbol q_\varepsilon,\varepsilon}>\frac{a_\varepsilon}{2\pi}\ln\frac{1}{\varepsilon}\}}-\frac{2}{s_\varepsilon q_\varepsilon}\phi_\varepsilon(r,\theta)\boldsymbol\delta_{|\boldsymbol x-\boldsymbol q_\varepsilon|=s_\varepsilon}\bigg).
	\end{equation*}
	By Lemma \ref{B4} in Appendix \ref{appB}, we see that the nonlinear term $R_\varepsilon(\phi_\varepsilon)$ is supported on an annulus such that
	\begin{equation*}
		R_\varepsilon(\phi_\varepsilon)=0, \ \ \ \text{in} \ \left(\mathbb R^2_+\setminus B_{2s_\varepsilon}(\boldsymbol q_\varepsilon)\right) \cup B_{s_\varepsilon/2}(\boldsymbol q_\varepsilon).
	\end{equation*}
	
	To derive a better estimate for $\phi_\varepsilon$, we will establish the following lemma concerning the linear operator $\mathbb L_\varepsilon$, which is similar in form as the coercive estimate in Lemma \ref{lem2-3} at first glance. However, here we use the condition $\partial_1\psi_\varepsilon(\boldsymbol p_\varepsilon)=0$ (since $\boldsymbol p_\varepsilon$ is a maximum point) to overcome the problem of degeneracy, instead of considering a projective equation.

	\begin{lemma}\label{lem3-11}
		Suppose that $\mathrm{supp} \,\mathbb L_\varepsilon\phi_\varepsilon\subset B_{2s_\varepsilon}(\boldsymbol q_\varepsilon)$. Then for any $p\in(2, +\infty]$, there are $\varepsilon_0>0$ and $c_0>0$ such that for any $\varepsilon\in (0,\varepsilon_0]$, it holds
		\begin{equation*}
			\quad \varepsilon^{1-\frac{2}{p}} ||\mathbb{L}_\varepsilon \phi_\varepsilon||_{W^{-1,p}(B_{Ls_\varepsilon}(\boldsymbol q_\varepsilon))}+\varepsilon^2||\mathbb{L}_\varepsilon \phi_\varepsilon||_{L^\infty(B_{s_\varepsilon/2}(\boldsymbol q_\varepsilon))}\ge c_0 \left(\varepsilon^{1-\frac{2}{p}} ||\nabla \phi_\varepsilon||_{L^{p}(B_{Ls_\varepsilon}(\boldsymbol q_\varepsilon))}+||\phi_\varepsilon||_*\right).
		\end{equation*}
	\end{lemma}
	\begin{proof}
		We will argue by contradiction. Suppose on the contrary that there exists $\varepsilon_n\to 0$ such that $\phi_n:=\phi_{\varepsilon_n}$ satisfies
		\begin{equation*}
			\varepsilon_n^{1-\frac{2}{p}} ||\mathbb{L}_{\varepsilon_n} \phi_n||_{W^{-1,p}(B_{Ls_{\varepsilon_n}}(\boldsymbol q_{\varepsilon_n}))}+\varepsilon_n^2||\mathbb{L}_{\varepsilon_n} \phi_n||_{L^\infty(B_{s_\varepsilon/2}(\boldsymbol q_{\varepsilon_n}))}\le \frac{1}{n},
		\end{equation*}
		and
		\begin{equation}\label{3-8}
			\varepsilon_n^{1-\frac{2}{p}} ||\nabla \phi_n||_{L^{p}(B_{Ls_{\varepsilon_n}}(\boldsymbol q_{\varepsilon_n}))}+||\phi_n||_*=1.
		\end{equation}
		By letting $f_n=\mathbb{L}_{\varepsilon_n}\phi_n$, we have
		\begin{equation*}
			-{\Delta^*} \phi_n=\frac{2}{s_{\varepsilon_n}q_{\varepsilon_n}}\phi_n(r,\theta)\boldsymbol\delta_{|\boldsymbol x-\boldsymbol q_{\varepsilon_n}|=s_{\varepsilon_n}}+f_n.
		\end{equation*}
		Here, we also denote $\tilde v(\boldsymbol y):=v(s_{\varepsilon_n} \boldsymbol y+\boldsymbol q_{\varepsilon_n})$ for an arbitrary function. Then the above equation has a weak form
		\begin{equation*}
			\int_{\mathbb R^2_+}\frac{1}{s_{\varepsilon_n}y_1+q_{\varepsilon_n}}\cdot\nabla\tilde\phi_n\cdot\nabla\varphi d\boldsymbol y=2\int_{|\boldsymbol y|=1}\frac{1}{q_{\varepsilon_n}}\tilde\phi_n\varphi+\langle\tilde f_n,\varphi\rangle, \ \ \  \forall \, \varphi\in C_0^\infty(\mathbb R^2).
		\end{equation*}
		
		Since the right hand side of the equation is bounded in $W_{\text{loc}}^{-1,p}(\mathbb{R}^2)$, $\tilde \phi_n$ is bounded in $W^{1,p}_{\text{loc}}(\mathbb{R}^2)$ and hence bounded in $C_{\text{loc}}^\alpha (\mathbb{R}^2)$ for some $\alpha>0$ by Sobolev embedding. We may assume that $\tilde \phi_n$ converges uniformly in any compact subset of $\mathbb{R}^2$ to $\phi^*\in L^\infty(\mathbb{R}^2)\cap C(\mathbb{R}^2)$, and the limiting function $\phi^*$ satisfies
		\begin{equation*}
			-\Delta \phi^* =2\phi^*(r,\theta)\boldsymbol \delta_{|\boldsymbol y|=1}, \quad \text{in} \ \mathbb{R}^2.
		\end{equation*}
		Therefore, we conclude from the nondegeneracy of limiting operator and symmetry with respect to $x_1$-axis that
		\begin{equation*}
			\phi^*=C_1\cdot\frac{\partial w}{\partial y_1}
		\end{equation*}
		with $C_1$ a constant, and
		\begin{equation*}
			w(\boldsymbol y)=\left\{
			\begin{array}{lll}
				\frac{1}{4}(1-|\boldsymbol y|^2), \ \ \  & \mathrm{if} \  |\boldsymbol y|\le 1,\\
				\frac{1}{2}\ln\frac{1}{|\boldsymbol y|}, &\mathrm{if} \  |\boldsymbol y|\ge 1.
			\end{array}
			\right.
		\end{equation*}
		On the other hand, since $\varepsilon_n^2|f_n|\leq 1/n$ in $B_{s_{\varepsilon_n}/2}(\boldsymbol q_{\varepsilon_n})$ and $|\tilde \phi_n|\leq 1$, we know that $\tilde \phi_n$ is bounded in $W^{2,p}({B_{1/4}(\boldsymbol 0)})$. Thus we may assume $\tilde \phi_n\to \phi^*$ in $C^1({B_{1/4}(\boldsymbol 0)})$. Since by condition $\partial_1\Phi_{\boldsymbol q_\varepsilon,\varepsilon}(\boldsymbol p_\varepsilon)=0$ and \eqref{3-4} we have
		$$\partial_1 \tilde \phi_n \left(\frac{\boldsymbol p_{\varepsilon_n}-\boldsymbol q_{\varepsilon_n}}{s_{\varepsilon_n}}\right)=\partial_1 \phi_n(\boldsymbol p_{\varepsilon_n})=0$$
		and $\frac{\boldsymbol p_{\varepsilon_n}-\boldsymbol q_{\varepsilon_n}}{s_{\varepsilon_n}}\to 0$, it holds $\partial_1 \phi^*(\boldsymbol 0)=0$. This implies $C_1=0$ and hence $\phi^*\equiv 0$.
		
		Therefore, we have proved $\phi_n=o_n(1)$ in $B_{Ls_{\varepsilon_n}}(\boldsymbol q_{\varepsilon_n})$. Then, using the strong maximum principle and a similar argument as in the proof of Lemma \ref{lem2-2}, we can derive
		\begin{equation}\label{3-9}
			||\phi_n||_*\le C||\phi_n||_{L^\infty(B_{Ls_{\varepsilon_n}}(\boldsymbol q_{\varepsilon_n}))}=o_n(1).
		\end{equation}
		Now we turn to consider the norm  $||\nabla \phi_\varepsilon||_{L^{p}(B_{Ls_{\varepsilon_n}}(\boldsymbol q_{\varepsilon_n}))}$. For any $\tilde\varphi\in C_0^\infty (B_{L}(\boldsymbol 0))$,  it holds
		\begin{equation}\label{3-10}
			\begin{split}
				&\left|\int_{D_n}\frac{1}{s_{\varepsilon_n}y_1+q_{\varepsilon_n}}\cdot\nabla\tilde\phi_n\cdot\nabla\tilde\varphi d\boldsymbol y\right|=\left|2\int_{|\boldsymbol y|=1}\frac{1}{q_{\varepsilon_n} }\tilde\phi_n\varphi+\langle\tilde f_n,\tilde\varphi\rangle \right|\\
				& \ \ \ \ \ \ \ \ \ =o_n(1)\cdot\|\tilde\varphi\|_{W^{1,1}(B_L(\boldsymbol 0))}+o_n(1)\cdot\|\tilde\varphi\|_{W^{1,p'}(B_L(\boldsymbol 0))}\\
				& \ \ \ \ \ \ \ \ \ =o_n(1)\cdot\left(\int_{B_L(\boldsymbol 0)}|\nabla\tilde\varphi|^{p'}\right)^{\frac{1}{p'}}.
			\end{split}
		\end{equation}
		Thus we have
		\begin{equation*}
			\varepsilon_n^{1-\frac{2}{p}}\|\nabla\phi_n\|_{L^p(B_{Ls_{\varepsilon_n}}(\boldsymbol q_{\varepsilon_n}))}\le C||\nabla\tilde\phi_n||_{L^p(B_L(\boldsymbol 0))}=o_n(1).
		\end{equation*}
		We see that \eqref{3-9} and \eqref{3-10} is a contradiction to \eqref{3-8}, and hence the proof of Lemma \ref{lem3-11} is finished.
	\end{proof}

\bigskip
	
    In the following, we always make the convention that $\tilde v(\boldsymbol y)=v(s_\varepsilon\boldsymbol y+\boldsymbol q_\varepsilon)$ for a general function $v$.	In Appendix \ref{appB}, we denote 
	$$\mathcal W_\varepsilon(s_\varepsilon)=-Wq_\varepsilon\ln\frac{1}{\varepsilon}+\frac{s_\varepsilon^2}{4\varepsilon^2}\cdot q_\varepsilon\ln\frac{8q_\varepsilon}{s_\varepsilon}-\frac{s_\varepsilon^2}{16\varepsilon^2}\cdot q_\varepsilon$$
	as the coefficient appearing in the linear part of the difference of $\tilde{\mathbf\Gamma}_\varepsilon:=\{\boldsymbol y\ | \ \tilde{\mathbf U}_{\boldsymbol q_\varepsilon,\varepsilon}(\boldsymbol y)=0\}$ with $\partial B_1(\boldsymbol 0)$. We also let 
	$$t_{\varepsilon,\mathcal W}(\theta)=\frac{\mathcal W_\varepsilon(s_\varepsilon)}{\mathcal N_\varepsilon}\cdot(\cos\theta,0) \quad \mathrm{and} \quad t_{\varepsilon,\tilde\phi_\varepsilon}=\frac{\tilde\phi_\varepsilon(\cos\theta,\sin\theta)}{s_\varepsilon\mathcal N_\varepsilon}\cdot (\cos\theta,\sin\theta), \quad \theta\in (0,2\pi]$$ 
	be the linear part of the difference of $\partial \tilde A_\varepsilon$ with $\partial B_1(\boldsymbol 0)$ induced by $\mathcal W_\varepsilon(s_\varepsilon)$ and $\phi_\varepsilon$ separately. Since $s_\varepsilon\mathcal N_\varepsilon$ is bounded as a perturbation near the gradient of stream function for Rankine vortex at the vortex boundary, we see $t_{\varepsilon,\mathcal W}(\theta)$ is of the same order as $s_\varepsilon\mathcal W_\varepsilon(s_\varepsilon)$, and $t_{\varepsilon,\tilde\phi_\varepsilon}(\theta)$ is of the same order as $\|\phi_\varepsilon\|_{L^{\infty}(B_{Ls_\varepsilon}(\boldsymbol q_\varepsilon))}$.

	By the coercive estimate obtained in Lemma \ref{lem3-11}, we are now in the position to improve the estimate for error $\phi_\varepsilon$ in terms of $t_{\varepsilon,\mathcal W}(\theta)$ and $t_{\varepsilon,\tilde\phi_\varepsilon}(\theta)$.
	\begin{lemma}\label{lem3-12}
		For $p\in (2,+\infty]$ and $\varepsilon\in (0,\varepsilon_0]$ small, it holds
		\begin{equation}\label{3-11}
			||\phi_\varepsilon||_*+\varepsilon^{1-\frac{2}{p}} ||\nabla \phi_\varepsilon||_{L^{p}(B_{Ls_\varepsilon}(\boldsymbol q_\varepsilon))}=O\left(t_{\varepsilon,\mathcal W}+\varepsilon |t_{\varepsilon,\mathcal W}+t_{\varepsilon,\tilde \phi_\varepsilon}|^{\frac{1}{2}+\frac{1}{p}}+\varepsilon^2|\ln\varepsilon|\right).
		\end{equation}
	\end{lemma}	
	\begin{proof}
		In view of Lemma \ref{lem3-11}, it is sufficient to verify that
\begin{equation*}
\begin{split}
				&\quad\varepsilon^{1-\frac{2}{p}} ||R_\varepsilon(\phi_\varepsilon)||_{W^{-1,p}(B_{Ls_\varepsilon}(\boldsymbol q_\varepsilon))}+\varepsilon^2||R_\varepsilon(\phi_\varepsilon)||_{L^\infty(B_{s_\varepsilon/2}(\boldsymbol q_\varepsilon))}\\
				&=O\left(t_{\varepsilon,\mathcal W}+\varepsilon |t_{\varepsilon,\mathcal W}+t_{\varepsilon,\tilde \phi_\varepsilon}|^{\frac{1}{2}+\frac{1}{p}}+\varepsilon^2|\ln\varepsilon|\right).
			\end{split}
\end{equation*}
		Notice that we have
		\begin{equation*}
			R_\varepsilon(\phi_\varepsilon)\equiv 0,\quad \text{in}\ B_{s_\varepsilon/2}(\boldsymbol q_\varepsilon).
		\end{equation*}
		So it remains to estimate $\varepsilon^{1-\frac{2}{p}} ||R_\varepsilon(\phi_\varepsilon)||_{W^{-1,p}(B_{Ls_\varepsilon}(\boldsymbol q_\varepsilon))}$.
		
		Then for each $\varphi(\boldsymbol x)\in C_0^1(B_{Ls_\varepsilon}(\boldsymbol q_\varepsilon))$, we have
		\begin{equation*}
			\begin{split}
           \langle R_\varepsilon(\phi_\varepsilon),\varphi \rangle&=\frac{s_\varepsilon^2}{\varepsilon^2}\int_{B_L(\boldsymbol 0)}(s_\varepsilon y_1+q_\varepsilon)\left(\boldsymbol1_{\{\tilde\psi_\varepsilon-\frac{W}{2}(s_\varepsilon y_1+q_\varepsilon)^2\ln\frac{1}{\varepsilon}>\mu_\varepsilon\}}-\boldsymbol 1_{\{\tilde V_{\boldsymbol q_\varepsilon,\varepsilon}>\frac{a_\varepsilon}{2\pi}\ln\frac{1}{\varepsilon}\}}\right)\tilde\varphi d\boldsymbol y\\
			&\quad-\frac{2}{q_\varepsilon}\int_0^{2\pi}\tilde\phi_\varepsilon\tilde\varphi(1,\theta)d\theta.
			\end{split}
\end{equation*}
		Denote
		\begin{equation*}
			\begin{split}
				\boldsymbol y_{\varepsilon}(\theta):&=(1+t_{\varepsilon}(\theta))(\cos\theta,\sin\theta)\\
				&=(1+t_{\varepsilon,\mathcal W}(\theta)+t_{\varepsilon,\tilde\phi_m}(\theta)+o(\varepsilon))(\cos\theta,\sin\theta)\\
				&\in\{\boldsymbol y \,\,\,|\,\, \ \tilde{\mathbf U}_{\boldsymbol q_\varepsilon,\varepsilon}(\boldsymbol y_{\varepsilon})+\tilde\phi_\varepsilon(\boldsymbol y_\varepsilon)=\mu_\varepsilon\}\cap B_{L}(\boldsymbol 0).
			\end{split}
		\end{equation*}
		as the notations given in Lemma \ref{B4}. We deduce that
		\begin{equation*}
			\begin{split}
				&\frac{s_\varepsilon^2}{\varepsilon^2}\int_{B_L(\boldsymbol 0)}(s_\varepsilon y_1+q_\varepsilon)\left(\boldsymbol1_{\{\tilde\psi_\varepsilon-\frac{W}{2}(s_\varepsilon y_1+q_\varepsilon)^2\ln\frac{1}{\varepsilon}>\mu_\varepsilon\}}-\boldsymbol 1_{\{\tilde V_{\boldsymbol q_\varepsilon,\varepsilon}>\frac{a_\varepsilon}{2\pi}\ln\frac{1}{\varepsilon}\}}\right)\tilde\varphi d\boldsymbol y\\
				&=\frac{s_\varepsilon^2}{\varepsilon^2}\int_0^{2\pi}\int_1^{1+t_\varepsilon(\theta)} q_\varepsilon r \tilde \varphi(r,\theta) dr d\theta+O(\varepsilon)\cdot |t_{\varepsilon,\mathcal W}+t_{\varepsilon,\tilde \phi_\varepsilon}|^{\frac{1}{q'}}\cdot\|\tilde \varphi\|_{L^q(B_L(\boldsymbol 0))}\\
				&=\frac{s_\varepsilon^2}{\varepsilon^2}\int_0^{2\pi}\int_1^{1+t_\varepsilon(\theta)} q_\varepsilon r \tilde \varphi(1,\theta) dr d\theta+\frac{s_\varepsilon^2}{\varepsilon^2}\int_0^{2\pi}\int_1^{1+t_\varepsilon(\theta)} q_\varepsilon t (\tilde \varphi(r,\theta)-\tilde \varphi(1,\theta)) dr d\theta\\
				& \ \ \ +O(\varepsilon)\cdot|t_{\varepsilon,\mathcal W}+t_{\varepsilon,\tilde \phi_\varepsilon}|^{\frac{1}{2}+\frac{1}{p}}\cdot\|\tilde\varphi\|_{W^{1,p'}(B_L(\boldsymbol 0))}\\
				&=I_1+I_2+O_\varepsilon\left(\varepsilon |t_{\varepsilon,\mathcal W}+t_{\varepsilon,\tilde \phi_\varepsilon}|^{\frac{1}{2}+\frac{1}{p}}\right)\cdot \|\tilde\varphi\|_{W^{1,p'}(B_L(\boldsymbol 0))},
			\end{split}
		\end{equation*}
		where we use Sobolev embedding and choose $q=\frac{2p'}{2-p'}$. It follows from Lemma \ref{lem3-10} and Lemma \ref{B4} that
		\begin{equation*}
			\begin{split}
					I_1&=\frac{s_\varepsilon^2}{\varepsilon^2}\int_0^{2\pi}\int_1^{1+t_\varepsilon(\theta)}q_\varepsilon r \tilde \varphi(1,\theta) dr d\theta\\
				&=\frac{2}{q_\varepsilon}\int_0^{2\pi}\left(\tilde\phi_\varepsilon(\boldsymbol y_\varepsilon(\theta))+O\left(t_{\varepsilon,\mathcal W}+\|\tilde\phi_\varepsilon\|_{L^\infty(B_L(\boldsymbol 0))}^2+\varepsilon^2|\ln\varepsilon|\right)\right) \tilde \varphi(1,\theta) d\theta\\
				&=\frac{2}{q_\varepsilon}\int_{|\boldsymbol y|=1} \tilde\phi_\varepsilon\tilde \varphi d\theta+\frac{2}{q_\varepsilon}\int_0^{2\pi} (\tilde\phi_\varepsilon(\boldsymbol y_\varepsilon(\theta))-\tilde\phi_\varepsilon(1,\theta))\tilde \varphi d\theta\\
				&\quad+O\left(t_{\varepsilon,\mathcal W}+o_\varepsilon(1)\cdot\|\tilde\phi_\varepsilon\|_{L^\infty(B_L(\boldsymbol 0))}+\varepsilon^2|\ln\varepsilon|\right) \cdot\int_{|\boldsymbol y|=1}\tilde \varphi(1,\theta) d\theta\\
				&=\frac{2}{q_\varepsilon}\int_{|\boldsymbol y|=1} \tilde\phi_\varepsilon\tilde \varphi d\theta +\frac{2}{q_\varepsilon}\int_0^{2\pi}\int_{1}^{1+t_\varepsilon(\theta)} \frac{\partial \tilde\phi_\varepsilon(r,\theta)}{\partial r}\tilde \varphi dr d\theta\\
				&\quad+O\left(t_{\varepsilon,\mathcal W}+o_\varepsilon(1)\cdot\|\tilde\phi_\varepsilon\|_{L^\infty(B_L(\boldsymbol 0))}+\varepsilon^2|\ln\varepsilon|\right)\cdot\int_{|\boldsymbol y|=1}\tilde \varphi(1,\theta) d\theta\\
				&=\frac{2}{q_\varepsilon}\int_{|\boldsymbol y|=1} \tilde\phi_\varepsilon\tilde \varphi d\theta\\
				&\quad+O\left(t_{\varepsilon,\mathcal W}+o_\varepsilon(1)\cdot\|\tilde\phi_\varepsilon\|_{L^\infty(B_L(\boldsymbol 0))}+o_\varepsilon(1)\cdot\|\nabla\tilde\phi_\varepsilon\|_{L^p(B_L(\boldsymbol 0))}+\varepsilon^2|\ln\varepsilon|\right)\cdot||\tilde \varphi||_{W^{1,p'}(B_L(\boldsymbol 0))}.
			\end{split}
		\end{equation*}
		Using Lemma \ref{B4}, we can  also deduce that
		\begin{equation*}
			\begin{split}
				I_2&=\frac{s_\varepsilon^2}{\varepsilon^2}\int_0^{2\pi}\int_1^{1+t_\varepsilon(\theta)} q_\varepsilon r (\tilde \varphi(r,\theta)-\tilde \varphi(1,\theta)) dr d\theta\\
				&=\frac{s^2}{\varepsilon^2}\int_0^{2\pi}\int_1^{1+t_\varepsilon(\theta)}q_\varepsilon r \int_1^t  \frac{\partial \tilde \varphi(r,\theta)}{\partial r} drdt d\theta\\
				&\leq \frac{s_\varepsilon^2}{\varepsilon^2}\int_0^{2\pi}q_\varepsilon|t_{\varepsilon,\mathcal W}(\theta)+t_{\varepsilon,\tilde \phi_\varepsilon}(\theta)|\int_1^{1+t_\varepsilon(\theta)} \left|\frac{\partial \tilde \varphi(r,\theta)}{\partial r} \right|dr d\theta\\
				&=O\left(t_{\varepsilon,\mathcal W}+\|\tilde\phi_\varepsilon\|_{L^\infty(B_L(\boldsymbol 0))}+\varepsilon^2|\ln\varepsilon|\right)\cdot\int_0^{2\pi}\int_1^{1+t_\varepsilon(\theta)}\left|\frac{\partial \tilde \varphi(r,\theta)}{\partial r} \right|dr d\theta\\
				&=o_\varepsilon(1)\cdot O\left(t_{\varepsilon,\mathcal W}+||\tilde\phi_\varepsilon||_{L^\infty}+\varepsilon^2|\ln\varepsilon|\right)\cdot||\tilde \varphi||_{W^{1,p'}(B_L(\boldsymbol 0))}.
			\end{split}
		\end{equation*}
		Combining above estimates, we arrive at
		\begin{equation*}
			\begin{array}{ll}
				\langle \mathcal{R_\varepsilon}(\phi_\varepsilon), \varphi\rangle
				&=O\left(t_{\varepsilon,\mathcal W}+\varepsilon |t_{\varepsilon,\mathcal W}+t_{\varepsilon,\tilde \phi_\varepsilon}|^{\frac{1}{2}+\frac{1}{p}}+\varepsilon^2|\ln\varepsilon|\right)\cdot||\tilde \varphi||_{W^{1,p'}(B_L(\boldsymbol 0))}\\
				&\quad+o_\varepsilon(1)\cdot\left(\|\tilde\phi_\varepsilon\|_{L^\infty(B_L(\boldsymbol 0))}+\|\nabla\tilde\phi_\varepsilon\|_{L^p(B_L(\boldsymbol 0))}\right)\cdot||\tilde \varphi||_{W^{1,p'}(B_L(\boldsymbol 0))},
			\end{array}
		\end{equation*}
		which implies
		\begin{equation*}
			\begin{array}{ll}
				\varepsilon^{1-\frac{2}{p}} ||\mathcal{R_\varepsilon}(\phi_\varepsilon)||_{W^{-1,p}(B_{Ls_\varepsilon}(\boldsymbol q_\varepsilon))}
				&=O\left(t_{\varepsilon,\mathcal W}+\varepsilon |t_{\varepsilon,\mathcal W}+t_{\varepsilon,\tilde \phi_\varepsilon}|^{\frac{1}{2}+\frac{1}{p}}+\varepsilon^2|\ln\varepsilon|\right)\\
				&\quad+o_\varepsilon(1)\cdot\left(\|\phi_\varepsilon\|_*+\varepsilon^{1-\frac{2}{p}}\|\nabla\phi_\varepsilon\|_{L^p(B_L(\boldsymbol 0))}\right).
			\end{array}
		\end{equation*}
		Thus from the above discussion, we finally obtain
\[	
||\phi_\varepsilon||_*+\varepsilon^{1-\frac{2}{p}} ||\nabla \phi_\varepsilon||_{L^{p}(B_{Ls_\varepsilon}(\boldsymbol q_\varepsilon))}
			=O\left(t_{\varepsilon,\mathcal W}+\varepsilon |t_{\varepsilon,\mathcal W}+t_{\varepsilon,\tilde \phi_\varepsilon}|^{\frac{1}{2}+\frac{1}{p}}+\varepsilon^2|\ln\varepsilon|\right),
\]
		which is exactly the result we desired.
	\end{proof}
    
    \bigskip
    
    By the discussion before Lemma \ref{lem3-12} and the fact $\|\phi_\varepsilon\|_{L^{\infty}(B_{Ls_\varepsilon}(\boldsymbol q_\varepsilon))}\le\|\phi_\varepsilon\|_*$, it holds
	$$t_{\varepsilon,\tilde \phi_\varepsilon}=O\left(t_{\varepsilon,\mathcal W}+\varepsilon |t_{\varepsilon,\mathcal W}+t_{\varepsilon,\tilde \phi_\varepsilon}|^{\frac{1}{2}+\frac{1}{p}}+\varepsilon^2|\ln\varepsilon|\right).$$
	Hence we can improve the estimate for $\tilde{\mathbf\Gamma}_{\varepsilon,\tilde \phi_\varepsilon}$ in Lemma \ref{B4} as follows.
	\begin{lemma}\label{lem3-13}
		The set
		$$\tilde{\mathbf\Gamma}_{\varepsilon,\tilde \phi_\varepsilon}:=\left\{\boldsymbol y \ | \ \psi_\varepsilon(s_\varepsilon\boldsymbol y+\boldsymbol q_\varepsilon)-\frac{W}{2}(s_\varepsilon y_1+ q_\varepsilon)^2\ln\frac{1}{\varepsilon}\cdot \mathbf{e}_1=\mu_\varepsilon\right\}$$
		is a continuous closed convex curve in $\mathbb{R}^2$, and it holds
		\begin{equation*}
			\begin{split}
				\tilde{\mathbf\Gamma}_{\varepsilon,\tilde \phi_\varepsilon}&=(1+t_\varepsilon(\theta)+t_{\varepsilon,\tilde \phi_\varepsilon}(\theta))(\cos\theta,\sin\theta)\\
				&=(\cos\theta,\sin\theta)+O\left(t_{\varepsilon,\mathcal W}+\varepsilon |t_{\varepsilon,\mathcal W}+t_{\varepsilon,\tilde \phi_\varepsilon}|^{\frac{1}{2}+\frac{1}{p}}+\varepsilon^2|\ln\varepsilon|\right)
			\end{split}
		\end{equation*}
		for each $\theta\in (0,2\pi]$ with $p\in(2,+\infty]$.
	\end{lemma}

\bigskip
	
	Using a bootstrap method, we can further improve the estimate for $\phi_\varepsilon$ and $|A_\varepsilon\bigtriangleup B_{s_\varepsilon}(\boldsymbol q_\varepsilon)|$ to our desired level. According to the assumption on $s_\varepsilon$ in \eqref{gradient}, we should always remember the facts that $s_\varepsilon=O(\varepsilon)$, and $s_\varepsilon \mathcal N_\varepsilon$ is near a constant.
	\begin{lemma}\label{lem3-14}
	    One has
	    $$t_{\varepsilon,\mathcal W}(\theta)=O(\varepsilon^3|\ln\varepsilon|) \quad \mathrm{and} \quad t_{\varepsilon,\tilde \phi_\varepsilon}(\theta)=O(\varepsilon^2|\ln\varepsilon|).$$
		As a result, it holds
		$$|A_\varepsilon\bigtriangleup B_{s_\varepsilon}(\boldsymbol q_\varepsilon)|=O(\varepsilon^4|\ln\varepsilon|).$$
	\end{lemma}
	\begin{proof}
		At the first stage, we have $t_{\varepsilon,\mathcal W}(\theta)=O(\varepsilon|\ln\varepsilon|)$ in hand since $s_\varepsilon\mathcal W_\varepsilon(s_\varepsilon)=O(\varepsilon|\ln\varepsilon|)$ and $s_\varepsilon\mathcal N_\varepsilon$ is bounded. Hence from \eqref{3-11}, we can deduce
		\begin{equation*}
			||\phi_\varepsilon||_*+\varepsilon^{1-\frac{2}{p}} ||\nabla \phi_\varepsilon||_{L^{p}(B_{Ls_\varepsilon}(\boldsymbol q_\varepsilon))}=O(\varepsilon|\ln\varepsilon|),
		\end{equation*}
		and hence $t_{\varepsilon,\tilde \phi_\varepsilon}(\theta)=O(\varepsilon|\ln\varepsilon|)$. Note that we set $s_\varepsilon^*=(\frac{\kappa}{q_\varepsilon\pi})^{1/2}\varepsilon$ in Appendix \ref{appC}. By the circulation constraint \eqref{3-2} and Lemma \ref{B3}, we have
		\begin{equation*}
			\begin{split}
				\frac{s_\varepsilon^{*2}}{\varepsilon^2}\cdot q_\varepsilon\pi&=\frac{s_\varepsilon^2}{2\varepsilon^2}\int_0^{2\pi}(q_\varepsilon+s_\varepsilon\cos\theta)\left(1+t_{\varepsilon,\mathcal W}(\theta)+t_{\varepsilon,\tilde \phi_\varepsilon}(\theta)+\varepsilon^2|\ln\varepsilon|\right)^2 d\theta\\
				&=\frac{s_\varepsilon^2}{\varepsilon^2}\cdot q_\varepsilon\pi+O(|t_{\varepsilon,\mathcal W}(\theta)+t_{\varepsilon,\tilde \phi_\varepsilon}(\theta)|+O(\varepsilon^2|\ln\varepsilon|)).
			\end{split}
		\end{equation*}
		Hence we can derive the difference of $s_\varepsilon^*$ and $s_\varepsilon$ as
		\begin{equation}\label{sestmate}
			\begin{split}
		|s_\varepsilon^*-s_\varepsilon|&=O\left(\varepsilon |t_{\varepsilon,\mathcal W}(\theta)+t_{\varepsilon,\tilde \phi_\varepsilon}(\theta)|+\varepsilon^3|\ln\varepsilon|\right)\\
		&=O(\varepsilon^2|\ln\varepsilon|)
		\end{split}
		\end{equation}
		and
		$$\mathbf{er}_{\varepsilon}:=|A_\varepsilon\bigtriangleup B_{s_\varepsilon^*}(\boldsymbol q_\varepsilon)|=O(\varepsilon^3|\ln\varepsilon|)$$
		by the estimate for $\partial A_\varepsilon$ in Lemma \ref{B3}. Then according to Lemma \ref{C5}, it holds
		\begin{equation}\label{West}
			\begin{array}{ll}
			    \mathcal W_\varepsilon(s_\varepsilon)&=Wq_\varepsilon\ln\frac{1}{\varepsilon}-\frac{\kappa}{4\pi}\ln\frac{8q_\varepsilon}{s_\varepsilon^*}+\frac{\kappa}{16\pi}
	+O_\varepsilon\left(\varepsilon^2|\ln\varepsilon|\right)\\
				&=O(\varepsilon|\ln\varepsilon|).
			\end{array}
		\end{equation}
		By the relation ship of $t_{\varepsilon,\mathcal W}(\theta)$ and $\mathcal W_\varepsilon(s_\varepsilon)$ in Lemma \ref{B2}, we have improved the estimate for $t_{\varepsilon,\mathcal W}(\theta)$ from $O(\varepsilon|\ln\varepsilon|)$ to $O(\varepsilon^2|\ln\varepsilon|)$.
		
		In the second step, we combine above estimates with \eqref{3-11} to obtain
		\begin{equation*}
			t_{\varepsilon,\tilde \phi_\varepsilon}(\theta)\le s_\varepsilon\mathcal N_\varepsilon\cdot||\phi_\varepsilon||_*=O\left(\varepsilon t_{\varepsilon,\tilde \phi_\varepsilon}^{\frac{1}{2}+\frac{1}{p}}+\varepsilon^2|\ln\varepsilon|\right), \quad \forall \, p\in (2,+\infty].
		\end{equation*}
		Now we claim
		\begin{equation}\label{phiest}
			t_{\varepsilon,\tilde \phi_\varepsilon}(\theta)=O(\varepsilon^2|\ln\varepsilon|).
		\end{equation}
		Suppose not. Then there exists a series $\{\varepsilon_n\}$ tends to $0$, such that $t_{\varepsilon_n,\tilde \phi_\varepsilon}(\theta)>n\varepsilon_n^2|\ln\varepsilon_n|$. Since it holds
		\begin{equation*}
			\begin{split}
				\varepsilon_nt_{\varepsilon_n,\tilde \phi_\varepsilon}^{\frac{1}{2}+\frac{1}{p}}&=\varepsilon_n\left(n\varepsilon_n^2|\ln\varepsilon_n|\right)^{\frac{1}{p}-\frac{1}{2}}\cdot\left(n\varepsilon_n^2|\ln\varepsilon_n|\right)^{\frac{1}{2}-\frac{1}{p}}t_{\varepsilon_n,\tilde \phi_\varepsilon}^{\frac{1}{2}+\frac{1}{p}}\\
				&\le \varepsilon_n\left(n\varepsilon_n^2|\ln\varepsilon_n|\right)^{\frac{1}{p}-\frac{1}{2}}t_{\varepsilon_n,\tilde \phi_\varepsilon},
			\end{split}
		\end{equation*}
		we can let $p>2$ be sufficiently close to $2$ and $\varepsilon_n\left(n\varepsilon_n^2|\ln\varepsilon_n|\right)^{\frac{1}{p}-\frac{1}{2}}=o_{\varepsilon_n}(1)$. According to \eqref{3-11}, we have
		\begin{equation*}
			t_{\varepsilon_n,\tilde \phi_\varepsilon}(\theta)=o_{\varepsilon_n}(1)\cdot t_{\varepsilon_n,\tilde \phi_\varepsilon}(\theta)+O(\varepsilon_n^2|\ln\varepsilon_n|),
		\end{equation*}
		which is a contradiction to our assumption $t_{\varepsilon_n,\tilde \phi_\varepsilon}(\theta)>n\varepsilon_n^2|\ln\varepsilon_n|$, and verifies \eqref{phiest}.
		
		In the last step, we use \eqref{3-11} again, and improve the estimate for $t_{\varepsilon,\tilde \phi_\varepsilon}(\theta)$ to
		\begin{equation*}
			t_{\varepsilon,\tilde \phi_\varepsilon}(\theta)=O\left(\varepsilon^2|\ln\varepsilon|+\varepsilon(\varepsilon^2|\ln\varepsilon|)^{\frac{1}{2}+\frac{1}{p}}\right)=O(\varepsilon^2|\ln\varepsilon|).
		\end{equation*}
		Note that we have obtained $t_{\varepsilon,\mathcal W}(\theta)=O(\varepsilon^2|\ln\varepsilon|)$ after \eqref{West}. Proceeding as the first step, we deduce
		$$|A_\varepsilon\bigtriangleup B_{s_\varepsilon}(\boldsymbol q_\varepsilon)|=O(\varepsilon^4|\ln\varepsilon|)$$
		and $t_{\varepsilon,\mathcal W}(\theta)=O(\varepsilon^3|\ln\varepsilon|)$. Hence the proof is complete.
	\end{proof}

   \bigskip
	
	Remember that $\sigma_\varepsilon=\frac{1}{2}\text{diam}\, A_\varepsilon$ is  defined as the cross-section parameter. Now we can obtain the Kelvin--Hicks formula in Proposition \ref{prop3-2}.
	
	\noindent{\bf Proof of Proposition \ref{prop3-2}:}
	By $|A_\varepsilon\bigtriangleup B_{s_\varepsilon}(\boldsymbol q_\varepsilon)|=O(\varepsilon^4|\ln\varepsilon|)$ obtained in Lemma \ref{lem3-14} and the estimate in Lemma \ref{C5}, we have
	\begin{equation}\label{3-12}
		Wq_\varepsilon\ln\frac{1}{\varepsilon}-\frac{\kappa}{4\pi}\ln\frac{8q_\varepsilon}{s_\varepsilon^*}+\frac{\kappa}{16\pi}=O(\varepsilon^2|\ln\varepsilon|).
	\end{equation}
	Recall that $s_\varepsilon^*=\sqrt{\varepsilon^2\kappa/q_\varepsilon\pi}$ and $\sigma_\varepsilon=\frac{1}{2}\text{diam} A_\varepsilon$, from \eqref{sestmate} we have
	\begin{equation}\label{3-13}
		\begin{split}
			|s_\varepsilon^*-\tau_\varepsilon|&\le|s_\varepsilon^*-s_\varepsilon|+|s_\varepsilon-\tau_\varepsilon|\\
			&=O\left(\varepsilon |t_{\varepsilon,\mathcal W}(\theta)+t_{\varepsilon,\tilde \phi_\varepsilon}(\theta)|+\varepsilon^3|\ln\varepsilon|\right)\\
			&=O(\varepsilon^3|\ln\varepsilon|).
		\end{split}
	\end{equation}
	By combining \eqref{3-12} and \eqref{3-13}, we have verified Proposition \ref{prop3-2}.
	\qed
	
	\bigskip
	
	\subsection{The uniqueness}
	
	To show the uniqueness of $\psi_\varepsilon$ satisfying \eqref{3-1} and \eqref{3-2}, we first refine the estimate for the cross-section $A_\varepsilon$ on the location $\boldsymbol q_\varepsilon=(q_\varepsilon,0)$. Noticing that the value of $s_\varepsilon$ depends on $\varepsilon$ and $q_\varepsilon$ by \eqref{gradient}, we can write $s_\varepsilon$ as $s_\varepsilon(q_\varepsilon)$. The following result is a direct consequence of Lemma \ref{lem3-14} and Proposition \ref{prop3-2}, which is essentially dependent on $q_\varepsilon$ near a non-degenerate zero point of $g(x)$ defined in the following lemma.  
	\begin{lemma}\label{lem3-15}
		For each $\varepsilon\in(0,\varepsilon_0]$ with $\varepsilon_0>0$ sufficiently small, let $q_0$ be the only zero point of
		$$ g(x)=Wx\ln\frac{1}{\varepsilon}-\frac{\kappa}{4\pi}\left(\ln\frac{8x}{s_\varepsilon^*(x)}-\frac{1}{4}\right), \quad x>0,$$
		with $s_\varepsilon^*(x)=(\frac{\kappa}{\pi x})^{1/2}\varepsilon$. Then one has
		\begin{equation*}
			|q_\varepsilon-q_0|=O(\varepsilon^2), \quad \mathrm{and} \quad |s_\varepsilon(q_\varepsilon)-s_\varepsilon(q_0)|=O(\varepsilon^3|\ln\varepsilon|).
		\end{equation*}
	\end{lemma}
	\begin{proof}
		Direct computation yields the estimate for derivative at $q_0$ as $g'(q_0)/|\ln\varepsilon|=W+o_\varepsilon(1)$. By \eqref{3-12}, we have
		\begin{equation*}
			|q_\varepsilon-q_0|=O(\varepsilon^2).
		\end{equation*}
		To derive the estimate for $s_\varepsilon(q_\varepsilon)$, we can use the definition $s_\varepsilon^*(x)$ and above estimate for $q_\varepsilon$ to obtain
		$$|s_\varepsilon^*(q_\varepsilon)-s_\varepsilon^*(q_0)|=O(\varepsilon^3).$$
		Since $|s_\varepsilon^*(q_0)-s_\varepsilon(q_0)|=O(\varepsilon^3|\ln\varepsilon|)$ and $|s_\varepsilon^*(q_\varepsilon)-s_\varepsilon(q_\varepsilon)|=O(\varepsilon^3|\ln\varepsilon|)$ the same as in \eqref{3-13}, we then conclude by triangle inequality that
		\begin{equation*}
			|s_\varepsilon(q_\varepsilon)-s_\varepsilon(q_0)|=O(\varepsilon^3|\ln\varepsilon|).
		\end{equation*}
\end{proof}

Using the estimates in previous lemma, $t_{\varepsilon,\mathcal W}(\theta)=O(\varepsilon^3|\ln\varepsilon|)$ from Lemma \ref{lem3-14} and the expansion in Remark \ref{remarkB}, we then have
$$\psi_\varepsilon-\frac{W}{2}x_1^2\ln\frac{1}{\varepsilon}-\mu_\varepsilon=V_{\boldsymbol q_\varepsilon,\varepsilon}+\phi_\varepsilon-\eta(\boldsymbol x)-\frac{a_\varepsilon}{2\pi}\ln\frac{1}{\varepsilon}+O(\varepsilon^2)$$
with
$$\eta(\boldsymbol x)=\frac{3\kappa}{16\pi q_\varepsilon}\ln\frac{1}{\varepsilon}\cdot |\boldsymbol x-\boldsymbol q_\varepsilon|^2+\frac{W}{2}\ln\frac{1}{\varepsilon}\cdot(x_1-q_\varepsilon)^2.$$ 
Moreover, by applying the equation \eqref{phieq} for $\phi_\varepsilon$, we derive a more precise estimate for $\phi_\varepsilon$ as the unique function satisfying
\begin{equation}\label{phiestimate}
	\begin{cases}
		\mathbb L_\varepsilon\phi_\varepsilon=x_1\boldsymbol1_{\{V_{\boldsymbol q_\varepsilon,\varepsilon}-\eta(\boldsymbol x)>\frac{a_\varepsilon}{2\pi}\ln\frac{1}{\varepsilon}\}}-x_1\boldsymbol1_{\{V_{\boldsymbol q_\varepsilon,\varepsilon}>\frac{a_\varepsilon}{2\pi}\ln\frac{1}{\varepsilon}\}}+O(\varepsilon^2),\\
		\partial_1\phi_\varepsilon(\boldsymbol p_\varepsilon)=0,\\
		\phi_\varepsilon=0, & \text{on} \ x_1=0,\\
		\phi_\varepsilon, \ |\nabla\phi_\varepsilon|/x_1\to0, &\text{as} \ |\boldsymbol x |\to \infty,
	\end{cases}
\end{equation}
from which we can obtain the expansion for the Stokes stream function \eqref{exp} by Remark \ref{remarkB}. (the uniqueness of $\phi_\varepsilon$ is derived by the condition $\partial_1\phi_\varepsilon(\boldsymbol p_\varepsilon)=0$, which eliminate the functions in kernel of $\mathbb L_\varepsilon$ as in the proof of Lemma \ref{lem3-11}.)

We also see that the $x_1$-coordinate $q_\varepsilon$ of limit location is unique by the implicit function theorem. However, this is not enough for the uniqueness of $\psi_\varepsilon$ for any $\varepsilon\in (0,\varepsilon_0]$ and $\varepsilon_0>0$ small. In the following, we will use a second order Pohozaev identity to obtain a contradiction of this non-vanishing property for $\mathbf c_{2,\varepsilon}$ in \eqref{exp} and multiplicity of solutions.

\bigskip
	
	To derive the uniqueness of $\psi_\varepsilon$, we suppose on the contrary that there are two different $\psi_\varepsilon^{(1)}$ and $\psi_\varepsilon^{(2)}$ even symmetric respect to $x_1$-axis and solving \eqref{3-1} \eqref{3-2}. Define
	\begin{equation*}
		\Theta_\varepsilon(\boldsymbol x):=\frac{\psi_\varepsilon^{(1)}(\boldsymbol x)-\psi_\varepsilon^{(2)}(\boldsymbol x)}{||\psi_\varepsilon^{(1)}-\psi_\varepsilon^{(2)}||_{L^\infty(\mathbb R^2_+)}}.
	\end{equation*}
	Then, $\Theta_\varepsilon$ satisfies $||\Theta_\varepsilon||_{L^\infty(\mathbb R^2_+)}=1$ and
	\begin{equation*}
		\begin{cases}
			-\varepsilon^2x_1{\Delta^*}\Theta_\varepsilon=f_\varepsilon(\boldsymbol x), & \text{in} \ \mathbb R^2_+,
			\\
			\Theta_\varepsilon=0, & \text{on} \ x_1=0,
			\\
			\Theta_\varepsilon, \ |\nabla\Theta_\varepsilon|/x_1\to0, &\text{as} \ |\boldsymbol x |\to \infty,
		\end{cases}
	\end{equation*}
	where
	\begin{equation*}
		f_\varepsilon(\boldsymbol x)=\frac{x_1\left(\boldsymbol 1_{\{\psi_\varepsilon^{(1)}-\frac{W}{2}x_1^2\ln\frac{1}{\varepsilon}>\mu_\varepsilon^{(1)}\}}-\boldsymbol 1_{\{\psi_\varepsilon^{(2)}-\frac{W}{2}x_1^2\ln\frac{1}{\varepsilon}>\mu_\varepsilon^{(2)}\}}\right)}{\varepsilon^2||\psi_\varepsilon^{(1)}-\psi_\varepsilon^{(2)}||_{L^\infty(\mathbb R^2_+)}}.
	\end{equation*}
	We see that $f_\varepsilon(\boldsymbol x)=0$ in $\mathbb R^2_+\setminus B_{Ls_\varepsilon^{(1)}}(\boldsymbol q_\varepsilon^{(1)})$ due to Lemma \ref{lem3-15}.
	
	\bigskip
	
	At this stage, we are to obtain a series of necessary estimates for $\Theta_\varepsilon$ and $f_\varepsilon$. Then we will use these estimates to a a local Pohozaev identity and obtain a contradiction provided $\psi_\varepsilon^{(1)} \not\equiv \psi_\varepsilon^{(2)}$. For simplicity, we always use $|\,\cdot\,|_\infty$ to denote $||\,\cdot\,||_{L^\infty(\mathbb R^2_+)}$, and abbreviate the parameters $s_\varepsilon^{(1)}$ as $s_\varepsilon$, $q_\varepsilon^{(1)}$ as $q_\varepsilon$.
	\begin{lemma}\label{lem3-16}
		For $p\in (2,\infty]$, $s_\varepsilon^2 f_\varepsilon(s_\varepsilon \boldsymbol y+\boldsymbol q_\varepsilon)$ is bounded in $W^{-1,p}(B_L(\boldsymbol 0))$. As $\varepsilon\to 0$, for all $\tilde\varphi(\boldsymbol y)\in \ C_0^\infty(\mathbb{R}^2)$ it holds
		\begin{equation*}
			\int_{\mathbb{R}^2} s_\varepsilon^2f_\varepsilon(s_\varepsilon\boldsymbol y+\boldsymbol q_\varepsilon)\tilde\varphi d\boldsymbol y= \frac{2}{q_\varepsilon}\int_{|\boldsymbol y|=1} \left(b_\varepsilon\cdot\frac{\partial w}{\partial y_1}+O(\varepsilon)\right) \tilde\varphi,
		\end{equation*}
		where $b_\varepsilon$ is bounded independent of $\varepsilon$, and $w$ is defined by
		\begin{equation*}
			w(\boldsymbol y)=\left\{
			\begin{array}{lll}
				\frac{1}{4}(1-|\boldsymbol y|^2), \ & \mathrm{if} \ |\boldsymbol y|\le 1,\\
				\frac{1}{2}\ln\frac{1}{|\boldsymbol y|}, & \mathrm{if} \ |\boldsymbol y|\ge 1.
			\end{array}
			\right.
		\end{equation*}
	\end{lemma}
	\begin{proof}
		Let
		$$\tilde{\mathbf\Gamma}_\varepsilon^{(i)}:=\left\{\boldsymbol y\,\,\, | \,\, \psi_\varepsilon^{(i)}(s_\varepsilon\boldsymbol y+\boldsymbol q_\varepsilon^{(i)})-\frac{W}{2}(s_\varepsilon y_1+ q_\varepsilon^{(i)})^2\ln\frac{1}{\varepsilon}=\mu_\varepsilon^{(i)}\right\}, \ \ \ i=1,2.$$
		We take
		$$\boldsymbol y_\varepsilon^{(1)}=\left(1+ t_\varepsilon^{(1)}(\theta)\right)(\cos\theta,\sin\theta)\in \tilde{\mathbf \Gamma}_\varepsilon^{(1)}$$
		with $| t_\varepsilon^{(1)}(\theta)|=O(\varepsilon^2|\ln\varepsilon|)$ by Lemma \ref{lem3-14}. Similarly, there is a $t_\varepsilon^{(2)}$ satisfying $| t_\varepsilon^{(2)}(\theta)|=O(\varepsilon^2|\ln\varepsilon|)$ such that
		$$\boldsymbol y_\varepsilon^{(2)}=\left(1+ t_\varepsilon^{(2)}(\theta)\right)(\cos\theta,\sin\theta)\in \tilde{\mathbf \Gamma}_\varepsilon^{(2)}.$$
		We will take $\boldsymbol q_\varepsilon^{(1)}$ and $\boldsymbol q_\varepsilon^{(2)}$ as a same point $\boldsymbol q_\varepsilon=\boldsymbol q_\varepsilon^{(1)}$ in the following. As a cost, this leads to some loss on the estimate of $t_\varepsilon^{(2)}(\theta)$: since $|q_\varepsilon^{(i)}-q_0|=O(\varepsilon^2)$ from Lemma \ref{lem3-15}, at this time we only have
		$$| t_\varepsilon^{(2)}(\theta)|=O(\varepsilon)$$
		by letting $\boldsymbol q_\varepsilon^{(2)}$ coincide with $\boldsymbol q_\varepsilon^{(1)}$.
		
		\bigskip
		
		Using the definition of $\tilde{\mathbf\Gamma}_\varepsilon^{(i)}$ and the estimate
		$$t_{\varepsilon,\mathcal W}+t_{\varepsilon,\tilde\phi_\varepsilon^{(i)}}=O(\varepsilon^2\ln|\varepsilon|)$$
		obtained from Lemma \ref{lem3-14}, we have
		\begin{equation*}
			\begin{split}
				&\quad \psi_\varepsilon^{(1)}\left(s_\varepsilon\boldsymbol y_\varepsilon^{(2)}+\boldsymbol q_\varepsilon\right)-\psi_\varepsilon^{(2)}\left(s_\varepsilon \boldsymbol y_\varepsilon^{(2)}+\boldsymbol q_\varepsilon\right)\\
				&=\psi_\varepsilon^{(1)}\left(s_\varepsilon\boldsymbol y_\varepsilon^{(2)}+\boldsymbol q_\varepsilon\right)-\psi_\varepsilon^{(1)}\left(s_\varepsilon \boldsymbol y_\varepsilon^{(1)}+\boldsymbol q_\varepsilon\right)+\psi_\varepsilon^{(1)}\left(s_\varepsilon\boldsymbol y_\varepsilon^{(1)}+\boldsymbol q_\varepsilon\right)-\psi_\varepsilon^{(2)}\left(s_\varepsilon \boldsymbol y_\varepsilon^{(2)}+\boldsymbol q_\varepsilon\right)\\
				&=\psi_\varepsilon^{(1)}\left(s_\varepsilon\boldsymbol y_\varepsilon^{(2)}+\boldsymbol q_\varepsilon\right)-\psi_\varepsilon^{(1)}\left(s_\varepsilon \boldsymbol y_\varepsilon^{(1)}+\boldsymbol q_\varepsilon\right)-\left(\mu_\varepsilon^{(2)}-\mu_\varepsilon^{(1)}\right)\\
				 &\quad-W\left(s_\varepsilon y_{1,\varepsilon}^{(2)}+q_\varepsilon\right)^2\ln\frac{1}{\varepsilon}+W\left(s_\varepsilon y_{1,\varepsilon}^{(1)}+q_\varepsilon\right)^2\ln\frac{1}{\varepsilon}\\
				&= (-s_\varepsilon\mathcal N_\varepsilon+O(\varepsilon^2|\ln\varepsilon|))\left( t_\varepsilon^{(2)}(\theta)- t_\varepsilon^{(1)}(\theta)\right)-\left(\mu_\varepsilon^{(2)}-\mu_\varepsilon^{(1)}\right),
			\end{split}
		\end{equation*}
		with
		$$\mathcal N_\varepsilon=\frac{s_\varepsilon}{2\varepsilon^2}\cdot q_\varepsilon^2$$
		the value of $|\nabla V_{\boldsymbol q_\varepsilon,\varepsilon}(\boldsymbol x)|$ at $|\boldsymbol x-\boldsymbol q_\varepsilon|=s_\varepsilon$. Thus it holds
		\begin{equation}\label{3-14}
			\begin{split}
				t_\varepsilon^{(2)}(\theta)- t_\varepsilon^{(1)}(\theta)&=(-s_\varepsilon\mathcal N_\varepsilon+O(\varepsilon^2|\ln\varepsilon|))\\
				&\quad\times\left(\psi_{1,\varepsilon}^{(1)}\left(s_\varepsilon\boldsymbol y_\varepsilon^{(2)}+\boldsymbol q_\varepsilon\right)-\psi_{1,\varepsilon}^{(2)}\left(s_\varepsilon \boldsymbol y_\varepsilon^{(1)}+\boldsymbol q_\varepsilon\right)-\left(\mu_\varepsilon^{(2)}-\mu_\varepsilon^{(2)}\right)\right).
			\end{split}
		\end{equation}
		On the other hand, the circulation constraint \eqref{3-2} yields
		\begin{equation*}
			\begin{split}
				\kappa&=\frac{s_\varepsilon^2}{2\varepsilon^2}\int_0^{2\pi}q_\varepsilon\left(1+t_\varepsilon^{(1)}(\theta)\right)^2 d\theta+\frac{s_\varepsilon^3}{3\varepsilon^2}\int_0^{2\pi}\left(1+ t_\varepsilon^{(1)}(\theta)\right)^3\cos\theta d\theta\\
				&=\frac{s_\varepsilon^2}{2\varepsilon}\int_0^{2\pi}q_\varepsilon\left(1+t_\varepsilon^{(2)}(\theta)\right)^2 d\theta+\frac{s_\varepsilon^3}{3\varepsilon^2}\int_0^{2\pi}\left(1+t_\varepsilon^{(2)}(\theta)\right)^3\cos\theta d\theta,
			\end{split}
		\end{equation*}
		and hence
		$$\int_0^{2\pi} q_\varepsilon\left( t_\varepsilon^{(2)}(\theta)- t_\varepsilon^{(1)}(\theta)\right) \left(1+\frac{1}{2} t_\varepsilon^{(1)}(\theta)+\frac{1}{2}t_\varepsilon^{(2)}(\theta)+O(\varepsilon)\right) d\theta=0.$$
		It follows that
		\begin{equation*}
			\begin{split}
				&\int_0^{2\pi}(s_\varepsilon\mathcal N_\varepsilon+O(\varepsilon^2|\ln\varepsilon|))\left(\psi_\varepsilon^{(1)}\left(s_\varepsilon\boldsymbol y_\varepsilon^{(2)}+\boldsymbol q_\varepsilon\right)-\psi_\varepsilon^{(2)}\left(s_\varepsilon \boldsymbol y_\varepsilon^{(2)}+\boldsymbol q_\varepsilon\right)\right) \\
				&\quad\times\left(2+t_\varepsilon^{(1)}(\theta)+ t_\varepsilon^{(2)}(\theta)+O(\varepsilon)\right) d\theta\\
				&=\left(\mu_\varepsilon^{(1)}-\mu_\varepsilon^{(2)}\right)\int_0^{2\pi}(s_\varepsilon\mathcal N_\varepsilon+O(\varepsilon^2|\ln\varepsilon|))\left(2+ t_\varepsilon^{(1)}(\theta)+t_\varepsilon^{(2)}(\theta)+O(\varepsilon)\right) d\theta,
			\end{split}
		\end{equation*}
		which implies 	
		\begin{equation*}
			\frac{|\mu_\varepsilon^{(1)}-\mu_\varepsilon^{(2)}|}{|\psi_\varepsilon^{(1)}-\psi_\varepsilon^{(2)}|_\infty}=O_\varepsilon(1),
		\end{equation*}
		and
		\begin{equation*}
			 \frac{|t_\varepsilon^{(2)}(\theta)-t_\varepsilon^{(1)}(\theta)|}{|\psi_\varepsilon^{(1)}-\psi_\varepsilon^{(2)}|_\infty}=O_\varepsilon(1)
		\end{equation*}
		by \eqref{3-14}.
		
		We then define the normalized difference of $\psi_\varepsilon^{(i)}-\mu_\varepsilon^{(i)}$ as
		\begin{equation*}
			 \Theta_{\varepsilon,\mu}:=\frac{\left(\psi_\varepsilon^{(1)}-\mu_\varepsilon^{(1)}\right)-\left(\psi_\varepsilon^{(2)}-\mu_\varepsilon^{(2)}\right)}{|\psi_\varepsilon^{(1)}-\psi_\varepsilon^{(2)}|_\infty}.
		\end{equation*}
		Recall that for a general function $v$, we denote $\tilde v(\boldsymbol y)=v(s_\varepsilon\boldsymbol y+\boldsymbol q_\varepsilon)$, and $D_s=\{\boldsymbol y \ | \ s_\varepsilon\boldsymbol y+\boldsymbol q_\varepsilon\in\mathbb R^2_+\}$. $\tilde\Theta_{\varepsilon,\mu}$ will satisfy the equation
		\begin{equation*}
			-\text{div}\left(\frac{\nabla\tilde\Theta_{\varepsilon,\mu}}{s_\varepsilon y_1+q_\varepsilon}\right)=\tilde f_\varepsilon(\boldsymbol y), \quad \text{in} \ D_s.
		\end{equation*}
		For any $\tilde\varphi(\boldsymbol y) \in C_0^\infty(B_{L}(\boldsymbol 0))$ and $p'\in [1,2)$, we have
		\begin{equation*}
			\begin{split}
				&\quad \int_{\mathbb{R}^2} s_\varepsilon^2 f_\varepsilon (s_\varepsilon \boldsymbol y+ \boldsymbol q_\varepsilon) \tilde\varphi d\boldsymbol y\\
				&=-\frac{s_\varepsilon^2}{\varepsilon^2|\psi_\varepsilon^{(1)}-\psi_\varepsilon^{(2)}|_\infty}\int_0^{2\pi} \int_{1+t_\varepsilon^{(1)}}^{1+ t_\varepsilon^{(2)}} (q_\varepsilon+t\cos\theta)t\tilde\varphi(t,\theta) dt d\theta\\
				&=-\frac{s_\varepsilon^2}{\varepsilon^2|\psi_\varepsilon^{(1)}-\psi_\varepsilon^{(2)}|_\infty}\int_{|\boldsymbol y|=1} q_\varepsilon \tilde\varphi(\boldsymbol y)( t_\varepsilon^{(2)}(\theta)- t_\varepsilon^{(1)}(\theta)) \left(1+\frac{1}{2} t_\varepsilon^{(1)}(\theta)+\frac{1}{2} t_\varepsilon^{(2)}(\theta)\right) d\boldsymbol y\\
				&\quad -\frac{s_\varepsilon^2}{\varepsilon^2|\psi_\varepsilon^{(1)}-\psi_\varepsilon^{(2)}|_\infty}\int_{|\boldsymbol y|=1} \int_{1+ t_\varepsilon^{(1)}}^{1+t_\varepsilon^{(2)}} (q_\varepsilon+t\cos\theta)t\left[\tilde\varphi(t\boldsymbol y)-\tilde\varphi(\boldsymbol y)\right] dt d\boldsymbol y\\
				&\quad+o_\varepsilon\left(\int_{|\boldsymbol y|=1}  |\tilde\varphi(\boldsymbol y)| d\boldsymbol y\right)\\
				&= -\frac{s_\varepsilon^2(1+o_\varepsilon(1))}{\varepsilon^2|\psi_\varepsilon^{(1)}-\psi_\varepsilon^{(2)}|_\infty}\int_{|\boldsymbol y|=1} \int_{1+ t_\varepsilon^{(1)}}^{1+ t_\varepsilon^{(2)}}\int_0^1 q_\varepsilon t(t-1)\nabla\tilde\varphi((1+\sigma (t-1))\boldsymbol y)\cdot \boldsymbol y d\sigma dt d\boldsymbol y\\
				&\quad+O_\varepsilon\left(\int_{|\boldsymbol y|=1}  |\tilde\varphi(\boldsymbol y)| d\boldsymbol y\right)\\
				&=o_\varepsilon\left(\frac{||\nabla\tilde\varphi||_{L^1(B_2(\boldsymbol 0))}}{|\psi_\varepsilon^{(1)}-\psi_\varepsilon^{(2)}|_\infty}\int_{1+ t_\varepsilon^{(1)}}^{1+ t_\varepsilon^{(2)}} tdt \right)+O_\varepsilon\left(\int_{|\boldsymbol y|=1}  |\tilde\varphi(\boldsymbol y)| d\boldsymbol y\right)\\
				&=O_\varepsilon\left(\int_{|\boldsymbol y|=1}  |\tilde\varphi(\boldsymbol y)| dy+||\nabla \tilde\varphi||_{L^1(B_2(\boldsymbol 0))}\right) \\
				&=O_\varepsilon(||\tilde\varphi||_{W^{1,p'}(B_2(\boldsymbol 0))}).
			\end{split}
		\end{equation*}
		So for $p\in(2,+\infty]$, we obtain
		$$||s_\varepsilon^2f_\varepsilon(s_\varepsilon \boldsymbol y+\boldsymbol q_\varepsilon)||_{W^{-1,p}(B_L(\boldsymbol 0))}=O_\varepsilon(1).$$
		By standard theory for elliptic equations, $\tilde\Theta_{\varepsilon,\mu}(\boldsymbol y)$ is bounded in $W^{1,p}_{\text{loc}}(\mathbb{R}^2)$ for $p\in[2,+\infty)$ and hence in  $C^\alpha_{\text{loc}}(\mathbb{R}^2)$. For further use, we let
		\begin{equation}\label{projection}
			\tilde\Theta_{\varepsilon,\mu}^*(\boldsymbol y):=\tilde\Theta_{\varepsilon,\mu}-b_\varepsilon\cdot\frac{\partial w}{\partial y_1}
		\end{equation}
		with $w$ defined in the statement of lemma, and
		\begin{equation*}
			b_\varepsilon=\left(\int_{B_L(\boldsymbol 0)} \tilde\Theta_{\varepsilon,\mu}(\boldsymbol y)\cdot(-\Delta)\frac{\partial w}{\partial y_1} d\boldsymbol y\right)\left(\int_{B_L(\boldsymbol 0)}\frac{\partial w}{\partial y_1}\cdot(-\Delta)\frac{\partial w}{\partial y_1} d\boldsymbol y\right)^{-1}
		\end{equation*}
		as the projection coefficient bounded independent of $\varepsilon$. Then for any $\tilde\varphi(\boldsymbol y) \in C_0^\infty(B_{L}(\boldsymbol 0))$,  $\tilde\Theta_{\varepsilon,\mu}^*$ satisfies
		\begin{equation}\label{3-15}
			\begin{split}
				&\quad\int_{B_L(\boldsymbol 0)} \frac{1}{s_\varepsilon y_1+q_\varepsilon}\cdot\nabla \tilde\Theta_{\varepsilon,\mu}^*\cdot\nabla\tilde\varphi d\boldsymbol y-\frac{2}{q_\varepsilon}\int_{|\boldsymbol y|=1}\tilde\Theta_{\varepsilon,\mu}^*\tilde\varphi\\
				&=-b_\varepsilon\left(\int_{B_L(\boldsymbol 0)}  \frac{1}{s_\varepsilon y_1+q_\varepsilon}\cdot\nabla\left(\frac{\partial w}{\partial y_1}\right)\cdot\nabla\tilde\varphi d\boldsymbol y-\frac{2}{q_\varepsilon}\int_{|\boldsymbol y|=1}\frac{\partial w}{\partial y_1}\tilde\varphi\right)\\
				&\quad+\left(\int_{B_L(\boldsymbol 0)}  s_\varepsilon^2\tilde f_\varepsilon\tilde\varphi d\boldsymbol y-\frac{2}{q_\varepsilon}\int_{|\boldsymbol y|=1}\tilde\Theta_{\varepsilon,\mu}\tilde\varphi\right)\\
				&=I_1+I_2.
			\end{split}
		\end{equation}
		Since the kernel of
		\begin{equation*}
			\tilde{\mathbb L}^*v=-\Delta v-2v(r,\theta)\boldsymbol{\delta}_{|\boldsymbol y|=1}, \ \ \  \text{in} \ \mathbb R^2
		\end{equation*}
		is spanned by
		$$\left\{\frac{\partial w}{\partial y_1},\frac{\partial w}{\partial y_2}\right\},$$
		we deduce $I_1=O(\varepsilon)\cdot \|\tilde\varphi\|_{W^{1,p'}(B_L(\boldsymbol 0))}$. For the term $I_2$, using \eqref{3-14} and the estimate $| t_\varepsilon^{(2)}(\theta)|=O(\varepsilon)$, we have
		\begin{equation*}
			\begin{split}
				I_2&=-\frac{s_\varepsilon^2}{\varepsilon^2|\psi_\varepsilon^{(1)}-\psi_\varepsilon^{(2)}|_\infty}\int_0^{2\pi} \int_{1+ t_\varepsilon^{(1)}}^{1+ t_\varepsilon^{(2)}} (q_\varepsilon+t\cos\theta)t\tilde\varphi(t,\theta) dt d\theta-\frac{2}{q_\varepsilon}\int_{|\boldsymbol y=1|}\tilde\Theta_{\varepsilon,\mu}\tilde\varphi\\
				&=-\frac{s_\varepsilon^2}{\varepsilon^2|\psi_\varepsilon^{(1)}-\psi_\varepsilon^{(2)}|_\infty}\int_0^{2\pi} \int_{1+ t_\varepsilon^{(1)}}^{1+ t_\varepsilon^{(2)}} (q_\varepsilon+t\cos\theta)t(\tilde\varphi(t,\theta)-\tilde\varphi(1,\theta)) dt d\theta\\
				&\quad -\frac{s_\varepsilon^2}{\varepsilon^2|\psi_\varepsilon^{(1)}-\psi_\varepsilon^{(2)}|_\infty}\int_0^{2\pi} \int_{1+ t_\varepsilon^{(1)}}^{1+ t_\varepsilon^{(2)}} (q_\varepsilon+t\cos\theta)t\tilde\varphi(1,\theta) dt d\theta-\frac{2}{q_\varepsilon}\int_{|\boldsymbol y=1|}\tilde\Theta_{\varepsilon,\mu}\tilde\varphi\\
				&=-\frac{s_\varepsilon^2(1+o_\varepsilon(1))}{\varepsilon^2|\psi_\varepsilon^{(1)}-\psi_\varepsilon^{(2)}|_\infty}\int_{1+ t_\varepsilon^{(1)}}^{1+ t_\varepsilon^{(2)}} q_\varepsilon t(t-1)\nabla\tilde\varphi((1+\sigma (t-1))\boldsymbol y)\cdot \boldsymbol y d\sigma dt d\boldsymbol y\\
				&\quad+\left(\frac{2}{q_\varepsilon}+O(\varepsilon^2|\ln\varepsilon|)\right)\int_{|\boldsymbol y|=1}\tilde\Theta_{\varepsilon,\mu}(1+O_\varepsilon( t_\varepsilon^{(2)}))\tilde\varphi\\
				&\quad-\frac{2}{q_\varepsilon}\int_{|\boldsymbol y|=1}\tilde\Theta_{\varepsilon,\mu}\tilde\varphi+O(\varepsilon)\cdot \|\tilde\varphi\|_{W^{1,p'}(B_L(\boldsymbol 0))}\\
				&=O(\varepsilon)\cdot \|\tilde\varphi\|_{W^{1,p'}(B_L(\boldsymbol 0))}.
			\end{split}
		\end{equation*}
		Actually, we can regard the left hand side of \eqref{3-15} as the weak form of linear operator
		$$\tilde{\mathbb L}_{s_\varepsilon}v=-\text{div}\left(\frac{\nabla v}{s_\varepsilon y_1+q_\varepsilon}\right)-\frac{2}{q_\varepsilon}v(r,\theta)\boldsymbol{\delta}_{|\boldsymbol y|=1}$$
		acting on $\tilde\Theta_{\varepsilon,\mu}^*$. Since both $\tilde\Theta_{\varepsilon,\mu}$ and $\tilde\Theta_{\varepsilon,\mu}^*$ are even with respect to $x_1$-axis, the kernel of $\tilde{\mathbb L}_{s_\varepsilon}$ is then approximated by $\partial w/\partial y_1$. Consequently, if a function $v^*\in W^{-1,p}(B_L(\boldsymbol 0))$ with $p\in(2,+\infty]$ satisfies orthogonality condition
		$$\int_{B_L(\boldsymbol 0)} v^*\cdot(-\Delta)\frac{\partial w}{\partial y_1} d\boldsymbol y=0,$$
		then it holds following local coercive estimate
		$$\|v^*\|_{L^\infty(B_L(\boldsymbol 0))}+\|\nabla v^*\|_{L^p(B_{L}(\boldsymbol 0))}\le C\|\tilde{\mathbb L}_{s_\varepsilon}v^*\|_{W^{-1,p}(B_L(\boldsymbol 0))}, \quad \forall \, p\in(2,+\infty],$$
		which is verified in the proof of Lemma \ref{lem2-2}. Since $\tilde\Theta_{\varepsilon,\mu}^*$ satisfy the orthogonality condition by projection \eqref{projection}, we deduce from the estimates for $I_1,I_2$ that
		\begin{equation*}
			\|\tilde\Theta_{\varepsilon,\mu}^*\|_{L^\infty(B_L(\boldsymbol 0))}+\|\nabla\tilde\Theta_{\varepsilon,\mu}^*\|_{L^p(B_{L}(\boldsymbol 0))}=O(\varepsilon), \quad \forall \, p\in(2,+\infty].
		\end{equation*}
		
		Now we arrive at a conclusion: by the definition of $\tilde\Theta_{\varepsilon,\mu}^*$ in \eqref{projection}, for each $p\in(2,+\infty]$, it holds
		\begin{equation*}
			\tilde\Theta_{\varepsilon,\mu}(\boldsymbol y)=b_\varepsilon\cdot\frac{\partial w}{\partial y_1}+O(\varepsilon), \quad \text{in} \ W^{1,p}(B_{L}(\boldsymbol 0)),
		\end{equation*}
		and for all $\tilde\varphi(\boldsymbol y)\in \ C_0^\infty(\mathbb{R}^2)$, it holds
		\begin{equation*}
			\int_{\mathbb{R}^2} s_\varepsilon^2f_\varepsilon(s_\varepsilon\boldsymbol y+\boldsymbol q_\varepsilon)\tilde\varphi d\boldsymbol y= \frac{2}{q_\varepsilon}\int_{|\boldsymbol y|=1} \left(b_\varepsilon\cdot\frac{\partial w}{\partial y_1}+O(\varepsilon)\right) \tilde\varphi,
		\end{equation*}
		where $b_\varepsilon$ is bounded independent of $\varepsilon$. So we have completed the proof of Lemma \ref{lem3-16}.
	\end{proof}

\bigskip
	
	To make use of the local Pohozaev identity in Appendix \ref{appC}, we will deal with the main part of $\Theta_\varepsilon$ related to $\psi_{\varepsilon,1}$ and let it be
	\begin{equation*}
		\xi_\varepsilon(\boldsymbol x):=\frac{\psi_{1,\varepsilon}^{(1)}(\boldsymbol x)-\psi_{1,\varepsilon}^{(2)}(\boldsymbol x)}{|\psi_\varepsilon^{(1)}-\psi_\varepsilon^{(2)}|_\infty}.
	\end{equation*}
	Then $\xi_\varepsilon$ has the following integral representation
	\begin{equation}\label{3-16}
		\xi_\varepsilon=q_\varepsilon^2\int_{\mathbb R^2_+} G(\boldsymbol x,\boldsymbol x')\cdot x_1'^{-1}f_\varepsilon(\boldsymbol x')d\boldsymbol x'.
	\end{equation}
	By the asymptotic estimate for $f_\varepsilon(s_\varepsilon\boldsymbol y+\boldsymbol q_\varepsilon)$ in Lemma \ref{lem3-16}, it holds
	\begin{equation*}
		\frac{\psi_{2,\varepsilon}^{(1)}(\boldsymbol x)-\psi_{2,\varepsilon}^{(2)}(\boldsymbol x)}{|\psi_\varepsilon^{(1)}-\psi_\varepsilon^{(2)}|_\infty}=\int_{\mathbb R^2_+} H(\boldsymbol x,\boldsymbol x')\cdot x_1'^{-1}f_\varepsilon(\boldsymbol x')d\boldsymbol x'=o_\varepsilon(1).
	\end{equation*}
	So we see that $\xi_\varepsilon$ is indeed the main part in $\Theta_\varepsilon$, and $||\xi_\varepsilon||_{L^\infty(\mathbb R^2_+)}=1-o_\varepsilon(1)$. To derive a contradiction and obtain uniqueness, we only have to show $||\xi_\varepsilon||_{L^\infty(\mathbb R^2_+)}=o_\varepsilon(1)$.

	For the purpose of dealing with boundary terms in the local Pohozaev identity (the left hand side of \eqref{C-1} in Appendix \ref{appC}), we need the following lemma concerning the behavior of $\xi_\varepsilon$ away from $\boldsymbol q_\varepsilon$.
	\begin{lemma}\label{lem3-17}
		It holds
		\begin{equation}\label{3-17}
			\begin{split}
			\xi_\varepsilon(\boldsymbol x)&=\mathbf B_\varepsilon\cdot \frac{s_\varepsilon q_\varepsilon^2}{2\pi}\frac{x_1-q_\varepsilon}{|\boldsymbol x-\boldsymbol q_\varepsilon|^2}+\mathbf B_\varepsilon\cdot \frac{s_\varepsilon q_\varepsilon^2}{2\pi}\frac{x_1+q_\varepsilon}{|\boldsymbol x-\boldsymbol {\bar q_\varepsilon}|^2}+\mathbf B_\varepsilon\cdot \frac{s_\varepsilon q_\varepsilon}{2\pi} \ln \frac{|\boldsymbol x-\boldsymbol {\bar q_\varepsilon}|}{|\boldsymbol x-\boldsymbol q_\varepsilon|}\\
			&\quad+O(\varepsilon^2), \quad \mathrm{in} \ \, C^1(\mathbb R^2_+\setminus B_{\delta_0}(\boldsymbol q_\varepsilon))
			\end{split}
		\end{equation}
		with $\delta_0$ a small positive constant, and
		\begin{equation*}
			\mathbf B_\varepsilon:=\frac{1}{s_\varepsilon}\int_{B_{2s_\varepsilon}(\boldsymbol q_\varepsilon)} (x_1-q_\varepsilon) x_1^{-1}f_\varepsilon(\boldsymbol x) d\boldsymbol x
		\end{equation*}
		bounded independent of $\varepsilon$.
	\end{lemma}
	\begin{proof}
		Since $\xi_\varepsilon$ is symmetric with respect to $x_1$-axis, for $\boldsymbol x\in \mathbb{R}^2_+\setminus B_{\delta_0}(\boldsymbol q_\varepsilon)$ we have
		\begin{equation*}
			\begin{split}
				\xi_\varepsilon(\boldsymbol x)&=\frac{q_\varepsilon^2}{2\pi}\int_{\mathbb{R}^2_+} x_1'^{-1}\ln \left(\frac{|\boldsymbol x-\boldsymbol {\bar x}'|}{|\boldsymbol x-\boldsymbol x'|} \right)f_\varepsilon(\boldsymbol x')d\boldsymbol x'=\frac{q_\varepsilon^2}{2\pi}\int_{B_{Ls_\varepsilon}(\boldsymbol z)} x_1'^{-1}\ln \left(\frac{|\boldsymbol x-\boldsymbol {\bar x}'|}{|\boldsymbol x-\boldsymbol x'|}\right) f_\varepsilon(\boldsymbol x')d\boldsymbol x'\\
				&=\frac{q_\varepsilon}{2\pi}\ln \frac{1}{|\boldsymbol x-\boldsymbol q_\varepsilon|} \int_{B_{Ls_\varepsilon}(\boldsymbol q_\varepsilon)}  f_\varepsilon(\boldsymbol x')d\boldsymbol x'+ \frac{q_\varepsilon^2}{2\pi}\int_{B_{Ls_\varepsilon}(\boldsymbol q_\varepsilon)}x_1'^{-1}\ln \left(\frac{|\boldsymbol x-\boldsymbol z|}{|\boldsymbol x-\boldsymbol x'|}\right) f_\varepsilon(\boldsymbol x')d\boldsymbol x'\\
				&\quad-\frac{q_\varepsilon}{2\pi}\ln \frac{1}{|\boldsymbol x-\boldsymbol {\bar q}_\varepsilon|} \int_{B_{Ls_\varepsilon}(\boldsymbol q_\varepsilon)}  f_\varepsilon(\boldsymbol x')d\boldsymbol x'-\frac{q_\varepsilon^2}{2\pi}\int_{B_{Ls_\varepsilon}(\boldsymbol q_\varepsilon)}x_1'^{-1}\ln \left(\frac{|\boldsymbol x-\boldsymbol {\bar q_\varepsilon}'|}{|\boldsymbol x-\boldsymbol {\bar x}'|} \right)f_\varepsilon(\boldsymbol x')d\boldsymbol x'\\
				&\quad-\frac{q_\varepsilon}{2\pi} \ln \frac{|\boldsymbol x-\boldsymbol {\bar q_\varepsilon}|}{|\boldsymbol x-\boldsymbol q_\varepsilon|} \int_{B_{Ls_\varepsilon}(\boldsymbol q_\varepsilon)}(x_1-q_\varepsilon) x_1^{-1}f_\varepsilon(\boldsymbol x)d\boldsymbol x\\
				&=-\frac{q_\varepsilon^2}{4\pi}\int_{B_{Ls_\varepsilon}(\boldsymbol q_\varepsilon)}x_1^{-1}\ln\left(1+ \frac{2(\boldsymbol x-\boldsymbol q_\varepsilon)\cdot (\boldsymbol q_\varepsilon-\boldsymbol x')}{|\boldsymbol x-\boldsymbol q_\varepsilon|^2} +\frac{ |\boldsymbol q_\varepsilon-\boldsymbol x'|^2}{|\boldsymbol x-\boldsymbol q_\varepsilon|^2}\right)f_\varepsilon(\boldsymbol x')d\boldsymbol x'\\
				&\quad +\frac{q_\varepsilon^2}{4\pi}\int_{B_{Ls_\varepsilon}(\boldsymbol q_\varepsilon)}x_1^{-1}\ln\left(1+ \frac{2(\boldsymbol x-\boldsymbol {\bar q}_\varepsilon)\cdot (\boldsymbol {\bar q}_\varepsilon-\boldsymbol {\bar x}')}{|\boldsymbol x-\boldsymbol {\bar q}_\varepsilon|^2} +\frac{ |\boldsymbol {\bar q}_\varepsilon-\boldsymbol {\bar x}'|^2}{|\boldsymbol x-\boldsymbol q_\varepsilon|^2}\right)f_\varepsilon(\boldsymbol x')d\boldsymbol x'\\
				&\quad-\frac{q_\varepsilon}{2\pi} \ln \frac{|\boldsymbol x-\boldsymbol {\bar q_\varepsilon}|}{|\boldsymbol x-\boldsymbol q_\varepsilon|} \int_{B_{Ls_\varepsilon}(\boldsymbol q_\varepsilon)}(x_1-q_\varepsilon) x_1^{-1}f_\varepsilon(\boldsymbol x)d\boldsymbol x\\
				&=\mathbf B_\varepsilon\cdot \frac{s_\varepsilon q_\varepsilon^2}{2\pi}\frac{x_1-q_\varepsilon}{|\boldsymbol x-\boldsymbol q_\varepsilon|^2}+\mathbf B_\varepsilon\cdot \frac{s_\varepsilon q_\varepsilon^2}{2\pi}\frac{x_1+q_\varepsilon}{|\boldsymbol x-\boldsymbol {\bar q_\varepsilon}|^2}+\mathbf B_\varepsilon\cdot \frac{s_\varepsilon q_\varepsilon}{2\pi} \ln \frac{|\boldsymbol x-\boldsymbol {\bar q_\varepsilon}|}{|\boldsymbol x-\boldsymbol q_\varepsilon|}+O(\varepsilon^2).
			\end{split}
		\end{equation*}
		Moreover, $\mathbf B_\varepsilon$ is bounded independent of $\varepsilon$ since $s_\varepsilon^2 f_\varepsilon(s_\varepsilon \boldsymbol y+\boldsymbol q_\varepsilon)$ is bounded in $W^{-1,p}(B_L(\boldsymbol 0))$ for $p\in (2,\infty]$. In the next step, we can verify the convergence in $C^1(\mathbb R^2_+\setminus B_{\delta_0}(\boldsymbol q_\varepsilon))$ by a same argument.
	\end{proof}

\bigskip
	
	If we apply \eqref{C-1} in Appendix \ref{appC} to $\psi_{1,\varepsilon}^{(1)}$ and $\psi_{1,\varepsilon}^{(2)}$ separately and calculate their difference, we can obtain the following local Pahozaev identity:
	\begin{equation}\label{3-18}
		\begin{split}
			&\quad-\int_{\partial B_\delta(\boldsymbol q_\varepsilon)}\frac{\partial\xi_\varepsilon}{\partial\nu}\frac{\partial\psi_{1,\varepsilon}^{(1)}}{\partial x_1}dS-\int_{\partial B_\delta(\boldsymbol q_\varepsilon)}\frac{\partial \psi_{1,\varepsilon}^{(2)}}{\partial\nu }\frac{\partial\xi_\varepsilon}{\partial x_1}dS+\frac{1}{2}\int_{\partial B_\delta(\boldsymbol q_\varepsilon)}\langle\nabla(\psi_{1,\varepsilon}^{(1)}+\psi_{1,\varepsilon}^{(2)}),\nabla\xi_\varepsilon\rangle\nu_1dS\\
			&=-\frac{q_\varepsilon^2}{\varepsilon^2|\psi_\varepsilon^{(1)}-\psi_\varepsilon^{(2)}|_\infty}\int_{B_\delta(\boldsymbol q_\varepsilon)} \left(\partial_1\psi_{2,\varepsilon}^{(1)}\cdot\boldsymbol 1_{A_\varepsilon^{(1)}}-\partial_1\psi_{2,\varepsilon}^{(2)}\cdot\boldsymbol 1_{A_\varepsilon^{(2)}}\right)d\boldsymbol x.
		\end{split}
	\end{equation}
	The proof of the uniqueness of a vortex ring with small cross-section is based on a careful estimate for each term in \eqref{3-18}. We will show that the left hand side of \eqref{3-18} is of order $O(\varepsilon)$, while the right hand side is equal to the quantity 
	$$\frac{3\kappa}{4q_\varepsilon}\cdot b_\varepsilon s_\varepsilon\ln\left(\frac{1}{s_\varepsilon}\right)=4\pi \mathbf c_{2,\varepsilon}\cdot b_\varepsilon \cdot s_\varepsilon|\ln s_\varepsilon|$$ 
	with $\mathbf c_{2,\varepsilon}$ in the expansion \eqref{exp} and $b_\varepsilon$ defined in Lemma \ref{lem3-16}. Thus we can conclude $b_\varepsilon=O(1/|\ln\varepsilon|)$, so that a contradiction can be derived by $||\xi_\varepsilon||_{L^\infty(\mathbb R^2_+)}=O(1/|\ln\varepsilon|)$ and the discussion before Lemma \ref{lem3-17}.
	
	\bigskip
	
	\noindent{\bf Proof of Proposition \ref{prop3-1}:}
	Using the asymptotic estimate for $\psi_{1,\varepsilon}$ in Lemma \ref{C2} and $\xi_\varepsilon$ in Lemma \ref{lem3-17}, we see that
	\begin{equation}\label{3-19}
		\begin{split}
			&\quad\int_{\partial B_\delta(\boldsymbol q_\varepsilon)}\frac{\partial\xi_\varepsilon}{\partial\nu}\frac{\partial\psi_{1,\varepsilon}^{(1)}}{\partial x_1}dS+\int_{\partial B_\delta(\boldsymbol q_\varepsilon)}\frac{\partial \psi_{1,\varepsilon}^{(2)}}{\partial\nu }\frac{\partial\xi_\varepsilon}{\partial x_1}dS-\frac{1}{2}\int_{\partial B_\delta(\boldsymbol q_\varepsilon)}\langle\nabla(\psi_{1,\varepsilon}^{(1)}+\psi_{1,\varepsilon}^{(2)}),\nabla\xi_\varepsilon\rangle\nu_1dS\\
			&=O(\varepsilon)\cdot \mathbf B_\varepsilon+O(\varepsilon^2).
		\end{split}
	\end{equation}
	To deal with the right hand side of \eqref{3-18}, we write
	\begin{equation*}
		\begin{split}
			&\quad\frac{q_\varepsilon^2}{\varepsilon^2|\psi_\varepsilon^{(1)}-\psi_\varepsilon^{(2)}|_\infty}\int_{B_\delta(\boldsymbol q_\varepsilon)} \left(\partial_1\psi_{2,\varepsilon}^{(1)}\cdot\boldsymbol 1_{A_\varepsilon^{(1)}}-\partial_1\psi_{2,\varepsilon}^{(2)}\cdot\boldsymbol 1_{A_\varepsilon^{(2)}}\right)d\boldsymbol x\\
			&=\frac{q_\varepsilon^2}{\varepsilon^2|\psi_\varepsilon^{(1)}-\psi_\varepsilon^{(2)}|_\infty}\int_{B_\delta(\boldsymbol q_\varepsilon)} \left(\partial_1\psi_{2,\varepsilon}^{(1)}(\boldsymbol 1_{A_\varepsilon^{(1)}}-\boldsymbol 1_{A_\varepsilon^{(2)}})+\boldsymbol 1_{A_\varepsilon^{(2)}}(\partial_1\psi_{2,\varepsilon}^{(1)}-\partial_1\psi_{2,\varepsilon}^{(2)})\right)d\boldsymbol x\\
			&=G_1+G_2,
		\end{split}
	\end{equation*}
	and
	\begin{equation*}
		G_1=\frac{q_\varepsilon^2}{\varepsilon^2}\int_{B_\delta(\boldsymbol q_\varepsilon)} x_1^{-1}f_\varepsilon(\boldsymbol x)\int_{B_\delta(\boldsymbol q_\varepsilon)} \partial_{x_1}H(\boldsymbol x,\boldsymbol x')\cdot\boldsymbol 1_{A_\varepsilon^{(1)}}d\boldsymbol x'd\boldsymbol x=G_{11}+G_{12}+G_{13}+G_{14},
	\end{equation*}
	where
	\begin{equation*}
		G_{11}=\frac{q_\varepsilon^2}{4\pi\varepsilon^2}\cdot\ln\left(\frac{1}{s_\varepsilon}\right)\cdot\int_{B_\delta(\boldsymbol q_\varepsilon)} x_1^{-3/2}f_\varepsilon(\boldsymbol x)\int_{A_\varepsilon^{(1)}}x_1'^{3/2} d\boldsymbol x'd\boldsymbol x,
	\end{equation*}
	\begin{equation*}
		G_{12}=\frac{q_\varepsilon^2}{4\pi\varepsilon^2}\cdot\int_{B_\delta(\boldsymbol q_\varepsilon)} x_1^{-3/2}f_\varepsilon(\boldsymbol x)\int_{A_\varepsilon^{(1)}}x_1'^{3/2} \ln\left(\frac{s_\varepsilon}{|\boldsymbol x-\boldsymbol x'|}\right)d\boldsymbol x'd\boldsymbol x,
	\end{equation*}
	\begin{equation*}
		G_{13}=-\frac{q_\varepsilon^2}{2\pi\varepsilon^2}\cdot\int_{B_\delta(\boldsymbol q_\varepsilon)} x_1^{-1}f_\varepsilon(\boldsymbol x)\int_{A_\varepsilon^{(1)}}\left(x_1^{1/2}x_1'^{3/2}-q_\varepsilon^2\right)\cdot\frac{x_1-x_1'}{|\boldsymbol x-\boldsymbol x'|^2}d\boldsymbol x'd\boldsymbol x,
	\end{equation*}
	and $G_{14}$ a regular term. Using the circulation constraint \eqref{3-2} and Lemma \ref{lem3-16}, we have
	\begin{equation*}
		\begin{split}
			G_{11}&=\frac{q_\varepsilon^2}{4\pi}\cdot\ln\left(\frac{1}{s_\varepsilon}\right)\cdot\int_{B_\delta(\boldsymbol q_\varepsilon)} x_1^{-3/2}f_\varepsilon\cdot \frac{1}{\varepsilon^2}\int_{\Omega_\varepsilon^{(1)}}x_1'\left(q_\varepsilon^{1/2}+O(\varepsilon)\right)d\boldsymbol x'd\boldsymbol x\\
			&=\frac{\kappa q_\varepsilon^2}{4\pi}\cdot \left(q_\varepsilon^{1/2}+O(\varepsilon)\right)\cdot\ln\left(\frac{1}{s_\varepsilon}\right)\cdot\int_{B_\delta(\boldsymbol q_\varepsilon)} x_1^{-3/2}f_\varepsilon(\boldsymbol x)d\boldsymbol x\\
			&=\frac{\kappa q_\varepsilon^2}{4\pi}\cdot \left(q_\varepsilon^{1/2}+O(\varepsilon)\right)\cdot\ln\left(\frac{1}{s_\varepsilon}\right)\cdot\int_{B_\delta(\boldsymbol q_\varepsilon)} f_\varepsilon\cdot\left(q_\varepsilon^{-3/2}-\frac{3}{2q_\varepsilon^{5/2}}\cdot(x_1-q_\varepsilon)+O(\varepsilon^2)\right)d\boldsymbol x\\
			&=\frac{\kappa q_\varepsilon^2}{4\pi}\cdot \left(q_\varepsilon^{1/2}+O(\varepsilon)\right)\cdot\ln\left(\frac{1}{s_\varepsilon}\right)\cdot\int_{\mathbb R^2}\left(-\frac{3}{2q_\varepsilon^{5/2}}\cdot s_\varepsilon y_1+O(\varepsilon^2)\right)s_\varepsilon^2f_\varepsilon(s_\varepsilon\boldsymbol y+\boldsymbol q_\varepsilon)d\boldsymbol y\\
			&=\frac{\kappa q_\varepsilon^2}{4\pi}\cdot \left(q_\varepsilon^{1/2}+O(\varepsilon)\right)\cdot\ln\left(\frac{1}{s_\varepsilon}\right)\cdot\int_{|\boldsymbol y|=1}\left(-\frac{3}{2q_\varepsilon^{5/2}}\cdot s_\varepsilon y_1+O(\varepsilon^2)\right)\left(b_\varepsilon\cdot\frac{y_1}{q_\varepsilon|\boldsymbol y|^2}+O(\varepsilon)\right)d\boldsymbol y\\
			&=-\frac{3\kappa}{8q_\varepsilon}\cdot b_\varepsilon s_\varepsilon\ln\left(\frac{1}{s_\varepsilon}\right)+O(\varepsilon).
		\end{split}
	\end{equation*}
	For the term $G_{12}$, it holds
	\begin{equation*}
		\begin{split}
			G_{12}&=\frac{q_\varepsilon^2}{4\pi\varepsilon^2}\int_{B_\delta(\boldsymbol q_\varepsilon)} \left(q_\varepsilon^{-3/2}+O(\varepsilon)\right)f_\varepsilon\int_{A_\varepsilon^{(1)}}\left(q_\varepsilon^{3/2}+O(\varepsilon)\right) \ln\left(\frac{s_\varepsilon}{|\boldsymbol x-\boldsymbol x'|}\right)d\boldsymbol x'd\boldsymbol x\\
			&=\frac{q_\varepsilon^2}{4\pi\varepsilon^2}\int_{B_\delta(\boldsymbol q_\varepsilon)} f_\varepsilon\int_{B_\delta(\boldsymbol q_\varepsilon)} \ln\left(\frac{s_\varepsilon}{|\boldsymbol x-\boldsymbol x'|}\right)d\boldsymbol x'd\boldsymbol x+O(\varepsilon)\\
			&=\frac{q_\varepsilon^2s_\varepsilon^2}{4\pi\varepsilon^2}\int_{|\boldsymbol y|=1}\left(b_\varepsilon\cdot\frac{y_1}{q_\varepsilon|\boldsymbol y|^2}+O(\varepsilon)\right)\left(\int_{B_1(\boldsymbol 0)} \ln\left(\frac{1}{|\boldsymbol y-\boldsymbol y'|}\right)d\boldsymbol y'\right)+O(\varepsilon)\\
			&=O(\varepsilon),
		\end{split}
	\end{equation*}
	where we have used the integral formula of Rankine vortex
	\begin{equation*}
		\frac{1}{2\pi}\int_{B_1(\boldsymbol 0)} \ln\left(\frac{1}{|\boldsymbol y-\boldsymbol y'|}\right)d\boldsymbol y'=\left\{
		\begin{array}{lll}
			\frac{1}{4}(1-|\boldsymbol y|^2), \    &\mathrm{if} \ |\boldsymbol y|\le 1,\\
			\frac{1}{2}\ln\frac{1}{|\boldsymbol y|}, &\mathrm{if} \ |\boldsymbol y|\ge 1.
		\end{array}
		\right.
	\end{equation*}
	Similarly, for $G_{13}$ we have
	\begin{equation*}
		\begin{split}
			G_{13}&=-\frac{q_\varepsilon^2}{4\pi\varepsilon^2}\int_{B_\delta(\boldsymbol q_\varepsilon)}  f_\varepsilon\int_{A_\varepsilon^{(1)}}\left((x_1-q_\varepsilon)+3(x_1'-q_\varepsilon)\right) \cdot\frac{x_1-x_1'}{|\boldsymbol x-\boldsymbol x'|^2}d\boldsymbol x'd\boldsymbol x+O(\varepsilon)\\
			&=-\frac{q_\varepsilon^2}{4\pi\varepsilon^2}\int_{B_\delta(\boldsymbol q_\varepsilon)}  f_\varepsilon\int_{B_\delta(\boldsymbol q_\varepsilon)}\left((x_1-q_\varepsilon)+3(x_1'-q_\varepsilon)\right) \cdot\frac{x_1-x_1'}{|\boldsymbol x-\boldsymbol x'|^2}d\boldsymbol x'd\boldsymbol x+O(\varepsilon).
		\end{split}
	\end{equation*}
	Notice that
	\begin{equation*}
		g(\boldsymbol y)=\int_{B_1(\boldsymbol 0)}\left(y_1+3y'_1\right) \cdot\frac{y_1-y_1'}{|\boldsymbol y-\boldsymbol y'|^2}d\boldsymbol y'
	\end{equation*}
	is a bounded function even symmetric with respect to $y_1=0$. While $\partial w/\partial y_1$ is odd symmetric with respect to $y_1=0$. Hence it holds
	\begin{equation*}
		G_{13}=-\frac{q_\varepsilon^2s_\varepsilon^2}{4\pi\varepsilon^2}\int_{|\boldsymbol y|=1}\left(\frac{2}{q_\varepsilon}\cdot b_\varepsilon\cdot\frac{\partial w}{\partial y_1}+O(\varepsilon)\right)g(\boldsymbol y)+O(\varepsilon)=O(\varepsilon).
	\end{equation*}
	For the regular term $G_{14}$, it is easy to verify that $G_{14}=O(\varepsilon)$. Summarizing all the estimates above, we get
	\begin{equation}\label{3-20}
		G_1=-\frac{3\kappa}{8q_\varepsilon}\cdot b_\varepsilon s_\varepsilon\ln\left(\frac{1}{\varepsilon}\right)+O(\varepsilon).
	\end{equation}
	
	Then we turn to deal with $G_2$. Using Fubini's theorem, we obtain
	\begin{equation*}
		\begin{split}
			 G_2&=\frac{q_\varepsilon^2}{\varepsilon^4|\psi_\varepsilon^{(1)}-\psi_\varepsilon^{(2)}|_\infty}\int_{A_\varepsilon^{(2)}}\left(\int_{A_\varepsilon^{(1)}}\partial_{x_1}H(\boldsymbol x,\boldsymbol x')d\boldsymbol x'-\int_{A_\varepsilon^{(2)}}\partial_{x_1}H(\boldsymbol x,\boldsymbol x')d\boldsymbol x'\right)d\boldsymbol x\\
			&=\frac{q_\varepsilon^2}{\varepsilon^4|\psi_\varepsilon^{(1)}-\psi_\varepsilon^{(2)}|_\infty}\int_{B_\delta(\boldsymbol q_\varepsilon)} \left(\boldsymbol 1_{A_\varepsilon^{(1)}}-\boldsymbol 1_{A_\varepsilon^{(2)}}\right)\int_{A_\varepsilon^{(2)}}\partial_{x_1}H(\boldsymbol x,\boldsymbol x')d\boldsymbol x'd\boldsymbol x\\
			&=\frac{q_\varepsilon^2}{\varepsilon^4|\psi_\varepsilon^{(1)}-\psi_\varepsilon^{(2)}|_\infty}\int_{B_\delta(\boldsymbol q_\varepsilon)} \left(\boldsymbol 1_{A_\varepsilon^{(1)}}-\boldsymbol 1_{A_\varepsilon^{(2)}}\right)\partial_1\psi_{2,\varepsilon}^{(2)}d\boldsymbol x.
		\end{split}
	\end{equation*}
	Since $G_2$ has a very similar formulation with $G_1$ (where $\partial_1\psi_{2,\varepsilon}^{(2)}$ takes the place of $\partial_1\psi_{2,\varepsilon}^{(1)}$), we claim
	\begin{equation}\label{3-21}
		G_2=-\frac{3\kappa}{8q_\varepsilon}\cdot b_\varepsilon s_\varepsilon\ln\left(\frac{1}{\varepsilon}\right)+O(\varepsilon).
	\end{equation}
	
	Now from \eqref{3-19} \eqref{3-20} \eqref{3-21}, we have
	\begin{equation}\label{3-22}
		\frac{3\kappa}{4q_\varepsilon}\cdot b_\varepsilon s_\varepsilon\ln\left(\frac{1}{\varepsilon}\right)=O(\varepsilon).
	\end{equation}
	Since $q_\varepsilon$ is near $q_0>0$ defined in Lemma \ref{lem3-15}, and $s_\varepsilon=O(\varepsilon)$, we can derive from \eqref{3-22} that
	\begin{equation*}
		b_\varepsilon =O\left(\frac{1}{|\ln\varepsilon|}\right).
	\end{equation*}
	According to Lemma \ref{lem3-16}, we can also use the fact that for fixed $\boldsymbol y\in\mathbb R^2$ it holds
	$$\frac{1}{2\pi}\ln\left(\frac{1}{|\boldsymbol y-\cdot|}\right)\in W_{\text{loc}}^{1,p'}(\mathbb R^2),\quad \forall \, p'\in [1,2),$$
	and deduce
	\begin{equation*}
		\begin{split}
			\tilde\xi_\varepsilon(\boldsymbol y)&=\frac{q_\varepsilon}{2\pi}\int_{\mathbb R^2_+} \ln\left(\frac{1}{s_\varepsilon|\boldsymbol y-\boldsymbol y'|}\right)\cdot \left(1-\frac{s_\varepsilon y_1'}{q_\varepsilon}\right)s_\varepsilon^2f_\varepsilon(s_\varepsilon\boldsymbol y'+\boldsymbol q_\varepsilon)d\boldsymbol y'+O\left(\frac{1}{|\ln\varepsilon|}\right)\\
			&=\frac{1}{\pi}\int_{|\boldsymbol y'=1|} \ln\left(\frac{1}{|\boldsymbol y-\boldsymbol y'|}\right)\cdot \left(1-\frac{s_\varepsilon y_1'}{q_\varepsilon}\right)\left(b_\varepsilon\cdot\frac{\partial w(\boldsymbol y')}{\partial y_1}+O(\varepsilon)\right)\\
			&\quad+\frac{1}{\pi}\ln\left(\frac{1}{s_\varepsilon}\right)\cdot\int_{|\boldsymbol y'|=1}\left(b_\varepsilon\cdot\frac{\partial w(\boldsymbol y')}{\partial y_1}+O(\varepsilon)\right)+O\left(\frac{1}{|\ln\varepsilon|}\right)\\
			&=O\left(\frac{1}{|\ln\varepsilon|}\right).
		\end{split}
	\end{equation*}
	Thus we conclude $||\xi_\varepsilon||_{L^\infty(\mathbb R^2_+)}=O(1/|\ln\varepsilon|)$, which is a contradiction to $||\xi_\varepsilon||_{L^\infty(\mathbb R^2_+)}=1-o_\varepsilon(1)$. By the discussion given before Lemma \ref{lem3-17}, we have verified the uniqueness of $\psi_\varepsilon$ for $\varepsilon>0$ small, which means the vortex ring $\zeta_\varepsilon$ with assumptions in Proposition \ref{prop3-1} is unique.\qed
	
	\bigskip
	
	As we discussed before, the value
	$$G_1=G_2=-\frac{3\kappa}{8q_\varepsilon}\cdot b_\varepsilon s_\varepsilon\ln\left(\frac{1}{\varepsilon}\right)+O(\varepsilon)$$
	in the right hand side of \eqref{3-18} extracts the coefficient $\mathbf c_{2,\varepsilon}\cdot s_\varepsilon^2|\ln\varepsilon|$ in expansion \eqref{exp} for $\psi_\varepsilon$, where there is a difference of multiplier $2\pi s_\varepsilon$ due to the integration around $\partial B_{s_\varepsilon}(\boldsymbol q_\varepsilon)$. In view of \eqref{3-18}, the non-vanishing property $\mathbf c_{2,\varepsilon}\neq 0$ makes $b_\varepsilon$ a quantity of order $O(1/|\ln\varepsilon|)$, and hence leads to a contradiction and the uniqueness of $\psi_\varepsilon$ for $\varepsilon\in(0,\varepsilon_0]$, which is consistent with our observation (b) for \eqref{exp}.
	
	\bigskip
	\bigskip
	
	\section{Stability}\label{sec4}
	In this section, we study nonlinear orbital stability of the steady vortex ring $\zeta_\varepsilon$ constructed in Theorem \ref{thm1}. We will provide the proof of Theorem \ref{thm4}. The key idea is to build a bridge between the existence result of \cite{Bad, CWZ} based on Arnol'd's dual variational principle and the uniqueness result established in the proceeding section in order to apply the concentration-compactness theorem of Lions \cite{Lions} to a maximizing sequence.
	
	\bigskip
	
	\subsection{Variational setting}
	Let $\kappa$ and $W$ be as in Theorem \ref{thm1}. We now show that ${\zeta}_\varepsilon$ enjoys a variational characteristic.
	We set the space of admissible functions
	\begin{equation*}
		\mathcal{F}_\varepsilon:=\left\{\zeta\in L^\infty(\mathbb R^3)\mid \zeta: \text{axi-symmetric},\ 0\le  \zeta\le 1/\varepsilon^2, \ \|\zeta\|_{L^1(\mathbb R^3)}\le2\pi\kappa \right\}.
	\end{equation*}
	The same as in the introduction, we also introduce the kinetic energy of the fluid
	\begin{equation*}
		E[\zeta]:=\frac{1}{2}\int_{\mathbb{R}^3}|\mathbf{v}(\boldsymbol x)|^2d\boldsymbol x,\ \ \ \mathbf{v}=\nabla\times(-\Delta)^{-1}\left(r\zeta \right),
	\end{equation*}
	and its impulse
	\begin{equation*}
		\mathcal{P}[\zeta]=\frac{1}{2}\int_{\mathbb{R}^3}r^2\zeta(\boldsymbol x)d \boldsymbol x=\pi\int_\Pi r^3 \zeta drdz.
	\end{equation*}
	We shall consider the maximization problem:
	\begin{equation}\label{4-1}
		\mathcal{E}_{\varepsilon}=\sup_{\zeta\in \mathcal{F}_\varepsilon}\left(E[\zeta]-W\ln\frac{1}{\varepsilon}\mathcal{P}[\zeta]   \right).
	\end{equation}
	Denote by $\mathcal{S}_{\varepsilon}\subset \mathcal{F}_{\varepsilon}$ the set of maximizers of \eqref{4-1}. Note that any $z$-directional translation of $\zeta\in \mathcal{S}_{\varepsilon}$ is still in $ \mathcal{S}_{\varepsilon}$.
	
	\bigskip
	
	The following result is essentially contained in \cite{Bad, CWZ}.
	\begin{proposition}\label{pro4-1}
		If $\varepsilon>0$ is sufficiently small, then $\mathcal{S}_{\varepsilon}\neq \emptyset$ and each maximizer $\hat{\zeta}_\varepsilon \in \mathcal{S}_{\varepsilon}$ is a steady vortex ring with circulation $\kappa$ and  translational velocity $W\ln \varepsilon\,\mathbf{e}_z$. Furthermore,
		\begin{itemize}
			\item [(i)]$\hat{\zeta}_\varepsilon=\varepsilon^{-2}\boldsymbol 1_{\hat{\Omega}_\varepsilon}$ for some axi-symmetric topological torus $\hat{\Omega}_\varepsilon\subset \R^3$.
			\item [(ii)]It holds $C_1\ep \le \sigma\left(\hat{\Omega}_\varepsilon\right)<C_2\ep$ for some constants $0<C_1<C_2$.
			\item [(iii)]As $\varepsilon\to 0$,  $\mathrm{dist}_{\mathcal C_{r^*}}(\hat{\Omega}_\varepsilon)\to0$ with $r^*={\kappa}/{4\pi W}$.
		\end{itemize}
	\end{proposition}
	
	\bigskip

    For easy comprehension, we will give a brief on the proof of Proposition \ref{pro4-1}, which follows the general idea of vorticity method. Firstly, by the $L^1\cap L^\infty$ constraint of $\mathcal{F}_\varepsilon$, we can show that $\mathcal{E}_{\varepsilon}$ is bounded from above and attainable over $\mathcal{F}_\varepsilon$. Then for each $\zeta\in\mathcal{F}_\varepsilon$, let
    $$\zeta_\tau=\hat{\zeta}_\varepsilon+\tau(\zeta-\hat{\zeta}_\varepsilon), \quad \tau\in[0,1].$$
    It is easy to see that $\zeta_\tau$ is still in $\mathcal{F}_\varepsilon$. Since $\hat{\zeta}_\varepsilon \in \mathcal{S}_{\varepsilon}$ is a maximizer, it must hold
    $$\frac{d}{d\tau}\left(E[\zeta_\tau]-W\ln\frac{1}{\varepsilon}\mathcal{P}[\zeta_\tau]   \right)\bigg|_{\tau=0^+}\le 0,$$
    that is
    $$\int_\Pi\hat{\zeta}_\varepsilon\left(\hat\psi_\varepsilon-\frac{W}{2}r^2\ln\frac{1}{\varepsilon}\right)rdrdz\ge \int_\Pi\zeta\left(\hat\psi_\varepsilon-\frac{W}{2}r^2\ln\frac{1}{\varepsilon}\right)rdrdz,$$
    with $\hat\psi_\varepsilon$ the Stokes stream function corresponding to the maximizer $\hat{\zeta}_\varepsilon$. Using the bathtub Lemma (Theorem 1.14 in \cite{LL}), we have
    \begin{equation}\label{4-2}
    	\begin{cases}
    		\hat\psi_\varepsilon-\frac{W}{2}r^2\ln\frac{1}{\varepsilon}> \mu_\varepsilon, & \text{if} \ \hat{\zeta}_\varepsilon=1/\varepsilon^2,
    		\\
    		\hat\psi_\varepsilon-\frac{W}{2}r^2\ln\frac{1}{\varepsilon}= \mu_\varepsilon, & \text{if} \ 0<\hat{\zeta}_\varepsilon<1/\varepsilon^2,
    		\\
    		\hat\psi_\varepsilon-\frac{W}{2}r^2\ln\frac{1}{\varepsilon}< \mu_\varepsilon, & \text{if} \ \hat{\zeta}_\varepsilon=0,
    	\end{cases}
    \end{equation}
    where the flux constant is determined by
    $$\mu_\varepsilon=\inf\left\{x\in\mathbb R \ \big| \ |\{\boldsymbol x\in\Pi \ | \ \hat\psi_\varepsilon-\frac{W}{2}r^2\ln\frac{1}{\varepsilon}>x \}|\le\kappa\varepsilon^2\right\}.$$
    Notice that the kinetic energy $E$ is increased by Steiner symmetrization with respect to $z=h$, while the impulse $\mathcal{P}$ will keep invariant. We claim that $\hat\psi_\varepsilon$ is strictly symmetric decreasing with respect to some horizontal plane, and every level set of $\hat\psi_\varepsilon$ has measure zero, which means the part $0<\hat{\zeta}_\varepsilon<1/\varepsilon^2$ in \eqref{4-2} does not exist. Thus we have obtained the explicit formulation of $\hat{\zeta}_\varepsilon$. To finish the proof of Proposition \ref{pro4-1}, we can investigate the asymptotic shape of $\hat{\zeta}_\varepsilon$ as a disk by energy analysis and Riesz rearrangement inequality, see \cite{CWZ}. 
    
    The last step is to determine location parameter $r^*$ by the variational characteristic, we let $\xi^0_\varepsilon=\boldsymbol 1_{\Omega^0_\varepsilon}\in \mathcal F_\varepsilon$. Here $\Omega^0_\varepsilon$ is a topological torus in $\mathbb R^3$, whose cross-section is a disk centered at $r=q_\varepsilon+\Delta r$ with radius $s_\varepsilon^*=\sqrt{\varepsilon^2\kappa/q_\varepsilon\pi}$. (The parameters $q_\varepsilon$ and $s_\varepsilon^*$ are defined in Section \ref{sec3}). Then according to the estimates in Section \ref{sec3}, it holds
    \begin{equation*}
    	\begin{split}
    		E[\zeta^0_\varepsilon]-W\ln\frac{1}{\varepsilon}\mathcal{P}[\zeta^0_\varepsilon]&=\left(\frac{\kappa}{4\pi}\cdot q_\varepsilon-\frac{W}{2}\cdot q_\varepsilon^{2}\right)\ln\frac{1}{\varepsilon}\\
    		&\quad-(\mathbf c_{2,\varepsilon}+\mathbf c_{2,W})\cdot\Delta r^2\ln\frac{1}{\varepsilon}-E_{\phi}+O(\varepsilon^2),
    	\end{split}
    \end{equation*}
    where $E_{\phi}>0$ is fixed and induced by $\phi_\varepsilon$. We see that $q_\varepsilon$ is the exact location for the energy attending its maximum, where $\mathbf c_{2,\varepsilon}+\mathbf c_{2,W}>0$ leads to the convexity as we have discussed in condition (c) for \eqref{exp}.  Hence we can derive the limit location $r^*=q_\varepsilon={\kappa}/{4\pi W}+o_\varepsilon(1)$ from the angle of energy maximization. This description is actually more precise than that in \cite{CWZ,VS}, where the authors claimed that $r^*$ is the maximum point of a quadratic function $-\frac{W}{2}\cdot x^2+\frac{\kappa}{4\pi}\cdot x$. 
         
	Since $\hat\zeta_\varepsilon\in \mathcal{S}_{\varepsilon}$ must be symmetric with respect to some horizontal line $z=h$ by former discussion, it can be centralized by a unique translation in the $z$-direction that makes it a centralized steady vortex ring. We shall denote its centralized version by $\zeta^\#$. We also set $\mathcal{S}^\#_{\varepsilon}:=\{\zeta^\#\mid \zeta\in \mathcal{S}_{\varepsilon}\}$. In view of the uniqueness result in Theorem \ref{thm2}, we see that $\mathcal{S}^\#_{\varepsilon}=\{{\zeta}_\varepsilon\}$ for all $\varepsilon>0$ small.
	
	\bigskip
	
	The following elementary estimates can be found in \cite{Choi20} (see Lemma 2.3 in \cite{Choi20}), which will be used as basic tools to prove the the nonlinear orbital stability.
	\begin{lemma}\label{le4-2}
		There exists a positive number $C$ such that
		\begin{equation*}
			\begin{split}
				&  |E[\zeta]|\le E[|\zeta|]\le C\left(\|r^2\zeta\|_{L^1(\mathbb R^3)}+\|\zeta\|_{L^1\cap L^2(\mathbb R^3)} \right)\|r^2\zeta\|_{L^1(\mathbb R^3)}^{1/2}\|r^2\zeta\|_{L^1(\mathbb R^3)}^{1/2}, \\
				& |E[\zeta_1]-E[\zeta_2]|\le C\left(\|r^2(\zeta_1+\zeta_2)\|_{L^1(\mathbb R^3)}+\|\zeta_1+\zeta_2\|_{L^1\cap L^2(\mathbb R^3)} \right)\\
				&\quad\quad\quad\quad\quad\quad\quad\quad\quad\quad\quad\quad\quad\quad\quad\quad \ \ \times\|r^2(\zeta_1-\zeta_2)\|_{L^1(\mathbb R^3)}^{1/2}\|r^2(\zeta_1-\zeta_2)\|_{L^1(\mathbb R^3)}^{1/2},
			\end{split}
		\end{equation*}
		for any axi-symmetric $\zeta$, $\zeta_1$, $\zeta_2\in \left(L^1\cap L^2\cap L^1_\mathrm{w}\right)(\mathbb R^3)$.
	\end{lemma}

    \bigskip
	
	\subsection{Reduction to absurdity}
	
	Since $\zeta_\varepsilon$ is the unique centralized maximizer of $\mathcal E_\varepsilon$, and the solution will conserve its kinetic energy $E$ and impulse $\mathcal P$, it is reasonable to conjecture that any small perturbation of $\zeta_\varepsilon$ can not leave $\zeta_\varepsilon$ too far as time evolves (possibly modulo a $z$-directional translation), namely, the nonlinear orbital stability. Following this intuitive idea, we are now in a position to prove Theorem \ref{thm4}.
	
	\bigskip
	
	\noindent{\bf Proof of Theorem \ref{thm4}:}
	We argue by contradiction. Suppose that there exist a positive number $\eta_0$, a sequence $\{\zeta_{0,n}\}_{n=1}^\infty$ of non-negative axi-symmetric functions, and a sequence $\{t_n\}_{n=1}^{\infty}$ of non-negative numbers such that, for each $n\ge 1$, we have $\zeta_{0,n}$, $(r\zeta_{0,n})\in L^\infty(\mathbb R^3)$,
	\begin{equation*}
		\|\zeta_{0,n}-\zeta_\varepsilon\|_{L^1\cap L^2(\mathbb R^3)}+\|r^2(\zeta_{0,n}-\zeta_\varepsilon)\|_{L^1(\mathbb R^3)}\le \frac{1}{n^2},
	\end{equation*}
	and
	\begin{equation*}
		\inf_{\tau\in \mathbb R} \|\zeta_n(\cdot-\tau \mathbf{e}_z,t_n)-\zeta_\varepsilon\|_{L^1\cap L^2(\mathbb R^3)}+\|r^2(\zeta_n(\cdot-\tau \mathbf{e}_z,t_n)-\zeta_\varepsilon)\|_{L^1(\mathbb R^3)}\ge \eta_0,
	\end{equation*}
	where $\zeta_n(\boldsymbol x,t)$ is the global-in-time weak solution of \eqref{1-7} for the initial data $\zeta_{0,n}$ obtained by Proposition \ref{Pro1}.  Using Lemma \ref{le4-2}, we get
	\begin{equation*}
		\lim_{n\to \infty}E[\zeta_{0,n}]=E[\zeta_\varepsilon].
	\end{equation*}
	Thus, we have
	\begin{equation*}
		\begin{split}
			& \lim_{n\to \infty} \mathcal{P}[\zeta_{0,n}]=\mathcal{P}[\zeta_\varepsilon],\ \  \lim_{n\to \infty}E[\zeta_{0,n}]= E[\zeta_\varepsilon],\\
			&  \lim_{n\to \infty}\|\zeta_{0,n}\|_{L^p(\mathbb R^3)}=\|\zeta_{\varepsilon}\|_{L^p(\mathbb R^3)},\ \ \forall\, 1\le p\le 2.
		\end{split}
	\end{equation*}
	Let us write $\zeta_n=\zeta_n(\cdot,t_n)$. By virtue of the conservations, we conclude that
	\begin{equation}\label{4-5}
		\begin{split}
			& \lim_{n\to \infty} \mathcal{P}[\zeta_{n}]=\mathcal{P}[\zeta_\varepsilon],\ \  \lim_{n\to \infty}E[\zeta_{n}]=E[\zeta_\varepsilon],\\
			&  \lim_{n\to \infty}\|\zeta_{n}\|_{L^p(\mathbb R^3)}=\|\zeta_{\varepsilon}\|_{L^p(\mathbb R^3)},\ \ \forall\, 1\le p\le 2.
		\end{split}
	\end{equation}
	Note that
	\begin{equation*}
		\int_{\left\{\boldsymbol x\in \mathbb R^3\,\,\mid \,\,|\zeta_n(\boldsymbol x)-1/\varepsilon^2|\ge 1/n\right\}}\zeta_n d\boldsymbol x=\int_{\left\{\boldsymbol x\in \mathbb R^3\,\,\mid\,\,\, |\zeta_{0,n}(\boldsymbol x)-1/\varepsilon^2|\ge 1/n\right\}}\zeta_{0,n}d\boldsymbol x.
	\end{equation*}
	Set $D(n):=\left\{\boldsymbol x\in \mathbb R^3\,\,\mid\,\, |\zeta_{0,n}(\boldsymbol x)-1/\varepsilon^2|\ge 1/n \right\}$ and $Q:=\text{supp}\, \zeta_\varepsilon$. We check that
	\begin{equation*}
		\begin{split}
			\int_{D(n)}\zeta_{0,n}d\boldsymbol x & =\|\zeta_{0,n} \|_{L^1\left(D(n)\cap Q\right)}+\|\zeta_{0,n} \|_{L^1\left(D(n)\backslash Q \right)} \\
			& \le \|\zeta_{0,n}-\zeta_\varepsilon \|_{L^1\left(D(n)\cap Q \right)}+\|\zeta_\varepsilon \|_{L^1\left(D(n)\cap Q \right)}+\|\zeta_{0,n}-\zeta_\varepsilon \|_{L^1\left(D(n)\backslash Q \right)}\\
			&\le \|\zeta_{0,n}-\zeta_\varepsilon \|_{L^1(\mathbb R^3)}+\|\zeta_\varepsilon \|_{L^1\left(D(n)\cap Q\right)}\\
			&\le \|\zeta_{0,n}-\zeta_\varepsilon \|_{L^1(\mathbb R^3)}+|D(n)\cap Q|\le (n+1)\|\zeta_{0,n}-\zeta_\varepsilon \|_{L^1(\mathbb R^3)}\le \frac{n+1}{n^2}\to 0
		\end{split}
	\end{equation*}
	as $n\to \infty$, where we used the fact that
	\begin{equation*}
		\frac{1}{n} |D(n)\cap Q|\le \|\zeta_{0,n}-\zeta_\varepsilon \|_{L^1\left(D(n)\cap Q \right)}\le \|\zeta_{0,n}-\zeta_\varepsilon \|_{L^1(\mathbb R^3)}.
	\end{equation*}
	Set
	\begin{equation*}
		\mathcal{F}^*_\varepsilon:=\left\{\zeta\in\mathcal{F}_\varepsilon\mid \mathcal{P}[\zeta]=\mathcal{P}[\zeta_\varepsilon] \right\}.
	\end{equation*}
	It is easy to see that
	\begin{equation*}
		E[\zeta_\varepsilon]=\max_{\zeta\in \mathcal{F}^*_\varepsilon}E[\zeta]\ \ \ \text{and}\ \ \ \mathcal{S}_{\varepsilon}=\left\{\zeta\in \mathcal{F}^*_\varepsilon \mid E[\zeta]
		=E[\zeta_\varepsilon]\right\}.
	\end{equation*}
	Therefore, we can now use Theorem 3.1 in \cite{Choi20} as a version of the concentration-compactness principle to obtain a subsequence (still using the same index $n$) and $\{\tau_n\}_{n=1}^\infty\subset \mathbb R$ such that
	\begin{equation*}
		\|r^2\left(\zeta_{n}(\cdot-\tau_n\mathbf{e}_z)-\zeta_\varepsilon\right)\|_{L^1(\mathbb R^3)}\to 0, \quad \text{as}\  n\to \infty.
	\end{equation*}
	Recalling \eqref{4-5}, we can further deduce that
	\begin{equation*}
		\|\zeta_{n}(\cdot-\tau_n\mathbf{e}_z)-\zeta_\varepsilon\|_{L^2(\mathbb R^3)}\to 0, \quad \text{as}\  n\to \infty.
	\end{equation*}
	By H$\ddot{\text{o}}$lder's inequality, we get
	\begin{equation*}
		\lim_{n\to\infty}\int_Q \zeta_{n}(\boldsymbol x-\tau_n\mathbf{e}_z)d\boldsymbol x=\int_Q \zeta_\varepsilon(\boldsymbol x) d\boldsymbol x,
	\end{equation*}
	which implies
	\begin{equation*}
		\lim_{n\to\infty}\int_{\mathbb R^3\backslash Q}\zeta_n(\boldsymbol x-\tau_n\mathbf{e}_z)d\boldsymbol x=\lim_{n\to\infty}\int_{\mathbb R^3}\zeta_n(\boldsymbol x-\tau_n\mathbf{e}_z)d\boldsymbol x-\lim_{n\to\infty}\int_{Q}\zeta_n(\boldsymbol x-\tau_n\mathbf{e}_z)d\boldsymbol x=0.
	\end{equation*}
	It follows that
	\begin{equation*}
		\begin{split}
			\|\zeta_{n}(\cdot-\tau_n\mathbf{e}_z)-\zeta_\varepsilon\|_{L^1(\mathbb R^3)} & =\|\zeta_{n}(\cdot-\tau_n\mathbf{e}_z)-\zeta_\varepsilon\|_{L^1(Q)}+\|\zeta_{n}(\cdot-\tau_n\mathbf{e}_z)-\zeta_\varepsilon\|_{L^1(\mathbb R^3\backslash Q)} \\
			& \le |Q|^{1/2}\|\zeta_{n}(\cdot-\tau_n\mathbf{e}_z)-\zeta_\varepsilon\|_{L^2(\mathbb R^3)}+\|\zeta_{n}(\cdot-\tau_n\mathbf{e}_z)\|_{L^1(\mathbb R^3\backslash Q)}\to 0
		\end{split}
	\end{equation*}
	as $n\to \infty$. In sum, we have
	\begin{equation*}
		0 =\lim_{n\to\infty}\|\zeta_n(\cdot-\tau_n \mathbf{e}_z,t_n)-\zeta_\varepsilon\|_{L^1\cap L^2(\mathbb R^3)}+\|r^2(\zeta_n(\cdot-\tau_n \mathbf{e}_z,t_n)-\zeta_\varepsilon)\|_{L^1(\mathbb R^3)}\ge \eta_0>0,
	\end{equation*}
	which is a contradiction. The proof is thus complete.
	\qed
	
\bigskip	
\bigskip

	\appendix
	
	\section{Method of moving planes}\label{appA}
	
	In this appendix, we will show that the cross-section $A_\varepsilon$ and Stokes stream function $\psi_\varepsilon$ are symmetric with respect to the line $\{x_2=h\}$ for some $h$ by the method of moving planes (see also Lemma 2.1 in \cite{AF88}). Though the proof is almost the same as that of Appendix A in \cite{CQZZ2}, we give it in detail here for readers' convenience.
	\begin{proposition}\label{A1}
		Suppose that a bounded set $A$ with $\bar A\subset \mathbb R^2_+$, satisfies
		$$A=\{\boldsymbol x\in B_R(\boldsymbol 0)\cap \{x_1>0\}\mid \psi(\boldsymbol x)+\frac{W}{2}x_1^2>\mu\}$$
		for some constants $W$ and $\mu$. Moreover, $\psi$ is the potential of $A$ in the sense
		\begin{equation*}	
			\psi(\boldsymbol x)=\frac{1}{4\pi}\int_{\mathbb R^2_+}{G_*}(\boldsymbol x,\boldsymbol x')\boldsymbol 1_A(\boldsymbol x')d\boldsymbol x'.	 
		\end{equation*}	
		Then, there is $h\in \mathbb{R}$ such that  $A$ is symmetric with respect to the line $\{x_2=h\}$.
	\end{proposition}
	\begin{proof}
		To prove this proposition, the key observation is that  ${G_*}(\boldsymbol x,\boldsymbol x')$ is a strictly decreasing function of $|x_2-x_2'|^2$ for fixed $x_1$ and $x_1'$. Namely, for any fixed $x_1$ and $x_1'$, if we denote $r_2:=|x_2-x_2'|^2$, then we have ${G_*}(\boldsymbol x,\boldsymbol x')=J_{x_1,x_1'}(r_2)$ for some strictly decreasing function $J_{x_1,x_1'}(\cdot)$.
		
		For $-R<t<R$, define
		$$A_t:= \{\boldsymbol   x\in A\mid x_2<t\},\quad A_t^*:=\{\boldsymbol x\in\mathbb{R}^2\mid (x_1,2t-x_2)\in A_t\}.$$ This is, $A_t^*$ is the reflection of $A_t$ with respect to the line $x_2=t$.
		Let $\mathtt d:=\inf_{\boldsymbol y\in A} y_2$.
		We will carry out the proof of Proposition  \ref{A1} by two steps.
		
		\medskip
		
		\emph{Step 1.}\,Let us first show that there exists $\epsilon>0$ small enough such that, for any $\mathtt d<t\leq \mathtt d+\epsilon$,
		$$A_t^*\subset A.$$
		For any $\boldsymbol x\in \{x_2=\mathtt d\}\cap \bar A$, we compute
		\begin{align*}
			\partial_{x_2} \psi (\boldsymbol  x)= \int_{A} 2\partial_{r_2}J_{x_1,x_1'}(|x_2-x_2'|^2)(x_2-x_2')d\boldsymbol x'\geq c_0>0,
		\end{align*}
		for some constant $c_0$ independent of $\boldsymbol x$. 
		We define the set $S_\epsilon:=\{\boldsymbol x\in A\mid \mathtt d<x_2<\mathtt d+\epsilon\}$. Arguing by contradiction, we can show that $\sup_{\boldsymbol  x\in S_\epsilon} \text{dist}(\boldsymbol  x, \{x_2=\mathtt d\}\cap \bar A)\to 0$ as $\epsilon\to 0$. Then, by the $C^1_{\text{loc}}$ continuity of $\psi$ in $\mathbb{R}^2_+$, there exists $\varepsilon_1>0$ small such that $\partial_{x_2}\psi(\boldsymbol  x)>c_0/2>0$ for all $\boldsymbol x\in S_\epsilon$ whenever $0<\epsilon<\epsilon_1$. Since $\psi\in C_{\text{loc}}^{1,\alpha}(\mathbb{R}^2_+)$ by the regularity theory and $A$ is far away from the boundary $x_1=0$, for $\mathtt d<t<\mathtt d+\epsilon_1$, we have for all $\boldsymbol  x\in A_t$,
		\begin{align*}
			\psi(x_1,2t-x_2)-\psi(x_1,x_2)&=2\partial_{x_2}\psi(\boldsymbol x)(t-x_2)+O((t-x_2)^{1+\alpha})\\
			&\geq c_0 (t-x_2)+O((t-x_2)^{1+\alpha}).
		\end{align*}
		Thus, there exists $0<\varepsilon_2\leq \epsilon_1$ small such that for any $\mathtt d<t<\mathtt d+\epsilon_2$, it holds $$\psi(x_1,2t-x_2)-\psi(x_1,x_2)\geq 0,\,\,\,\forall\, \boldsymbol x\in A_t,$$
		which implies $A_t^*\subset A$.
		
		\medskip
		
		\emph{Step 2.}\,We move the line continuously until its limiting position. Step 1 provides a starting point for us to move lines. Define the limiting position
		$$h:=\sup\{t\mid A_\tau^*\subset A,\,\,\forall \,\mathtt d<\tau\leq t\}.$$
		We will show that $A$ is symmetric with respect to the line $\{x_2=h\}$.
		In fact, we are going to prove that $$|N|=0, \quad \text{for} \ N=A\setminus (A_h\cup A_h^*).$$
		Suppose that $|N|>0$, we will get a contradiction.
		
		By step 1, we have $\mathtt d<h<\sup_{\boldsymbol x\in A} x_2$. By the definition of $h$, we have $A^*_h\subset \bar A$. We first claim that  $ \partial  A\cap \partial A_h^*\not=\emptyset$. Indeed, suppose on the contrary that $\bar A_h^*\subset A$. This means that $A_h$ is far away from the line $\{x_2=h\}$ and the set $A$ is divided into disjoint sets by $\{x_2=h\}$. Then, it is easy to see that there exists a $\mathtt d<t<h$ such that $A_t^*\not\subset A$, which contradicts the definition of $h$. Therefore, we must have $ \partial  A\cap \partial A_h^*\not=\emptyset$.
		
		Suppose that there exists a point $\boldsymbol x^*\in \partial  A\cap \partial A_h^*$ such that $x_2^*>h$. We write $\boldsymbol  x=(x_1^*,2h-x_2^*)$. Then, we calculate
		\begin{align*}
			0&=\psi(\boldsymbol x)-\psi(\boldsymbol x^*)\\
			&=\int_{N} \left({G_*}(\boldsymbol x,\boldsymbol x')-{G_*}(\boldsymbol x^*,\boldsymbol x')\right)d\boldsymbol x'<0,
		\end{align*}
		if $|N|>0$. Here, we have used the fact that $|x_2-x_2'|>|x_2^*-x_2'|$ for any $\boldsymbol x' \in N$. This is a contradiction and thus we must have $|N|=0$ in this case.
		
		Now, we consider the remaining case, where for any $\boldsymbol x^*\in \partial  A\cap \partial A_h^*$, it must holds $x_2^*=h$ and thus $\boldsymbol x=\boldsymbol x^*$. However, for any $\boldsymbol x\in \bar  A\cap\{x_2=h\}$, it holds
		\begin{align*}
			\partial_{x_2} \psi (\boldsymbol x)= \int_{N} 2\partial_{r_2}J_{x_1,x_1'}(|x_2-x_2'|^2)(x_2-x_2')d\boldsymbol x'\geq c_0>0,
		\end{align*}
		for some constant $c_0$ independent of $\boldsymbol x$ provided that $|N|>0$.  We can take $\varepsilon_3>0$ small such that $\partial_{x_2}\psi(\boldsymbol  x)\geq c_0/2>0$ for all $\boldsymbol  x$ lies in the strip $\{\boldsymbol  x\in A\mid h-\epsilon_3<x_2<h+\varepsilon_3\}$.
		We denote $A_b^{*,c}$ as the reflection of the set $A_b$ with respect to line $x_2=c$ for any $b,c \in \mathbb{R}$. We first have $\text{dist} (A_{h-\epsilon_3}^{*,h},\partial A)\geq c_{\epsilon_3}$ for some constant $c_{\epsilon_3}>0$. Otherwise, we will obtain a point $\boldsymbol x^*\in \partial A_{h}^{*}\cap \partial A$ with $x^*_2\geq h+\epsilon>h$, which has already been considered. Therefore, if we take $\epsilon_4:=\min\{\epsilon_3, c_{\epsilon_3}\}$, then  for all $h<t<h+\epsilon_4$, it holds $$A_{h-\epsilon_3}^{*,t}\subset A.$$ For $\boldsymbol x$ in the strip $A\cap \{h-\epsilon_3\leq x_2<t\}$, we have
		\begin{align*}
			\psi(x_1,2t-x_2)-\psi(x_1,x_2)&=2\partial_{x_2}\psi(\boldsymbol x)(t-x_2)+O((t-x_2)^{1+\alpha})\\
			&\geq c_0 (t-x_2)+O((t-x_2)^{1+\alpha}).
		\end{align*}
		Thus, there exists $0<\epsilon_5\leq \epsilon_4$ small such that for any $h<t<h+\epsilon_5$, it holds $$\psi(x_1,2t-x_2)-\psi(x_1,x_2)\geq 0,\quad \forall\,\,x\in A\cap \{s-\epsilon_3\leq x_2<t\},$$
		which implies $A_t^*\subset A$. This contradicts the definition of $h$ and hence we must have $|N|=0$, which means that $A$ is symmetric with respect to some  line $\{x_2=h\}$.
		
		The proof is thus finished.	
	\end{proof}
	
	\bigskip
	
	\section{Expansion for stream function and cross-section boundary}\label{appB}
	
	In this appendix, we will expand the stream function $\psi_\varepsilon$ and cross-section boundary $\partial A_\varepsilon$ by Taylor's formula, where $A_\varepsilon$ is given by the level set 
	$$A_\varepsilon=\left\{\boldsymbol x\in \mathbb R^2_+ \ \big| \ \psi_\varepsilon-\frac{W}{2}x_1^2\ln\frac{1}{\varepsilon}>\mu_\varepsilon\right\}$$
	for some flux constant $\mu_\varepsilon>0$. Recall that we let
	\begin{equation*}
		\begin{split}
			\mathbf U_{\boldsymbol q_\varepsilon,\varepsilon}(\boldsymbol x)&=\psi_\varepsilon(\boldsymbol x)-\phi_\varepsilon(\boldsymbol x)-\frac{W}{2}x_1^2\ln\frac{1}{\varepsilon}-\mu_\varepsilon\\
			&=V_{\boldsymbol q_\varepsilon,\varepsilon}(\boldsymbol x)-V_{\bar{\boldsymbol q}_\varepsilon,\varepsilon}(\boldsymbol x)+\mathcal H_{\boldsymbol q_\varepsilon,\varepsilon}(\boldsymbol x)-\frac{W}{2}x_1^2\ln\frac{1}{\varepsilon}-\mu_\varepsilon.
		\end{split}
	\end{equation*}
	be an approximation of $V_{\boldsymbol q_\varepsilon,\varepsilon}(\boldsymbol x)-\frac{a_\varepsilon}{2\pi}\ln\frac{1}{\varepsilon}$, where $\phi_\varepsilon(\boldsymbol x)$ is the error function, functions
	\begin{equation*}
		V_{\boldsymbol q_\varepsilon,\varepsilon}(\boldsymbol x)=\left\{
		\begin{array}{lll}
			\frac{a_\varepsilon}{2\pi}\ln\frac{1}{\varepsilon}+\frac{q_\varepsilon^2}{4\varepsilon^2}(s_\varepsilon^2-|\boldsymbol x-\boldsymbol q_\varepsilon|^2), \ \ \  &\mathrm{if} \ |\boldsymbol x-\boldsymbol q_\varepsilon|\le s_\varepsilon,\\
			\frac{a_\varepsilon}{2\pi}\ln\frac{1}{\varepsilon}\cdot\frac{\ln|\boldsymbol x-\boldsymbol q_\varepsilon|}{\ln s_\varepsilon},&\mathrm{if} \ |\boldsymbol x-\boldsymbol q_\varepsilon|\ge s_\varepsilon,
		\end{array}
		\right.
	\end{equation*}
	$V_{\bar{\boldsymbol q}_\varepsilon,\varepsilon}(\boldsymbol x)$ with $\bar{\boldsymbol q}_\varepsilon=(-q_\varepsilon,0)$, and
	$$\mathcal H_{\boldsymbol q_\varepsilon,\varepsilon}(\boldsymbol x)=\frac{1}{\varepsilon^2}\int_{\mathbb R^2_+}H(\boldsymbol x,\boldsymbol x') \boldsymbol 1_{B_{s_\varepsilon}(\boldsymbol q_\varepsilon)}(\boldsymbol x')d\boldsymbol x'$$
	correspond to the vortex, mirror vortex and remaining regular part in the stream function respectively. Moreover, $s_\varepsilon$ as the radius of $\{\boldsymbol x\in \mathbb R^2_+ \mid V_{\boldsymbol q_\varepsilon,\varepsilon}(\boldsymbol x)-\frac{a_\varepsilon}{2\pi}\ln\frac{1}{\varepsilon}>0\}$ will satisfy the regularity condition \eqref{2-14}, namely, 
	\begin{equation}\label{gradient}
		\mathcal N_\varepsilon:=\frac{a_\varepsilon}{2\pi}\ln\frac{1}{\varepsilon}\cdot\frac{1}{s_\varepsilon|\ln s_\varepsilon|}=\frac{s_\varepsilon}{2\varepsilon^2}\cdot q_\varepsilon^2,
	\end{equation}
	and the parameter $a_\varepsilon$ is determined by
	\begin{equation}\label{acondition}
		\frac{a_\varepsilon}{2\pi}\ln\frac{1}{\varepsilon}=\mu_\varepsilon+\frac{W}{2}q_\varepsilon^2\ln\frac{1}{\varepsilon}-\mathcal H_{\boldsymbol q_\varepsilon,\varepsilon}(\boldsymbol q_\varepsilon)+V_{\boldsymbol{\bar q_\varepsilon},\varepsilon}(\boldsymbol q_\varepsilon),
	\end{equation}
	All these assumptions make 
	$$\Phi_{\boldsymbol q_\varepsilon,\varepsilon}(\boldsymbol x)=V_{\boldsymbol q_\varepsilon,\varepsilon}(\boldsymbol x)-V_{\bar{\boldsymbol q}_\varepsilon,\varepsilon}(\boldsymbol x)+\mathcal H_{\boldsymbol q_\varepsilon,\varepsilon}(\boldsymbol x)$$ 
	a good approximation to $\psi_\varepsilon$, where the difference is perturbation term $\phi_\varepsilon$.
	Notice that the set 
	$$\{\boldsymbol x\in \mathbb R^2_+ \mid V_{\boldsymbol q_\varepsilon,\varepsilon}(\boldsymbol x)-\frac{a_\varepsilon}{2\pi}\ln\frac{1}{\varepsilon}>0\}$$ 
	equals the disk $B_{s_\varepsilon}(\boldsymbol q_\varepsilon)$. To give a precise description for $\partial A_\varepsilon$, our strategy is to compare the two sets  $\{\boldsymbol x\in \mathbb R^2_+ \mid \mathbf U_{\boldsymbol q_\varepsilon,\varepsilon}(\boldsymbol x)>0\}$ and $B_{s_\varepsilon}(\boldsymbol q_\varepsilon)$ first, and then obtain the estimates of $A_\varepsilon$ by calculate its difference with $\{\boldsymbol x\in \mathbb R^2_+ \mid \mathbf U_{\boldsymbol q_\varepsilon,\varepsilon}(\boldsymbol x)>0\}$ considering the influence of error term $\phi_\varepsilon$. At the first stage, we will compare the value of two functions $\mathbf U_{\boldsymbol q_\varepsilon,\varepsilon}(\boldsymbol x)$ and $V_{\boldsymbol q_\varepsilon,\varepsilon}(\boldsymbol x)-\frac{a_\varepsilon}{2\pi}\ln\frac{1}{\varepsilon}$ near the cross-section $A_\varepsilon$. For this purpose, we let the function $\mathcal W_\varepsilon(x):\mathbb R_+\to\mathbb R$ be defined as 
	\begin{equation}\label{B-1}
		\begin{split}
			\mathcal W_\varepsilon(x)&=\frac{s_\varepsilon^2}{4\varepsilon^2}\cdot q_\varepsilon\ln\frac{1}{s}-Wq_\varepsilon\ln\frac{1}{\varepsilon}\\
			&\quad +\frac{1}{8\varepsilon^2}\cdot q_\varepsilon\left\{
			\begin{array}{lll}
				(s_\varepsilon^2-x^2), \ &  0<x<s_\varepsilon\\
				2\ln(s_\varepsilon/x), \ &  x\ge s_\varepsilon
			\end{array}
			\right.
			+\frac{3}{16\varepsilon^2}\cdot q_\varepsilon\left\{
			\begin{array}{lll}
				2s_\varepsilon^2-x^2, \ & 0<x<s_\varepsilon\\
				s_\varepsilon^4/x^2, \ &  x\ge s_\varepsilon
			\end{array}
			\right.\\
			&\quad+\frac{s_\varepsilon^2}{4\varepsilon^2}\cdot q_\varepsilon(\ln(8q_\varepsilon)-1),
		\end{split}
	\end{equation}
	which will be shown as the coefficient in linear part of the difference. It should be noticed that $s_\varepsilon^2q_\varepsilon\pi/\varepsilon^2\to \kappa$ as $\varepsilon\to 0$ by our assumptions on $s_\varepsilon,\mu_\varepsilon, a_\varepsilon$ in \eqref{gradient}-\eqref{acondition}. As a result, as $\varepsilon\to 0$ we have
	\begin{equation*}
		\begin{split}
	-\mathcal W_\varepsilon(s_\varepsilon)&=Wq_\varepsilon\ln\frac{1}{\varepsilon}-\frac{s_\varepsilon^2}{4\varepsilon^2}\cdot q_\varepsilon\ln\frac{8q_\varepsilon}{s_\varepsilon}+\frac{s_\varepsilon^2}{16\varepsilon^2}\cdot q_\varepsilon\\
	&=Wq_\varepsilon\ln\frac{1}{\varepsilon}-\frac{\kappa}{4\pi}\ln\frac{8q_\varepsilon}{s_\varepsilon}+\frac{\kappa}{16\pi}+o_\varepsilon(1),
	\end{split}
	\end{equation*}
	which is exactly the leading order term of Kelvin-Hicks formula in Proposition \ref{prop3-2} or $\mathbf c_1$ in \eqref{exp}. We will see it once again at the end of Appendix \ref{appC}.
	
	\bigskip
	
	 In this appendix, we always assume that $L>0$ is a large fixed constant, and denote $\tilde v(\boldsymbol y)=v(s_\varepsilon\boldsymbol y+\boldsymbol q_\varepsilon)$ for a general function $v$. At the first stage, the Taylor's expansion of $\mathbf U_{\boldsymbol q_\varepsilon,\varepsilon}(\boldsymbol x)$ at $V_{\boldsymbol q_\varepsilon,\varepsilon}(\boldsymbol x)-\frac{a_\varepsilon}{2\pi}\ln\frac{1}{\varepsilon}$ is given in the following lemma.
	
	\begin{lemma}\label{B1}
		For every $\boldsymbol y\in D_\varepsilon=\{\boldsymbol y \ | \ s_\varepsilon\boldsymbol y+\boldsymbol q_\varepsilon\in\mathbb R^2_+\}$ bounded, it holds
		$$\tilde{\mathbf U}_{\boldsymbol q_\varepsilon,\varepsilon}(\boldsymbol y)=\tilde V_{\boldsymbol q_\varepsilon,\varepsilon}(\boldsymbol y)-\frac{a_\varepsilon}{2\pi}\ln\frac{1}{\varepsilon}+s_\varepsilon y_1\cdot\mathcal W_\varepsilon(|s_\varepsilon\boldsymbol y|)+O(\varepsilon^2|\ln\varepsilon|).$$
	\end{lemma}
	\begin{proof}
		By the definition of $\mathbf U_{\boldsymbol q_\varepsilon,\varepsilon}(\boldsymbol x)$ and asymptotic estimate for $G_*$ in \eqref{2-12}, it holds
		\begin{equation*}
			\begin{split}
				\mathbf U_{\boldsymbol q_\varepsilon,\varepsilon}(\boldsymbol x)&=\frac{1}{2\pi\varepsilon^2}\int_{ B_{s_\varepsilon}(\boldsymbol q_\varepsilon)}x_1^{1/2}x_1'^{3/2}\ln\left(\frac{1}{|\boldsymbol x-\boldsymbol x'|}\right)d\boldsymbol x'-\frac{W}{2}x_1^2\ln\frac{1}{\varepsilon}-\mu_\varepsilon\\
				&\quad+\frac{1}{4\pi\varepsilon^2}\int_{ B_{s_\varepsilon}(\boldsymbol q_\varepsilon)}x_1^{1/2}x_1'^{3/2}\left(\ln(x_1x_1')+2\ln 8-4+O\left(\rho\ln\frac{1}{\rho}\right)\right)d\boldsymbol x'\\
				&=\frac{q_\varepsilon^2}{2\pi\varepsilon^2}\int_{ B_{s_\varepsilon}(\boldsymbol q_\varepsilon)}\ln\left(\frac{1}{|\boldsymbol x-\boldsymbol x'|}\right)d\boldsymbol x'+\frac{1}{2\pi\varepsilon^2}\int_{ B_{s_\varepsilon}(\boldsymbol q_\varepsilon)}(x_1^{1/2}x_1'^{3/2}-q_\varepsilon^2)\ln\left(\frac{1}{|\boldsymbol x-\boldsymbol x'|}\right)d\boldsymbol x'\\
				&\quad+\frac{1}{4\pi\varepsilon^2}\int_{ B_{s_\varepsilon}(\boldsymbol q_\varepsilon)}x_1^{1/2}x_1'^{3/2}\left(\ln(x_1x_1')+2\ln 8-4+O\left(\rho\ln\frac{1}{ \rho}\right)\right)d\boldsymbol x'\\
				&\quad-\frac{W}{2}x_1^2\ln\frac{1}{\varepsilon}-\mu_\varepsilon,
			\end{split}
		\end{equation*}
		with $\rho$ the distance defined in \eqref{rho}. According to the Taylor's formula, we have
		\begin{equation*}
			\begin{split}
				&\quad\frac{1}{2\pi\varepsilon^2}\int_{ B_{s_\varepsilon}(\boldsymbol q_\varepsilon)}(x_1^{1/2}x_1'^{3/2}-q_\varepsilon^2)\ln\left(\frac{1}{|\boldsymbol x-\boldsymbol x'|}\right)d\boldsymbol x'\\
				&=\frac{1}{2\pi\varepsilon^2}\int_{ B_{s_\varepsilon}(\boldsymbol q_\varepsilon)}\left(\left(q_\varepsilon^{1/2}+\frac{1}{2q_\varepsilon^{1/2}}(x_1-q_\varepsilon)+O(\varepsilon^2)\right)\left(q_\varepsilon^{3/2}+\frac{3q_\varepsilon^{1/2}}{2}(x_1'-q_\varepsilon)+O(\varepsilon^2)\right)-q_\varepsilon^2\right)\\
				&\quad\times\ln\left(\frac{1}{|\boldsymbol x-\boldsymbol x'|}\right)d\boldsymbol x'\\
				&=\frac{q_\varepsilon}{2\pi\varepsilon^2}\int_{ B_{s_\varepsilon}(\boldsymbol q_\varepsilon)}\left(\frac{x_1-q_\varepsilon}{2}+\frac{3(x_1'-q_\varepsilon)}{2}\right)\ln\left(\frac{1}{|\boldsymbol x-\boldsymbol x'|}\right)d\boldsymbol x'+O(\varepsilon^2|\ln\varepsilon|)\\
				&=\frac{s_\varepsilon^2}{4\varepsilon^2}\cdot q_\varepsilon(x_1-q_\varepsilon)\ln\frac{1}{s_\varepsilon}+\frac{(x_1-q_\varepsilon)}{8\varepsilon^2}\cdot q_\varepsilon\left\{
				\begin{array}{lll}
					(s_\varepsilon^2-|\boldsymbol x-\boldsymbol q_\varepsilon|^2), \ &  |\boldsymbol x-\boldsymbol q_\varepsilon|<s_\varepsilon\\
					2\ln\frac{s_\varepsilon}{|\boldsymbol x-\boldsymbol q_\varepsilon|}, \ &  |\boldsymbol x-\boldsymbol q_\varepsilon|\ge s_\varepsilon
				\end{array}
				\right.\\
				&\quad+\frac{3(x_1-q_\varepsilon)}{16\varepsilon^2}\cdot q_\varepsilon\left\{
				\begin{array}{lll}
					2s_\varepsilon^2-|\boldsymbol x-\boldsymbol q_\varepsilon|^2, \ &  |\boldsymbol x-\boldsymbol q_\varepsilon|<s_\varepsilon\\
					\frac{s_\varepsilon^4}{|\boldsymbol x-\boldsymbol q_\varepsilon|^2}, \ &  |\boldsymbol x-\boldsymbol q_\varepsilon|\ge s_\varepsilon
				\end{array}
				\right.+O(\varepsilon^2|\ln\varepsilon|),
			\end{split}
		\end{equation*}
		where we have used the formula of planar Rankine vortex and integral
		\begin{equation*}
			\frac{1}{2\pi}\int_{B_1(\boldsymbol 0)}y_1'\ln\frac{1}{|\boldsymbol y-\boldsymbol y'|}d\boldsymbol y'=\left\{
			\begin{array}{lll}
				\frac{y_1}{4}-\frac{|\boldsymbol y|^2y_1}{8}, \ \ \ &\mathrm{if} \  |\boldsymbol y|<1,\\
				\frac{y_1}{8|\boldsymbol y|^2},  &\mathrm{if} \  |\boldsymbol y|\ge 1.
			\end{array}
			\right.\\
		\end{equation*}
		Let
		\begin{equation*}
			\begin{split}
				\mathcal R_\varepsilon(\boldsymbol x)&=\frac{1}{4\pi\varepsilon^2}\int_{ B_{s_\varepsilon}(\boldsymbol q_\varepsilon)}x_1^{1/2}x_1'^{3/2}\left(\ln(x_1x_1')+2\ln 8-4+O( \rho\ln(1/\rho))\right)d\boldsymbol x'\\
				&\quad-\frac{W}{2}x_1^2\ln\frac{1}{\varepsilon}-\mu_\varepsilon.
			\end{split}
		\end{equation*}
		By the choice of $a_\varepsilon$ in \eqref{acondition}, it holds
		\begin{equation*}
			\mathcal R_\varepsilon(\boldsymbol x)=\mathcal R_\varepsilon(\boldsymbol q_\varepsilon)+(x_1-q_\varepsilon)\cdot\partial_1\mathcal R_\varepsilon(\boldsymbol q_\varepsilon)+O(\varepsilon^2|\ln\varepsilon|)
		\end{equation*}
		with
		\begin{equation*}
			\mathcal R_\varepsilon(\boldsymbol q_\varepsilon)=-\frac{a_\varepsilon}{2\pi}\ln\frac{1}{\varepsilon},
		\end{equation*}
		and
		\begin{equation*}
			\begin{split}
				\partial_1\mathcal R_\varepsilon(\boldsymbol q_\varepsilon)&=\frac{1}{4\pi\varepsilon^2}\int_{ B_{s_\varepsilon}(\boldsymbol q_\varepsilon)}\left(\frac{x_1'^{3/2}}{2q_\varepsilon^{1/2}}\left(\ln(q_\varepsilon x_1')+2\ln 8-4\right)+\frac{x_1'^{3/2}}{q_\varepsilon^{1/2}}\right)d\boldsymbol x'-Wq_\varepsilon\ln\frac{1}{\varepsilon}\\
				&=\frac{s^2}{4\varepsilon^2} \cdot q_\varepsilon(\ln 8q_\varepsilon-1)-Wq_\varepsilon\ln\frac{1}{\varepsilon}+O(\varepsilon|\ln\varepsilon|).
			\end{split}
		\end{equation*}
		Combining all the facts above, we have
		$$\mathbf U_{\boldsymbol q_\varepsilon,\varepsilon}(\boldsymbol x)=V_{\boldsymbol q_\varepsilon,\varepsilon}(\boldsymbol x)-\frac{a_\varepsilon}{2\pi}\ln\frac{1}{\varepsilon}+(x_1-q_\varepsilon)\cdot \mathcal W(|\boldsymbol x-\boldsymbol q_\varepsilon|)+O(\varepsilon^2|\ln\varepsilon|).$$
		By letting $\boldsymbol x=s_\varepsilon\boldsymbol y+\boldsymbol q_\varepsilon$, the proof of Lemma \ref{B1} is then complete.
	\end{proof}
	
\begin{remark}\label{remarkB}
	Using the asymptotic property for $G_*$ in \eqref{2-12} and a similar technique as in the previous Lemma, a more careful computation can be carried out to show that the approximate function $\Phi_{\boldsymbol q_\varepsilon,\varepsilon}-Wx_1^2|\ln\varepsilon|/2$ for Stokes stream function $\psi_\varepsilon-Wx_1^2|\ln\varepsilon|/2$ can be rescaled and expanded as
	\begin{equation*}
		\begin{split}
			&\quad\Phi_{\boldsymbol q_\varepsilon,\varepsilon}(s_\varepsilon \boldsymbol y+q_\varepsilon)-\frac{W}{2}\cdot(s_\varepsilon y_1+q_\varepsilon)^2\ln\frac{1}{\varepsilon}\\
			&=\tilde V_{\boldsymbol q_\varepsilon,\varepsilon}(y_1,0)-\left(Wq_\varepsilon\ln\frac{1}{\varepsilon}-\frac{\kappa}{4\pi}\ln\frac{8q_\varepsilon}{s_\varepsilon}+\frac{\kappa}{16\pi}\right)\cdot s_\varepsilon y_1\\
			&\quad-\frac{3\kappa}{16\pi q_\varepsilon}\ln\frac{1}{\varepsilon}\cdot (s_\varepsilon |\boldsymbol y|)^2-\frac{W}{2}\ln\frac{1}{\varepsilon}\cdot(s_\varepsilon y_1)^2+O(\varepsilon^2).
		\end{split}
	\end{equation*}
    Using this fact together with the estimate for the $O(\varepsilon^2|\ln\varepsilon|)$-perturbation term $\phi_\varepsilon$ in \eqref{phiestimate} we can obtain \eqref{exp} in the introduction, which yields the uniqueness and nonlinear orbital stability of vortex ring by conditions (b) (c).
\end{remark}

\bigskip
	
Since $\nabla \tilde V_{\boldsymbol q_\varepsilon, \varepsilon}(\boldsymbol y)=s_\varepsilon\mathcal N_\varepsilon$ on $\partial B_{1}(\boldsymbol 0)$, and the difference of $\tilde{\mathbf U}_{\boldsymbol q_\varepsilon,\varepsilon}$ and $\tilde V_{\boldsymbol q_\varepsilon,\varepsilon}-\frac{a_\varepsilon}{2\pi}\ln\frac{1}{\varepsilon}$ is mainly given by the small term $s_\varepsilon y_1\cdot\mathcal W(|s_\varepsilon\boldsymbol y|)$ by the previous lemma, we can use the implicit function theorem to show that the level set $\{\boldsymbol y\ | \ \tilde{\mathbf U}_{\boldsymbol q_\varepsilon,\varepsilon}(\boldsymbol y)=0\}$ is almost a $t_{\varepsilon,\mathcal W}$-perturbation around $\{\boldsymbol y\ | \ \tilde V_{\boldsymbol q_\varepsilon,\varepsilon}(\boldsymbol y)-\frac{a_\varepsilon}{2\pi}\ln\frac{1}{\varepsilon}=0\}$, where the subscript $\mathcal W$ in $t_{\varepsilon,\mathcal W}$ means that the perturbation is induced by $\mathcal W_\varepsilon(s_\varepsilon)$.
	\begin{lemma}\label{B2}
		The set
		$$\tilde{\mathbf\Gamma}_\varepsilon:=\{\boldsymbol y\ | \ \tilde{\mathbf U}_{\boldsymbol q_\varepsilon,\varepsilon}(\boldsymbol y)=0\}$$
		is a closed convex curve in $\mathbb{R}^2$ parameterized as
		\begin{equation}\label{B-2}
			\begin{split}
				\tilde{\mathbf\Gamma}_\varepsilon (\theta)&=(1+t^*_{\varepsilon}(\theta))(\cos\theta,\sin\theta)\\
				&=(1+t_{\varepsilon,\mathcal W}(\theta)+O(\varepsilon^2|\ln\varepsilon|))(\cos\theta,\sin\theta), \,\,\,\theta\in [0,2\pi],
			\end{split}
         \end{equation}
where
\begin{equation*}
				t_{\varepsilon,\mathcal W}(\theta)(\cos\theta,\sin\theta)=\frac{\mathcal W_\varepsilon(s_\varepsilon)}{\mathcal N_\varepsilon}\cdot(\cos\theta,0), \quad \theta\in (0,2\pi]
\end{equation*}
with $\mathcal N_\varepsilon$ in \eqref{gradient}.  Moreover, one has
		\begin{equation*}
			\tilde{\mathbf{U}}_{\boldsymbol q_\varepsilon,\varepsilon}((1+t)(\cos\theta,\sin\theta))\left\{
			\begin{array}{lll}
				>0, \ & \mathrm{if} \  t<t^*_{\varepsilon}(\theta),\\
				<0, \ & \mathrm{if} \  t>t^*_{\varepsilon}(\theta),
			\end{array}
			\right.
		\end{equation*}
	and
 $$\|t_{\varepsilon,\mathcal W}(\theta)\|_{C^1((0,2\pi])}=O(\varepsilon|\ln\varepsilon|).$$ 
	\end{lemma}
	\begin{proof}
		In view of Lemma \ref{B1}, for every $\boldsymbol y\in D_\varepsilon=\{\boldsymbol y\ | \ s_\varepsilon\boldsymbol y+\boldsymbol q_\varepsilon\in\mathbb R^2_+\}$ bounded, it holds
		$$\tilde{\mathbf U}_{\boldsymbol q_\varepsilon,\varepsilon}(\boldsymbol y)=\tilde V_{\boldsymbol q_\varepsilon,\varepsilon}(\boldsymbol y)-\frac{a_\varepsilon}{2\pi}\ln\frac{1}{\varepsilon}+s_\varepsilon y_1\cdot \mathcal W_\varepsilon(|s_\varepsilon\boldsymbol y|)+O(\varepsilon^2|\ln\varepsilon|).$$
		Notice that
		\begin{equation*}
			\tilde V_{\boldsymbol q_\varepsilon,\varepsilon}(\boldsymbol y)=\left\{
			\begin{array}{lll}
				\frac{a_\varepsilon}{2\pi}\ln\frac{1}{\varepsilon}+\frac{q_\varepsilon^2s_\varepsilon^2}{4\varepsilon^2}(1-|\boldsymbol y|^2), \ \ \ & |\boldsymbol y|\le 1,\\
				\frac{a_\varepsilon}{2\pi}\ln\frac{1}{\varepsilon}\left(1+\frac{\ln|\boldsymbol y|}{\ln s_\varepsilon}\right),    & |\boldsymbol y|\ge 1,
			\end{array}
			\right.
		\end{equation*}
		and
		$$s_\varepsilon\mathcal W_\varepsilon(|s_\varepsilon\boldsymbol y|)=O(\varepsilon|\ln\varepsilon|).$$
		If $|\boldsymbol y|<1-L_1\varepsilon|\ln\varepsilon|$ for some large $L_1>0$, then
		\begin{equation*}
			\tilde{\mathbf{U}}_{\boldsymbol q_\varepsilon,\varepsilon}(\boldsymbol y)\ge \frac{q_\varepsilon^2s_\varepsilon^2}{4\varepsilon^2}(1-|1-L_1\varepsilon|\ln\varepsilon||^2)+O(\varepsilon|\ln\varepsilon|)>0.
		\end{equation*}
		If $|\boldsymbol y|>1+L_2\varepsilon|\ln\varepsilon|$ for some large $L_2>0$, then
		\begin{equation*}
			\tilde{\mathbf{U}}_{\boldsymbol q_\varepsilon,\varepsilon}(\boldsymbol y)\le
			\frac{a_\varepsilon}{2\pi}\ln\frac{1}{\varepsilon}\cdot\frac{\ln|1+L_2\varepsilon|\ln\varepsilon||}{\ln s_\varepsilon}<0.
		\end{equation*}
		So we have proved that for any $(\cos\theta,\sin\theta)$, there exists a $t^*_\varepsilon(\theta)$, such that $|t^*_\varepsilon(\theta)|=O(\varepsilon|\ln\varepsilon|)$, and
		$$(1+t^*_\varepsilon)(\cos\theta,\sin\theta)\in\tilde{\mathbf\Gamma}_\varepsilon(\theta).$$
		On the other hand, it holds
		\begin{equation*}
			\frac{\partial\tilde{\mathbf{U}}_{\boldsymbol q_\varepsilon,\varepsilon}((1+t^*_\varepsilon)(\cos\theta,\sin\theta))}{\partial t}\bigg|_{t=0}=-s_\varepsilon\mathcal N_\varepsilon+O(\varepsilon|\ln\varepsilon|)=-\frac{s_\varepsilon^2q_\varepsilon^2}{2\varepsilon^2}+O(\varepsilon|\ln\varepsilon|)<0.
		\end{equation*}
		By the implicit function theorem, we see that $t^*_\varepsilon(\theta)$ is unique, and satisfies
		\begin{equation*}
			t^*_\varepsilon(\theta)=\frac{\cos\theta\cdot s_\varepsilon\mathcal W_\varepsilon(s_\varepsilon)+O(\varepsilon^2|\ln\varepsilon|)}{s_\varepsilon\mathcal N_\varepsilon+O_\varepsilon(\varepsilon|\ln\varepsilon|)}.
		\end{equation*}
		Hence it holds
		$$t^*_\varepsilon(\theta)=\frac{\cos\theta}{\mathcal N_\varepsilon}\cdot\mathcal W_\varepsilon(s)+O(\varepsilon^2|\ln\varepsilon|),$$
		and \eqref{B-2} is verified.
		
		To obtain an estimate for $t^{*'}_\varepsilon(\theta)$, we differentiate $\tilde{\mathbf U}_{\boldsymbol q_\varepsilon,\varepsilon}((1+t^*_\varepsilon)(\cos\theta,\sin\theta))=0$ with respect to $\theta$ and derive
		\begin{equation*}
			\frac{\partial\tilde{\mathbf{U}}_{\boldsymbol q_\varepsilon,\varepsilon}((1+t)(\cos\theta,\sin\theta))}{\partial \theta}\bigg|_{t=t^*_\varepsilon}=O(\varepsilon)\cdot |t^{*'}_\varepsilon(\theta)|+O(\varepsilon|\ln\varepsilon|).
		\end{equation*}
		Using the implicit function theorem again, we have
		\begin{equation*}
			\frac{\partial\tilde{\mathbf{U}}_{\boldsymbol q_\varepsilon,\varepsilon}((1+t)(\cos\theta,\sin\theta))}{\partial \theta}\bigg|_{t=t^*_\varepsilon}=(s_\varepsilon\mathcal N_\varepsilon+O(\varepsilon|\ln\varepsilon|))\cdot t^{*'}_\varepsilon(\theta).
		\end{equation*}
		We conclude that $|t^{*'}_\varepsilon(\theta)|=O(\varepsilon|\ln\varepsilon|)$, and hence $\tilde{\mathbf\Gamma}_\varepsilon$ is a closed convex curve.
	\end{proof}

\bigskip
	
	Noting that $\mathbf U_{\boldsymbol q_\varepsilon,\varepsilon}(\boldsymbol x)+\phi_\varepsilon(\boldsymbol x)=\psi_\varepsilon(\boldsymbol x)-\frac{W}{2}x_1^2\ln\frac{1}{\varepsilon}-\mu_\varepsilon$, now we can also estimate the cross-section boundary 
	$$\partial \tilde A_\varepsilon=\left\{\boldsymbol y\ | \ \tilde\psi_\varepsilon(\boldsymbol y)-\frac{W}{2}(s_\varepsilon y_1+q_\varepsilon)^2\ln\frac{1}{\varepsilon}-\mu_\varepsilon=0\right\}$$
	as a almost $t_{\varepsilon,\mathcal W}+t_{\varepsilon,\tilde\phi_\varepsilon}$ perturbation near $B_1(\boldsymbol 0)$, where $t_{\varepsilon,\mathcal W}$ is given in former lemma, and $t_{\varepsilon,\tilde\phi_\varepsilon}$ is induced by scaled error function $\tilde\phi_\varepsilon$. Without loss of generality, here we use a general small $\tilde\phi$ instead of $\tilde\phi_\varepsilon$.
	\begin{lemma}\label{B3}
		Suppose that $\tilde\phi$ is a function satisfying
		\begin{equation}\label{B-3}
			\|\nabla\tilde\phi\|_{L^\infty(B_L(\boldsymbol 0))}=O(\varepsilon|\ln\varepsilon|), \ \ \ \ \ \ \|\tilde\phi\|_{L^\infty(B_L(\boldsymbol 0))}=O(\varepsilon|\ln\varepsilon|).
		\end{equation}
		Then the set
		$$\tilde{\mathbf\Gamma}_{\varepsilon,\tilde\phi}:=\{\boldsymbol y\ | \ \tilde{\mathbf U}_{\boldsymbol q_\varepsilon,\varepsilon}(\boldsymbol y)+\tilde\phi=0\}$$
		is a closed convex curve in $\mathbb{R}^2$, and
		\begin{equation}\label{B-4}
			\begin{split}
				\tilde{\mathbf\Gamma}_{\varepsilon,\tilde\phi}(\theta)&=(1+t_\varepsilon(\theta))(\cos\theta,\sin\theta)\\
				&=(1+t_{\varepsilon,\tilde\phi}+t_{\varepsilon,\mathcal W}+O(\varepsilon^2|\ln\varepsilon|))(\cos\theta,\sin\theta)\\
				&=\left(1+\frac{1}{s_\varepsilon\mathcal N_\varepsilon}\cdot\tilde\phi(\cos\theta,\sin\theta)\right)(\cos\theta,\sin\theta)+\frac{\mathcal W_\varepsilon(s_\varepsilon)}{\mathcal N_\varepsilon}\cdot(\cos\theta,0)\\
				&\quad+O(\varepsilon^2|\ln\varepsilon|), \quad \theta\in (0,2\pi]
			\end{split}
		\end{equation}
		with $\mathcal N_\varepsilon$ in \eqref{gradient}. Moreover, one has
		\begin{equation*}
			(\tilde{\mathbf{U}}_{\boldsymbol q_\varepsilon,\varepsilon}+\tilde\phi)((1+t_\varepsilon)(\cos\theta,\sin\theta))\left\{
			\begin{array}{lll}
				>0, \ & \mathrm{if} \ t<t_\varepsilon(\theta),\\
				<0, \ &\mathrm{if} \ t>t_\varepsilon(\theta),
			\end{array}
			\right.
		\end{equation*}
		and
		\begin{equation}\label{B-5}
			\left|\tilde{\mathbf\Gamma}_{\varepsilon,\tilde\phi_1}-\tilde{\mathbf\Gamma}_{\varepsilon,\tilde\phi_2}\right|=\left(\frac{1}{s_\varepsilon\mathcal N_\varepsilon}+O(\varepsilon|\ln\varepsilon|^2)\right)\|\tilde\phi_1-\tilde\phi_2\|_{L^\infty(B_L(\boldsymbol 0))}
		\end{equation}
		for functions $\tilde\phi_1,\tilde\phi_2$ satisfying \eqref{B-3}.
	\end{lemma}
	\begin{proof}
		From Lemma \ref{B1}, we have
		\begin{equation*}
			\tilde{\mathbf U}_{\boldsymbol q_\varepsilon,\varepsilon}(\boldsymbol y)+\tilde\phi=\tilde V_{\boldsymbol q_\varepsilon,\varepsilon}-\frac{a_\varepsilon}{2\pi}\ln\frac{1}{\varepsilon}+s_\varepsilon y_1\cdot \mathcal W_\varepsilon (|s_\varepsilon\boldsymbol y|)+\tilde\phi+O(\varepsilon^2|\ln\varepsilon|).
		\end{equation*}
		Hence it holds
		\begin{equation*}
			1+t_\varepsilon\in(1-L_1\varepsilon|\ln\varepsilon|^2,1+L_2\varepsilon|\ln\varepsilon|^2)
		\end{equation*}
		in a similar way as Lemma \ref{B2}. Using the fact
		\begin{equation*}
			\frac{(\partial\tilde{\mathbf U}_{\boldsymbol q_\varepsilon,\varepsilon}+\partial\tilde\phi)((1+t_{\varepsilon,\mathcal W}+t)(\cos\theta,\sin\theta))}{\partial t}\bigg|_{t=0}=-s_\varepsilon\mathcal N_\varepsilon+O(\varepsilon|\ln\varepsilon|)<0,
		\end{equation*}
		we see that $t_{\varepsilon}$ is unique, and $\tilde{\mathbf\Gamma}_{\varepsilon,\tilde\phi}$ is a continuous closed curve in $\mathbb{R}^2$. Then we let
		\begin{equation*}
			\boldsymbol y_\varepsilon=(1+t_\varepsilon)(\cos\theta,\sin\theta)\in \tilde{\mathbf\Gamma}_{\varepsilon,\tilde\phi}(\theta).
		\end{equation*}
		By the implicit function theorem, it holds
		\begin{equation*}
			|\boldsymbol y_\varepsilon|-1=\frac{\cos\theta\cdot s_\varepsilon\mathcal W_\varepsilon(s_\varepsilon)+\tilde\phi(\boldsymbol y_\varepsilon)+(t_{\varepsilon,\mathcal W}+t_{\varepsilon,\tilde\phi})\cdot O(\varepsilon|\ln\varepsilon|)+O(\varepsilon^2|\ln\varepsilon|)}{s_\varepsilon\mathcal N_\varepsilon+(t_{\varepsilon,\mathcal W}+t_{\varepsilon,\tilde\phi})\cdot O_\varepsilon(1)}.
		\end{equation*}
		While for $\tilde\phi(\boldsymbol y_\varepsilon)$, it holds
		\begin{equation*}
			|\tilde\phi(\boldsymbol y_\varepsilon)-\tilde\phi(\cos\theta,\sin\theta)|\le\|\nabla\tilde\phi\|_{L^\infty(B_L(\boldsymbol 0))}\cdot |t_{\varepsilon,\mathcal W}(\theta)+t_{\varepsilon,\tilde\phi}(\theta)|,
		\end{equation*}
		from which we can verify \eqref{B-4}. Moreover, we can obtain $|t_{\varepsilon,\mathcal W}'(\theta)+t_{\varepsilon,\tilde\phi}'(\theta)|=O(\varepsilon|\ln\varepsilon|^2)$ as in Lemma \ref{B2}. So $\tilde{\mathbf\Gamma}_{\varepsilon,\tilde\phi}$ is also convex.
		
		Denote $\boldsymbol y_{\varepsilon,m}$ as the coordinate corresponding to $\tilde\phi_m$ ($m=1,2$). Then according to the definition of $\boldsymbol y_{\varepsilon,m}$, we have
		\begin{equation*}
			\begin{split}
				\tilde{\mathbf U}_{\boldsymbol q_\varepsilon,\varepsilon}(\boldsymbol y_{\varepsilon,1})-\tilde{\mathbf U}_{\boldsymbol q_\varepsilon,\varepsilon}(\boldsymbol y_{\varepsilon,2})&=\tilde\phi_1(y_{\varepsilon,1})-\tilde\phi_2(y_{\varepsilon,1})+\tilde\phi_2(y_{\varepsilon,1})-\tilde\phi_2(y_{\varepsilon,2})\\
				&=\|\tilde\phi_1-\tilde\phi_2\|_{L^\infty(B_L(\boldsymbol 0))}+\|\nabla\tilde\phi\|_{L^\infty(B_L(\boldsymbol 0))}\cdot|\boldsymbol y_{\varepsilon,1}-\boldsymbol y_{\varepsilon,2}|\\
				&=\|\tilde\phi_1-\tilde\phi_2\|_{L^\infty(B_L(\boldsymbol 0))}+O(\varepsilon|\ln\varepsilon|^2)\cdot|\boldsymbol y_{\varepsilon,1}-\boldsymbol y_{\varepsilon,2}|.
			\end{split}
		\end{equation*}
		Since
		\begin{equation*}
			\frac{\partial\tilde{\mathbf U}_{\boldsymbol q_\varepsilon,\varepsilon}((1+t_\varepsilon+t)(\cos\theta,\sin\theta))}{\partial t}\bigg|_{t=0}=-s_\varepsilon\mathcal N_\varepsilon+O(\varepsilon|\ln\varepsilon|),
		\end{equation*}
		we conclude \eqref{B-5} and finish our proof.
	\end{proof}

\bigskip
	
	In Section 3 for uniqueness of steady vortex rings, we have used a coarse version of Lemma \ref{B3} to show the nonlinear part $R_\varepsilon(\phi_\varepsilon)$ satisfies
	\begin{equation*}
		R_\varepsilon(\phi_\varepsilon)=0, \ \ \ \text{in} \ \left(\mathbb R^2_+\setminus B_{2s_\varepsilon}(\boldsymbol q_\varepsilon)\right) \cup B_{s_\varepsilon/2}(\boldsymbol q_\varepsilon),
	\end{equation*}
	and to derive the estimate for $\phi_\varepsilon$ in Lemma \ref{lem3-12}, which is summarized as follows. Since the proof is similar to Lemma \ref{B3}, we omit the detail here therefore.
	\begin{lemma}\label{B4}
		Suppose that $\tilde\phi$ is a function satisfying
		\begin{equation*}
			\|\nabla\tilde\phi\|_{L^\infty(B_L(\boldsymbol 0))}= o_{\varepsilon}(1), \ \ \ \ \ \ \|\tilde\phi\|_{L^\infty(B_L(\boldsymbol 0))}= o_\varepsilon(1).
		\end{equation*}
		Then the set
		$$\tilde{\mathbf\Gamma}_{\varepsilon,\tilde\phi}:=\{\boldsymbol y\ | \ \tilde{\mathbf U}_{\boldsymbol q_\varepsilon,\varepsilon}(\boldsymbol y)+\tilde\phi=0\}$$
		is a closed convex curve in $\mathbb{R}^2$, and
		\begin{equation}\label{B-6}
			\begin{split}
				\tilde{\mathbf\Gamma}_{\varepsilon,\tilde\phi}(\theta)&=(1+t_\varepsilon(\theta))(\cos\theta,\sin\theta)\\
				&=(1+t_{\varepsilon,\tilde\phi}+t_{\varepsilon,\mathcal W}+o(\varepsilon))(\cos\theta,\sin\theta)\\
				&=\left(1+\frac{1}{s_\varepsilon\mathcal N_\varepsilon}\cdot\tilde\phi(\cos\theta,\sin\theta)\right)(\cos\theta,\sin\theta)+\frac{\mathcal W_\varepsilon(s_\varepsilon)}{\mathcal N_\varepsilon}\cdot(\cos\theta,0)\\
				&\quad+o(\varepsilon), \quad \theta\in (0,2\pi]
			\end{split}
		\end{equation}
		with $\mathcal N_\varepsilon$ in \eqref{gradient}.
	\end{lemma}

    \bigskip
	
	\section{Estimates for the Pohozaev identity}\label{appC}
	
	This appendix is devoted to the estimates for $\mathcal W_\varepsilon(s_\varepsilon)$ or $t_{\varepsilon,\mathcal W}(\theta)$ that have been used in obtaining the uniqueness of steady vortex rings in Section 3.  To this aim, we will apply the following local Pohozaev identity, which corresponds to the translation transformation of semilinear elliptic equations, and contains the information of first order derivatives for the approximation in Appendix \ref{appB}, namely, the coefficient $\mathcal W_\varepsilon(s_\varepsilon)$ in linear part of the difference.
	\begin{lemma}\label{C1}
		Suppose that $u\in H^1(\mathbb R^2_+)\cap C^1(\mathbb R^2_+)$ is a weak solution to
		\begin{equation*}
			-\Delta u=f(\boldsymbol x,u), \ \ \ \mathrm{in} \ \mathbb R^2_+,
		\end{equation*}
		where $f(\boldsymbol x,u)$ is a function continuous in $\boldsymbol x$, and nondecreasing with respect to $u$. Then for any bounded smooth domain $D\subset \mathbb R^2_+$, it holds
		\begin{equation*}
			\int_{\partial D}\frac{\partial u}{\partial x_i}\frac{\partial u}{\partial \nu}dS-\frac{1}{2}\int_{\partial D}|\nabla u|^2\nu_idS+\int_{\partial D} F(\boldsymbol x,u)dS=\int_D F_{x_i}(\boldsymbol x,u)d\boldsymbol x, \quad i=1,2,
		\end{equation*}
		where $\nu$ the unit outward normal to the boundary $\partial D$, and $F(\boldsymbol x, u)$ is defined by
		$$F(\boldsymbol x, u):=\int_0^u f(\boldsymbol x,t)dt.$$
	\end{lemma}
	
	\bigskip
	
	The proof of Lemma \ref{C1} can be found in \cite{CPYb} (see Theorem 6.2.1 in \cite{CPYb}) together with an approximation procedure. In our case, we let the domain $D\subset \mathbb R^2_+$ be $B_\delta(\boldsymbol q_\varepsilon)$ with a small constant $\delta>0$, let the function $u$ be $\psi_{1,\varepsilon}$, and let the nonlinearity $f$ be
	$$f(\boldsymbol x, \psi_{1,\varepsilon})=\frac{q_\varepsilon^2}{\varepsilon^2}\cdot \boldsymbol 1_{\{\psi_\varepsilon-\frac{W}{2}x_1^2\ln\frac{1}{\varepsilon}>\mu_\varepsilon\}}.$$
	Thus the primitive function for $f$ is
	$$F(\boldsymbol x, \psi_{1,\varepsilon})=\frac{q_\varepsilon^2}{\varepsilon^2}\cdot\left(\psi_\varepsilon-\frac{W}{2}x_1^2\ln\frac{1}{\varepsilon}-\mu_\varepsilon\right)_+,$$
	and the local Pohozaev identity in Lemma \ref{C1} with $i=1$ turns to be
	\begin{equation}\label{C-1}
		\begin{split}
			&\quad-\int_{\partial B_\delta(\boldsymbol q_\varepsilon)} \frac{\partial\psi_{1,\varepsilon}}{\partial \nu}\frac{\partial\psi_{1,\varepsilon}}{\partial x_1}dS+ \frac{1}{2}\int_{\partial B_\delta(\boldsymbol q_\varepsilon)} |\nabla\psi_{1,\varepsilon}|^2 \nu_1dS\\
			&=-\frac{q_\varepsilon^2}{\varepsilon^2}\int_{B_\delta(\boldsymbol q_\varepsilon)} \partial_1\psi_{2,\varepsilon}(\boldsymbol x)\cdot\boldsymbol 1_{A_\varepsilon}(\boldsymbol x)d\boldsymbol x+\frac{q_\varepsilon^2}{\varepsilon^2}\int_{B_\delta(\boldsymbol q_\varepsilon)} Wx_1\ln\frac{1}{\varepsilon}\cdot\boldsymbol 1_{A_\varepsilon}(\boldsymbol x)d\boldsymbol x
		\end{split}
	\end{equation}
	with the cross-section 
	$$A_\varepsilon:=\left\{\boldsymbol x=(x_1,x_2)\in \mathbb R^2_+ \ \big| \ \psi_\varepsilon-\frac{W}{2}x_1^2\ln\frac{1}{\varepsilon}>\mu_\varepsilon\right\}.$$
	
	According to the estimates obtained in Section 3, we see that $A_\varepsilon$ is an area tending to the disk $B_{s_\varepsilon^*}(\boldsymbol q_\varepsilon)$ with radius
	$$s_\varepsilon^*=\sqrt{\varepsilon^2\kappa/q_\varepsilon\pi}.$$
	By denoting the symmetry difference
	$$A_\varepsilon\bigtriangleup B_{s_\varepsilon^*}(\boldsymbol q_\varepsilon):=\left(A_\varepsilon\setminus B_{s_\varepsilon^*}(\boldsymbol q_\varepsilon)\right)\cup\left(B_{s_\varepsilon^*}(\boldsymbol q_\varepsilon)\setminus A_\varepsilon\right),$$
	and the error
	$$\mathbf {er}_\varepsilon:=|A_\varepsilon\bigtriangleup B_{s_\varepsilon^*}(\boldsymbol q_\varepsilon)|,$$
	we will proceed a series of lemmas to compute each term in \eqref{C-1}, and obtain a relationship of $\mathbf {er}_\varepsilon$ with $\mathcal W_\varepsilon(s_\varepsilon)$ or $t_{\varepsilon,\mathcal W}(\theta)$ in Appendix \ref{appB}. 
	
	\bigskip
	
	To estimate the left hand side of \eqref{C-1}, we have following lemma concerning the asymptotic behavior of $\psi_\varepsilon$ away from $A_\varepsilon$. 
	\begin{lemma}\label{C2}
		For every $\boldsymbol x\in \mathbb{R}^2_+\setminus \{\boldsymbol x \  | \ \mathrm{dist}(\boldsymbol x,A_\varepsilon)\le L\varepsilon\}$, we have
		\begin{equation*}
			\psi_{1,\varepsilon}(\boldsymbol x)=\frac{\kappa}{2\pi}\cdot q_\varepsilon\ln \frac{|\boldsymbol x-\boldsymbol {\bar q}_\varepsilon|}{|\boldsymbol x-\boldsymbol q_\varepsilon|}+O\left(\frac{\mathbf {er}_\varepsilon}{\varepsilon|\boldsymbol x-\boldsymbol q_\varepsilon|}\right),
		\end{equation*}
		and
		\begin{equation*}
			\nabla \psi_{1,\varepsilon}(\boldsymbol x)=-\frac{\kappa}{2\pi}\cdot q_\varepsilon\frac{\boldsymbol x-\boldsymbol q_\varepsilon}{|\boldsymbol x-\boldsymbol q_\varepsilon|^2}+\frac{\kappa}{2\pi}\cdot q_\varepsilon\frac{\boldsymbol x-\boldsymbol {\bar q}_\varepsilon}{|\boldsymbol x-\boldsymbol {\bar q}_\varepsilon|^2}+O\left(\frac{\mathbf {er}_\varepsilon}{\varepsilon|\boldsymbol x-\boldsymbol q_\varepsilon|^2}\right).
		\end{equation*}
	\end{lemma}
	\begin{proof}
		For each $\boldsymbol x\in \mathbb{R}^2\setminus \{\boldsymbol x\ | \ \mathrm{dist}(\boldsymbol x,A_\varepsilon)\le L\varepsilon\}$, it must hold $\boldsymbol x \notin A_\varepsilon$. Then, using  Taylor's formula
		\begin{equation*}
			|\boldsymbol x-\boldsymbol x'|=|\boldsymbol x-\boldsymbol q_\varepsilon|-\langle\frac{\boldsymbol x-\boldsymbol q_\varepsilon}{|\boldsymbol x-\boldsymbol q_\varepsilon|},\boldsymbol x'-\boldsymbol q_\varepsilon\rangle+O\left(\frac{|\boldsymbol x'-\boldsymbol q_\varepsilon|^2}{|\boldsymbol x-\boldsymbol q_\varepsilon|}\right), \quad \forall\, \boldsymbol x'\in A_\varepsilon,
		\end{equation*}
		we obtain
		\begin{equation*}
			\begin{split}
				\psi_{1,\varepsilon}(\boldsymbol x)&=\frac{q_\varepsilon^2}{2\pi\varepsilon^2}\int_{A_\varepsilon}\ln \left(\frac{|\boldsymbol x-\boldsymbol {\bar x}'|}{|\boldsymbol x-\boldsymbol x'|}\right)d\boldsymbol x'\\
				&=\frac{\kappa}{2\pi}\cdot q_\varepsilon\ln \frac{|\boldsymbol x-\boldsymbol {\bar q}_\varepsilon|}{|\boldsymbol x-\boldsymbol q_\varepsilon|}+\frac{q_\varepsilon^2}{2\pi\varepsilon^2}\int_{A_\varepsilon}\ln \left(\frac{|\boldsymbol x-\boldsymbol q_\varepsilon|}{|\boldsymbol x-\boldsymbol x'|}\right)d\boldsymbol x'\\
				&\quad\quad\quad\quad\quad\quad-\frac{q_\varepsilon^2}{2\pi\varepsilon^2}\int_{A_\varepsilon}\ln \left(\frac{|\boldsymbol x-\boldsymbol {\bar q}_\varepsilon|}{|\boldsymbol x-\boldsymbol {\bar x}'|}\right)d\boldsymbol x'+O\left(\frac{\mathbf {er}_\varepsilon}{\varepsilon|\boldsymbol x-\boldsymbol q_\varepsilon|}\right)\\
				&=\frac{\kappa}{2\pi}\cdot q_\varepsilon\ln \frac{|\boldsymbol x-\boldsymbol {\bar q_\varepsilon}|}{|\boldsymbol x-\boldsymbol q_\varepsilon|}-\frac{q_\varepsilon^2}{2\pi\varepsilon^2}\int_{A_\varepsilon}\frac{(\boldsymbol x-\boldsymbol q_\varepsilon)\cdot(\boldsymbol q_\varepsilon-\boldsymbol x')}{|\boldsymbol x-\boldsymbol q_\varepsilon|^2}d\boldsymbol x'\\
				&\quad\quad\quad\quad\quad\quad +\frac{q_\varepsilon^2}{2\pi\varepsilon^2}\int_{A_\varepsilon}\frac{(\boldsymbol x-\boldsymbol {\bar q}_\varepsilon)\cdot(\boldsymbol {\bar q}_\varepsilon-\boldsymbol {\bar x}')}{|\boldsymbol x-\boldsymbol {\bar q}_\varepsilon|^2}d\boldsymbol x'+O\left(\frac{\mathbf {er}_\varepsilon}{\varepsilon|\boldsymbol x-\boldsymbol q_\varepsilon|}\right).
			\end{split}
		\end{equation*}
		Using the odd symmetry, we have
		\begin{align*}
			& \ \ \ \int_{A_\varepsilon} \frac{(\boldsymbol x-\boldsymbol q_\varepsilon)\cdot(\boldsymbol q_\varepsilon-\boldsymbol x')}{|\boldsymbol x-\boldsymbol q_\varepsilon|^2}d\boldsymbol x'\\
			&=\left(\int_{A_\varepsilon\setminus B_{s_\varepsilon^*}(\boldsymbol q_\varepsilon)}-\int_{B_{s_\varepsilon^*}(\boldsymbol q_\varepsilon)\setminus A_\varepsilon}\right)\frac{(\boldsymbol x-\boldsymbol q_\varepsilon)\cdot(\boldsymbol q_\varepsilon-\boldsymbol x')}{|\boldsymbol x-\boldsymbol q_\varepsilon|^2}d\boldsymbol x'+ \int_{ B_{s_\varepsilon^*}(\boldsymbol q_\varepsilon)} \frac{(\boldsymbol x-\boldsymbol q_\varepsilon)\cdot(\boldsymbol q_\varepsilon-\boldsymbol x')}{|\boldsymbol x-\boldsymbol q_\varepsilon|^2}d\boldsymbol x'\\
			&=\left(\int_{A_\varepsilon\setminus B_{s_\varepsilon^*}(\boldsymbol q_\varepsilon)}-\int_{B_{s_\varepsilon^*}(\boldsymbol q_\varepsilon)\setminus A_\varepsilon}\right)\frac{(\boldsymbol x-\boldsymbol q_\varepsilon)\cdot(\boldsymbol q_\varepsilon-\boldsymbol x')}{|\boldsymbol x-\boldsymbol q_\varepsilon|^2}d\boldsymbol x'\\
			&=O\left(\frac{\varepsilon}{|\boldsymbol x-\boldsymbol q_\varepsilon|}\right)\cdot |A_\varepsilon\bigtriangleup B_{s_\varepsilon^*}(\boldsymbol q_\varepsilon)|=O\left(\frac{\varepsilon\cdot \mathbf {er}_\varepsilon}{|\boldsymbol x-\boldsymbol q_\varepsilon|}\right).
		\end{align*}
		While, for the other terms, we can use a same argument to deduce
		\begin{equation*}
			\int_{A_\varepsilon}\frac{(\boldsymbol x-\boldsymbol {\bar q}_\varepsilon)\cdot(\boldsymbol {\bar q}_\varepsilon-\boldsymbol {\bar x}')}{|\boldsymbol x-\boldsymbol {\bar q}_\varepsilon|^2}d\boldsymbol x'=O\left(\frac{\varepsilon\cdot\mathbf {er}_\varepsilon}{|\boldsymbol x-\boldsymbol {\bar q}_\varepsilon|}\right)=O(\varepsilon\cdot\mathbf {er}_\varepsilon).
		\end{equation*}
		Hence we have verified the first part of this lemma. The second part can be verified by similar procedure.
	\end{proof}
	
	Using Lemma \ref{C2}, we can compute the left hand side of \eqref{C-1} as follows.
	\begin{lemma}\label{C3}
		It holds
		\begin{equation*}
			-\int_{\partial B_\delta(\boldsymbol q_\varepsilon)} \frac{\partial\psi_{1,\varepsilon}}{\partial \nu}\frac{\partial\psi_{1,\varepsilon}}{\partial x_1}dS+ \frac{1}{2}\int_{\partial B_\delta(\boldsymbol q_\varepsilon)} |\nabla\psi_{1,\varepsilon}|^2 \nu_1dS=\kappa\cdot\frac{s_\varepsilon^{*2}}{4\varepsilon^2}\cdot q_\varepsilon^2+O\left(\frac{\mathbf {er}_\varepsilon}{\varepsilon}\right).
		\end{equation*}
	\end{lemma}
	\begin{proof}
		Using the identity
		\begin{equation*}
			-\int_{\partial B_\delta(\boldsymbol q_\varepsilon)}\frac{G(\boldsymbol x,\boldsymbol x')}{\partial \nu}\frac{G(\boldsymbol x,\boldsymbol x')}{\partial x_1}dS+\frac{1}{2}\int_{\partial B_\delta(\boldsymbol q_\varepsilon)} |\nabla G(\boldsymbol x,\boldsymbol x')|^2 \nu_1dS
			=-\partial_1\left(\frac{1}{2\pi}\ln\frac{1}{|\boldsymbol x-\boldsymbol {\bar q}_\varepsilon|}\right)\bigg|_{\boldsymbol x=\boldsymbol q_\varepsilon},
		\end{equation*}
		and the asymptotic estimate in Lemma \ref{C2}, this lemma can be verified by direct computation.
	\end{proof}

\bigskip

	Using the circulation constraint \eqref{3-2}, it is obvious that
	\begin{equation}\label{C-2}
		\frac{q_\varepsilon^2}{\varepsilon^2}\int_{B_\delta(\boldsymbol q_\varepsilon)} W x_1\ln\frac{1}{\varepsilon}\cdot\boldsymbol 1_{A_\varepsilon}(\boldsymbol x)d\boldsymbol x=\kappa \cdot Wq_\varepsilon^2\ln\frac{1}{\varepsilon}.
	\end{equation}
	Thus we will focus on the first term in the right hand side of \eqref{C-1} relevant to $\partial_1\psi_{2,\varepsilon}$.
	\begin{lemma}\label{C4}
		It holds
		\begin{equation*}
			-\frac{q_\varepsilon^2}{\varepsilon^2}\int_{B_\delta(\boldsymbol q_\varepsilon)} \partial_1\psi_{2,\varepsilon}(\boldsymbol x)\cdot\boldsymbol 1_{A_\varepsilon}(\boldsymbol x)d\boldsymbol x=-\kappa\cdot\frac{s_\varepsilon^{*2}}{4\varepsilon^2}\cdot q_\varepsilon^2\left(\ln\frac{8q_\varepsilon}{s_\varepsilon^*}-\frac{5}{4}\right)+O\left(\frac{\mathbf {er}_\varepsilon}{\varepsilon^2}+\varepsilon^2|\ln\varepsilon|\right).
		\end{equation*}
	\end{lemma}
	\begin{proof}
		By the definition of $\partial_1\psi_{2,\varepsilon}$, it holds
		$$\partial_1\psi_{2,\varepsilon}=\frac{1}{\varepsilon^2}\int_{\mathbb R^2_+}\partial_{x_1}H(\boldsymbol x,\boldsymbol x')\boldsymbol 1_{A_\varepsilon}(\boldsymbol x')d\boldsymbol x',$$
		where
		\begin{equation*}
			\begin{split}
				H(\boldsymbol x,\boldsymbol x')&=\left(\frac{x_1^{1/2}x_1'^{3/2}}{2\pi}-\frac{q_\varepsilon^2}{2\pi}\right)\ln\frac{1}{|\boldsymbol x-\boldsymbol x'|}+\frac{q_\varepsilon^2}{2\pi}\ln\frac{1}{|\boldsymbol x-\boldsymbol {\bar x}'|}\\
				&\quad+\frac{x_1^{1/2}x_1'^{3/2}}{4\pi}\left(\ln(x_1x_1')+2\ln 8-4+\boldsymbol \rho\right),
			\end{split}
		\end{equation*}
		with $\boldsymbol \rho=O\left(\rho\ln(1/\rho)\right)$ a regular remainder and $\rho$ defined in \eqref{rho}. For simplicity, we let
		\begin{equation*}
			-\frac{q_\varepsilon^2}{\varepsilon^2}\int_{B_\delta(\boldsymbol q_\varepsilon)} \partial_1\psi_{2,\varepsilon}(\boldsymbol x)\cdot\boldsymbol 1_{A_\varepsilon}(\boldsymbol x)d\boldsymbol x=I_1+I_2+I_3+I_4,
		\end{equation*}
		where
		\begin{equation*}
			 I_1=-\frac{q_\varepsilon^2}{4\pi\varepsilon^4}\int_{A_\varepsilon}x_1^{-1/2}\int_{A_\varepsilon}x_1'^{3/2}\ln\left(\frac{1}{s_\varepsilon^*}\right)d\boldsymbol x'd\boldsymbol x,
		\end{equation*}
		\begin{equation*}
			I_2=-\frac{q_\varepsilon^2}{4\pi\varepsilon^4}\int_{A_\varepsilon}x_1^{-1/2}\int_{A_\varepsilon}x_1'^{3/2}\ln\left(\frac{s_\varepsilon^*}{|\boldsymbol x-\boldsymbol x'|}\right)d\boldsymbol x'd\boldsymbol x,
		\end{equation*}
		\begin{equation*}
			 I_3=\frac{q_\varepsilon^2}{2\pi\varepsilon^4}\int_{A_\varepsilon}\int_{A_\varepsilon}\left(x_1^{1/2}x_1'^{3/2}-q_\varepsilon^2\right)\cdot\frac{x_1-x_1'}{|\boldsymbol x-\boldsymbol x'|^2}d\boldsymbol x'd\boldsymbol x,
		\end{equation*}
		and $I_4$ the remaining regular terms.
		
		Let us consider $I_1$ first. Using Taylor's expansion, $I_1$ can be rewritten as
		\begin{equation*}
			\begin{split}
				I_1&=-\frac{q_\varepsilon^2}{4\pi\varepsilon^4}\cdot\ln\frac{1}{s_\varepsilon^*}\cdot\int_{A_\varepsilon}x_1 \left (q_\varepsilon^{-3/2}-\frac{3}{2q_\varepsilon^{5/2}}\cdot(x_1-q_\varepsilon)+O(|x_1-q_\varepsilon|^2)\right)d\boldsymbol x\\
				&\quad \times\int_{A_\varepsilon}x_1'\left(q_\varepsilon^{1/2}+\frac{1}{2q_\varepsilon^{1/2}}\cdot(x'_1-q_\varepsilon)+O(|x'_1-q_\varepsilon|^2)\right)d\boldsymbol x'.
			\end{split}
		\end{equation*}
		Then, we are to estimate each terms in the product. Using circulation constraint \eqref{3-2}, we have
		\begin{equation*}
			\quad\frac{q_\varepsilon}{4\pi\varepsilon^4}\cdot\ln\frac{1}{s_\varepsilon^*}\cdot\int_{A_\varepsilon}x_1d\boldsymbol x\int_{A_\varepsilon}x_1'd\boldsymbol x'=\kappa\cdot\frac{s_\varepsilon^{*2}}{4\varepsilon^2}\cdot q_\varepsilon^2\ln\frac{1}{s_\varepsilon^*}.
		\end{equation*}
		By the odd symmetry of $x_1-q_\varepsilon$ on $x_1=q_\varepsilon$, it holds
		\begin{equation*}
			\begin{split}
				&\quad\frac{1}{\varepsilon^2}\int_{A_\varepsilon}x_1(x_1-q_\varepsilon)d\boldsymbol x=\frac{1}{\varepsilon^2}\int_{A_\varepsilon}x_1'(x_1'-q_\varepsilon)d\boldsymbol x'\\
				&=\frac{1}{\varepsilon^2}\int_{B_{s_\varepsilon^*}(\boldsymbol q_\varepsilon)}q_\varepsilon(x_1-q_\varepsilon)d\boldsymbol x+\frac{1}{\varepsilon^2}\int_{B_{s_\varepsilon^*}(\boldsymbol q_\varepsilon)}(x_1-q_\varepsilon)^2d\boldsymbol x\\
				&\quad\quad\quad\quad+\frac{1}{\varepsilon^2}\left(\int_{A_\varepsilon}x_1(x_1-q_\varepsilon)d\boldsymbol x-\int_{B_{s_\varepsilon^*}(\boldsymbol q_\varepsilon)}x_1(x_1-q_\varepsilon)d\boldsymbol x\right)\\
				&=O(\varepsilon^2)+O\left(\frac{1}{\varepsilon}\right)\cdot|A_\varepsilon\bigtriangleup B_{s_\varepsilon^*}(\boldsymbol q_\varepsilon)|\\
				&=O\left(\varepsilon^2+\frac{\mathbf {er}_\varepsilon}{\varepsilon}\right)
			\end{split}
		\end{equation*}
		Notice that the remaining terms in the product have a higher order on $\varepsilon$. Thus we have shown
		\begin{equation}\label{C-3}
			I_1=\kappa\cdot\frac{s_\varepsilon^{*2}}{4\varepsilon^2}\cdot q_\varepsilon^2\ln\frac{1}{s_\varepsilon^*}+O\left(\varepsilon^2|\ln\varepsilon|+\frac{\mathbf {er}_\varepsilon}{\varepsilon^2}\right).
		\end{equation}
		For the second term $I_2$, we also expand it as
		\begin{equation*}
			\begin{split}
				I_2&=-\frac{q_\varepsilon^2}{4\pi\varepsilon^4}\int_{A_\varepsilon} \left (q_\varepsilon^{-1/2}-\frac{1}{2q_\varepsilon^{3/2}}\cdot(x_1-q_\varepsilon)+O(|x_1-q_\varepsilon|^2)\right)\\
				&\quad \times\int_{A_\varepsilon}\left(q_\varepsilon^{3/2}+\frac{3q_\varepsilon^{1/2}}{2}\cdot(x'_1-q_\varepsilon)+O(|x'_1-q_\varepsilon|^2)\right)\ln\left(\frac{s_\varepsilon^*}{|\boldsymbol x-\boldsymbol x'|}\right)d\boldsymbol x'd\boldsymbol x.
			\end{split}
		\end{equation*}
		Using a similar method as we deal with $I_1$, it holds
		\begin{equation}\label{C-4}
			\begin{split}
				I_2&=-\frac{q_\varepsilon^2}{4\pi\varepsilon^4}\int_{B_{s_\varepsilon^*}(\boldsymbol q_\varepsilon)}x_1^{-1/2}\int_{B_{s_\varepsilon^*}(\boldsymbol q_\varepsilon)}x_1'^{3/2}\ln\left(\frac{s_\varepsilon^*}{|\boldsymbol x-\boldsymbol x'|}\right)d\boldsymbol x'd\boldsymbol x+O\left(\frac{\mathbf {er}_\varepsilon}{\varepsilon^2}\right)\\
				&=-\frac{q_\varepsilon^3}{8\varepsilon^4}\int_{B_{s_\varepsilon^*}(\boldsymbol q_\varepsilon)}(s_\varepsilon^{*2}-|\boldsymbol x-\boldsymbol q_\varepsilon|^2)d\boldsymbol x+O\left(\frac{\mathbf {er}_\varepsilon}{\varepsilon^2}\right)\\
				&=-\frac{s_\varepsilon^{*4}}{16\varepsilon^4}\cdot q_\varepsilon^3\pi+O\left(\frac{\mathbf {er}_\varepsilon}{\varepsilon^2}\right)=-\kappa\cdot\frac{s_\varepsilon^{*2}}{16\varepsilon^2}\cdot q_\varepsilon^2+O\left(\frac{\mathbf {er}_\varepsilon}{\varepsilon^2}\right).
			\end{split}
		\end{equation}
		Now we turn to $I_3$ and obtain
		\begin{equation}\label{C-5}
			\begin{split}
				 I_3&=\frac{q_\varepsilon^2}{2\pi\varepsilon^4}\int_{A_\varepsilon}\int_{A_\varepsilon}\bigg(\big(q_\varepsilon^{1/2}+\frac{1}{2q_\varepsilon^{1/2}}\cdot(x_1-q_\varepsilon)+O(|x_1-q_\varepsilon|^2)\big)\\
				&\quad\quad\quad\quad\times \big(q_\varepsilon^{3/2}+\frac{3q_\varepsilon^{1/2}}{2}\cdot(x_1'-q_\varepsilon)+O(|x_1'-q_\varepsilon|^2)\big)-q_\varepsilon^2\bigg)\cdot\frac{x_1-x_1'}{|\boldsymbol x-\boldsymbol x'|^2}d\boldsymbol x'd\boldsymbol x\\
				&=\frac{q_\varepsilon^3}{4\pi\varepsilon^4}\int_{B_{s_\varepsilon^*}(\boldsymbol q_\varepsilon)}\int_{B_{s_\varepsilon^*}(\boldsymbol q_\varepsilon)}((x_1-q_\varepsilon)+3(x_1'-q_\varepsilon))\cdot\frac{x_1-x_1'}{|\boldsymbol x-\boldsymbol x'|^2}d\boldsymbol x'd\boldsymbol x+O\left(\frac{\mathbf {er}_\varepsilon}{\varepsilon^2}\right)\\
				&=-\frac{q_\varepsilon^3}{2\pi\varepsilon^4}\int_{B_{s_\varepsilon^*}(\boldsymbol q_\varepsilon)}\partial_1\left(\int_{B_{s_\varepsilon^*}(\boldsymbol q_\varepsilon)}(x_1'-q_\varepsilon)\ln\left(\frac{s_\varepsilon^*}{|\boldsymbol x-\boldsymbol x'|}\right)d\boldsymbol x'\right)d\boldsymbol x+O\left(\frac{\mathbf {er}_\varepsilon}{\varepsilon^2}\right)\\
				&=-\frac{q_\varepsilon^3}{\varepsilon^4}\int_{B_{s_\varepsilon^*}(\boldsymbol q_\varepsilon)}\partial_1\left(\frac{s_\varepsilon^{*2}(x_1-q_\varepsilon)}{4}-\frac{|\boldsymbol x-\boldsymbol q_\varepsilon|^2(x_1-q_\varepsilon)}{8}\right)d\boldsymbol x+O\left(\frac{\mathbf {er}_\varepsilon}{\varepsilon^2}\right)\\
				&=-\frac{q_\varepsilon^3}{\varepsilon^4}\int_{B_{s_\varepsilon^*}(\boldsymbol q_\varepsilon)}\left(\frac{s_\varepsilon^{*2}}{4}-\frac{(x_1-q_\varepsilon)^2}{4}-\frac{|\boldsymbol x-\boldsymbol q_\varepsilon|^2}{8}\right)d\boldsymbol x+O\left(\frac{\mathbf {er}_\varepsilon}{\varepsilon^2}\right)\\
				&=-\kappa\cdot\frac{s_\varepsilon^{*2}}{8\varepsilon^2}\cdot q_\varepsilon^2+O\left(\frac{\mathbf {er}_\varepsilon}{\varepsilon^2}\right).
			\end{split}
		\end{equation}
		For the last term $I_{\boldsymbol \rho}$, it is easy to verify that
		\begin{equation}\label{C-6}
			\begin{split}
				I_4&=-\frac{q_\varepsilon^2}{4\pi\varepsilon^4}\int_{A_\varepsilon}x_1^{-1/2}\int_{A_\varepsilon}\left(\frac{x_1'^{3/2}}{2}\cdot(\ln(x_1x_1')+2\ln 8-4)+x_1^{3/2}\right)d\boldsymbol x'd\boldsymbol x\\
				&\quad\quad\quad+\kappa\cdot\frac{s_\varepsilon^{*2}}{4\varepsilon^2}\cdot q_\varepsilon^2+O\left(\frac{\mathbf {er}_\varepsilon}{\varepsilon^2}+\varepsilon^2|\ln\varepsilon|\right)\\
				&=-\kappa\cdot\frac{s_\varepsilon^{*2}}{4\varepsilon^2}\cdot q_\varepsilon^2\left(\ln(8q_\varepsilon)-2\right)+O\left(\frac{\mathbf {er}_\varepsilon}{\varepsilon^2}+\varepsilon^2|\ln\varepsilon|\right).
			\end{split}
		\end{equation}
		
		Combining \eqref{C-3} \eqref{C-4} \eqref{C-5} \eqref{C-6}, we finally obtain
		\begin{equation*}
			-\frac{1}{\varepsilon^2}\int_{B_\delta(\boldsymbol q_\varepsilon)} x_1^2\partial_1\psi_{2,\varepsilon}(\boldsymbol x)\cdot\boldsymbol 1_{A_\varepsilon}(\boldsymbol x)d\boldsymbol x=-\kappa\cdot\frac{s_\varepsilon^{*2}}{4\varepsilon^2}\cdot q_\varepsilon^2\left(\ln\frac{8q_\varepsilon}{s_\varepsilon^*}-\frac{5}{4}\right)+O\left(\frac{\mathbf {er}_\varepsilon}{\varepsilon^2}+\varepsilon^2|\ln\varepsilon|\right),
		\end{equation*}
		which is the desired result.
	\end{proof}

\bigskip
	
	Notice that $s_\varepsilon^{*2}q_\varepsilon\pi=\varepsilon^2\kappa$ by the definition of $s_\varepsilon^*$. From \eqref{C-2}, Lemma \ref{C3} and Lemma \ref{C4}, we obtain a relation of circulation $\kappa$, traveling speed $W$, radius $s_\varepsilon^*$ and location $q_\varepsilon$, which has been used to derive Kelvin--Hicks formula in Proposition \ref{prop3-2}. We summarize this result as follows.
	\begin{lemma}\label{C5}
		It holds
		\begin{equation}\label{C-7}
			Wq_\varepsilon\ln\frac{1}{\varepsilon}-\frac{\kappa}{4\pi}\ln\frac{8q_\varepsilon}{s_\varepsilon^*}+\frac{\kappa}{16\pi}=O\left(\frac{\mathbf{er}_\varepsilon}{\varepsilon^2}+\varepsilon^2|\ln\varepsilon|\right).
		\end{equation}
	\end{lemma}

\bigskip

    In view of Lemma \ref{B3}, we have 
    $$\mathbf {er}_\varepsilon=|A_\varepsilon\bigtriangleup B_{s_\varepsilon^*}(\boldsymbol q_\varepsilon)|=O(\varepsilon^2)\cdot O(|t_{\varepsilon,\mathcal W}+t_{\varepsilon,\tilde\phi_\varepsilon}|+\varepsilon^2|\ln\varepsilon|).$$
    On the other hand, in Section \ref{sec3} we will see that the left hand of \eqref{C-7} is a $O(\varepsilon^2|\ln\varepsilon|)$-perturbation of $\mathcal W_\varepsilon(s_\varepsilon)$. Noting that $$t_{\varepsilon,\mathcal W}(\theta)(\cos\theta,\sin\theta)=\frac{\mathcal W_\varepsilon(s_\varepsilon)}{\mathcal N_\varepsilon}\cdot(\cos\theta,0), \quad \theta\in (0,2\pi]$$ 
    is of the same order as $s_\varepsilon\mathcal W_\varepsilon(s_\varepsilon)$ since $s_\varepsilon\mathcal N_\varepsilon$ is near a constant (the gradient of the stream function for Rankine vortex on vortex boundary). Hence it will hold
    $$t_{\varepsilon,\mathcal W}=O(\varepsilon|t_{\varepsilon,\mathcal W}+t_{\varepsilon,\tilde\phi_\varepsilon}|+\varepsilon^2|\ln\varepsilon|),$$
    which is the crucial fact that we will use in Lemma \ref{lem3-14} for the bootstrap argument.
    
    \bigskip
    
    \textbf{Data available statement:} Data will be made available on reasonable request.
	
	\bigskip
	\bigskip
	
	\phantom{s}
	\thispagestyle{empty}

\end{document}